\documentstyle[12pt]{article}

\makeatletter
\ifx\@auxdefsloaded\relax\endinput
\else\let\@auxdefsloaded\relax\fi

\def\newenvironment{%
   \@ifnextchar *{\@@newenv{\global\@ignoretrue}}{\@@newenv{}*}}

\def\@@newenv#1*#2{%
   \@ifnextchar [{\@newenv{#1}{#2}}{\@newenv{#1}{#2}[0]}}

\long\def\@newenv#1#2[#3]#4#5{%
   \expandafter\newcommand\csname#2\endcsname[#3]{#4}%
   \expandafter\long\expandafter\def\csname end#2\endcsname{#5#1}}

\def\renewenvironment{%
   \@ifnextchar *{\@@renewenv{\global\@ignoretrue}}{\@@renewenv{}*}}

\def\@@renewenv#1*#2{%
   \@ifnextchar [{\@renewenv{#1}{#2}}{\@renewenv{#1}{#2}[0]}}

\long\def\@renewenv#1#2[#3]#4#5{%
   \expandafter\renewcommand\csname#2\endcsname[#3]{#4}%
   \expandafter\long\expandafter\def\csname end#2\endcsname{#5#1}}

\def\newoptcommand#1#2{%
   \@ifnextchar [{\@optargdef#1#2}{\@optargdef#1#2[1]}}

\def\renewoptcommand#1#2{%
   \edef\@tempa{\expandafter\@cdr\string#1\@nil}%
   \@ifundefined{\@tempa}{%
      \@latexerr{\string#1\space undefined}\@ehc}{}%
   \@ifnextchar [{\@reoptargdef#1#2}{\@reoptargdef#1#2[1]}}

\long\def\@optargdef#1#2[#3]#4{%
   \@ifdefinable #1{\@reoptargdef#1#2[#3]{#4}}}

\long\def\@reoptargdef#1#2[#3]#4{%
   \@tempcnta#3\relax \@tempcntb \@ne
   \let#1\relax \let\@tempa\relax
   \edef\@tempb{\long\def\csname\string#1\endcsname
      [\@tempa\the\@tempcntb]}%
   \advance\@tempcntb \@ne \advance\@tempcnta \m@ne
   \@whilenum\@tempcnta>0\do{%
      \edef\@tempb{\@tempb\@tempa\the\@tempcntb}%
      \advance\@tempcntb \@ne \advance\@tempcnta \m@ne}%
   \let\@tempa=##\@tempb{#4}%
   \def#1{\@ifnextchar [{\csname\string#1\endcsname}{%
      \csname\string#1\endcsname[#2]}}}

\def\newoptenvironment{%
   \@ifnextchar *{\@@newoptenv{\global\@ignoretrue}}{%
      \@@newoptenv{}*}}

\def\@@newoptenv#1*#2#3{%
   \@ifnextchar [{\@newoptenv{#1}{#2}{#3}}{%
      \@newoptenv{#1}{#2}{#3}[0]}}

\long\def\@newoptenv#1#2#3[#4]#5#6{%
   \expandafter\newoptcommand\csname#2\endcsname{#3}[#4]{#5}%
   \expandafter\long\expandafter\def\csname end#2\endcsname{#6#1}}
\def\renewoptenvironment{%
   \@ifnextchar *{\@@renewoptenv{\global\@ignoretrue}}{%
      \@@renewoptenv{}*}}

\def\@@renewoptenv#1*#2#3{%
   \@ifnextchar [{\@renewoptenv{#1}{#2}{#3}}{%
      \@renewoptenv{#1}{#2}{#3}[0]}}

\long\def\@renewoptenv#1#2#3[#4]#5#6{%
   \expandafter\renewoptcommand\csname#2\endcsname{#3}[#4]{#5}%
   \expandafter\long\expandafter\def\csname end#2\endcsname{#6#1}}

\newcounter{keepoptional}
\newcounter{optnestctr}

\newoptenvironment*{optional}{1}[1]{%
   \ifnum#1>\value{keepoptional}\relax
      \setcounter{optnestctr}{0}\@powerup\expandafter\@endopt\fi
   \ignorespaces}{}

\def\@powerup{\catcode`\{=12 \catcode`\}=12 \catcode`\\=12 \relax}
\def\@powerdown{\catcode`\{=1 \catcode`\}=2 \catcode`\\=0 \relax}
\begingroup \catcode`|=0 \catcode`[=1 \catcode`]=2 \@powerup
   |long|gdef|@endopt#1\end{optional}[%
      |@beginopt#1\begin{optional}*]%
   |long|gdef|@beginopt#1\begin{optional}[%
      |@ifstar[%
         |ifnum|value[optnestctr]>0|relax
            |addtocounter[optnestctr][-1]|let|@zzzap=|@endopt
         |else
            |@powerdown|end[optional]|let|@zzzap=|relax|fi
	 |@zzzap]%
      [
         |addtocounter[optnestctr][1]|@beginopt]]%
|endgroup

\makeatother
 \makeatletter
\newskip\dgARROWLENGTH  \dgARROWLENGTH=2.5em\relax
\newskip\dgHORIZPAD     \dgHORIZPAD=1em\relax
\newskip\dgVERTPAD      \dgVERTPAD=2ex\relax
\newskip\dgLABELOFFSET  \dgLABELOFFSET=.7ex\relax
\newcount\dgARROWPARTS  \dgARROWPARTS=4\relax
\newcommand{\dgeverynode}{\displaystyle}
\newcommand{\dgeverylabel}{\scriptstyle}
\newskip\dgDOTSPACING   \dgDOTSPACING=0.35em
\newskip\dgDOTSIZE      \dgDOTSIZE=1.5\fontdimen8\tenln
\newcount\dgMAXSQUARE   \dgMAXSQUARE=4\relax
\newcount\dgMINDBLSQ    \dgMINDBLSQ=10\relax
\newskip\dgCOLUMNWIDTH  \dgCOLUMNWIDTH=2em\relax
\chardef\f@ur=4
\@namedef{dgo@}{\let\dg@VECTOR=\vector
   \dg@LBLPOS=\dgARROWPARTS \divide\dg@LBLPOS\tw@}%
\@namedef{dgo@-}{\let\dg@VECTOR=\line}%
\@namedef{dgo@!}{\let\dg@VECTOR=\dg@novector}%
\@namedef{dgo@1}{\dg@LBLPOS=\@ne\relax}%
\@namedef{dgo@2}{\dg@LBLPOS=\tw@\relax}%
\@namedef{dgo@3}{\dg@LBLPOS=\thr@@\relax}%
\@namedef{dgo@4}{\dg@LBLPOS=\f@ur\relax}%
\@namedef{dgo@5}{\dg@LBLPOS=5\relax}%
\@namedef{dgo@6}{\dg@LBLPOS=6\relax}%
\@namedef{dgo@7}{\dg@LBLPOS=7\relax}%
\@namedef{dgo@8}{\dg@LBLPOS=8\relax}%
\@namedef{dgo@9}{\dg@LBLPOS=9\relax}%
\def\dgt@e{\dg@DX=\@ne \dg@DY=\z@ \dg@SIZE=\@ne}%
\def\dgt@w{\dg@DX=\m@ne \dg@DY=\z@ \dg@SIZE=\@ne}%
\def\dgt@n{\dg@DX=\z@ \dg@DY=\@ne \dg@SIZE=\@ne}%
\def\dgt@s{\dg@DX=\z@ \dg@DY=\m@ne \dg@SIZE=\@ne}%
\def\dgt@ne{\dg@DX=\@ne \dg@DY=\@ne \dg@SIZE=\@ne}%
\def\dgt@se{\dg@DX=\@ne \dg@DY=\m@ne \dg@SIZE=\@ne}%
\def\dgt@nw{\dg@DX=\m@ne \dg@DY=\@ne \dg@SIZE=\@ne}%
\def\dgt@sw{\dg@DX=\m@ne \dg@DY=\m@ne \dg@SIZE=\@ne}%
\def\dgt@nne{\dg@DX=\@ne \dg@DY=\tw@ \dg@SIZE=\@ne}%
\def\dgt@nnw{\dg@DX=\m@ne \dg@DY=\tw@ \dg@SIZE=\@ne}%
\def\dgt@sse{\dg@DX=\@ne \dg@DY=-\tw@ \dg@SIZE=\@ne}%
\def\dgt@ssw{\dg@DX=\m@ne \dg@DY=-\tw@ \dg@SIZE=\@ne}%
\def\dgt@ene{\dg@DX=\tw@ \dg@DY=\@ne \dg@SIZE=\tw@}%
\def\dgt@ese{\dg@DX=\tw@ \dg@DY=\m@ne \dg@SIZE=\tw@}%
\def\dgt@wnw{\dg@DX=-\tw@ \dg@DY=\@ne \dg@SIZE=\tw@}%
\def\dgt@wsw{\dg@DX=-\tw@ \dg@DY=\m@ne \dg@SIZE=\tw@}%
\def\dggeometry{
   \dg@ZTEMP=\dg@XGRID \multiply\dg@ZTEMP\tw@
   \ifnum\dg@YGRID=\z@ \dg@ZTEMP=\tw@
   \else \divide\dg@ZTEMP\dg@YGRID \fi
   \ifnum\dg@ZTEMP>\f@ur \dg@ZTEMP=\f@ur \fi
   \ifnum\dg@ZTEMP<\@ne \dg@ZTEMP=\@ne \fi
   \unitlength=2sp\relax
   \ifnum\dg@ZTEMP<\tw@
      \advance\dg@ZTEMP\@ne
      \multiply\unitlength\dg@YGRID
   \else
      \multiply\unitlength\dg@XGRID \divide\unitlength\dg@ZTEMP
   \fi
   \dg@XGRID=\dg@ZTEMP \dg@YGRID=\tw@
   \dg@rmcommondiv\tw@\dg@XGRID\dg@YGRID
   \divide\unitlength\dg@YGRID \divide\unitlength\@m\relax}
\@namedef{dgo@..}{\let\dg@VECTOR=\dg@dotvector}%

\def\dg@dotvector(#1,#2)#3{%
   \begingroup
   \dg@XTEMP=#1\relax \dg@YTEMP=#2\relax
   \let\dg@NDOTS=\dg@XEND \let\dg@DOTDIAM=\dg@WEND
   \dg@NDOTS=\unitlength \multiply\dg@NDOTS #3\relax
   \dg@ZTEMP=\dg@YTEMP \dg@changesign\dg@YTEMP\dg@ZTEMP
   \ifnum\dg@XTEMP>\z@
      \ifnum\dg@YTEMP>\dg@XTEMP
         \multiply\dg@NDOTS\dg@YTEMP \divide\dg@NDOTS\dg@XTEMP \fi
   \else\ifnum\dg@XTEMP<\z@
      \ifnum\dg@YTEMP>-\dg@XTEMP
         \multiply\dg@NDOTS\dg@YTEMP \divide\dg@NDOTS-\dg@XTEMP \fi
   \fi\fi
   \dg@YTEMP=\dg@ZTEMP
   \divide\dg@NDOTS\dgDOTSPACING
   \ifnum\dg@NDOTS>\z@\else \dg@NDOTS=\@ne \fi
   \dg@ZTEMP=\unitlength \multiply\dg@ZTEMP #3\relax
   \divide\dg@ZTEMP\dg@NDOTS
   \ifnum\dg@XTEMP=\z@
      \dg@changesign\dg@ZTEMP\dg@YTEMP \dg@YTEMP=\dg@ZTEMP
   \else
      \dg@changesign\dg@ZTEMP\dg@XTEMP
      \multiply\dg@YTEMP\dg@ZTEMP \divide\dg@YTEMP\dg@XTEMP
      \dg@XTEMP=\dg@ZTEMP
   \fi
   \divide\dg@XTEMP\unitlength \divide\dg@YTEMP\unitlength
   \begin{picture}(0,0)
      \dg@DOTDIAM=\dgDOTSIZE \divide\dg@DOTDIAM\unitlength
      \multiput(0,0)(\dg@XTEMP,\dg@YTEMP){\dg@NDOTS}{%
         \circle*{\dg@DOTDIAM}}%
      \multiply\dg@XTEMP\dg@NDOTS \multiply\dg@YTEMP\dg@NDOTS
      \put(\dg@XTEMP,\dg@YTEMP){\vector(#1,#2){0}}%
   \end{picture}%
   \endgroup}%
\newif\ifdg@LATEXGEOM \dg@LATEXGEOMfalse

\@ifundefined{lamsvector}{}{%
   \@namedef{dgo@}{%
      \let\dg@VECTOR=\lamsvector
      \lamsshaft{?}\lamssource{?}\lamstarget{?}%
      \dg@LBLPOS=\dgARROWPARTS \divide\dg@LBLPOS\tw@\relax}%
   \@namedef{dgo@..}{\lamsshaft{.}}%
   \@namedef{dgo@=}{\lamsshaft{=}}%
   \@namedef{dgo@-}{\lamstarget{0}}%
   \@namedef{dgo@'}{\lamstarget{'}}%
   \@namedef{dgo@`}{\lamstarget{`}}%
   \@namedef{dgo@A}{\lamstarget{A}}%
   \@namedef{dgo@V}{\lamssource{V}}%
   \@namedef{dgo@J}{\lamssource{J}}%
   \@namedef{dgo@L}{\lamssource{L}}%
   \@namedef{dgo@S}{\lamssource{S}}%
   \def\dggeometry{
      \dg@ZTEMP=\dg@XGRID \multiply\dg@ZTEMP\tw@
      \ifnum\dg@YGRID=\z@ \dg@ZTEMP=\tw@
      \else \divide\dg@ZTEMP\dg@YGRID \fi
      \ifnum\dg@ZTEMP>6\relax \dg@ZTEMP=6\relax \fi
      \ifdg@LATEXGEOM\ifnum\dg@ZTEMP>\f@ur \dg@ZTEMP=\f@ur \fi\fi
      \ifnum\dg@ZTEMP<\@ne \dg@ZTEMP=\@ne \fi
      \unitlength=2sp\relax
      \ifnum\dg@ZTEMP<\tw@
         \advance\dg@ZTEMP\@ne
         \multiply\unitlength\dg@YGRID
      \else
         \multiply\unitlength\dg@XGRID \divide\unitlength\dg@ZTEMP
      \fi
      \dg@XGRID=\dg@ZTEMP \dg@YGRID=\tw@
      \dg@rmcommondiv\tw@\dg@XGRID\dg@YGRID
      \divide\unitlength\dg@YGRID \divide\unitlength\@m
      \dg@LATEXGEOMfalse}
   \def\dgt@nee{\dg@DX=\tw@ \dg@DY=\@ne \dg@SIZE=\tw@}%
   \def\dgt@see{\dg@DX=\tw@ \dg@DY=\m@ne \dg@SIZE=\tw@}%
   \def\dgt@nww{\dg@DX=-\tw@ \dg@DY=\@ne \dg@SIZE=\tw@}%
   \def\dgt@sww{\dg@DX=-\tw@ \dg@DY=\m@ne \dg@SIZE=\tw@}%
   \def\dgt@nnne{\dg@DX=\@ne \dg@DY=\thr@@ \dg@SIZE=\@ne}%
   \def\dgt@nnnw{\dg@DX=\m@ne \dg@DY=\thr@@ \dg@SIZE=\@ne}%
   \def\dgt@sssw{\dg@DX=\m@ne \dg@DY=-\thr@@ \dg@SIZE=\@ne}%
   \def\dgt@ssse{\dg@DX=\@ne \dg@DY=-\thr@@ \dg@SIZE=\@ne}%
   \def\dgt@nnnee{\dg@DX=\tw@ \dg@DY=\thr@@ \dg@SIZE=\tw@}%
   \def\dgt@nnnww{\dg@DX=-\tw@ \dg@DY=\thr@@ \dg@SIZE=\tw@}%
   \def\dgt@sssww{\dg@DX=-\tw@ \dg@DY=-\thr@@ \dg@SIZE=\tw@}%
   \def\dgt@sssee{\dg@DX=\tw@ \dg@DY=-\thr@@ \dg@SIZE=\tw@}%
   \def\dgt@nneee{\dg@DX=\thr@@ \dg@DY=\tw@ \dg@SIZE=\thr@@}%
   \def\dgt@nnwww{\dg@DX=-\thr@@ \dg@DY=\tw@ \dg@SIZE=\thr@@}%
   \def\dgt@sswww{\dg@DX=-\thr@@ \dg@DY=-\tw@ \dg@SIZE=\thr@@}%
   \def\dgt@sseee{\dg@DX=\thr@@ \dg@DY=-\tw@ \dg@SIZE=\thr@@}%
   \def\dgt@neee{\dg@DX=\thr@@ \dg@DY=\@ne \dg@SIZE=\thr@@
      \global\dg@LATEXGEOMtrue}%
   \def\dgt@nwww{\dg@DX=-\thr@@ \dg@DY=\@ne \dg@SIZE=\thr@@
      \global\dg@LATEXGEOMtrue}%
   \def\dgt@swww{\dg@DX=-\thr@@ \dg@DY=\m@ne \dg@SIZE=\thr@@
      \global\dg@LATEXGEOMtrue}%
   \def\dgt@seee{\dg@DX=\thr@@ \dg@DY=\m@ne \dg@SIZE=\thr@@
      \global\dg@LATEXGEOMtrue}%
}
\newcount\dg@HORIZ      \newcount\dg@VERT
\newcount\dg@XLPAD      \newcount\dg@YBPAD
\newcount\dg@XRPAD      \newcount\dg@YTPAD
\newcount\dg@X          \newcount\dg@Y
\newcount\dg@XGRID      \newcount\dg@YGRID
\newcount\dg@SIZE
\newcount\dg@USERSIZE
\newcount\dg@DX         \newcount\dg@DY
\newcount\dg@XLBL       \newcount\dg@YLBL
\newcount\dg@XOFFSET    \newcount\dg@YOFFSET
\newcount\dg@LBLOFF
\newcount\dg@XLBLOFF    \newcount\dg@YLBLOFF
\newcount\dg@LBLPOS
\newcount\dg@XTEMP      \newcount\dg@YTEMP
\newcount\dg@ZTEMP
\newcount\dg@XNODE      \newcount\dg@YNODE
\newcount\dg@XEND       \newcount\dg@YEND
\newcount\dg@WEND       \newcount\dg@HEND
\newcount\dg@COUNT
\newbox\dg@NODEBOX
\newoptenvironment*{diagram}{}[1]{%
   \global\advance\dg@COUNT\@ne \typeout{[diagram \the\dg@COUNT:}%
   \let\node=\dg@node \let\\=\dg@cr \let\arrow=\dg@arrow
   \def\dg@BIGNODE{#1}%
   \ifx\dg@BIGNODE\@empty\else
      \setbox\dg@NODEBOX=\hbox{$\dgeverynode{#1}$}%
      \dg@XTEMP=\wd\dg@NODEBOX
      \dg@YTEMP=\ht\dg@NODEBOX \advance\dg@YTEMP\dp\dg@NODEBOX
      \ifnum\dg@YTEMP=\z@ \dg@ZTEMP=\z@
      \else \dg@ZTEMP=\dg@XTEMP \divide\dg@ZTEMP\dg@YTEMP\fi
      \ifnum\dg@ZTEMP>\dgMAXSQUARE
         \ifnum\dg@ZTEMP<\dgMINDBLSQ \dg@XGRID=\thr@@ \dg@YGRID=\tw@
         \else \dg@XGRID=\tw@ \dg@YGRID=\@ne \fi
      \else \dg@XGRID=\@ne \dg@YGRID=\@ne \fi
      \advance\dg@XTEMP\dgHORIZPAD \advance\dg@YTEMP\dgVERTPAD
      \dg@XLPAD=\dg@XTEMP \divide\dg@XLPAD\tw@
      \dg@XRPAD=\dg@XTEMP \divide\dg@XRPAD\tw@
      \dg@YBPAD=\dg@YTEMP \divide\dg@YBPAD\tw@
      \dg@YTPAD=\dg@YTEMP \divide\dg@YTPAD\tw@
      \advance\dg@XTEMP\dgARROWLENGTH
      \divide\dg@XTEMP\dg@XGRID \divide\dg@XTEMP\@m
      \unitlength=1sp\relax \multiply\unitlength\dg@XTEMP
   \fi
   \typeout{saving...}%
   \dg@X=-\@ne \dg@Y=\z@ \dg@HORIZ=\z@\relax
   \let\dg@SLIST=\@empty
   \let\dg@NLIST=\@empty \let\dg@ALIST=\@empty
   \let\dg@PASS=\dg@savepass
}{
   \dg@VERT=-\dg@Y\relax
   \typeout{calculating...}%
   \ifx\dg@BIGNODE\@empty
      \dg@XGRID=\z@ \dg@YGRID=\z@
      \dg@XLPAD=\z@ \dg@XRPAD=\z@ \dg@YBPAD=\z@ \dg@YTPAD=\z@
      \let\dg@PASS=\dg@geompass
      \dg@NLIST \dg@ALIST
      \ifnum\dg@XGRID>\z@\else \dg@XGRID=\quad\relax \fi
      \ifnum\dg@YGRID>\z@\else \dg@YGRID=\quad\relax \fi
      \dggeometry
      \ifdim\unitlength>\z@\else \unitlength=\quad \fi
      \ifnum\dg@XGRID>\z@\else \dg@XGRID=\@ne \fi
      \ifnum\dg@YGRID>\z@\else \dg@YGRID=\@ne \fi
   \fi
   \dg@LBLOFF=\dgLABELOFFSET \divide\dg@LBLOFF\unitlength
   \multiply\dg@HORIZ\@m \multiply\dg@HORIZ\dg@XGRID
   \multiply\dg@VERT\@m \multiply\dg@VERT\dg@YGRID
   \divide\dg@XLPAD\unitlength \divide\dg@XRPAD\unitlength
   \divide\dg@YBPAD\unitlength \divide\dg@YTPAD\unitlength
   \typeout{drawing...}%
   \let\dg@PASS=\dg@drawpass
   \hspace*{\dg@XLPAD\unitlength}%
   \vcenter{%
      \hsize=0pt\relax\parindent=0pt\relax
      \parskip=0pt\relax\baselineskip=0pt\relax
      \vspace*{\dg@YTPAD\unitlength}%
      \begin{picture}(0,0)(0,0)\dg@NLIST\dg@ALIST\end{picture}%
      \raisebox{\z@}[\z@][\dg@VERT\unitlength]{}%
      \vspace*{\dg@YBPAD\unitlength}%
   }
   \hspace*{\dg@HORIZ\unitlength}%
   \hspace*{\dg@XRPAD\unitlength}%
   \typeout{done]}}%

\def\dg@savepass{s}
\def\dg@geompass{g}
\def\dg@drawpass{d}

\newoptcommand{\dg@node}{\@ne}[2]{%
   \ifx\dg@PASS\dg@savepass
      \dg@XTEMP=#1\relax
      \ifnum\dg@XTEMP<\@ne \dg@XTEMP=\@ne\fi
      \advance\dg@X\dg@XTEMP
      \ifnum\dg@HORIZ<\dg@X \dg@HORIZ=\dg@X \fi
      \setbox\dg@NODEBOX=\hbox{$\dgeverynode{#2}$}%
      \dg@XTEMP=\wd\dg@NODEBOX \advance\dg@XTEMP\dgHORIZPAD
      \dg@YTEMP=\ht\dg@NODEBOX \advance\dg@YTEMP\dp\dg@NODEBOX
      \advance\dg@YTEMP\dgVERTPAD
      \toks\z@=\expandafter{\dg@SLIST}%
      \edef\dg@SLIST{\the\toks\z@
         ,\noexpand\dg@XNODE=\number\dg@X\noexpand\relax
         \noexpand\dg@YNODE=\number\dg@Y\noexpand\relax
         \noexpand\dg@XTEMP=\number\dg@XTEMP\noexpand\relax
         \noexpand\dg@YTEMP=\number\dg@YTEMP\noexpand\relax}%
      \toks\z@=\expandafter{\dg@NLIST}%
      \toks\tw@={\dg@node{#2}}%
      \edef\dg@NLIST{\the\toks\z@
         \noexpand\dg@X=\number\dg@X\noexpand\relax
         \noexpand\dg@Y=\number\dg@Y\noexpand\relax
         \the\toks\tw@}%
   \else\ifx\dg@PASS\dg@geompass
      \ifnum\dg@X=\z@
         \dg@getnodesize
            {\dg@SLIST}{\dg@X}{\dg@Y}{\dg@WEND}{\dg@HEND}%
         \divide\dg@WEND\tw@
         \ifnum\dg@XLPAD<\dg@WEND \dg@XLPAD=\dg@WEND \fi\fi
      \ifnum\dg@X=\dg@HORIZ
         \dg@getnodesize
            {\dg@SLIST}{\dg@X}{\dg@Y}{\dg@WEND}{\dg@HEND}%
         \divide\dg@WEND\tw@
         \ifnum\dg@XRPAD<\dg@WEND \dg@XRPAD=\dg@WEND \fi\fi
      \ifnum\dg@Y=\z@
         \dg@getnodesize
            {\dg@SLIST}{\dg@X}{\dg@Y}{\dg@WEND}{\dg@HEND}%
         \divide\dg@HEND\tw@
         \ifnum\dg@YTPAD<\dg@HEND \dg@YTPAD=\dg@HEND \fi\fi
      \ifnum\dg@Y=-\dg@VERT\relax
         \dg@getnodesize
            {\dg@SLIST}{\dg@X}{\dg@Y}{\dg@WEND}{\dg@HEND}%
         \divide\dg@HEND\tw@
         \ifnum\dg@YBPAD<\dg@HEND \dg@YBPAD=\dg@HEND \fi\fi
   \else\ifx\dg@PASS\dg@drawpass
      \dg@XNODE=\dg@X \multiply\dg@XNODE\@m
      \multiply\dg@XNODE\dg@XGRID
      \dg@YNODE=\dg@Y \multiply\dg@YNODE\@m
      \multiply\dg@YNODE\dg@YGRID
      \setbox\dg@NODEBOX=\hbox{$\dgeverynode{#2}$}%
      \put(\dg@XNODE,\dg@YNODE){\dg@makebox{\box\dg@NODEBOX}}%
   \fi\fi\fi}%

\newoptcommand{\dg@cr}{\@ne}[1]{%
   \ifx\dg@PASS\dg@savepass
      \dg@YTEMP=#1\relax
      \ifnum\dg@YTEMP<\@ne \dg@YTEMP=\@ne \fi
      \advance\dg@Y -\dg@YTEMP\relax
      \dg@X=-\@ne\relax\fi}%
\newoptcommand{\dg@arrow}{\@ne}[2]{%
   \begingroup
   \dg@USERSIZE=#1\relax
   \ifnum\dg@USERSIZE<\@ne \dg@USERSIZE=\@ne \fi
   \dg@parse{#2}%
   \ifx\dg@PASS\dg@savepass
      \ifx\dg@label\dgl@b \let\dg@label=\dgl@t \fi
      \ifx\dg@label\dgl@r \let\dg@label=\dgl@l \fi
      \let\dg@process=\dg@save
   \else\ifx\dg@PASS\dg@geompass
      \let\dg@process=\dg@ignore
      \dg@geomcalc
   \else\ifx\dg@PASS\dg@drawpass
      \let\dg@process=\dg@draw
      \dg@drawcalc
   \fi\fi\fi
   \dg@label{\dg@process{#1}{#2}}}%

\newoptcommand{\arrow}{\@ne}[2]{%
   \dg@parse{#2}%
   \ifx\dg@label\dgl@b \let\dg@label=\dgl@t \fi
   \ifx\dg@label\dgl@r \let\dg@label=\dgl@l \fi
   \dg@label{\dg@textarrow{#1}{#2}}}%

\def\dg@textarrow#1#2#3#4{%
   \mathop{{\dgHORIZPAD=0pt\relax\dgVERTPAD=0pt\relax
      \begin{diagram}
         \node{}\arrow[#1]{#2}{#3}{#4}\node{}
      \end{diagram}}}}

\def\dg@parse#1{%
   \let\dg@label=\dgl@ \dgo@
   \let\dg@type=\@empty \let\dg@lbltype=\@empty
   \@for\dg@list:=#1\do{%
      \ifx\dg@type\@empty \let\dg@type=\dg@list
      \else\ifx\dg@lbltype\@empty \let\dg@lbltype=\dg@list
         \@ifundefined{dgo@\dg@list}{}{\@nameuse{dgo@\dg@list}}%
      \else
         \@ifundefined{dgo@\dg@list}{}{\@nameuse{dgo@\dg@list}}%
      \fi\fi}%
   \@ifundefined{dgt@\dg@type}{\dgt@e}{\@nameuse{dgt@\dg@type}}%
   \@ifundefined{dgl@\dg@lbltype}{}{%
      \dg@letname\dg@label{dgl@\dg@lbltype}}}

\def\dg@draw#1#2#3#4{%
   \put(\dg@X,\dg@Y){\dg@makebox{%
      \begin{picture}(0,0)%
         \thinlines
         \put(\dg@XOFFSET,\dg@YOFFSET){%
            \dg@VECTOR(\dg@DX,\dg@DY){\dg@SIZE}}%
         \put(\dg@XLBL,\dg@YLBL){\dg@makebox{%
            \begin{picture}(0,0)%
               \put(\dg@XLBLOFF,\dg@YLBLOFF){%
                  \dg@makebox[\dg@LBLONE]{$\dgeverylabel{#3}$}}%
               \put(-\dg@XLBLOFF,-\dg@YLBLOFF){%
                  \dg@makebox[\dg@LBLTWO]{$\dgeverylabel{#4}$}}%
            \end{picture}}}%
      \end{picture}}}%
   \endgroup}%

\def\dg@save#1#2#3#4{%
   \endgroup 
   \toks\z@=\expandafter{\dg@ALIST}%
   \toks\tw@={\dg@arrow[#1]{#2}{#3}{#4}}%
   \edef\dg@ALIST{\the\toks\z@%
      \noexpand\dg@X=\number\dg@X\noexpand\relax
      \noexpand\dg@Y=\number\dg@Y\noexpand\relax
      \the\toks\tw@}}%

\def\dg@ignore#1#2#3#4{\endgroup}

\def\dg@geomcalc{%
   \dg@XEND=\dg@SIZE \multiply\dg@XEND\dg@USERSIZE
   \ifnum\dg@DX=\z@
      \dg@YEND=\dg@XEND \dg@XEND=\z@
      \dg@changesign\dg@YEND\dg@DY
   \else
      \dg@changesign\dg@XEND\dg@DX \dg@YEND=\dg@XEND
      \multiply\dg@YEND\dg@DY \divide\dg@YEND\dg@DX
   \fi
   \advance\dg@XEND\dg@X \advance\dg@YEND\dg@Y
   \dg@getnodesize
      {\dg@SLIST}{\dg@XEND}{\dg@YEND}{\dg@WEND}{\dg@HEND}%
   \dg@XOFFSET=\dg@WEND \dg@YOFFSET=\dg@HEND
   \dg@getnodesize
      {\dg@SLIST}{\dg@X}{\dg@Y}{\dg@WEND}{\dg@HEND}%
   \advance\dg@XOFFSET\dg@WEND \divide\dg@XOFFSET\tw@
   \advance\dg@YOFFSET\dg@HEND \divide\dg@YOFFSET\tw@
   \dg@XTEMP=\dgARROWLENGTH \dg@YTEMP=\dgARROWLENGTH
   \ifnum\dg@DX>\z@
      \dg@ZTEMP=\dg@DX \multiply\dg@XTEMP\dg@DX
   \else \dg@ZTEMP=-\dg@DX \multiply\dg@XTEMP -\dg@DX \fi
   \ifnum\dg@DY>\z@
      \advance\dg@ZTEMP\dg@DY \multiply\dg@YTEMP\dg@DY
   \else \advance\dg@ZTEMP -\dg@DY \multiply\dg@YTEMP -\dg@DY\fi
   \ifnum\dg@ZTEMP=\z@\else
      \divide\dg@XTEMP\dg@ZTEMP \divide\dg@YTEMP\dg@ZTEMP
      \advance\dg@XOFFSET\dg@XTEMP \advance\dg@YOFFSET\dg@YTEMP
   \fi
   \divide\dg@XOFFSET\dg@SIZE \divide\dg@XOFFSET\dg@USERSIZE
   \divide\dg@YOFFSET\dg@SIZE \divide\dg@YOFFSET\dg@USERSIZE
   \ifnum\dg@DX=\z@ \dg@XOFFSET=\z@ \fi
   \ifnum\dg@DY=\z@ \dg@YOFFSET=\z@ \fi
   \ifnum\dg@XGRID<\dg@XOFFSET \global\dg@XGRID=\dg@XOFFSET\fi
   \ifnum\dg@YGRID<\dg@YOFFSET \global\dg@YGRID=\dg@YOFFSET\fi
   \relax}

\def\dg@drawcalc{%
   \dg@XEND=\dg@SIZE \multiply\dg@XEND\dg@USERSIZE
   \ifnum\dg@DX=\z@
      \dg@YEND=\dg@XEND \dg@XEND=\z@
      \dg@changesign\dg@YEND\dg@DY
   \else
      \dg@changesign\dg@XEND\dg@DX \dg@YEND=\dg@XEND
      \multiply\dg@YEND\dg@DY \divide\dg@YEND\dg@DX
   \fi
   \advance\dg@XEND\dg@X \advance\dg@YEND\dg@Y
   \dg@getnodesize
      {\dg@SLIST}{\dg@XEND}{\dg@YEND}{\dg@WEND}{\dg@HEND}%
   \divide\dg@WEND\unitlength \divide\dg@HEND\unitlength
   \multiply\dg@DX\dg@XGRID \multiply\dg@DY\dg@YGRID
   \dg@rmcommondiv\tw@\dg@DX\dg@DY
   \dg@rmcommondiv\tw@\dg@DX\dg@DY 
   \dg@rmcommondiv\thr@@\dg@DX\dg@DY
   \multiply\dg@SIZE\dg@USERSIZE \multiply\dg@SIZE\@m
   \ifnum\dg@DX=\z@
      \multiply\dg@SIZE\dg@YGRID
      \divide\dg@HEND\tw@ \advance\dg@SIZE -\dg@HEND
      \dg@getnodesize
         {\dg@SLIST}{\dg@X}{\dg@Y}{\dg@WEND}{\dg@YOFFSET}%
      \divide\dg@YOFFSET\unitlength \divide\dg@YOFFSET\tw@
      \advance\dg@SIZE -\dg@YOFFSET
      \dg@XOFFSET=\z@
      \def\dg@LBLONE{r}\def\dg@LBLTWO{l}%
      \dg@XLBL=\z@ \dg@YLBL=\dg@SIZE
      \multiply\dg@YLBL\dg@LBLPOS
      \divide\dg@YLBL\dgARROWPARTS\relax
      \advance\dg@YLBL\dg@YOFFSET
      \dg@changesign\dg@YLBL\dg@DY
      \dg@changesign\dg@YOFFSET\dg@DY
   \else
      \multiply\dg@SIZE\dg@XGRID
      \ifnum\dg@DY=\z@
         \divide\dg@WEND\tw@ \advance\dg@SIZE -\dg@WEND
         \dg@getnodesize
            {\dg@SLIST}{\dg@X}{\dg@Y}{\dg@XOFFSET}{\dg@HEND}%
         \divide\dg@XOFFSET\unitlength \divide\dg@XOFFSET\tw@
         \advance\dg@SIZE -\dg@XOFFSET
         \dg@YOFFSET=\z@
         \def\dg@LBLONE{b}\def\dg@LBLTWO{t}%
         \dg@YLBL=\z@ \dg@XLBL=\dg@SIZE
         \multiply\dg@XLBL\dg@LBLPOS
         \divide\dg@XLBL\dgARROWPARTS\relax
         \advance\dg@XLBL\dg@XOFFSET
         \dg@changesign\dg@XLBL\dg@DX
         \dg@changesign\dg@XOFFSET\dg@DX
      \else
         \divide\dg@WEND\tw@ \divide\dg@HEND\tw@
         \multiply\dg@HEND\dg@DX \divide\dg@HEND\dg@DY
         \ifnum\dg@HEND<\z@ \multiply\dg@HEND\m@ne \fi
         \ifnum\dg@WEND<\dg@HEND \advance\dg@SIZE -\dg@WEND
         \else \advance\dg@SIZE -\dg@HEND \fi
         \dg@getnodesize
            {\dg@SLIST}{\dg@X}{\dg@Y}{\dg@WEND}{\dg@HEND}%
         \divide\dg@WEND\unitlength \divide\dg@WEND\tw@
         \divide\dg@HEND\unitlength \divide\dg@HEND\tw@
         \multiply\dg@HEND\dg@DX \divide\dg@HEND\dg@DY
         \ifnum\dg@HEND<\z@ \multiply\dg@HEND\m@ne \fi
         \ifnum\dg@WEND<\dg@HEND \dg@XOFFSET=\dg@WEND
         \else \dg@XOFFSET=\dg@HEND \fi
         \advance\dg@SIZE -\dg@XOFFSET
         \dg@changesign\dg@XOFFSET\dg@DX
         \dg@YOFFSET=\dg@XOFFSET
         \multiply\dg@YOFFSET\dg@DY \divide\dg@YOFFSET\dg@DX
         \def\dg@LBLONE{br}\def\dg@LBLTWO{tl}%
         \ifnum\dg@DX<\z@ \ifnum\dg@DY>\z@
            \def\dg@LBLONE{bl}\def\dg@LBLTWO{tr}\fi\fi
         \ifnum\dg@DX>\z@ \ifnum\dg@DY<\z@
            \def\dg@LBLONE{bl}\def\dg@LBLTWO{tr}\fi\fi
         \dg@XLBL=\dg@SIZE
         \multiply\dg@XLBL\dg@LBLPOS
         \divide\dg@XLBL\dgARROWPARTS\relax
         \dg@changesign\dg@XLBL\dg@DX
         \dg@YLBL=\dg@XLBL
         \multiply\dg@YLBL\dg@DY \divide\dg@YLBL\dg@DX
         \advance\dg@XLBL\dg@XOFFSET
         \advance\dg@YLBL\dg@YOFFSET
      \fi
   \fi
   \dg@XLBLOFF=-\dg@DY \dg@changesign\dg@XLBLOFF\dg@DX
   \dg@YLBLOFF=\dg@DX \dg@changesign\dg@YLBLOFF\dg@DX
   \ifnum\dg@DX=\z@ \dg@XLBLOFF=\m@ne \fi
   \dg@XTEMP=\dg@DX \dg@changesign\dg@XTEMP\dg@DX
   \dg@YTEMP=\dg@DY \dg@changesign\dg@YTEMP\dg@DY
   \ifnum\dg@YTEMP>\dg@XTEMP \dg@XTEMP=\dg@YTEMP \fi
   \ifnum\dg@XTEMP=\z@ \dg@XTEMP=\@ne \fi
   \multiply\dg@XLBLOFF\dg@LBLOFF \divide\dg@XLBLOFF\dg@XTEMP
   \multiply\dg@YLBLOFF\dg@LBLOFF \divide\dg@YLBLOFF\dg@XTEMP
   \multiply\dg@X\@m \multiply\dg@X\dg@XGRID
   \multiply\dg@Y\@m \multiply\dg@Y\dg@YGRID
   \relax}%

\def\dg@rmcommondiv#1#2#3{%
   \dg@XTEMP=#2\relax
   \divide\dg@XTEMP #1\relax \multiply\dg@XTEMP #1\relax
   \dg@YTEMP=#3\relax
   \divide\dg@YTEMP #1\relax \multiply\dg@YTEMP #1\relax
   \ifnum\dg@XTEMP=#2\relax \ifnum\dg@YTEMP=#3\relax
      \divide#2#1\relax \divide#3#1\relax \fi\fi}%

\def\dg@changesign#1#2{%
   \ifnum #2<\z@ \multiply#1\m@ne
   \else\ifnum #2=\z@ #1=\z@ \fi\fi}%

\def\dg@getnodesize#1#2#3#4#5{%
   #4=\z@\relax #5=\z@\relax
   \expandafter\@for\expandafter\dg@trynode
   \expandafter:\expandafter=#1\do{%
      \dg@XNODE=\m@ne 
      \dg@trynode
      \ifnum #2=\dg@XNODE \ifnum #3=\dg@YNODE
         #4=\dg@XTEMP\relax #5=\dg@YTEMP\relax\fi\fi}}%

\newoptcommand{\dg@makebox}{}[2]{%
   \expandafter\makebox\expandafter(\expandafter
      0\expandafter,\expandafter0\expandafter)\expandafter
      [#1]{#2}}%

\def\dg@novector(#1,#2)#3{}%

\def\dg@letname#1#2{%
   \relax\expandafter
   \let\expandafter #1\csname #2\endcsname\relax}%

\def\dgl@#1{#1{}{}}%
\def\dgl@t#1#2{#1{#2}{}}%
\def\dgl@b#1#2{#1{}{#2}}%
\def\dgl@tb#1#2#3{#1{#2}{#3}}%
\def\dgl@l#1#2{#1{#2}{}}%
\def\dgl@r#1#2{#1{}{#2}}%
\def\dgl@lr#1#2#3{#1{#2}{#3}}%
\makeatother

\normalbaselines
\def\rank{\mathop{\rm rank}}

\def\rdots{\mathinner{\mkern1mu\raise1pt\vbox{\kern1pt\hbox{.}}\mkern2mu
   \raise4pt\hbox{.}\mkern2mu\raise7pt\hbox{.}\mkern1mu}}
\newcommand{\Z}{{\Bbb Z}}
\newcommand{\bP}{{\Bbb P}}
\newcommand{\Q}{{\Bbb Q}}
\newcommand{\A}{{\Bbb A}}
\newcommand{\R}{{\Bbb R}}
\newcommand{\h}{{\Bbb H}}
\newcommand{\C}{{\Bbb C}}
\newcommand{\T}{{\cal T}} 
\newcommand{\goth}[1]{#1}
\newcommand{\M}{{\cal M}}
\newtheorem{theorem}{Theorem}[section]
\newtheorem{example}[theorem]{Example}
\newtheorem{examples}[theorem]{Examples}
\newtheorem{rem}[theorem]{Remark}
\newtheorem{rems}[theorem]{Remarks}
\newtheorem{lem}[theorem]{Lemma}
\newtheorem{proposition}[theorem]{Proposition}
\newtheorem{corollary}[theorem]{Corollary}
\newtheorem{condition}[theorem]{Condition}
\newtheorem{convention}[theorem]{Convention}
\newtheorem{definition}[theorem]{Definition}
\newtheorem{question}[theorem]{Question}
\newcounter{bean}
\newcounter{bacon}

\def\J{{\cal J}}
\def\A{{\cal A}} 
\def\O{{\cal O}}
\def\L{{\cal L}}
\def\E{{\cal E}}

\def\P{{\Bbb P}} 
\def\K{{\widetilde{K}}}

\def\proclaim #1. #2\par{\medbreak{\bf#1.\enspace}{\it#2}\par
  \ifdim\lastskip<\medskipamount
  \removelastskip\penalty55\medskip\fi}
\newskip\Bigskipamount
\newcommand{\EndProof}{
\begin{flushright}
{$\Box$}
\end{flushright}
}
\newcommand{\RelPic}{{\cal P}ic}
\newcommand{\DisjointUnion}{\cup}
\newcommand{\RealNumbers}{{\Bbb R}}
\newcommand{\Integers}{{\Bbb Z}}
\newcommand{\ComplexNumbers}{{\Bbb C}}

\newcommand{\PiOne}{{\Bbb P}^1}

\newcommand{\IsomRightArrow}{\stackrel{\cong}{\rightarrow}}
\newcommand{\RightArrowOf}[1]{\stackrel{#1}{\rightarrow}}
\newcommand{\LongRightArrowOf}[1]{\stackrel{#1}{\longrightarrow}}

\newcommand{\HookRightArrowOf}[1]{\stackrel{#1}{\hookrightarrow}}

\newcommand{\SurjectiveRightArrow}{\twoheadrightarrow}
\newcommand{\DotRightArrow}{\rightarrow}
\newcommand{\StructureSheaf}[1]{{\cal O}_{#1}}
\newcommand{\CompletedSheafOfAt}[2]{\widehat{{\cal O}_{#1,#2}}}
\newcommand{\restricted}[2]{#1_{\mid_{#2}}}
\newcommand{\Normal}[2]{N_{#1\mid#2}}

\newcommand{\Sym}{\rm Sym}
\newcommand{\Ext}{\rm Ext}
\newcommand{\Hom}{\rm Hom}
\newcommand{\End}{\rm End}
\newcommand{\Pic}{\rm Pic}

\newcommand{\Hyper}{{\Bbb H}}
\newcommand{\Wedge}[1]{\stackrel{#1}{\wedge}}
\newcommand{\LieAlg}[1]{{\bf #1}}
\newcommand{\LoopGroup}{G_{\infty}}
\newcommand{\LoopAlg}{\LieAlg{g}_{\infty}}
\newcommand{\LevelInfinityGroup}{G^{+}_{\infty}}
\newcommand{\LevelInfinityAlg}{\LieAlg{g}^{+}_{\infty}}
\newcommand{\LevelInfinitySubgroup}[1]{G^{>#1}_{\infty}}
\newcommand{\LevelInfinitySubalg}[1]{\LieAlg{g}^{>#1}_{\infty}}
\newcommand{\LoopAlgSubtorus}[1]{\LieAlg{t}^{>#1}_{\infty}}

\newcommand{\HeisN}{{\rm Heis}_{\underline{n}}}
\newcommand{\heisN}{{\bf heis}_{\underline{n}}}
\newcommand{\HiggsModuli}{{\rm Higgs}}
\newcommand{\CoverHiggsModuli}{\widetilde{{\rm Higgs}}}
\newcommand{\Higgsm}{{\rm Higgsm}}
\newcommand{\HiggsmN}{{\rm Higgsm}^{\underline{n}}}
\newcommand{\CoverHiggsmN}{\widetilde{{\rm Higgsm}}^{\underline{n}}}
\newcommand{\CoverHiggsm}{\widetilde{{\rm Higgsm}}}
\newcommand{\BsmN}{Bsm^{\underline{n}}}
\newcommand{\CoverBsmN}{\widetilde{Bsm}^{\underline{n}}}
\newcommand{\HiggsComponent}{\HiggsModuli^{sm}}
\newcommand{\ModuliVB}{{\cal U}}
\newcommand{\ModuliVBLevels}{{\cal U}}
\newcommand{\Dual}{\vee}
\newcommand{\SkipAfterTitle}
{
\hspace{1ex}

\noindent
}
 \makeatletter
\expandafter\ifx\csname amssym.def\endcsname\relax \else\endinput\fi
\expandafter\edef\csname amssym.def\endcsname{%
       \catcode`\noexpand\@=\the\catcode`\@\space}
\catcode`\@=11
\def\undefine#1{\let#1\undefined}
\def\newsymbol#1#2#3#4#5{\let\next@\relax
 \ifnum#2=\@ne\let\next@\msafam@\else
 \ifnum#2=\tw@\let\next@\msbfam@\fi\fi
 \mathchardef#1="#3\next@#4#5}
\def\mathhexbox@#1#2#3{\relax
 \ifmmode\mathpalette{}{\m@th\mathchar"#1#2#3}%
 \else\leavevmode\hbox{$\m@th\mathchar"#1#2#3$}\fi}
\def\hexnumber@#1{\ifcase#1 0\or 1\or 2\or 3\or 4\or 5\or 6\or 7\or 8\or
 9\or A\or B\or C\or D\or E\or F\fi}

\font\tenmsa=msam10
\font\sevenmsa=msam7
\font\fivemsa=msam5
\newfam\msafam
\textfont\msafam=\tenmsa
\scriptfont\msafam=\sevenmsa
\scriptscriptfont\msafam=\fivemsa
\edef\msafam@{\hexnumber@\msafam}
\mathchardef\dabar@"0\msafam@39
\def\dashrightarrow{\mathrel{\dabar@\dabar@\mathchar"0\msafam@4B}}
\def\dashleftarrow{\mathrel{\mathchar"0\msafam@4C\dabar@\dabar@}}

\def\ulcorner{\delimiter"4\msafam@70\msafam@70 }
\def\urcorner{\delimiter"5\msafam@71\msafam@71 }
\def\llcorner{\delimiter"4\msafam@78\msafam@78 }
\def\lrcorner{\delimiter"5\msafam@79\msafam@79 }
\def\yen{{\mathhexbox@\msafam@55 }}
\def\checkmark{{\mathhexbox@\msafam@58 }}
\def\circledR{{\mathhexbox@\msafam@72 }}
\def\maltese{{\mathhexbox@\msafam@7A }}

\font\tenmsb=msbm10
\font\sevenmsb=msbm7
\font\fivemsb=msbm5
\newfam\msbfam
\textfont\msbfam=\tenmsb
\scriptfont\msbfam=\sevenmsb
\scriptscriptfont\msbfam=\fivemsb
\edef\msbfam@{\hexnumber@\msbfam}
\def\Bbb#1{{\fam\msbfam\relax#1}}
\def\widehat#1{\setbox\z@\hbox{$\m@th#1$}%
 \ifdim\wd\z@>\tw@ em\mathaccent"0\msbfam@5B{#1}%
 \else\mathaccent"0362{#1}\fi}
\def\widetilde#1{\setbox\z@\hbox{$\m@th#1$}%
 \ifdim\wd\z@>\tw@ em\mathaccent"0\msbfam@5D{#1}%
 \else\mathaccent"0365{#1}\fi}
\font\teneufm=eufm10
\font\seveneufm=eufm7
\font\fiveeufm=eufm5
\newfam\eufmfam
\textfont\eufmfam=\teneufm
\scriptfont\eufmfam=\seveneufm
\scriptscriptfont\eufmfam=\fiveeufm
\def\frak#1{{\fam\eufmfam\relax#1}}
\let\goth\frak

\csname amssym.def\endcsname
\makeatother
 \makeatletter
\expandafter\ifx\csname pre amssym.tex at\endcsname\relax \else \endinput\fi
\expandafter\chardef\csname pre amssym.tex at\endcsname=\the\catcode`\@
\catcode`\@=11
\newsymbol\boxdot 1200
\newsymbol\boxplus 1201
\newsymbol\boxtimes 1202
\newsymbol\square 1003
\newsymbol\blacksquare 1004
\newsymbol\centerdot 1205
\newsymbol\lozenge 1006
\newsymbol\blacklozenge 1007
\newsymbol\circlearrowright 1308
\newsymbol\circlearrowleft 1309
\undefine\rightleftharpoons
\newsymbol\rightleftharpoons 130A
\newsymbol\leftrightharpoons 130B
\newsymbol\boxminus 120C
\newsymbol\Vdash 130D
\newsymbol\Vvdash 130E
\newsymbol\vDash 130F
\newsymbol\twoheadrightarrow 1310
\newsymbol\twoheadleftarrow 1311
\newsymbol\leftleftarrows 1312
\newsymbol\rightrightarrows 1313
\newsymbol\upuparrows 1314
\newsymbol\downdownarrows 1315
\newsymbol\upharpoonright 1316
 
\newsymbol\downharpoonright 1317
\newsymbol\upharpoonleft 1318
\newsymbol\downharpoonleft 1319
\newsymbol\rightarrowtail 131A
\newsymbol\leftarrowtail 131B
\newsymbol\leftrightarrows 131C
\newsymbol\rightleftarrows 131D
\newsymbol\Lsh 131E
\newsymbol\Rsh 131F
\newsymbol\rightsquigarrow 1320
\newsymbol\leftrightsquigarrow 1321
\newsymbol\looparrowleft 1322
\newsymbol\looparrowright 1323
\newsymbol\circeq 1324
\newsymbol\succsim 1325
\newsymbol\gtrsim 1326
\newsymbol\gtrapprox 1327
\newsymbol\multimap 1328
\newsymbol\therefore 1329
\newsymbol\because 132A
\newsymbol\doteqdot 132B
 
\newsymbol\triangleq 132C
\newsymbol\precsim 132D
\newsymbol\lesssim 132E
\newsymbol\lessapprox 132F
\newsymbol\eqslantless 1330
\newsymbol\eqslantgtr 1331
\newsymbol\curlyeqprec 1332
\newsymbol\curlyeqsucc 1333
\newsymbol\preccurlyeq 1334
\newsymbol\leqq 1335
\newsymbol\leqslant 1336
\newsymbol\lessgtr 1337
\newsymbol\backprime 1038
\newsymbol\risingdotseq 133A
\newsymbol\fallingdotseq 133B
\newsymbol\succcurlyeq 133C
\newsymbol\geqq 133D
\newsymbol\geqslant 133E
\newsymbol\gtrless 133F
\newsymbol\sqsubset 1340
\newsymbol\sqsupset 1341
\newsymbol\vartriangleright 1342
\newsymbol\vartriangleleft 1343
\newsymbol\trianglerighteq 1344
\newsymbol\trianglelefteq 1345
\newsymbol\bigstar 1046
\newsymbol\between 1347
\newsymbol\blacktriangledown 1048
\newsymbol\blacktriangleright 1349
\newsymbol\blacktriangleleft 134A
\newsymbol\vartriangle 134D
\newsymbol\blacktriangle 104E
\newsymbol\triangledown 104F
\newsymbol\eqcirc 1350
\newsymbol\lesseqgtr 1351
\newsymbol\gtreqless 1352
\newsymbol\lesseqqgtr 1353
\newsymbol\gtreqqless 1354
\newsymbol\Rrightarrow 1356
\newsymbol\Lleftarrow 1357
\newsymbol\veebar 1259
\newsymbol\barwedge 125A
\newsymbol\doublebarwedge 125B
\undefine\angle
\newsymbol\angle 105C
\newsymbol\measuredangle 105D
\newsymbol\sphericalangle 105E
\newsymbol\varpropto 135F
\newsymbol\smallsmile 1360
\newsymbol\smallfrown 1361
\newsymbol\Subset 1362
\newsymbol\Supset 1363
\newsymbol\Cup 1264
 
\newsymbol\Cap 1265
 
\newsymbol\curlywedge 1266
\newsymbol\curlyvee 1267
\newsymbol\leftthreetimes 1268
\newsymbol\rightthreetimes 1269
\newsymbol\subseteqq 136A
\newsymbol\supseteqq 136B
\newsymbol\bumpeq 136C
\newsymbol\Bumpeq 136D
\newsymbol\lll 136E
 
\newsymbol\ggg 136F
 
\newsymbol\circledS 1073
\newsymbol\pitchfork 1374
\newsymbol\dotplus 1275
\newsymbol\backsim 1376
\newsymbol\backsimeq 1377
\newsymbol\complement 107B
\newsymbol\intercal 127C
\newsymbol\circledcirc 127D
\newsymbol\circledast 127E
\newsymbol\circleddash 127F
\newsymbol\lvertneqq 2300
\newsymbol\gvertneqq 2301
\newsymbol\nleq 2302
\newsymbol\ngeq 2303
\newsymbol\nless 2304
\newsymbol\ngtr 2305
\newsymbol\nprec 2306
\newsymbol\nsucc 2307
\newsymbol\lneqq 2308
\newsymbol\gneqq 2309
\newsymbol\nleqslant 230A
\newsymbol\ngeqslant 230B
\newsymbol\lneq 230C
\newsymbol\gneq 230D
\newsymbol\npreceq 230E
\newsymbol\nsucceq 230F
\newsymbol\precnsim 2310
\newsymbol\succnsim 2311
\newsymbol\lnsim 2312
\newsymbol\gnsim 2313
\newsymbol\nleqq 2314
\newsymbol\ngeqq 2315
\newsymbol\precneqq 2316
\newsymbol\succneqq 2317
\newsymbol\precnapprox 2318
\newsymbol\succnapprox 2319
\newsymbol\lnapprox 231A
\newsymbol\gnapprox 231B
\newsymbol\nsim 231C
\newsymbol\ncong 231D
\newsymbol\diagup 231E
\newsymbol\diagdown 231F
\newsymbol\varsubsetneq 2320
\newsymbol\varsupsetneq 2321
\newsymbol\nsubseteqq 2322
\newsymbol\nsupseteqq 2323
\newsymbol\subsetneqq 2324
\newsymbol\supsetneqq 2325
\newsymbol\varsubsetneqq 2326
\newsymbol\varsupsetneqq 2327
\newsymbol\subsetneq 2328
\newsymbol\supsetneq 2329
\newsymbol\nsubseteq 232A
\newsymbol\nsupseteq 232B
\newsymbol\nparallel 232C
\newsymbol\nmid 232D
\newsymbol\nshortmid 232E
\newsymbol\nshortparallel 232F
\newsymbol\nvdash 2330
\newsymbol\nVdash 2331
\newsymbol\nvDash 2332
\newsymbol\nVDash 2333
\newsymbol\ntrianglerighteq 2334
\newsymbol\ntrianglelefteq 2335
\newsymbol\ntriangleleft 2336
\newsymbol\ntriangleright 2337
\newsymbol\nleftarrow 2338
\newsymbol\nrightarrow 2339
\newsymbol\nLeftarrow 233A
\newsymbol\nRightarrow 233B
\newsymbol\nLeftrightarrow 233C
\newsymbol\nleftrightarrow 233D
\newsymbol\divideontimes 223E
\newsymbol\varnothing 203F
\newsymbol\nexists 2040
\newsymbol\Finv 2060
\newsymbol\Game 2061
\newsymbol\mho 2066
\newsymbol\eth 2067
\newsymbol\eqsim 2368
\newsymbol\beth 2069
\newsymbol\gimel 206A
\newsymbol\daleth 206B
\newsymbol\lessdot 236C
\newsymbol\gtrdot 236D
\newsymbol\ltimes 226E
\newsymbol\rtimes 226F
\newsymbol\shortmid 2370
\newsymbol\shortparallel 2371
\newsymbol\smallsetminus 2272
\newsymbol\thicksim 2373
\newsymbol\thickapprox 2374
\newsymbol\approxeq 2375
\newsymbol\succapprox 2376
\newsymbol\precapprox 2377
\newsymbol\curvearrowleft 2378
\newsymbol\curvearrowright 2379
\newsymbol\digamma 207A
\newsymbol\varkappa 207B
\newsymbol\Bbbk 207C
\newsymbol\hslash 207D
\undefine\hbar
\newsymbol\hbar 207E
\newsymbol\backepsilon 237F
\catcode`\@=\csname pre amssym.tex at\endcsname
\makeatother

\font\teneurm=eurm10
\font\seveneurm=eurm7
\font\fiveeurm=eurm5
\newfam\eurmfam
\textfont\eurmfam=\teneurm
\scriptfont\eurmfam=\seveneurm
\scriptscriptfont\eurmfam=\fiveeurm
\def\eurm#1{{\fam\eurmfam\relax#1}}
\newcommand{\bdl}[1]{{\eurm #1}}

\newtheorem{thm}[theorem]{Theorem}
\newtheorem{prop}[theorem]{Proposition}
\newtheorem{eg}[theorem]{Example}
\newtheorem{defn}[theorem]{Definition}

\newcommand{\UnderlinedEnd}[1]{\underline{End}( #1)}

\begin{document}
\begin{Large}
\centerline{\bf Spectral covers, }
\smallskip
\centerline{\bf algebraically completely integrable,}
\smallskip
\centerline{\bf Hamiltonian systems,}
\smallskip
\centerline{\bf and moduli of bundles}
\end{Large}

\vspace{0.3in}

\centerline{Ron Donagi
\footnote{
Partially supported by  NSA
Grant MDA904-92-H3047 and NSF Grant DMS 95-03249
}}
\vspace{0.1in}
\centerline{University of Pennsylvania}
\vspace{0.2in}
\centerline{and}
\vspace{0.2in}
\centerline{Eyal Markman
\footnote{Partially supported by a Rackham
Fellowship, University of Michigan, 1993}
}
\vspace{0.1in}
\centerline{University of Michigan}

\tableofcontents
\newpage

\section{Introduction} \label{ch1}

     The purpose of these notes is to present an
algebro-geometric point of view on several interrelated topics,
all involving integrable systems in symplectic-algebro-geometric
settings.  These systems range from some very old examples, such
as the geodesic flow on an ellipsoid, through the classical
hierarchies of $KP-$ and $KdV$-types, to some new systems which
are often based on moduli problems in algebraic geometry.

     The interplay between algebraic geometry and integrable
systems goes back quite a way.  It has been known at least since
Jacobi that many integrable systems can be solved explicitly in
terms of {\it theta functions}.  (There are numerous examples,
starting with various {\it spinning tops} and the {\it geodesic
flow on an ellipsoid}.)  Geometrically, this often means that the
system can be mapped to the total space of a family of Jacobians
of some curves, in such a way that the flows of the system are
mapped to linear flows along the Jacobians.  In practice, these
curves tend to arise as the spectrum (hence the name {\it
`spectral'} curves) of some parameter-dependent operator; they can
therefore be represented as branched covers of the parameter
space, which in early examples tended to be the Riemann sphere
${\bf CP}^1$.

     In {\it Hitchin's system}, the base ${\bf CP}^1$ is replaced
by an arbitrary (compact, non-singular) Riemann surface $\Sigma$.
The cotangent bundle $T^*{\cal U}_\Sigma$ to the moduli space
${\cal U}_\Sigma$ of stable vector bundles on $\Sigma$ admits two
very different interpretations:\ on the one hand, it parametrizes
certain {\it Higgs bundles}, or vector bundles with a
(canonically) twisted endomorphism; on the other, it parametrizes
certain {\it spectral data}, consisting of torsion-free sheaves
(generically, line bundles) on spectral curves which are branched
covers of $\Sigma$.  In our three central chapters
(\ref{ch4},\ref{ch5},\ref{ch6}) we
study this important system, its extensions and variants.  All
these systems are linearized on Jacobians of spectral curves.

     We also study some systems in which the spectral curve is
replaced by a higher-dimensional geometric object: \ a {\it
spectral variety} in Chapter \ref{ch9}, an algebraic {\it Lagrangian
subvariety} in Chapter \ref{ch8}, and a {\it Calabi-Yau manifold} in
Chapter \ref{ch7}.  Our understanding of some of these wild systems is
much less complete than in the case of the curve-based ones.  We
try to explain what we know and to point out some of what we do
not.  The Calabi-Yau systems seem particularly intriguing.  Not
only are the tori (on which these systems are linearized) not
Jacobians of curves, they are in general not even abelian
varieties.  There are some suggestive relations between these
systems and the conjectural mirror-symmetry for Calabi-Yaus.

     The first three chapters are introductory.  In Chapter \ref{ch2} we
collect the basic notions of {\it symplectic geometry} and {\it
integrable systems} which will be needed, including some
information about {\it symplectic reduction}. (An excellent further
reference is \cite{AG}.)  In Chapter \ref{ch3} we work out in some detail the
classical theory of geodesic flow on an ellipsoid, which is
integrable via hyperelliptic theta functions.  We think of this
both as a beautiful elementary and explicit example and as an
important special case of the much more powerful results which
follow.  (Our presentation follows \cite{knorrer,reid,donagi-group-law}).
Some of our main algebro-geometric objects of study are
introduced in Chapter \ref{ch4}:\ vector bundles and their moduli
spaces, spectral curves, and the {\it `spectral systems'}
constructed from them.  In particular, we consider the {\it
polynomial matrix system} \cite{AHH,B} (which contains the
geodesic flow on an ellipsoid as special case) and {\it Hitchin's
system} \cite{hitchin,hitchin-integrable-system}.

     Each of the remaining five chapters presents in some detail
a recent or current research topic.  Chapter \ref{ch5} outlines
constructions (from \cite{markman-higgs,botachin,tyurin-symplectic})
of the Poisson structure on
the spectral system of curves.  This is possible whenever the
twisting line bundle $K$ is a non-negative twist $\omega_\Sigma(D)$ of
the canonical bundle $\omega_\Sigma$, and produces an algebraically
completely integrable Hamiltonian system.  Following
\cite{markman-higgs} we emphasize the
deformation-theoretic construction, in which the Poisson
structure on an open subset of the system is obtained via
symplectic reduction from the cotangent bundle
$T^*{\cal U}_{\Sigma ,D}$ of the moduli
space ${\cal U}_{\Sigma ,D}$ of stable bundles with a {\it level-D structure}.

     In Chapter \ref{ch6} we explore the relation between these spectral
systems and the $KP$-hierarchy and its variants (multi-component
$KP$, Heisenberg flows,  and their $KdV$-type subhierarchies).  These
hierarchies are, of course, a rich source of geometry:\ The
Krichever construction (e.g. \cite{segal-wilson-loop-groups-and-kp})
shows that any Jacobian can be
embedded in $KP$-space, and these are the only finite-dimensional
orbits \cite{mulase-cohomological-structure, AdC, Sh}.
Following \cite{adams-bergvelt,li-mulase-category}
we describe some
``multi-Krichever'' constructions which take spectral data to the
spaces of the $KP$, $mcKP$ and Heisenberg systems.  Our
main new result is that the flows on the spectral system which
are obtained by pulling back the $mcKP$ or Heisenberg
flows via the corresponding Krichever maps are {\it Hamiltonian}
with respect to the Poisson structure constructed in Chapter \ref{ch5}.
In fact, we write down explicitly the Hamiltonians for these $KP$
flows on the spectral system, as residues of traces of
meromorphic matrices.  (Some related results have also been
obtained recently in \cite{li-mulase-compatibility}.)

     The starting point for Chapter \ref{ch7} is an attempt to
understand the condition for a given family of complex tori to
admit a symplectic structure and thus become an ACIHS.  We find
that the condition is a symmetry on the derivatives of the period
map, which essentially says that the periods are obtained as
partials of some field of symmetric cubic tensors on the base.
In the rest of this Chapter we apply this idea to an analytically
(not algebraically) integrable system constructed from any family
of Calabi-Yau $3$-folds. Some properties of this system suggest that it may
be relevant to a purely hodge-theoretic reformulation of the
mirror-symmetry conjectures.
(This chapter is based on \cite{cubics-calabi-yaus}.)

     Chapter \ref{ch8} is devoted to the construction of symplectic
and Poisson structures in some inherently non-linear
situations, vastly extending the results of Chapter \ref{ch5}.  The basic
space considered here is the moduli space parametrizing
line-bundle-like sheaves supported on (variable) subvarieties of a
given symplectic space $X$.  It is shown that when the
subvarieties are Lagrangian, the moduli space itself becomes
symplectic.  The spectral systems considered in Chapter \ref{ch5} can be
recovered as the case where $X$ is the total space of $T^*\Sigma$
and the Lagrangian subvarieties are the spectral curves.
(A fuller version of these results will appear in
\cite{markman-lagrangian-sheaves}.)

     In the final chapter we consider extensions of the spectral
system to allow a higher-dimensional base variety $S$, an arbitrary
reductive group $G$, an arbitrary representation $\rho: G \to Aut
V$, and values in an arbitrary vector bundle  $K$. (Arbitrary reductive
groups  $G$ were considered, over  a curve $S = \Sigma$ with
$K = \omega_\Sigma$, by Hitchin \cite{hitchin-integrable-system},
while the case
$K = \Omega_S$  over arbitrary base $S$  is Simpson's
\cite{simpson-moduli}).  We replace
spectral curves by various kinds of spectral covers, and introduce the
cameral cover, a version of the Galois-closure of a spectral cover which is
independent of $K$ and $\rho$.
It comes with an action of $W$, the Weyl group of $G$.
We analyze the decomposition, under the
action of $W$, of the cameral and spectral Picard varieties, and
identify the distinguished Prym in there.  This is shown to
correspond, up to certain shifts and twists, to the fiber of the
Hitchin map in this general setting, i.e. to moduli of Higgs
bundles with a  given $\widetilde{S}$. Combining this with
some obvious remarks about  existence of Poisson structures,
we find that  the moduli spaces of  K-valued Higgs bundles support
algebraically completely integrable systems. Our presentation closely follows
that of \cite{MSRI}

     It is a pleasure to express our gratitude to the organizers,
Mauro Francaviglia and Silvio Greco, for the opportunity to participate in the
CIME meeting and to publish
these notes here.
During the preparation of this long work we benefited from many enjoyable
conversations with M. Adams, M. Adler, A. Beauville, R. Bryant, C. L. Chai,
I. Dolgachev, L. Ein, B. van Geemen, A. Givental, M. Green, P. Griffiths,
N. Hitchin, Y. Hu, S. Katz, V. Kanev, L. Katzarkov, R. Lazarsfeld, P. van
Moerbeke,
D. Morrison,  T. Pantev,  E. Previato and E. Witten.

\newpage


\section{Basic Notions} \label{ch2}
\label{sec-basic-notions}

   We gather here those basic concepts and elementary results from symplectic
and Poisson geometry, completely integrable systems, and symplectic reduction
which will be helpful throughout these notes.  Included are a few useful
examples and only occasional proofs or sketches.  To the reader
unfamiliar with this material we were hoping to impart just as much of a
feeling for it as might be needed in the following chapters.  For more
details, we recommend the excellent survey \cite{AG}.

\subsection{Symplectic Geometry} \label{subsec-symplectic-geometry}

\noindent
\underline{{\bf Symplectic structure}}

A symplectic structure on a differentiable manifold $M$ of even dimension
$2n$ is given by a non-degenerate closed 2-form  $\sigma$.   The non
degeneracy means that either of the following equivalent conditions holds.
\begin{itemize}
\item $\sigma^n$ is a nowhere vanishing volume form.
\item Contraction with $\sigma$ induces an isomorphism $\rfloor \sigma :
TM \rightarrow T^*M$
\item For any non-zero tangent vector $v \in T_mM$ at $m \in M$, there is
some $v' \in T_mM$  such that $\sigma(v,v') \ne 0$.
\end{itemize}

\begin{examples} \label{examples-symplectic-varieties}
{\rm
\SkipAfterTitle
\begin{enumerate}
\item
\underline{Euclidean space}

The standard example of a symplectic manifold is Euclidean space
$\R^{2n}$ with $\sigma = \Sigma dp_i \wedge dq_i$, where $p_1,\cdots,p_n,
\ q_1, \cdots,q_n$ are linear coordinates.  Darboux's theorem says that any
symplectic manifold is locally equivalent to this example (or to any other).

\item
\underline{Cotangent bundles}

For any manifold $X$, the cotangent bundle $M := T^*X$ has a natural
symplectic structure.  First, $M$ has the tautological 1-form $\alpha$,
whose value at $(x,\theta) \in T^*X$ is $\theta$ pulled back to $T^*M$.
If $q_1 ,\cdots, q_n$ are local coordinates on $X$, then locally $\alpha =
\Sigma p_i dq_i$ where the $p_i$ are the fiber coordinates given by
$\partial / \partial q_i$.  The differential
$$
\sigma := d\alpha
$$
is then a globally defined closed (even exact) 2-form on $M$.  It is
given in local coordinates by $\Sigma dp_i \wedge dq_i$, hence is
non-degenerate.

\item
\underline{Coadjoint orbits}

Any Lie group $G$ acts on its Lie algebra $\LieAlg{g}$ (adjoint representation)
and hence on the dual vector space $\LieAlg{g}^*$ (coadjoint
representation).  Kostant and Kirillov noted that
for any $\xi \in \LieAlg{g}^*$, the coadjoint orbit
${\cal O}
= G \xi \subset \LieAlg{g}^*$ has
a natural symplectic structure.  The tangent space to $\cal O$ at $\xi$
is given by $\LieAlg{g}/\LieAlg{g}_\xi$, where $\LieAlg{g}_\xi$ is the
stabilizer of $\xi$:
$$
\LieAlg{g}_\xi :=
\{ x \in \LieAlg{g} \; |\; ad^*_x \xi = 0 \} =
\{x \in \LieAlg{g} \; | \; (\xi,[x,y]) = 0 \quad \forall \  y \in \LieAlg{g}
\}.
$$
Now $\xi$ determines an alternating bilinear form on $\LieAlg{g}$
$$
x,y \longmapsto (\xi, [x,y]),
$$
which clearly descends to $\LieAlg{g} / \LieAlg{g}_\xi$  and is non-degenerate
there.  Varying $\xi$ we get a non-degenerate 2-form $\sigma$ on $\cal O$.
The Jacobi identity on $g$ translates immediately into closedness of
$\sigma$.
\end{enumerate}
}
\end{examples}

\medskip
\noindent
\underline{{\bf Hamiltonians}}

To a function $f$ on a symplectic manifold $(M,\sigma)$ we associate its
{\it Hamiltonian vector field} $v_f$, uniquely determined by

$$
v_f \; \rfloor \; \sigma = df.
$$

A vector field $v$ on $M$ is Hamiltonian if and only if the 1-form $v \,
\rfloor \, \sigma$ is exact.  We say $v$ is {\it locally Hamiltonian}
if $v \; \rfloor \; \sigma$ is closed.  This is equivalent to saying that the
flow generated by $v$ preserves $\sigma$.  Thus on a symplectic surface
$(n=1)$, the locally Hamiltonian vector fields are the area-preserving
ones.

\medskip
\noindent
{\bf Example:} (Geodesic flow)

A Riemannian metric on a manifold $X$ determines an isomorphism of $M :=
TX$ with $T^*X$; hence we get on $M$ a natural symplectic structure
together with a $C^\infty$ function $f =$ (squared length).  The geodesic
flow on $X$ is the differential equation, on $M$, given by the
Hamiltonian vector field $v_f$.  Its integral curves are the geodesics on
$M$.

\medskip
\noindent
\underline{{\bf Poisson structures}}

The association $f \mapsto v_f$ gives a map of sheaves
\begin{equation} \label{eq-functions-to-hamiltonian-vectorfields}
v : C^\infty (M) \longrightarrow V(M)
\end{equation}
from $C^\infty$ functions on the symplectic manifold $M$ to vector
fields.  Now  $V(M)$ always has the structure of a Lie algebra, under
commutation of vector fields.  The symplectic structure on $M$ determines
a Lie algebra structure on $C^\infty(M)$ such that $v$ becomes a morphism
of (sheaves of) Lie algebras.  The operation on $C^\infty(M),$ called {\it
Poisson bracket}, is
$$
\{ f,g \} := (df, v_g) = -(dg, v_f)  = {{n df \wedge dg \wedge
\sigma^{n-1}} \over {\sigma^n}}.
$$
More generally, a {\it Poisson structure} on a manifold $M$ is a
Lie algebra bracket $\{\, ,\, \}$ on $C^\infty(M)$ which acts as a derivation
in each variable:
$$
\{f,gh\} = \{f,g\} h + \{f,h\}g, \ \ \  f,g,h \in C^\infty(M).
$$
Since the value at a point $m$ of a given derivation acting on a function
$g$ is a linear function of $d_mg$, we see that a Poisson structure on
$M$ determines a global 2-vector
$$
\psi \in H^0(M, \stackrel{2}{\wedge} TM).
$$
or equivalently a skew-symmetric homomorphism
$$
\Psi : T^*M \longrightarrow TM.
$$

Conversely, any 2-vector $\psi$ on $M$ determines an alternating bilinear
bracket on $C^\infty (M)$, by

$$
\{f,g\} := (df \wedge dg, \psi),
$$
and this acts as a derivation in each variable.  An equivalent way of
specifying a Poisson structure is thus to give a global 2-vector $\psi$
satisfying an integrability condition (saying that the above bracket
satisfies the Jacobi identity, hence gives a Lie algebra).

We saw that a symplectic structure $\sigma$ determines a Poisson bracket
$\{\ ,\ \}$.  The corresponding homomorphism $\Psi$ is just $(\rfloor
\sigma)^{-1}$; the closedness of $\sigma$ is equivalent to integrability
of $\psi$.  Thus, a Poisson structure which is (i.e. whose 2-vector is)
everywhere non-degenerate, comes from a symplectic structure.

A general Poisson structure can be degenerate in two ways:  first, there
may exist non-constant functions $f \in C^\infty(M)$, called {\it
Casimirs}, satisfying
$$
0 = df \rfloor \psi = \Psi(df),
$$
i.e.
$$
\{f,g\} = 0 \  \mbox{for all} \  g \in C^\infty(M).
$$
This implies that the rank of $\Psi$ is less than maximal everywhere.  In
addition, or instead, rank $\Psi$ could drop along some strata in $M$.
For even $r$, let
$$
M_r := \{ m \in M | rank(\Psi) = r \}.
$$
Then a basic result \cite{We} asserts that the $M_r$ are submanifolds, and
they are canonically foliated into {\it symplectic leaves}, i.e.
$r$-dimensional submanifolds $Z \subset M_r$ which inherit a symplectic
structure.  (This means that the restriction $\psi_{\mid_Z}$ is the image,
under the inclusion $ Z \hookrightarrow M_r$, of a two-vector $\psi_Z$
on $Z$ which is everywhere nondegenerate, hence comes from a symplectic
structure on $Z$.)  These leaves can be described in several ways:

\begin{itemize}
\item The image $\Psi(T^*M_r)$ is an involutive subbundle of rank $r$ in
$TM_r$;
the $Z$ are its integral leaves.
\item The leaf $Z$ through $m \in M_r$ is $Z = \{ z \in M_r | f(m) = f(z)\
\mbox{for all Casimirs} \  f  \    \mbox{on} \ M_r \}$.
\item
Say that two points of $M$ are $\psi$-connected if there is an integral
curve of some Hamiltonian vector field passing through both.  The leaves
are the equivalence classes for the equivalence relation generated by
$\psi$-connectedness.
\end{itemize}

\noindent
\begin{example}\label{example-coadjoint-orbits}

{\rm The Kostant-Kirillov symplectic structures on coadjoint orbits
of a Lie algebra $\LieAlg{g}$ extend to a Poisson structure
on the dual vector space $\LieAlg{g}^*$.
For a function $F \in C^\infty(\LieAlg{g}^*)$
we identify its differential $d_\xi F$ at $\xi \in
\LieAlg{g}*$ with an element of $\LieAlg{g} = \LieAlg{g}^{**}$.  We then set:
$$
\{F,G\}(\xi) := (\xi, [d_\xi F,d_\xi G]).
$$
This is a Poisson structure, whose symplectic leaves are precisely the
coadjoint orbits.  The rank of $\LieAlg{g}$ is, by definition, the smallest
codimension $\ell$ of a coadjoint orbit.  The Casimirs are the ad-invariant
functions on $\LieAlg{g}^*$.  Their restrictions to the largest stratum
$\LieAlg{g}^*_{\dim \LieAlg{g} - \ell}$ foliate this stratum, the leaves being
the {\it regular} (i.e. largest dimensional) coadjoint orbits.}
\end{example}

\subsection{Integrable Systems}

We say that two functions $h_1,h_2$ on a Poisson manifold $(M,\psi)$
{\em Poisson commute} if their Poisson
bracket $\{ h_1,h_2 \}$ is zero.  In this case the integral flow of the
Hamiltonian vector field of each function $h_i, \  i = 1,2$ is tangent to
the level sets of the other.  In other words, $h_2$ is a conservation law
for the Hamiltonian $h_1$ and the Hamiltonian flow of $h_2$ is a symmetry
of the Hamiltonian system associated with $(M, \psi, h_1)$
(the flow of the Hamiltonian vector field $v_{h_{1}}$  on $M$).

A map $f : M \rightarrow B$ between two Poisson manifolds is a {\it
Poisson map} if pullback of functions is a Lie algebra homomorphism with
respect to the Poisson bracket
$$
f^*\{F,G\}_B \, = \, \{f^*F,f^*G\}_M.
$$
Equivalently, if $df(\psi_M)$ equals $f^*(\psi_B)$ as sections of
$f^*(\stackrel{2}{\wedge} T_B)$.  If $H:M \rightarrow B$ is a Poisson map
with respect to the trivial (zero) Poisson structure on $B$ we will call
$H$ a {\em Hamiltonian map}.  Equivalently, $H$ is  Hamiltonian if the
Poisson structure $\psi$ vanishes on the
pullback $H^*(T^*B)$ of the cotangent bundle of $B$
(regarding the latter as a
subbundle of $(T^*M,\psi)$). In particular, the rank of the differential
$dH$ is less than or equal to $\dim M - {1 \over 2} \rank (\psi)$ at
every point.  A Hamiltonian map pulls back the algebra of functions on
$B$ to a commutative Poisson subalgebra of the algebra of functions on $M$.

The study of a Hamiltonian system $(M,\psi,h)$ simplifies tremendously if
one can extend the Hamiltonian function $h$ to a Hamiltonian map $H : M
\rightarrow B$ of maximal rank $\dim M - {1 \over 2} \rank(\psi)$.
Such a system is called a completely integrable Hamiltonian system.  The
Hamiltonian flow of a completely integrable system can often be realized
as a linear flow on tori embedded in $M$.  The fundamental theorem in
this case is Liouville's theorem (stated below).

\noindent
\begin{definition} {\rm
\begin{enumerate}
\item
\ Let $V$ be a vector space, $\sigma \in
\stackrel{2}{\wedge} V^*$ a (possibly degenerate) two form.  A
subspace $Z \subset V$ is called {\em isotropic (coisotropic)} if it is
contained in (contains) its symplectic complement.  Equivalently, $Z$ is
isotropic if $\sigma$ restricts to zero on $Z$.  If $\sigma$ is
nondegenerate,  a subspace $Z \subset V$ is called {\em Lagrangian} if it is
both isotropic and coisotropic.  In this case $V$ is even (say $2n$)
dimensional and the Lagrangian subspaces are the $n$ dimensional
isotropic subspaces.
\item
Let $(M,\sigma)$ be a symplectic manifold.  A submanifold $Z$ is
{\em isotropic} (respectively {\em coisotropic, Lagrangian})
if the tangent subspaces
$T_zZ$ are, for all $z \in Z$.
\end{enumerate}
}
\end{definition}

\begin{example} {\rm
For every manifold $X$,
the fibers of the cotangent bundle $T^*X$ over points of $X$
are Lagrangian submanifolds with respect to the
standard symplectic structure. A section of $T^*X$ over  $X$
is Lagrangian if and only if the corresponding $1$-form on $X$ is closed.
}
\end{example}

We will extend the above definition to Poisson geometry:

\noindent
\begin{definition} {\rm
\begin{enumerate}
\item
 Let $U$ be a vector space, $\psi$ an element of
$\stackrel{2}{\wedge}U$.  Let $V \subset U$ be the image of the
contraction $\rfloor \; \psi : U^* \rightarrow U$.  Let $W \subset U^*$
be its kernel.  $W$ is called the null space of $\psi$.  $\psi$ is in
fact a nondegenerate element of $\stackrel{2}{\wedge} V$ giving rise to a
symplectic form $\sigma \in \stackrel{2}{\wedge} V^*$ (its inverse).  A
subspace $Z \subset U$ is {\em Lagrangian} with respect to $\psi$ if $Z$ is a
Lagrangian subspace of $V \subset U$ with respect to $\sigma$.
Equivalently, $Z$ is Lagrangian if $(U/Z)^*$ is both an isotropic and a
coisotropic subspace of $U^*$ with respect to $\psi \in
\stackrel{2}{\wedge} U \cong \stackrel{2}{\wedge} (U^*)^*$.
\item
Let $(M,\psi)$ be a Poisson manifold, assume that $\psi$ has constant
rank (this condition will be relaxed in the complex analytic or algebraic
case).  A submanifold $Z \subset M$ is {\em Lagrangian} if the tangent
subspaces $T_zZ$ are, for all $z \in Z$.  Notice that the constant rank
assumption implies that each connected component of $Z$ is contained in a
single symplectic leaf.
\end{enumerate}
}
\end{definition}

\medskip
\noindent
\begin{theorem} \ (Liouville). \ Let $M$ be an m-dimensional Poisson
manifold with Poisson structure $\psi$ of constant rank $2g$.  Suppose
that $H : M \rightarrow B$ is a proper submersive Hamiltonian map of
maximal rank, i.e, dim $B = m - g$.  Then

\begin{description}
\item[i)] \ The null foliation of $M$ is induced locally by a foliation of
$B$ (globally if $H$ has connected fibers).
\item[ii)] \ The connected components of fibers of $H$ are Lagrangian
compact tori with a natural affine structure.
\item[iii)] \ The Hamiltonian vector fields of the pullback of functions
on $B$ by $H$ are tangent to the level tori and are translation invariant
(linear).
\end{description}
\end{theorem}

\medskip
\noindent
\begin {rem}:
{\rm
If $H$ is not proper, but the Hamiltonian flows are complete,
then the fibers of $H$ are generalized tori (quotients of a vector space
by a discrete subgroup, not necessarily of maximal rank).
}
\end{rem}

\medskip
\noindent
\underline{Sketch of proof of Liouville's theorem:}

\begin{description}
\item[i)] Since $H$ is a proper submersion the connected components of the
fibers of $H$ are smooth compact submanifolds.  Since $H$ is a Hamiltonian map
of maximal rank $m - g$, the pullback $H^*(T^*_B)$ is isotropic and
coisotropic and hence $H$ is a Lagrangian fibration.  In particular, each
connected component of a fiber of $H$ is contained in a single symplectic
leaf.
\item[ii),iii)]  Let $A_b$ be a connected component of the fiber
$H^{-1}(b)$.  Let $0 \rightarrow T_{A_b} \rightarrow T_{{M|}_{A_b}}
\stackrel{dH}{\longrightarrow} (T_bB) \otimes {\cal O}_{A_b} \rightarrow 0$
be the exact sequence of the differential of $H$.  Part i) implies that the
null subbundle $W_{{|A}_b} := Ker [\Psi : T^*M \rightarrow TM]_{|A_b}$ is the
pullback of a subspace $W_b$ of $T^*_bB$.  Since $H$ is a Lagrangian
fibration, the Poisson structure induces a surjective homomorphism
$\phi_b : H^*(T^*_bB) \rightarrow T_{A_b}$ inducing a trivialization
$\bar{\phi}_b : (T^*_bB / W_b) \otimes {\cal O}_{A_b}
\stackrel{\sim}{\longrightarrow} \ T_{A_b}$.

A basis of the vector space $T^*_bB/W_b$ corresponds to a frame of global
independent vector fields on the fiber $A_b$ which commute since the map
$H$ is Hamiltonian.  Hence $A_b$ is a compact torus.
\end{description}
\EndProof
\subsection{Algebraically Completely Integrable Hamiltonian Systems}

All the definitions and most of the results stated in this chapter for
$C^\infty$-manifolds translate verbatim and hold in the complex analytic
and complex algebro-geometric categories replacing the real symplectic
form by a holomorphic or algebraic $(2,0)$-form (similarly for Poisson
structures).  The (main) exception listed below is due to the differences
between the Zariski topology and the complex or $C^\infty$ topologies.  A
Zariski open subset is the complement of the zero locus of a system of
polynomial equations.  It is hence always a dense open subset.

The (local) foliation by symplectic leaves exists only local analytically.
For example, a rank 2 translation invariant section $\psi \in H^0(A,
\stackrel{2}{\wedge} TA)$ on a 3 dimensional abelian variety $A$
which is simple
(does not contain any abelian subvariety) is an algebraic Poisson
structure with a non algebraic null foliation.

We will relax the definitions of a Lagrangian subvariety and integrable
system in the algebro-geometric category:

\begin{definition} {\rm
 Let $(M,\psi)$ be a Poisson smooth algebraic variety.
An irreducible and reduced subvariety $Z \subset M$ is {\em Lagrangian} if the
tangent subspace $T_zZ \subset T_zM$ is Lagrangian for a generic point
$z \in Z$.
}
\end{definition}

\begin{definition} {\rm
An {\em algebraically completely integrable Hamiltonian system} consists
of a proper flat morphism $H:M \rightarrow B$ where
$(M,\psi)$ is a smooth Poisson variety and $B$ is a smooth variety such that,
over the complement $B \smallsetminus \Delta$ of some proper closed subvariety
$\Delta \subset B$, $H$ is a Lagrangian fibration whose fibers are
isomorphic to abelian varieties.
}
\end{definition}

Multiples of a theta line bundle embed an abelian variety in projective
spaces with the coordinates being theta functions.  Thus, a priori, the
solutions of an algebraically completely integrable Hamiltonian system can
be expressed in terms of theta functions.  Finding explicit formulas is
usually hard.  In the next chapter we will study one example, the geodesic
flow on ellipsoids, in some detail. Later we will encounter
certain equations of $Kdv$ type, the Hitchin system, and a few
other examples.  Other classical integrable systems include various
Euler-Arnold systems, spinning tops, the Neumann system of evolution of a point
on the sphere subject to a quadratic potential.

Most of these systems are the complexification of real algebraic systems.
Given a real algebraic symplectic variety $(M,\sigma)$ and an algebraic
Hamiltonian $h$ on $M$ we say that the system is {\em algebraically completely
integrable} if its complexification $(M_{\C},\sigma_{\C}, \, h_{\C})$
is.  A real completely integrable system $(M,\sigma,h)$ need not be
algebraically completely integrable even if $(M,\sigma,h)$ are algebraic:

\medskip
\noindent
{\bf A Counter Example:} \ \ Let $(M,\sigma)$ be $(\R^2, \, dx \wedge
dy)$ and $h:\R^2 \rightarrow \R$ a polynomial of degree $d$ whose level sets
are nonsingular.  The system is
trivially completely integrable, but it is algebraically completely
integrable if and only if $d=3$ because in all other cases the
generic fiber of the complexification is a complex affine plane curve of
genus $\frac{(d-1)(d-2)}{2} \neq 1$.

\medskip
\noindent
\underline{{\bf Action Angle Coordinates}}:

Let $(M,\sigma)$ be a $2n$-dimensional symplectic manifold, $H:M
\rightarrow B$ a Lagrangian fibration by compact connected tori.

\noindent
\begin{theorem} (real action angle coordinates).

\noindent
In a neighborhood of a fiber of $H:M \rightarrow B$ one can introduce the
structure of a direct product $(\R^n / \Z^n) \times \R^n$ with action
coordinates $(I_1 \, \cdots \, I_n)$ on the factor $\R^n$ and angular
coordinates $(\phi_1, \, \cdots \, \phi_n)$ on the torus $(\R^n/\Z^n)$ in
which the symplectic structure has the form $\sum^n_{k=1} dI_k \wedge
d\phi_k$.
\end{theorem}

The Local action coordinates on $B$ are canonical up to affine
transformation on $\R^n$ with differential in $SL(n,\Z)$.  The angle
coordinates depend canonically on the action coordinates and a choice of a
Lagrangian section of $H:M \rightarrow B$.

\noindent
{\bf Remarks:}

\begin{itemize}
\begin{enumerate}
\item  In action angle coordinates the equations of the Hamiltonian flow of
a function $h$ on $B$ becomes: $\dot{I}_k = 0$, $\dot{\varphi}_k =
c_k(I_1,\cdots,I_n)$ where the slopes $c_k$ are $c_k = {{\partial h} \over
{\partial I_k}}$.
\item  In the polarized complex analytic case, we still have local
holomorphic action coordinates.  They depend further on a choice of a
Lagrangian subspace of the integral homology $H_1(A_b,\Z)$ with respect to
the polarization (a section of $\stackrel{2}{\wedge}H^1(A_b,\Z))$.
\end{enumerate}
\end{itemize}

\subsection{Moment Maps and Symplectic Reduction}
\label{subsec-moment-maps}

\noindent
\underline{{\bf Poisson Actions}}

An action $\rho$ of a connected Lie group $G$ on a manifold $M$ determines an
{\it infinitesimal action}
$$
d \rho : \LieAlg{g} \longrightarrow V(M),
$$
which is a homomorphism from the Lie algebra of $G$ to the Lie algebra of
$C^\infty$ vector fields on $M$.  When $(M,\sigma)$ is symplectic, we say
that the action $\rho$ is {\it symplectic} if
$$
(\rho(g))^* \sigma = \sigma, \quad \quad {\rm all} \quad g \in G,
$$
or equivalently if the image of $d\rho$ consists of locally Hamiltonian
vector fields.

We say that the action $\rho$ is {\it Poisson} if it factors through
the Lie algebra homomorphism (\ref{eq-functions-to-hamiltonian-vectorfields})
$v : C^\infty(M) \rightarrow V(M)$
and a Lie algebra homomorphism
$$
H:\LieAlg{g} \longrightarrow C^\infty(M).
$$
This imposes two requirements on $\rho$, each of a cohomological nature:  the
locally Hamiltonian fields $d\rho(X)$ should be globally Hamiltonian,
$d\rho(X) = v(H(X))$; and it must be possible to
choose the $H(X)$ consistently so that
$$
H([X,Y]) \ = \ \{H(X), H(Y)\}.
$$
(a priori the difference between the two terms is a constant function,
since its $v$ is zero, so the condition is that it should be possible to
make all these constants vanish simultaneously.)

\medskip
\noindent
\underline{{\bf Moment Maps}}

Instead of specifying the Hamiltonian lift
$$
H : \LieAlg{g} \longrightarrow C^\infty(M)
$$
for a Poisson action of $G$ on $(M,\sigma)$, it is convenient to consider
the equivalent data of the {\it moment map}
$$
\mu : M \longrightarrow \LieAlg{g}^*
$$
defined by
$$
(\mu(m),X) \, := \,H(X)(m).
$$
It is a Poisson map with respect to the Kostant-Kirillov Poisson structure on
$\LieAlg{g}$ (example \ref{example-coadjoint-orbits}), and is $G$-equivariant.

\begin{examples} \label{examples-moment-maps}
{\rm
\begin{enumerate}
\item
Any action of $G$ on a manifold $X$ lifts to an action on $M := T^*X$.
This action is Poisson.  The corresponding moment map $T^*X \rightarrow
\LieAlg{g}^*$ is the dual of the infinitesimal action $\LieAlg{g} \rightarrow
\Gamma (TX)$.  It can be identified with the pullback of differential forms
from $X$ to $G$ via the action.
\item
The coadjoint action of $G$ on $\LieAlg{g}^*$ is Poisson, with the identity as
moment map.
\end{enumerate}
}
\end{examples}

\medskip
\noindent
\underline{{\bf Symplectic Reduction}}

Consider a Poisson action of $G$ on $(M,\sigma)$ for which a reasonable
quotient $G/M$ exists.  (We will remain vague about this for now, and discuss
the properties of the quotient on a case-by-case basis.  A general
sufficient condition for the quotient to be a manifold is that the action
is proper and free.)  The Poisson bracket on $M$ then descends to give a
Poisson structure on $M/G$.  The moment map,
$$\mu : M \longrightarrow \LieAlg{g}^*,
$$
determines the symplectic leaves of this Poisson structure: \ let $\xi =
\mu(m)$, let $\cal O$ be the coadjoint orbit through $\xi$ and let $G_\xi$
be the stabilizer of $\xi$. Assume for simplicity that $\mu^{-1}(\xi)$ is
connected and
$\mu$ is submersive at $\mu^{-1}(\xi)$.
Then, the leaf through $m$ is
$$
\mu^{-1}({\cal O}_\xi) / G \approx \mu^{-1}(\xi)/G_\xi.
$$
These symplectic leaves are often called the Marsden-Weinstein reductions
$M_{red}$ of $M$.

As an example, consider a situation where $G$ acts on $X$ with nice
quotient $X/G$.  The lifted action of $G$ on $M = T^*X$ is Poisson, and has
a quotient $M/G$ which is a vector bundle over $X/G$.  The cotangent
$T^*(X/G)$ sits inside $(T^*X)/G$  as the symplectic leaf over the trivial
orbit ${\cal O}_0 = \{0\} \subset \LieAlg{g}^*.$

In contrast, the action of $G$ on $\LieAlg{g}^*$ does not in general admit a
reasonable quotient.  Its action on the dense open subset $\LieAlg{g}^*_{reg}$
of regular elements (cf. example \ref{example-coadjoint-orbits})
does have a quotient,
which is a manifold.  The Poisson structure on the quotient is trivial, so
the symplectic leaves are points, in one-to-one correspondence with the
regular orbits.  We refer to this quotient simply as $\LieAlg{g}^*/G$.  The
map $\pi_{reg} \; : \; \LieAlg{g}^*_{reg} \rightarrow \LieAlg{g}^*/G$ extends
to $\pi \, : \, \LieAlg{g}^* \rightarrow \LieAlg{g}^*/G$, and there is a sense
in which $\LieAlg{g}^*/G$ really is the quotient of all $\LieAlg{g}^*$.  Each
coadjoint orbit $\cal O$ is contained in the closure of a unique regular
orbit ${\cal O}'$ and $\pi({\cal O}) = \pi_{reg}({\cal O}')$.

\medskip
\noindent
\underline{{\bf A Diagram of Quotients}} \nopagebreak[3]

In the general situation of Poisson action (with a nice quotient $\pi$) of
$G$ on a symplectic manifold $(M,\sigma)$, there is another, larger, Poisson
manifold $\bar{M}$, which can also be considered as a reduction of $M$ by
$G$.  Everything fits together in the commutative diagram of Poisson maps:

\begin{equation}\label{diagram-of-quotients}
{\divide\dgARROWLENGTH by 2
\begin{diagram}
\node[2]{M}
\arrow{ssw,l}{\pi}
\arrow{s}
\arrow{sse,t}{\mu}
\\
\node[2]{\bar{M}}
\arrow{sw,r}{\bar{\pi}}
\arrow{se,b}{\bar{\mu}}
\\
\node{M/G}
\arrow{sse}
\arrow{se}
\node[2]{\LieAlg{g}^{*}}
\arrow{sw}
\arrow{ssw}
\\
\node[2]{\LieAlg{g}^{*}/G}
\arrow{s}
\\
\node[2]{(0)}
\end{diagram}
}
\end{equation}

$\bar{M}$ may be described in several ways:
\begin{itemize}
\item $\bar{M}$ is the quotient of  $M$ by the equivalence relation $m \sim
gm$ if $g \in G_{\mu(m)}$, i.e., if $g(\mu(m)) = \mu(m)$.
\item $\bar{M}$ is the fiber product $\bar{M}$ = $(M/G)
\times_{(\LieAlg{g}^*/G)} \LieAlg{g}^*$.
\item $\bar{M}$  is the dual realization to the realization $M \rightarrow
\LieAlg{g}^*/G.$
\end{itemize}

A {\em realization}  of a Poisson manifold $P$ is defined to be a Poisson map
from
a symplectic manifold $M$ to $P$ (see \cite{We}).  The realization will be
called {\it full} if it is submersive.  A pair of realizations $P_2
\stackrel{f_2}{\longleftarrow} M \stackrel{f_1}{\longrightarrow} P_1$ is
called a dual pair if functions on one induce vector fields along the
fibers of the other (i.e., the two opposite foliations are symplectic
complements of each other).

We note that in the diagram of quotients,
any two opposite spaces are a dual pair of realizations.

Given a full dual pair with connected fibers, the symplectic leaf foliations
on $P_1$ and $P_2$ induce the same foliation on $M$ ($P_1$ and $P_2$ have
the ``same'' Casimir functions).  The bijection between
symplectic leaves on $P_1$ and $P_2$ is given by
$$
P_1 \supset S_1 \mapsto f_2(f^{-1}_1(S_1)) = f_2(f^{-1}_1(x)) \ \ \
\forall \;
x \in S_1.
$$

Returning to moment maps, we have over a coadjoint orbit
${\cal O} \subset \LieAlg{g}^*$:

\begin{itemize}
\item $\mu^{-1}(\cal{O})$ is coisotropic in $M$
\item $\pi(\mu^{-1}(\cal{O}))$ is a symplectic leaf $M_{red}$ in $M/G$
\item $\bar{\mu}^{-1}(\cal{O})$ is also a symplectic leaf in $\bar{M}$.  It is
isomorphic to $\mu^{-1}(\cal{O})$/(null), or to $\mu^{-1}(\cal{O})/ \sim$,
or to $M_{red} \times \cal O$.
\end{itemize}

\medskip
\noindent
\begin{example}  {\rm Take $M$ to be the cotangent bundle $T^*G$ of a Lie
group $G$.  Denote by $\mu_L : T^*G \rightarrow \LieAlg{g}^*$ the moment map
for the lifted left action of $G$.  The quotient $\pi : M \rightarrow M/G$
is just the moment map $\mu_R : T^*G \rightarrow \LieAlg{g}^*$ for
the lifted right action, and $\bar{M}$ is the fiber product $\LieAlg{g}^*
\times_{(\LieAlg{g}^*/G)} \LieAlg{g}^*$.}
\end{example}

\noindent
\begin{example}  {\rm If $G$ is a connected commutative group $T$, the
pair of nodes $\LieAlg{t}^*$ and $\LieAlg{t}^*/T$ coincide.
Consequently, so do $M/T$ and $\bar{M}$.  The diagram of quotients
degenerates to}

\begin{equation}
{\divide\dgARROWLENGTH by 4
\begin{diagram}[M]
\node{M}
\arrow{s,r}{\pi}
\\
\node{M/T}
\arrow{s,r}{\bar{\mu}}
\\
\node{\LieAlg{t}^{*}}
\arrow{s}
\\
\node{(0)}
\end{diagram}
}
\end{equation}

\end{example}
\noindent
\begin{example}\label{diagram-two-quotients}  {\rm Consider two Poisson
actions on $(M,\sigma)$ of two
groups $G,T$ with moment maps $\mu_G, \mu_T$ with connected fibers.  Assume
that

\begin{description}
\item[i)]   The actions of $G$ and $T$ commute.
\end{description}
\noindent
It follows that $\mu_T : M \rightarrow \LieAlg{t}^*$ factors through
$M/G$ and $\mu_G : M \rightarrow \LieAlg{g}^*$ factors through $M/T$.
Assume moreover
\begin{description}
\item[ii)]  $T$ is commutative,
\item[iii)] $M \rightarrow \LieAlg{g}^*/G$ factors through $\LieAlg{t}^*$
\end{description}

Then $\bar{\mu}_G : \bar{M}_G \rightarrow \LieAlg{g}^*$ factors through
$M/T$ and the two quotient diagrams fit nicely together:}

\begin{equation}
{\divide\dgARROWLENGTH by 4
\begin{diagram}
\node[3]{M}
\arrow{s}
\\
\node[3]{\bar{M}_{G}}
\arrow[2]{sw}
\arrow{se}
\\
\node[4]{M/T}
\arrow[2]{sw,b}{\bar{\mu}_{T}}
\arrow{se}
\\
\node{M/G}
\arrow{se}
\node[4]{\LieAlg{g}^*}
\arrow[2]{sw}
\\
\node[2]{\LieAlg{t}^*}
\arrow{se}
\\
\node[3]{\LieAlg{g}^*/G}
\arrow{s}
\\
\node[3]{(0)}
\end{diagram}
}
\end{equation}
\end{example}

\begin{rem} \label{rem-acihs-implies-maximal-commutative-subalgebra}
{\rm
Note that condition iii in example \ref{diagram-two-quotients}
holds whenever $M/G \rightarrow \LieAlg{t}^*$ is a completely
integrable system (with connected fibers). In that case the map
$M/G \rightarrow \LieAlg{t}^*$ pulls back
$C^\infty(\LieAlg{t}^*)$ to a
{\em maximal} commutative subalgebra ${\cal I}_T$ of $(C^\infty(M/G),\{,\})$.
The map $M/G \rightarrow \LieAlg{g}^*/G$ pulls back
$C^\infty(\LieAlg{g}^*/G)$ to a commutative Lie subalgebra ${\cal I}_G$ of
$(C^\infty(M/G),\{,\})$.
As the two group actions commute so do
the subalgebras ${\cal I}_G$ and ${\cal I}_T$.
By maximality, ${\cal I}_T$ contains ${\cal I}_G$ and consequently
$M/G \rightarrow \LieAlg{g}^*/G$ factors through $\LieAlg{t}^*$.
}
\end{rem}

The diagram of quotients for a Poisson action
(diagram \ref{diagram-of-quotients}) generalizes to an
analogous diagram for any full dual pair of realizations
$P_2\stackrel{f_2}{\longleftarrow} M \stackrel{f_1}{\longrightarrow} P_1$.
Denote by $\bar{M}$ the image of $M$ in the Poisson manifold
$P_1 \times P_2$ under the diagonal Poisson map $f_1 \times f_2 : M
\rightarrow P_1 \times P_2$.
The realization dual to $f_1 \times f_2 : M \rightarrow \bar{M}$
is the pullback of the symplectic leaf foliations on
$P_1$ or $P_2$ (they pull back to the same foliation of $M$).

The following is the analogue of Example \ref{diagram-two-quotients}
 replacing the commutative $T$-action by a realization:

\begin{example} \label{example-diagram-hexagon-plus-realization}
{\rm
\ Let $M/G \stackrel{\pi}{\longleftarrow} M
\stackrel{\mu}{\longrightarrow} \LieAlg{g}^*$ be the full dual pair
associated to a
Poisson action of $G$ on $M$ and $N \stackrel{\ell}{\longleftarrow} M
\stackrel{h}{\longrightarrow} B$ a full dual pair of realizations with
connected fibers where:

\begin{description}
\item [(i)] $h$ is $G$-invariant
\item [(ii)] $h : M \rightarrow B$ is a Hamiltonian map ($B$ is endowed
with the trivial Poisson structure) and
\item [(iii)] The composition $M \stackrel{\mu}{\longrightarrow} \LieAlg{g}^*
\longrightarrow \LieAlg{g}^*/G$ factors through $h : M \rightarrow B$.
\end{description}

\noindent
Then we get a diagram analogous to the one in example
\ref{diagram-two-quotients}:

\begin{equation}
{\divide\dgARROWLENGTH by 4
\begin{diagram}
\node[3]{M}
\arrow{s}
\\
\node[3]{\bar{M}_{G}}
\arrow[2]{sw,t}{\bar{\pi}}
\arrow{se}
\\
\node[4]{N}
\arrow[2]{sw,b}{\bar{h}}
\arrow{se}
\\
\node{M/G}
\arrow{se}
\node[4]{\LieAlg{g}^*}
\arrow[2]{sw}
\\
\node[2]{B}
\arrow{se}
\\
\node[3]{\LieAlg{g}^*/G}
\arrow{s}
\\
\node[3]{(0)}
\end{diagram}
}
\end{equation}

It follows that the Poisson map
$M \stackrel{h \times \mu}{\longrightarrow}
B \times_{(\LieAlg{g}^*/G)} \LieAlg{g}^*$
into the fiber product space factors through
the realization $M \stackrel{\ell}{\rightarrow} N$ dual to
$h : M \rightarrow B$.
If, moreover, $M/G \rightarrow B$ is a Lagrangian fibration,
then $M \stackrel{h \times \mu}{\longrightarrow}
B \times_{(\LieAlg{g}^*/G)} \LieAlg{g}^*$ is itself a realization dual to
$h : M \rightarrow B$.
}
\end{example}



\subsection{Finite dimensional Poisson loop group actions}
\label{sec-finite-dim-loop-group-actions}

We present in this section two elementary constructions related to finite
dimensional symplectic leaves in the
Poisson quotient $Q_\infty$ of an infinite dimensional
symplectic space $M$ by subgroups of loop groups. The material in this
section will only be used in section \ref{sec-compatibility-of-heirarchies}
so the reader may prefer to read it in conjunction with that section.

We will not construct the quotient $Q_\infty$. The spaces involved are
constructed independently. Rather, we will analyze the relationship between
the Poisson action of the loop group on the infinite dimensional
spaces and its descent to the finite dimensional symplectic leaves of
$Q_\infty$. In fact,
our main purpose in this section is to provide the terminology needed
in order to study the Poisson loop group action in the finite dimensional
setting (convention
\ref{convention-abused-hamiltonian-language}
and corollary
\ref{cor-hamiltonians-on-the-base}).

In section \ref{sec-finite-dim-approaximations} we note that the infinitesimal
Hamiltonian actions of elements of the loop group descend to Hamiltonian
vector fields on finite dimensional symplectic approximations $M_{(l,l)}$. The
$M_{(l,l)}$'s dominate finite dimensional Poisson subvarieties $Q_l$ of
$Q_{\infty}$ with positive dimensional fibers. In section
\ref{sec-type-loci} the action
of certain maximal tori in the loop group further descends to
finite Galois covers of certain (type) loci in $Q$ and we examine the sense
in which it is Hamiltonian.

\subsubsection{Finite dimensional approximations}
\label{sec-finite-dim-approaximations}

The loop group $\LoopGroup$ is the group
$GL(n,\ComplexNumbers((z)))$.
The level infinity group
$\LevelInfinityGroup$ is its positive part
$GL(n,\ComplexNumbers[[z]])$.
Let $(M,\sigma)$ be a symplectic variety with a Poisson loop group
action whose moment map is

\[ \mu : M \rightarrow \LoopAlg^*. \]

In section \ref{sec-compatibility-of-heirarchies} $M$ will be the
cotangent bundle of a projective (inverse) limit of finite dimensional
smooth algebraic varieties
(the cotangent bundles of
the moduli spaces of vector bundles with level structure).
It is thus the inductive
(direct) limit of projective limits of finite dimensional varieties.
All constructions (morphisms, group actions, symplectic structures
etc ...) can be made precise as limits of the standard constructions on
finite dimensional approximations. We will omit the technical details as
our point is to transfer the discussion back to the finite dimensional
symplectic leaves of the Poisson quotient
$Q_\infty := M/\LevelInfinityGroup$.

Let $\LevelInfinitySubgroup{l}$, $l \geq -1$, be the subgroup of
$\LevelInfinityGroup$ of elements equal to $1$ up to order $l$.
Denote by $\mu_{\LevelInfinitySubgroup{l}}$ its moment map.
We assume that the  subquotients
\[
M_{(l,k)} := \mu_{\LevelInfinitySubgroup{l}}^{-1}(0)/\LevelInfinitySubgroup{k},
\ \ \ k \geq l,
\]
are smooth, finite dimensional and that they approximate $M$:
\[M = \lim_{l \rightarrow \infty} \lim_{\infty \leftarrow k}  M_{(l,k)}.\]
Notice that $M_{(l,l)}$ is a symplectic reduction, hence symplectic.

Let $a$ be an element of the loop algebra $\LoopAlg$
with poles of order at most $l_{0}$.
The Hamiltonian vector field $\xi_a$ on $M$ is an infinite double
sequence of Hamiltonian vector fields on $M_{(l,k)}$, $l \geq 0$,
$k\geq \max\{l,l_0\}$
compatible with respect to projections and inclusions (by a Hamiltonian
vector field on $M_{(l,k)} \subset M_{(k,k)}$ we mean, the restriction
of a Hamiltonian vector field on $M_{(k,k)}$
which is tangent to $M_{(l,k)}$).

The quotient $Q_\infty := M/\LevelInfinityGroup$ is the direct limit
$\lim_{l \rightarrow \infty}Q_l$ of the finite dimensional Poisson varieties
\[
Q_l :=  \mu_{\LevelInfinitySubgroup{l}}^{-1}(0)/\LevelInfinityGroup
      = M_{(l,k)}/ G_k
      = M_{(l,l)}/ G_l
\]
where $G_k:=\LevelInfinityGroup/\LevelInfinitySubgroup{k}$ is the
finite dimensional level-$k$ group
(we assume that the quotients $Q_l$ are smooth).

\begin{example}
{\rm
The homogeneous $\LevelInfinityGroup$-space
${\cal U}_{\infty}:=\LevelInfinityGroup/ GL(n,\ComplexNumbers)$
is endowed with a canonical infinitesimal $\LoopGroup$-action
via its embedding as the degree-$0$ component of the homogeneous
$\LoopGroup$-space $\LoopGroup / GL(n,\ComplexNumbers[[z^{-1}]])$
\[\LevelInfinityGroup/ GL(n,\ComplexNumbers) \hookrightarrow
\LoopGroup / GL(n,\ComplexNumbers[[z^{-1}]])\]
(the degree of $a\in \LoopGroup$ is the signed order of the pole/zero
of $\det(a)$).
Let $M$ be an open subset of the cotangent bundle
$T^*{\cal U}_{\infty}$ for which
the regularity assumptions on the approximating quotients $M_{(l,k)}$
hold. This will be made precise in section \ref{sec-polynomial-matrices}
and the quotients $Q_l$ will be the spaces of conjugacy classes
of polynomial matrices studied in that section.
}
\end{example}

\bigskip
Unfortunately,
the action of $a \in \LoopAlg$ above is not defined on $Q_{l}$.
It is well defined only
when we retain at least the $l_0$-level structure, i.e., on $M_{(l,k)}$,
$k \geq l_0$. In section \ref{sec-type-loci} we will see that the action
of certain maximal tori in $\LoopGroup$ descends to
{\em finite} Galois covers of certain loci in $Q$.

\subsubsection{Type loci}  \label{sec-type-loci}

Let $(M,\sigma)$ be a smooth symplectic variety endowed with an
infinitesimal Poisson action $\mu_{G}^*: \LieAlg{g} \rightarrow
[\Gamma(M,\StructureSheaf{M}),\{,\}]$
of a group $G$. Consider a subgroup $G^+ \subset G$, a
commutative subgroup $T \subset G$, and their intersection
$T^+ := T \cap G^+$. Assume further that the following conditions hold:

\smallskip
\noindent
i) The infinitesimal $G^+$-action integrates to a free action on $M$,\\
ii) $T^+$ is a maximal commutative subgroup whose Weyl group
$W_{T^+}:= N_{G^+}(T^+)/T^+$ is finite.


\begin{definition} \label{def-group-theoretic-definition-type-loci}
{\rm
The {\em type} $\tau$ of $T$ is the class of all commutative subgroups
$T'$ of $G$ which are conjugate to $T$ via an element of $G^+$.
}
\end{definition}
Let $W := [N_{G^+}(T^+)\cap N(T)]/T^+$ be the corresponding subgroup of both
$W_{T^+}$ and $W_{T}$. Denote by
\[
\LieAlg{g}^*_\tau \subset \LieAlg{g}^*
\]
(respectively,
$\LieAlg{g}^*_T \subset \LieAlg{g}^*$) the subset of elements whose
stabilizer (with respect to the coadjoint action) is a torus of type
$\tau$ (respectively, precisely $T$).

\begin{example} \label{example-loop-group-level-infinity-group}
{\rm
Let $G$ be the loop group, $G^+$ the level infinity group
and $T \subset G$ a maximal torus of type $\underline{n}$ determined
by a partition of the integer $n$ (see section \ref{sec-the-heirarchies}).
In this case $G^+$ and $T$ generate $G$. It follows that $W=W_T=W_{T^+}$ and
the type $\tau$ is invariant throughout a coadjoint orbit in $\LieAlg{g}^*$.
}
\end{example}

Assume that a ``nice'' (Poisson) quotient $Q := M/G^+$ exists. Let
\[M^\tau := \mu_G^{-1}(\LieAlg{g}^*_\tau), \ \ \ \mbox{and} \ \ \
Q^\tau := M^\tau/G^+ \subset Q\]
be the loci of type $\tau$. Note that for each $T$ of type $\tau$ there is a
canonical isomorphism
\[
\mu^{-1}_G(\LieAlg{g}^*_T)/ [N_{G^+}(T^+) \cap N(T)] \IsomRightArrow
Q^{\tau} \subset Q.
\]
In particular, a choice of $T$ of type $\tau$ determines a canonical
$W$-Galois cover of $Q^\tau$
\begin{equation} \label{eq-the-galois-cover}
\tilde{Q}^T := \mu^{-1}_{G}(\LieAlg{g}^*_T)/T^+.
\end{equation}
All the $\tilde{Q}^T$ of type $\tau$ are isomorphic (not canonically) to a
fixed abstract $W$-cover $\tilde{Q}^\tau$. Note that $\tilde{Q}^T$ is a
subset of $M/T^+$.  We get a canonical ``section'' (the inclusion)
\begin{equation} \label{eq-the-section-from-the-galois-cover}
s_T : \tilde{Q}^T \hookrightarrow M/T^+
\end{equation}
into a $T$-invariant subset. Consequently, we get an induced
$T$-action on the Galois cover $\tilde{Q}^T$. The moment map
$\mu_T$ is $T$-invariant, hence, descends to $M/T^+$.
Restriction to $s^T(\tilde{Q}^T)$ gives rise to a canonical map
\begin{equation} \label{eq-loop-group-moment-map-on-galois-covers}
\bar{\mu}_T : \tilde{Q}^T \rightarrow \LieAlg{t}^*.
\end{equation}

\bigskip
The purpose of this section is to examine {\em the extent to which
$\bar{\mu}_T$ is the moment map of the
$T$-action with respect to the Poisson structure on $Q$}.
In general, the $G^+$-equivariant projection
\[
j: \LieAlg{g}^* \SurjectiveRightArrow (\LieAlg{g}^+)^*
\]
might {\em forget the type}. Coadjoint orbits
$S \subset (\LieAlg{g}^+)^*$ may intersect nontrivially the images
$j(\LieAlg{g}^*_\tau)$ of several types (e.g., take $S=0$ in example
\ref{example-loop-group-level-infinity-group} and observe that
the kernel of $j$ intersects coadjoint orbits of all types).
Consequently,
symplectic leaves $Q_S$ of $Q$ would intersect nontrivially several type loci
$Q^\tau$.
If $Q_S^{\tau^{open}}$ is an open subvariety of $Q_S$ of type $\tau$
(e.g., if $Q$ is the disjoint union of finitely many type loci
and $\tau$ is a {\em generic type})
then the corresponding open subvariety $\tilde{Q}_S^{T^{open}}$
of $\tilde{Q}_S^T$ will be a symplectic variety.
In this case the $T$-action
on $\tilde{Q}_S^{T^{open}}$ is Poisson whose moment map $\bar{\mu}_T$
is given by (\ref{eq-loop-group-moment-map-on-galois-covers}).

The Galois $W$-covers $\tilde{Q}_S^T$ of the nongeneric type loci
in $Q_S$ are not symplectic.
Nevertheless, motivated by the fact that $\bar{\mu}_T$ can be extended
canonically to $M/T^+$
\begin{equation}\label{diag-extending-the-moment-map-from-the-galois-cover}
{\divide\dgARROWLENGTH by 2
\begin{diagram}
\node{M}
\arrow{s}
\arrow[2]{e,t}{\mu_T}
\node[2]{\LieAlg{t}^*}
\\
\node{M/T^+}
\arrow{s}
\arrow{ene}
\node{\tilde{Q}^T_S}
\arrow{s}
\arrow{w,t}{\supset}
\arrow{ne,b}{\bar{\mu}_T}
\\
\node{Q}
\node{Q^\tau_S}
\arrow{w,t}{\supset}
\end{diagram}
}
\end{equation}
we will adopt the:

\begin{convention} \label{convention-abused-hamiltonian-language}
{\rm
i)
Given an element $h$ of $\LieAlg{t}$ we will say that the corresponding
vector field $\bar{\xi}_h$ on $\tilde{Q}_S^T$ is the
{\em Hamiltonian vector field of $h$}
(even if the type $\tau$ of $T$ is not generic in $Q_S$).
ii)
We will refer to the pair $(\bar{\mu}_T,\mu_T)$ as the moment map
of the $T$-action on $\tilde{Q}_S^T$.
}
\end{convention}

\begin{rems}
{\rm
Let $G$ be the loop group and $M,G^+,T$ as in section
\ref{sec-finite-dim-approaximations},
\begin{enumerate}
\item
Diagram (\ref{diag-extending-the-moment-map-from-the-galois-cover})
has an obvious finite dimensional approximation in which
$Q^\tau_S$, $\tilde{Q}^T_S$ and $T$ stay the same
but with $M$ replaced
by $M_{(l,l)}$ and $Q$ by $Q_{l}$. By $\mu_{T}$ we mean in this
context a linear homomorphism
$\mu_T^*:\LieAlg{t} \rightarrow \LieAlg{t}/\LoopAlgSubtorus{l}
\rightarrow \Gamma(M_{(l,l)},\StructureSheaf{M_{(l,l)}}).$
\item
(Relation with the diagram of quotients (\ref{diagram-of-quotients}))
Let $S$ be a coadjoint orbit of level $l$, i.e.,
$S \subset \LieAlg{g}_l^* :=
 (\LevelInfinityAlg/\LevelInfinitySubalg{l})^* \subset (\LevelInfinityAlg)^*$.
There is a rather subtle relationship between the Galois cover
$\tilde{Q}_S^T \rightarrow Q^\tau_S$ and the space
$\bar{M}_{(l,l)}$ dual to $\LieAlg{g}_l^*/G_{l}$ from the diagram of
quotients (\ref{diagram-of-quotients}) of level $l$.
The Galois cover
$\tilde{Q}_S^T \rightarrow Q^\tau_S$ factors canonically through
an intermediate subspace $\tilde{Q}_S^T/\sim$
of $\bar{M}_{(l,l)}$.
Note that the loop group moment map $\mu_{G_\infty}$
descends to a map
\[
\bar{\mu}_{G_\infty}:\tilde{Q}_S^T \rightarrow (\LoopAlg^*)_T
\subset \LoopAlg^*.
\]
Two points $\tilde{x}_1, \tilde{x}_2 \in \tilde{Q}_S^T$
in a fiber over $x \in Q^\tau_S$ are
identified in $\tilde{Q}_S^T/\sim$
if and only if
$\bar{\mu}_{G_\infty}(\tilde{x}_1)$ and $\bar{\mu}_{G_\infty}(\tilde{x}_1)$
project to the same point in $S \subset (\LevelInfinityAlg)^*$.
The relation $\sim$ is a geometric realization of the partial
type-forgetfullness of the projection
$j:\LoopAlg^* \rightarrow (\LevelInfinityAlg)^*$.
The loci in $\bar{M}_{(l,l)}$ at which the type is not forgotten
are precisely the loci to which the moment map of the infinitesimal
{\em loop group}
action descends
(from $M^T:= \mu_{G_{\infty}}^{-1}((\LoopAlg^*)_T)$. Note that
the moment map of the level infinity subgroup descends
by definition of the quotient $\bar{M}_{(l,l)}$).
In particular, the infinitesimal action of the maximal
torus $\LieAlg{t}$ integrates to a Poisson action in these loci.
(See section \ref{sec-compatibility-of-stratifications}
for examples of such loci.)
\end{enumerate}
}
\end{rems}

\bigskip
Assume further that we have a ``nice'' quotient $B := M/G$ and that
$G$ is generated by $T$ and $G^+$.
We get the type loci $B^\tau := M^\tau/G$.
Fixing $T$ of type $\tau$ we get the $W$-cover $\tilde{B}^T := M^T/T$.
The restriction of the moment map $\mu_T$ to $M^T$ descends further to
$\phi_T:\tilde{B}^T \rightarrow \LieAlg{t}^*$
and we get the commutative diagram

\begin{equation}\label{diag-factoring-the-moment-map-through-base}
{\divide\dgARROWLENGTH by 2
\begin{diagram}
\node[3]{M^T}
\arrow{sw}
\arrow{se,t}{\mu_G}
\\
\node[2]{\tilde{Q}^T}
\arrow{sw,t}{\tilde{h}}
\node[2]{\LieAlg{g}^*_T}
\arrow{se}
\\
\node{\tilde{B}^T}
\arrow[4]{e,t}{\phi_T}
\node[4]{\LieAlg{t}^*}
\end{diagram}
}
\end{equation}

\begin{corollary}
\label{cor-hamiltonians-on-the-base}
i) The $T$-action on $M$ descends to a canonical action on the
$W$-cover $\tilde{Q}^T$ of the type locus $Q^\tau \subset Q$.
ii) Its moment map,
in the sense of convention \ref{convention-abused-hamiltonian-language},
is $(\phi_T\circ \tilde{h},\mu_T)$.
iii) If the type $\tau$ of $T$ is the generic type in
a symplectic leaf $Q_S\subset Q$
and $Q_S^{\tau^{open}} \subset Q_S$  is an open subvariety,
then the corresponding $W$-cover $\tilde{Q}_S^{T^{open}}$
is symplectic and $\phi_T\circ \tilde{h}$ is the moment map of
the $T$-action in the usual sense.
\end{corollary}
\newpage
\section{Geodesic flow on an ellipsoid} \label{ch3}

     Consider the geodesic flow on an ellipsoid
$$
E =
\{(x_1,\cdots,x_{n+1})|\sum^{n+1}_{i=1} \frac{1}{a_i} x^2_i = 1\}
\subset \R^{n+1},
$$
where the metric is induced from the standard
one on $\R^{n+1}$, and where the $a_i$ are distinct positive
numbers, say
$$
0 < a_1 < \cdots < a_{n+1}.
$$

     For $n=1$, the problem is to compute arc length on an
ellipse.  It amounts to computing the integral
$$
s = \int{
\sqrt{\frac{a^2_1 + (a_2-a_1)x^2}{a_1(a_1 - x^2)}}\ \ dx}.
$$
(Hence the
name {\it elliptic} for this and similar integrals.)

     For $n=2$, the problem was solved by Jacobi.  Each geodesic
$\gamma$ on $E$ determines a hyperboloid $E'$, intersecting $E$
in a pair of ovals.  The geodesic $\gamma$ oscillates in the band
between these ovals, meeting them tangentially.  In fact, each
tangent line of $\gamma$ is also tangent to the hyperboloid $E'$.
 The solutions can be parametrized explicitly in terms of
hyperelliptic theta functions.

     The geodesic flow on an $n$-dimensional ellipsoid is
integrable, in fact algebraically integrable.  We will see this,
first using some elementary geometric techniques to describe the
geodesics concretely, and then again using the algebraic
description of hyperelliptic jacobians which will be extended
later to all spectral curves.

\subsection{Integrability}
     The geodesic flow on the $2n$-dimensional symplectic
manifold $TE \approx T^*E$ is given by the Hamiltonian function
$h$=length square.  We need $n-1$ farther, commuting,
independent, Hamiltonians.

     Consider the family of quadrics confocal to $E$:
$$
E_\lambda: \; \; \sum^{n+1}_{i=1} \frac{x^2_i}{a_i-\lambda} = 1,
$$
depending on a parameter $\lambda$.  (The name makes sense only
when $n=1$: we get the family of ellipses $(\lambda < a_1)$,
hyperbolas $(a_1 < \lambda < a_2)$, and empty (real) conics
$(\lambda > a_2)$, with fixed foci.)

     Here is an intrinsic way to think of this family.  Start
with a linear pencil $$Q_\lambda = Q_0 + \lambda Q_\infty,\ \ \ \
\ (\lambda \in \P^1)$$ of quadrics in general position in
projective space $\P^{n+1}$.  By ``general position'' we mean
that there are exactly $n+2$ values of $\lambda \in \P^1$ such
that $Q_\lambda$ is singular, and for those $\lambda$,
$Q_\lambda$ is a cone (i.e. its singular locus, or vertex, is a
single point).

\proclaim{Lemma}.  A generic linear subspace $L \approx \P^{k-1}$
in $\P^{n+1}$ is tangent to $Q_\lambda$ for $k$ values of
$\lambda$.  The points of tangency $p_\lambda$ are pairwise
harmonic with respect to each of the quadrics $Q_\mu$.\par

\noindent
{\bf Proof:}
     Four points of $\P^1$ are harmonic if their
cross ratio is $-1$; e.g. $0,\infty,a,-a$.  Two points $p_1,p_2
\in \P^1$ are harmonic with respect to a quadric $Q$ if the set
$\{p_1,p_2\} \cup (Q \cap \P^1)$ is harmonic.  For example, two
points on the line at infinity in $\P^2$ (i.e. two directions in
the affine plane) are harmonic with respect to some (hence
every) circle, iff the directions are perpendicular.

     Since $Q$ is tangent to $L$ if and only if $Q \cap L$ is
singular, the first part of the lemma follows by restriction of
the pencil to $L$.  The second part follows by restricting to the
line $p_{\lambda_1},p_{\lambda_2}$, where in appropriate
coordinates $Q_{\lambda_1} = x^2$ and $Q_{\lambda_2} = y^2$, so
the points of tangency are $0,\infty$ and the quadric $Q_\lambda$
vanishes at $\pm a$, where $a^2 = \lambda$.
\EndProof

     We choose the parameter $\lambda$ so that $Q_\infty$ is one
of the singular quadrics.  The dual $Q^*_\lambda$ of a
non-singular $Q_\lambda$ is a non-singular quadric in $(\P^{n+1})^*$.
The dual of $Q_{\infty}$ is a hyperplane $H_{\infty} \subset (\P^{n+1})^*$
(corresponding to the vertex of $Q_\infty$), with a
non-singular quadric $Q^*_\infty \subset H_{\infty}$.  We get a family of
confocal quadrics by restriction to the affine space $\R^{n+1}:=
(\P^{n+1})^*\smallsetminus H_{\infty}$:
$$
E_\lambda:=
Q^*_\lambda|_{\R^{n+1}}.
$$
If we choose coordinates so that
$$
\begin{array}{lcl}Q_\infty & = & \sum^{n+1}_{i=1} x^2_i\\Q_0 &
= & \sum^{n}_{i=0} a_i x^2_i\end{array}
$$
(where $a_0 = -1$ and
the other $a_i$ are as above), we retrieve the original
$E_\lambda$. (Euclidean geometry in $\R^{n+1}$ is equivalent, in
the sense of Klein's program, to the geometry of $\P^{n+1}$ with
a distinguished ``light-cone'' $Q_\infty$.  In this equivalence,
$E = E_0$ corresponds to $Q_0$, which determines the pencil
$\{Q_\lambda\}$, which corresponds to the confocal family
$\{E_\lambda\}$.)

     Dualizing the lemma, for $k = n+1,n$, gives the following
properties of the confocal family.
(The reader is invited to amuse herself by drawing the case
$n=1$ in the plane.)
\begin{list}{{\rm(\arabic{bean})}}{\usecounter{bean}}
\item Through a generic point $x$ of $\R^{n+1}$ pass $n+1$ of the
$E_\lambda$.
\item These $n+1$ quadrics intersect perpendicularly at $x$.
\item A generic line $\ell$ in $\R^{n+1}$ is tangent to $n$ of
the $E_\lambda$.
\item The tangent hyperplanes to these $n$ quadrics (at their
respective points of tangency to $\ell$) are perpendicular.
\end{list}

     By property (3) we can associate to a generic line $\ell$ in
$\R^{n+1}$ an unordered set of $n$ values $\lambda_i$, $1 \leq i
\leq n$, such that $Q_{\lambda_i}$ is tangent to $\ell$.  When
$\ell$ comes from a point of $TE$, one of these, say $\lambda_n$,
equals $0$.  The remaining $n-1$ values $\lambda_i$ (or rather,
their symmetric functions) give $n-1$ independent functions on
$TE$; in fact, they can take an arbitrary $(n-1)$-tuple of
values.  These functions descend to the projectivized tangent
bundle $\P(TE)$; so together with the original Hamiltonian $h$ (=
length squared) they give $n$ independent functions on $TE$.  The
key to integrability is:

\proclaim{Chasles' Theorem}.  The $\lambda_i$ are flow
invariants, i.e. they are constant along a geodesic $\gamma =
\gamma(t)$.\par

\noindent
{\bf Proof}.  For any curve $\gamma(t)$ in $\R^{n+1}$, the family
of tangent lines $\ell(t)$ gives a curve $\Lambda$ in the
Grassmannian $Gr(1,\R^{n+1})$ of affine lines in $\R^{n+1}$.
This curve is developable, i.e. its tangent line
$T_{\ell(t)}\Lambda$ is given by the pencil of lines through
$\gamma(t)$ in the osculating plane of $\gamma$ at $t$.  When
$\gamma$ is a geodesic, this plane is the span of $\ell(t)$ and
the normal vector $n(t)$ to $E$ at $\gamma(t)$. Write $\lambda_i(t)$
for the value of $\lambda_i$ at $\ell(t)$.

     Let $Z_i$ be the hypersurface in $Gr(1,\R^{n+1})$
parametrizing lines tangent to $E_{\lambda_i(t)}$, for some fixed
$t$.  The tangent space $T_{\ell(t)}Z_i$ contains all lines
through $\gamma(t)$ in the tangent hyperplane
$T_{p_i(t)}E_{\lambda_{i}(t)}$, and this hyperplane contains the
normal $n(t)$, by property (4).  Hence:
$$
T_{\ell(t)}\Lambda
\subset \{\mbox{lines\ through}\ \gamma(t),\ \mbox{in}\
T_{p_i(t)}E_{\lambda_{i(t)}}\} \subset T_{\ell(t)}Z_i\ \ \ \ \ i =
1,\cdots,n-1.
$$
If the family $\Lambda$ of  tangent lines to a geodesic meets $Z_i$,
it must therefore stay in it.
\EndProof


\subsection{Algebraic integrability}

     Since a line $\ell$ determines two (opposite) tangent
vectors of given non-zero length, we have identified the fiber of
the geodesic flow as a double cover $\widetilde{K}$ of
$$K :=
\{\ell \in Gr(1,\R^{n+1})|\ell \ \mbox{is\ tangent\ to}\
E_{\lambda_1},\cdots, E_{\lambda_{n}}\}.
$$
Next we want to
interpret this in terms of the real points of a complex abelian
variety.  We follow Kn\"{o}rrer's approach \cite{knorrer}, which in turn is
based on \cite{moser},\cite{reid} and \cite{donagi-group-law}.

     Start with the pencil of quadrics in $\P^{2n+1}$ (over $\C$):
$$
Y_\lambda := Y_0 + \lambda Y_\infty
$$
with
$$
\begin{array}{lcl}
Y_0 & = & \sum^{2n+1}_{i=1} a_i x^2_i - x^2_0\\
Y_\infty & = & \sum^{2n+1}_{i=1} x^2_i.\end{array}
$$
The
base locus $X = Y_0 \cap Y_\infty$ is non-singular if the $a_i$
are distinct.  We set $a_0 = \infty$.  The family of linear
subspaces $\P^n$ contained in a fixed quadric $Y_\lambda$
consists of two connected components, or rulings, for the
non-singular $Y_\lambda$
$(\lambda \not\in \{a_i\})$, and of a single
ruling for $\lambda = a_i$.  We thus get a double cover $$\pi: C
\to \P^1$$ of the $\lambda$-line, parametrizing the rulings.
(More precisely, one considers the variety of pairs
$$
{\cal P} =
\{(A,\lambda)|A\ \ \mbox{is\ a}\ \P^n\ \ \mbox{contained\ in}\
Q_\lambda\},
$$
and takes the Stein factorization of the second
projection.)  Explicitly, $C$ is the hyperelliptic curve of genus
$n$:
$$C: \;  s^2 = \Pi^{2n+1}_{i=1} (t-a_i).$$
Miles Reid \cite{reid} showed
that the Jacobian $J(C)$ is isomorphic to the variety
$$F:= \{A
\in Gr(n-1,\P^{2n+1})|A \subset X\}
$$
of linear subspaces in the
base locus.  An explicit group law on $F$ is given in
\cite{donagi-group-law}, and
corresponding results for rank $2$ vector bundles on $C$ are in
\cite{DR}.
Since we are interested in a family of varieties $F$
with varying parameters, we
need some information about the isomorphism.  Let $Pic^d(C)$
denote the variety parametrizing isomorphism classes of
degree-$d$ line bundles on $C$.  Then $Pic^0(C) = J(C)$ is a group,
$Pic^1(C)$ is a torser (= principal homogeneous space) over it,
but, these two have no natural identification; while $Pic^2(C)
\approx J(C)$ canonically, using the hyperelliptic bundle on
$C$.  It turns out that $F$ is isomorphic to $Pic^0(C)$ and
to $Pic^1(C)$, but neither isomorphism is canonical.  Rather, we
may think of $F$ as ``$Pic^{\frac{1}{2}}(C)$'': it is a torser
over $J(C)$, and has a natural torser map
$$
F \times F \to
Pic^1(C).
$$
All of this is based on the existence of a natural
morphism
$$
j: F \times C \to F .
$$
The ruling $p$ (on the quadric
$Y_{\pi(p)}$) contains a unique subspace $\P^n$ which contains a
given $\P^{n-1}$-subspace $A \in F$.  This $\P^n$ intersects $X$
in the union of $A$ and another element of $F$, which we call
$j(A,p)$.  We can also think of $j$ as a family of involutions of
$F$, indexed by $p \in C$.  This extends to a map
$$
F \times
J(C) \to F
$$ which gives the torser structure on $F$.  Once $F$
is thus identified with $J(C)$, the map $j$ becomes $$j(A,p) = p
- A,$$ up to an additive constant.  Since this is well defined
globally, points $A \in F$ must behave as line bundles on $C$ of
``degree $\frac{1}{2}$''.  In particular, we have for $0 \leq i
\leq 2n+1$ the involution
$$\begin{array}{cl}j_i: & F \to F\\& A
\mapsto j(A,a_i),\\ \end{array}$$ where we set $a_i = \infty$ and
identify the $a_i \in \P^1$ with the $2n+2$ Weierstrass points
$\pi^{-1}(a_i) \in C$.  Explicitly, each $j_i$ is induced by the
linear involution $$\overline{j}_i: \P^{2n+1} \to \P^{2n+1}$$
flipping the sign of the $i$-th coordinate.

     Consider the linear projection
$$
\begin{array}{rcl}\rho:
\P^{2n+1} & \to & \P^{n+1}\\
(x_0,\cdots,x_{2n+1}) & \mapsto &
(x_0,\cdots,x_{n+1}),\\
\end{array}
$$
which commutes with the
$\overline{j}_i,\ n + 2 \leq i$.  Recall that in $\P^{n+1}$ we
have the pencil of quadrics $Q_\lambda$, with dual quadrics
$E_\lambda$ in $(\P^{n+1})^*$.

\proclaim{Proposition}.
\begin{list}{{\rm(\roman{bacon})}}{\usecounter{bacon}}
\item The projection $\rho$ maps $F$ to $$F' := \{B \in Gr(n-
1,\P^{n+1})|B\ \mbox{is\ tangent\ to}\
Q_{n+2},\cdots,Q_{2n+1}\}.$$
\item The induced $\rho: F \to F'$ is a finite morphism of degree
$2^n$, and can be identified with the quotient of $F$ by the
group $G \approx (\Z/2\Z)^n$ generated by the involutions $j_i$,
$n + 2 \leq i \in 2n+1$.
\item Duality takes $F'$ isomorphically to the variety $K$ of
lines in $({\Bbb P}^{n+1})^*$ tangent to $E_{\lambda_i}$, $\lambda_i =
a_{n+1+i}$, $1 \leq i \leq n$.
\end{list}
\par

     We omit the straightforward proof.  Let $\widetilde{G}
\subset G$ be the index-$2$ subgroup generated by the products
$j_{i_1} \circ j_{i_2}$, and set
$$
\widetilde{K} :=
F/\widetilde{G}.
$$
We obtain natural commuting maps, whose
degrees are indicated next to the arrows:
$$
\begin{array}{ccccc}F
& \stackrel{2^{n-1}}{\longrightarrow} & \widetilde{K} &
\stackrel{2^{n+1}}{\longrightarrow} & Pic^{1}(C)\\
& & \downarrow^2 & & ^2\downarrow\\
& & K & \stackrel{2^{n+1}}{\longrightarrow}
& \mbox{Kummer}^{1}(C).\\
\end{array}
$$
Here $\mbox{Kummer}^d(C)$
stands for the quotient of $Pic^dC$ by the involution
$$
L \mapsto dH-L,
$$
where $H$ is the hyperelliptic bundle $\in Pic^2(C)$.
The composition of the maps in the top row is multiplication by
$2$:
$$
F \approx Pic^{\frac{1}{2}} (C)
\stackrel{\cdot 2}{\longrightarrow} Pic^{1}(C).
$$
In conclusion, the fiber of the geodesic flow on $E =
E_0$ with invariants $h = 1$ (say) and $\lambda_i = a_{n+1+i}$,
$1 \leq i \leq n$ can be identified with the real locus in
$\widetilde{K} = \widetilde{K}(\lambda_1,\cdots,\lambda_n)$.
The latter is a $2^{n+1}$-sheeted cover of $Pic^1(C)$, so up to translation by
some points of order $2$, it is an abelian variety, isomorphic to a
$2^{n+1}$ sheeted cover of the hyperelliptic Jacobian $J(C)$.


\subsection{The flows}

     Two details of the above story are somewhat unsatisfactory:
First, the asymmetry between the $n-1$ Hamiltonians $\lambda_i$
and the remaining Hamiltonian $H$ (length squared).
And second, the fact that the complexified
total space $TE$ of the system is not quite symplectic.  Indeed,
for an arbitrary algebraic hypersurface
$M \subset \C^{n+1}$, given by $f = 0$, the complexified metric on
$\ComplexNumbers^{n+1}$ induces bundle maps
$$
TM \hookrightarrow T\C^{n+1}|_M
\stackrel{\sim}{\rightarrow} T^*\C^{n+1}|_M
\rightarrow\!\!\!\rightarrow T^{*}M,
$$
but the composition is not an
isomorphism; rather, it is degenerate at points where
$$0 =
(\bigtriangledown f)^2 = \sum^{n+1}_{k=1}
\left(\frac{\partial f}{\partial
x_{k}}\right)^2.
$$
For an ellipsoid $\sum x_k^2/a_k = 1$
(other than a sphere) there will be
an empty real, but non-empty complex degeneracy locus, given by the equation
$\sum(x_k/a_k)^2 = 0$.

     Both of these annoyances disappear if we replace the total
space by the tangent bundle $TS$ of the sphere
$$S =
\{(x_1,\cdots,x_{n+1}) \in \C^{n+1}|\sum x_{k}^2 = 1\},
$$
i.e.
$$TS = \{(x,y) \in \C^{2n+2}| \sum x_{k}^2 = 1, \sum
x_k y_k = 0\}.
$$
This is globally symplectic, and the $n$
(unordered) Commuting Hamiltonians can be taken to be the values
$\lambda_i$, $1 \leq i \leq n$, such that the line $$\ell_{x,y}
:= (\mbox{line\ through}\ y\ \mbox{in\ direction}\ x\ )$$ is
tangent to $E_{\lambda_i}$.  The original system $TE$ can be
recovered as a $\C^*$-bundle (where $\C^*$ acts by
rescaling the tangent direction $x$) over the hypersurface
$\lambda = 0$ in $TS$.

     Here is the explicit equation of the hypersurface:
$$\begin{array}{rcl}
\lambda = 0 & \Leftrightarrow & \ell_{x,y} \ \
\mbox{is\ tangent\ to}\ E = \{\sum \frac{x_{k}^2}{a_k} = 1\}
\\
& \Leftrightarrow & -1+\sum \frac{(y_k + tx_k)^2}{a_k} = 0\
\mbox{has\ a\ unique\ solution}\ t
\\
& \Leftrightarrow &  (-
1+\sum \frac{y_{k}^2}{a_k} + 2t (\sum \frac{x_k y_k}{a_k}) +
t^2(\sum \frac{x_{k}^2}{a_k}) = 0\ \mbox{has\ a\ unique\
solution}
\\
& \Leftrightarrow & 0 = (\sum \frac{x_k
y_k}{a_k})^2 - (\sum \frac{x_{k}^2}{a_k})(-1 + \sum
\frac{y_{k}^2}{a_k}).
\end{array}
$$
More generally, this computation shows that $\ell_{x,y}$ is tangent
to $E_{\lambda}$ if and only if
\[
0 = (\sum \frac{x_k y_k}{a_{k^{-\lambda}}})^2 - (\sum
\frac{x_{k}^2}{a_{k^{-\lambda}}})(-1 + \sum \frac{y_{k}^2}{a_{k^-
\lambda}})
=
\sum_k \frac{x_k^2}{a_k-\lambda} +
\sum_{k\neq l}\frac{x_k y_k x_l y_l- x_k^2 y_l^2}{(a_k-\lambda)(a_l-\lambda)}.
\]
As a function of $\lambda$, the last expression has first order poles at
$\lambda=a_k$, $1\leq k \leq n+1$, so it can be rewritten as
\[
\sum \frac{1}{a_k - \lambda} F_k(x,y)
\]
where the $F_k$ are found by taking residue at $\lambda=a_k$:
\[
F_k (x,y) := x_{k}^2 + \sum_{\ell
\neq k} \frac{(x_k y_\ell - x_\ell y_k)^2}{a_k - a_\ell}.
\]
We see that fixing the $n+1$ values $F_k(x,y)$, $1
\leq k \leq n+1$, subject to the condition
$$
\sum^{n+1}_{k=1} F_k
= 1,
$$
is equivalent to fixing the $n$ (unordered) values
$\lambda_i$, $1 \leq i \leq n$.

     This determines the hyperelliptic curve
$$
C : \; s^2 =
\prod^{n+1}_{k=1} (t - a_k) \cdot \prod^{n}_{i=1} (t -
\lambda_i),
$$
and the corresponding abelian variety
$$
\widetilde{K} = \widetilde{K} (\lambda_1,\cdots,\lambda_n) =
J(C)/\widetilde{G} \approx \{(x,y)| \ell_{x,y}\ \mbox{is\ tangent\ to} \
E_{\lambda_i}, \ \ \ 1 \leq i \leq n\}.
$$

\vspace{0.1in}

\proclaim{Theorem}.
\begin{list}{{\rm(\arabic{bean})}}{\usecounter{bean}}
\item Geodesic flow on the quadric $E_{\lambda_i}$ is the
Hamiltonian vector field on $TS$ given by the (local) Hamiltonian
$\lambda_i$.  On $\widetilde{K}$ it is a constant vector field in
the direction of the Weierstrass point $\lambda_i \in C$.
\item The Hamiltonian vector field on $TS$ with Hamiltonian $F_k$
is constant on $\widetilde{K}$, in the direction of the
Weierstrass point $a_k \in C$.
\end{list}
\par

     The direction at $\ell \in \K$ of geodesic flow on
$E_{\lambda_i}$ was described in the proof of Chasles' theorem.
The direction given by the Weierstrass point $\lambda_i$ is given
at $A \in F$ as the tangent vector at $\lambda_i$ to the curve
$$
p \mapsto j(A_{i},p)\  \ (\mbox{where}\ A_i =
j(A,\lambda_i)).
$$
The proof of (1) amounts to unwinding the
definitions to see that these two directions agree.  (For
details, see \cite{knorrer} and
\cite{donagi-group-law}.)  Since the level sets of the
$\lambda_i$ and the $F_k$ are the same, the Hamiltonian vector
field of $F_k$ evolves on the same $\K$, and is constant there.
A monodromy argument on the family of hyperelliptic curves then
shows that its direction must agree with $a_k$.  Mumford gives an
explicit computation for this in \cite{Mum}, Theorem 4.7, following
Moser  \cite{moser}.

     We have identified the flows corresponding to $2n+1$ of the
Weierstrass points.  The remaining one, at $\lambda = \infty$,
corresponds to the Hamiltonian
$$
H = \frac{1}{2} \sum^{n+1}_{k=1}
a_k F_k = \frac{1}{2} \sum a_k x_{k}^2 + \frac{1}{2} \sum
y_{k}^2,
$$
giving Neumann's system, which is the starting point
for the analysis in \cite{moser} and
\cite{Mum}.


\subsection{Explicit parametrization}

     Fix the hyperelliptic curve of genus $n$
$$
C: \; s^2 = f(t) :=
\prod^{2n+1}_{i=1} (t - a_i),
$$
with projection
$$
\begin{array}{rcl}\pi:\ C & \rightarrow & \P^1\\(t,s) & \mapsto
& t\\ \end{array}
$$
and involution
$$
\begin{array}{rcl}i:\ C &
\rightarrow & C\\(t,s) & \mapsto & (t, - s).\\
\end{array}
$$
We identify the various components $Pic^d(C)$ by means of the base
point $\infty$, which is a Weierstrass point.
The affine open subset $J(C) \smallsetminus \Theta$ can
be described geometrically, by Riemann's theorem:
$$
J(C) \smallsetminus \Theta \approx
$$
$$\{L \in Pic^{n-1}(C)|h^0(L) = 0\}
\approx
$$
$$
\{(p_1,\cdots,p_n) \in Sym^nC| p_i \neq \infty, p_i
\neq i(p_j)\},
$$
where the last identification sends $L$ to the
unique effective divisor of $L(\infty)$.  Mumford \cite{Mum} gives an
explicit algebraic parametrization of the same open set, which he
attributes to Jacobi: to the $n$-tuple $D = (p_1,\cdots,p_n)$ he
associates three polynomials of a single variable $t$:
\begin{list}{{\rm(\roman{bacon})}}{\usecounter{bacon}}
\item $U(t) := \prod^{n}_{i=1} (t - t(p_i))$.
\item $V(t)$ is the unique polynomial of degree $\leq n-1$ such
that the meromorphic function
$$
V \circ \pi - s \ \ : \ \ C \to \P^1
$$
vanishes on the divisor $D \subset C$.  It is obtained by
Lagrange interpolation of the expansions of $s$ at the $p_i$,
e.g. when the $p_i$ are all distinct,
$$
V(t) = \sum^{n}_{i=1}
s(p_i) \prod_{j \neq i} \frac{t - t(p_j)}{t(p_i) - t(p_j)}.
$$
\item $W(t) = \frac{f(t) - V(t)^2}{U(t)}$; the definition of $V$
and the equation $s^2 = f(t)$ guarantee that this is a monic
polynomial of degree $n+1$.
\end{list}

     Conversely, the polynomials $U,V,W$ determine the values
$t(p_i), s(p_i)$, hence the divisor $D$.  By reading off the
coefficients, we obtain an embedding: $$(U,V,W): \J(C)\backslash
\Theta \hookrightarrow \C^{3n+1}.$$ The image is $$\{(U,V,W)| V^2
+ UW = f\}.$$

     This description fits beautifully with the integrable system
on $TS$ representing geodesic flow on the $E_\lambda$.  We can
rephrase our previous computation as: $$\ell_{x,y} \ \mbox{is\
tangent\ to}\ E_{\lambda_1},\cdots,E_{\lambda_n} \Leftrightarrow
f_1(t)f_2(t) = UW + V^2,$$ where, for $(x,y) \in TS$, we set:
$$
\begin{array}{lcl}
f_1(t) = \prod^{n+1}_{k=1}(t-a_k) & \ \ \ \ & \mbox{(this\ is\
independent\ of}\ x,y)
\\
f_2(t) = f_1(t) \cdot \sum_k \frac{F_k(x,y)}{t - a_k}
& \ \ \ \ & (\mbox{this\ varies\ with}\ x,y;\ \mbox{the\ roots\
are\ the}\ \lambda_i)
\\
U(t) = f_1(t) (\sum_k \frac{x_{k}^2}{t - a_k}) & &
\\
W(t) = f_1(t) (1 + \sum_k \frac{y_{k}^2}{t - a_k}) & &
\\
V(t) = \sqrt{-1} \cdot f_1(t) \cdot (\sum_k \frac{x_k y_k}{t -
a_k}) & &\\
\end{array}
$$

     The entire system $TS$ is thus mapped to $\C^{3n+1}$.
Each abelian variety $\K = \K(\lambda_1,\cdots,\lambda_n)$ is
mapped to $\J(C) \backslash \Theta$ embedded in $\C^{3n+1}$ as
before, where $C$ is defined by $s^2 = f(t)$, and $f = f_1 \cdot
f_2$, with $f_1$ fixed (of degree $n+1$) and $f_2$ variable (of
degree $n$).  On each $\K$ the map is of degree $2^{n+1}$; the
group $(\Z/2\Z)^{n+1}$ operates by sending
$$
(x_k,y_k)
\mapsto (\epsilon_k x_k, \epsilon_k y_k),\ \ \ \ \epsilon_k =
\pm 1.
$$
\newpage


\section{Spectral curves and vector bundles} \label{ch4}

We review in this chapter a general construction
of an integrable system on the moduli space of Higgs
pairs $(E,\varphi)$ consisting of a vector bundle $E$ on a curve and
a meromorphic $1$-form
valued endomorphism $\varphi$ (theorem \ref{thm-markman-botachin}).
These moduli
spaces admit a natural foliation by Jacobians of spectral curves.
The spectral curves are branched covers of the base
curve arising from the eigenvalues of the endomorphisms $\varphi$.

We concentrate on two examples:

\smallskip
\noindent
- The Hitchin system supported on the cotangent bundle of
the moduli space of vector bundles on a curve (section
\ref{sec-spectral-curves-and-the-hitchin-system}), and

\noindent
- An integrable system on the moduli space of conjugacy classes of
polynomial matrices (section \ref{sec-polynomial-matrices}).

\smallskip
The latter is then used to retrieve the Jacobi-Moser-Mumford system
which arose in chapter \ref{ch3} out of the geodesic flow on an ellipsoid.

Both examples are endowed with a natural symplectic or Poisson structure.
The general construction of the Poisson structure on the moduli spaces
of Higgs pairs is postponed to chapter \ref{ch5}.

We begin with a short survey of some basic facts about vector bundles on
a curve.

\subsection{Vector Bundles on a Curve}
\label{sec-vector-bundles-on-a-curve}

We fix a (compact, non-singular) curve $\Sigma$ of genus $g$.  A basic
object in these lectures will be the moduli space of stable (or semistable)
vector bundles on $\Sigma$ of given rank $r$ and degree $d$.  To motivate
the introduction of this object, let us try to describe a ``general''
vector bundle on $\Sigma$.  One simple operation which produces vector
bundles from line bundles is the direct image:  start with an $r$-sheeted
branched covering $\pi : C \rightarrow \Sigma$, ramified at points of some
divisor $R$ in the non-singular curve $C$.  Then any line bundle $L \in
Pic\; C$ determines a rank-$r$ vector bundle $E := \pi_*L$ on $\Sigma$.  As a
locally free sheaf of ${\cal O} _\Sigma$-modules of rank $r$, this is easy to
describe
$$
\Gamma({\cal U},\pi_*L) \ := \ \Gamma(\pi^{-1} {\cal U}, L),
$$
for open subsets ${\cal U} \subset \Sigma$.  As a vector bundle, the
description is clear only at unbranched points of $\Sigma$: \ if
$\pi^{-1}(p)$ consists of $r$ distinct points $p_1, \cdots , p_r$ then the
fiber of $E$ at $p$ is naturally isomorphic to the direct sum of the fibers
of $L$:
$$
E_p \; \approx \; \bigoplus^r_{i=1} \; L_{p_i}.
$$
At branch points of $\pi$, $E_p$ does not admit a natural decomposition,
but only a filtration.  This is reflected in a drop in the degree.  Indeed,
the Grothendieck-Riemann-Roch theorem says in our (rather trivial) case that
$$
\chi(\pi_*L) = \chi(L),
$$
where $\chi$ is the holomorphic Euler characteristic,
$$
\chi(E) := \deg E - (g-1) \rank \; E = \deg E + \chi({\cal O}) \cdot \rank
\; E.
$$
Using Hurwitz' formula:
$$
\chi({\cal O}_C) = r \, \chi({\cal O}_\Sigma) - {1 \over 2} \deg \; R,
$$
this becomes:
$$
\deg \; E = \deg \; L - {1 \over 2} \deg \; R.
$$
\begin{example} \label{example-direct-image}
{\rm Consider the double cover
\begin{eqnarray*}
\pi : \bP^1 &\longrightarrow& \bP^1 \\
w &\longrightarrow& z = w^2 \\
\end{eqnarray*}
branched over $0,\infty$.  The direct image of the structure sheaf is:
$$
\pi_* {\cal O} \approx {\cal O} \oplus {\cal O}(-1).
$$
We can think of this as sending a regular function $f = f(w)$ (on some
invariant open set upstairs) to the pair $(f_+(z), f_-(z))$ downstairs,
where
$$
f(w) = f_+(w^2) \; + \; wf_-(w^2).
$$
In the image we get all pairs with $f_+$ regular (i.e. a section of $\cal O$)
and $f_-$ regular and vanishing at $\infty$ (i.e. a section of ${\cal
O}(-1))$. (Similar considerations show that
$$
\pi_* {\cal O}(-1) \; \approx \; {\cal O}(-1) \oplus  {\cal O}(-1)
$$
and more generally:
$$
\pi_* {\cal O}(d) \approx {\cal O} ( [ {d \over 2} ]) \oplus {\cal O}([ {d-1
\over 2} ]).
$$
Note that this has degree $d-1$, as expected).
The structure of $\pi_*L$ near the branch point $z = 0$ can be described,
in this case, by the action on the local basis $a,b$ (of even, odd
sections) of multiplication by the section $w$ upstairs:
$$
a \mapsto b , \quad  \quad b \mapsto za,
$$
i.e. $w$ is represented by the matrix
$$
\left(
\begin{array}{cc}
0 & z \\
1 & 0 \\
\end{array}
\right)
$$
whose square is  $z \cdot I$.  At a branch point where $k$ sheets come
together, the corresponding action (in terms of a basis indexed by the
$k-th$ roots of unity) is given by the matrix:
\begin{equation}\label{eq-ramification-matrix}
P_k : =
\left(
\begin{array}{cccccc}
0      &        &          &          & 0      & z      \\
1      &        &          &          & \cdot  & \cdot  \\
\cdot  & \cdot  &          &          & \cdot  & \cdot  \\
\cdot  &        & \cdot    &          & \cdot  & \cdot  \\
\cdot  &        &          & \cdot    & 0      & 0      \\
0      &        &          &          & 1      & 0
\end{array}
\right)
\end{equation}
whose $k$-th power is $z \cdot I$.}
\end{example}

\begin{example} \ \ {\rm Now consider a 2-sheeted branched cover $\pi : C
\rightarrow \bP^1$, where $g(C) > 0$.  If we take $\chi(L) = 0$, i.e. $\deg
L = g - 1$, we get $\deg(\pi_*L) = -2$.  The equality
$$
\ell := h^0(C,L) = h^0(\bP^1, \pi_*L)
$$
implies
$$
\pi_*L \approx {\cal O}(\ell - 1) \; \oplus \; {\cal O}(-\ell - 1).
$$
In particular, we discover a very disturbing phenomenon: \ \ as the line
bundle $L$ varies continuously, in $Pic^{g-1}C$, so should presumably
$\pi_*L$; but if we consider a 1-parameter family of line bundles $L_t$
such that
\begin{eqnarray*}
L_0 &\in& \Theta \\
L_t &\notin& \Theta, \ \ \   t \not= 0 , \\
\end{eqnarray*}
we see that the vector bundle $\pi_*L_t$ jumps from its generic value,
${\cal O}(-1) \oplus {\cal O}(-1)$ to ${\cal O} \oplus {\cal O}(-2)$ at
$t=0$.  Similar jumps can clearly be forced on a rank-$r$ bundle by
considering $r$-sheeted branched covers.}
\end{example}
The moral of these examples is that if we want a moduli space parametrizing
the ``general'' vector bundle on a curve and having a reasonable (say,
separated) topology, we cannot consider {\it all} bundles.  In the case of
$\bP^1$, we will end up with only the balanced bundles such as ${\cal
O}(-1) \oplus {\cal O}(-1)$, thus avoiding the possibility of a
discontinuous jump.

The slope $\mu(E)$ of a vector bundle $E$ is defined by:
$$
\mu(E) := {\deg E \over \rank E}.
$$
A bundle $E$ is called stable (resp., semistable) if for every subbundle $F
\subset E$ (other than $0,E$),
$$
\mu(F) < \mu(E), \quad \quad (resp. \ \mu(F) \le \mu(E)).
$$
The basic result due to Mumford and Seshadri
\cite{seshadri-construction-moduli-vb},
is that reasonable (coarse) moduli
spaces ${\cal U}^s_\Sigma(r,d) \subset {\cal U}_\Sigma(r,d)$ exist, with
the following properties:

\begin{itemize}
\item ${\cal U}^s_\Sigma(r,d)$ is smooth; its points parametrize
isomorphism classes of stable bundles of rank $r$ and degree $d$ on
$\Sigma$; it is an open subset of ${\cal U}_\Sigma(r,d)$.
\item ${\cal U}_\Sigma(r,d)$ is projective; its points parametrize
equivalence classes of semistable bundles, where two bundles are
equivalent, roughly, if they admit filtrations by semistable subbundles
(of constant slope) with isomorphic graded pieces.
\item Both are coarse moduli spaces; this means that any ``family'', i.e.
vector bundle on a product $S \times \Sigma$, where $S$ is any scheme,
whose restrictions $E_s$ to copies $s \times \Sigma$ of $\Sigma$ are (semi)
stable of rank $r$ and degree $d$, determines a unique morphism of $S$ to
${\cal U}^s_\Sigma(r,d)$ (respectively, ${\cal U}_\Sigma(r,d))$ which
sends each $s \in S$ to the isomorphism (resp. equivalence) class of
$E_s$, and has the obvious functoriality properties.
\end{itemize}

\noindent
\underline{Examples}

\noindent
\underline{$g = 0$}. \ The stable bundles are the line bundles ${\cal
O}(d)$. The semi-stable bundles are the balanced vector bundles, ${\cal
O}(d)^{\oplus r}$.  Thus ${\cal U}_{P^1}(r,d)$ is a point if $r\mid d$, empty
otherwise, while the stable subset is empty when $r \ne 1$.

\smallskip
\noindent
\underline{$g = 1$}. \ Let $h := gcd(r,d)$. Atiyah
\cite{atiyah-vb-on-elliptic-curves}
shows that ${\cal U}_\Sigma(r,d)$
is isomorphic to the symmetric product $S^h\Sigma$, and
that each semistable equivalence class contains a unique decomposable
bundle $E = \oplus^h_{i=1} \; E_i$, where each $E_i$ is stable of rank
$r/h$ and degree $d/h$.  (Other bundles in this equivalence class are
filtered, with the $E_i$ as subquotients.)  Thus when $h=1, \  {\cal
U}^s_\Sigma = {\cal U}_\Sigma$, and when $h > 1,\  {\cal U}^s_\Sigma$
is empty.

The possibilities for semistable bundles are illustrated in the case $r=2$,
$d=0$:  given two line bundles $L_1,L_2 \in Pic^0 \  \Sigma$, the possible
extensions are determined, up to non zero scalars, by elements of
$$
Ext^1_{{\cal O}_\Sigma}(L_1,L_2) \approx H^1(L_2 \otimes L^{-1}_1).
$$
The direct sum is thus the only extension when $L_1 \not\approx L_2$, while
if $L_1 \approx L_2 \approx L$ there is, up to isomorphism, also a unique
non-trivial extension, say $E_L$.  There is, again, a jump phenomenon: \ \
by rescaling the extension class we get a family of vector bundles with
generic member isomorphic to $E_L$ and special member $L \oplus L$.  This
explains why there cannot exist a moduli space parametrizing {\it
isomorphism} classes of semistable bundles; neither $L \oplus L$ nor $E_L$
is excluded, and the point representing the former is in the closure of the
latter, so they must be identified, i.e. $E_L$ and $L \oplus L$ must be
declared to be equivalent.

\smallskip
\noindent
\underline{Higher Genus}. \ \  The only other cases where an explicit
description of ${\cal U}_\Sigma(r,d)$  is known are when $r=2$ and $g=3$
\cite{narasimhan-ramanan-rk2-genus3} or $r=2$ and $\Sigma$ is hyperelliptic of
any genus \cite{DR}.
In the latter case, the moduli space ${\cal U}_\Sigma(2,\xi)$ of rank $2$
vector bundles with a fixed determinant line bundle $\xi$ of odd degree
is isomorphic to the family of linear
spaces $\bP^{g-2}$ in the intersection of the two quadrics in $\bP^{2g+1}$
used in Chapter \ref{ch3}.  In the even degree case,  ${\cal U}_\Sigma(2,\xi)$
can also be described in terms
of the same two quadrics; when $g = 2$, it turns out to be isomorphic to
$\bP^3$, in which the locus of semistable but non-stable points is the
Kummer surface  $K := {\cal J}(\Sigma) / \pm 1$, with its classical
embedding in $\bP^3$ as a quadric with 16 nodes.

\bigskip
Elementary deformation theory lets us make some general statements about
${\cal U} := {\cal U}_\Sigma(r,d)$ and ${\cal U}^s :=
{\cal U}^s_\Sigma(r,d)$:

\begin{lem}  For $g \ge 2$:
\begin{itemize}
\begin{enumerate}
\item dim ${\cal U} = 1 + r^2(g-1)$, and ${\cal U}^s$ is a dense open
subset.
\item Stable bundles $E$ are simple, i.e. the only (global) endomorphisms
of $E$ are scalars.
\item Stable bundles are non-singular points of ${\cal U}$.
\item At points of ${\cal U}^s$ there are canonical
identifications
\end{enumerate}
\end{itemize}
\end{lem}
\begin{eqnarray*}
T_E{\cal U}^s &\approx& H^1(End \; E) \\
T^*_E{\cal U}^s &\approx& H^0(w_\Sigma \otimes \; End \; E). \\
\end{eqnarray*}

The proof of (2) is based on the observation that any nonzero $\alpha : E
\rightarrow E$ must be invertible, otherwise either ker $\alpha$ or im
$\alpha$ would violate stability.  Therefore $H^0(End \; E)$ is a finite
dimensional division algebra containing $\C$, hence equal to it.  Since a
vector bundle $E$ on $\Sigma$ is determined by a $1$-cocycle with values in
$GL(r,{\cal O}_\Sigma)$ (= transition matrices), a first order deformation
of $E$ is given by a 1-cocycle with values in the associated bundle of Lie
algebras, i.e. (up to isomorphism) by a class in $H^1(End \; E)$.  The
functoriality property of $\cal U$ (``coarse moduli space'') implies that
this is the Zariski tangent space, $T_E{\cal U}$.  By Riemann-Roch
$$
h^1 ({\rm End} \; E) = r^2(g-1) + h^0({\rm End} \; E).
$$
so the minimal value is obtained at the simple points, and equals $1 +
r^2(g-1)$ as claimed in (1).  The identification of $T^*_E{\cal U}^s$ follows
from
that of $T_E{\cal U}^s$ by Serre duality.
\EndProof

\subsection{Spectral Curves and the Hitchin System}
\label{sec-spectral-curves-and-the-hitchin-system}

The relation between vector bundles and finite dimensional integrable systems
arises from Hitchin's amazing result.

\begin{theorem}  \label{thm-hitchins-integrable-system}
\cite{hitchin,hitchin-integrable-system}
The cotangent bundle to the moduli space
of semistable vector bundles supports a natural ACIHS.
\end{theorem}

%

At the heart of Hitchin's theorem
is a construction of a spectral curve
associated to a $1$-form valued endomorphism of a vector bundle.  The
spectral construction allows a uniform treatment of a wide variety of
algebraically completely integrable Hamiltonian systems.  We will
concentrate in this section on the algebro-geometric aspects of these
systems leaving their symplectic geometry to Chapter \ref{ch8}.  We work with
vector bundles over curves, other structure
groups will be treated in Chapter \ref{ch9}. The
reader is referred to \cite{B-N-R} and \cite{hitchin-integrable-system}
for more details.

The total space of the cotangent bundle $T^*{\cal U}^s_\Sigma(r,d)$ of the
moduli space of stable vector bundles parametrizes pairs $(E,\varphi)$
consisting of a stable vector bundle $E$ and a covector $\varphi$ in
$H^1(\Sigma, {\rm End} \; E)^* \simeq H^0(\Sigma, {\rm End} E \otimes
\omega_\Sigma)$, i.e., a $1$-form valued endomorphism of $E$.

Consider more generally a pair $(E,\varphi)$ of a rank $r$ vector bundle $E$
and a section $\varphi \in {\rm Hom}(E, E \otimes K)$  where  $K$ is a line
bundle on $\Sigma$.  The $i$-th coefficient $b_i$ of the characteristic
polynomial  of $(E,\varphi)$ is a homogeneous polynomial of degree
$i$ on $K^{-1}$, hence a section of $H^0(\Sigma, K^{\otimes i})$.  In fact
$b_i = (-1)^i\cdot {\rm trace}(\stackrel{i}{\wedge}\varphi)$.

The Hamiltonian map of the Hitchin system is the characteristic polynomial
map
$$
H : T^* {\cal U}_\Sigma(r,d) \longrightarrow B_\omega := \bigoplus^r_{i=1}
H^0(\Sigma, \omega^{\otimes i}).
$$
The fibers of the Hitchin map $H$ turn out to be Jacobians of curves
associated canonically to characteristic polynomials.

Going back to the general $K$-valued pair $(E,\varphi)$, notice that a
characteristic polynomial \  char $(\varphi) = y^r - {\rm tr}(\varphi)
y^{r-1} + \cdots + (-1)^r\det \varphi$ in $B_K := \oplus^r_{i=1} H^0(\Sigma,
K^{\otimes i})$ defines a morphism from the line bundle $K$ to $K^{\otimes
r}$.  The inverse image $C$ of the zero section in $K^{\otimes r}$ under a
polynomial $P$ in $B_K$ is called a spectral curve.  If $P$ is the
characteristic polynomial of a pair $(E,\varphi)$ then indeed the fibers
of $\pi : C \rightarrow \Sigma$ consist of eigenvalues of $\varphi$.  If
$K^{\otimes r}$ has a section without multiple zeroes (e.g., if it is very
ample) then the generic spectral curve is smooth.

Lagrange interpolation extends a function on the inverse image $\pi^{-1}(U)
\subset C$ of an open set $U$ in $\Sigma$ to a unique function on the inverse
image of $U$ in the surface $K$ which is a polynomial of degree $\le r-1$ on
each
fiber.  It follows that the direct image $\pi_*{\cal O}_C$ is isomorphic
to ${\cal O}_\Sigma \oplus K^{-1} \oplus \cdots \oplus K^{1-r}$.  Assuming
that $K^{\otimes i}$
has no sections for $i < 0$,
the genus $h^1(C,{\cal O}_C) = h^1(\Sigma, \pi_*{\cal O}_C)$
of $C$ is equal to $\deg(K) \cdot r(r-1)/2 + r(g-1) + 1$.
In particular, when $K = \omega_\Sigma$, the
genus of $C$ is equal to half the dimension of the cotangent bundle.  The
data $(E,\varphi)$ determines moreover a sheaf $L$ on the spectral curve
which is a line bundle if the curve is smooth.  Away from the ramification
divisor $R$ in $C$, $L$ is the tautological eigenline subbundle of the
pullback $\pi^*E$.  More precisely, the homomorphism $(\pi^*(\varphi) - y
\cdot I) : \pi^*E \rightarrow \pi^*(E \otimes K)$, where  $y \in
H^0(C, \pi^*K)$ is the tautological eigenvalue section, has kernel
$L(-R)$.

Conversely, given a spectral curve $C$ and a line bundle $L$ on it we get a
pair $(\pi_*L, \pi_*(\otimes y))$ of a rank $r$ vector bundle on $\Sigma$
and a $K$ valued endomorphism (see example
\ref{example-direct-image}).
The two constructions are the inverse of
each other.

\begin{proposition} \label{prop-ordinary-spectral-construction-higgs-pairs}
\cite{hitchin-integrable-system,B-N-R}
If $C$ is an irreducible and reduced
spectral curve there is a bijection between isomorphism classes of
\begin{description}
\item [-] Pairs $(E,\varphi)$ with spectral curve $C$.
\item [-] Rank 1 torsion free sheaves $L$ on $C$.
\end{description}
\end{proposition}
Under this correspondence, line bundles on $C$ correspond to endomorphisms
$\varphi$ which are regular in every fiber, i.e., whose centralizer in
each fiber is an $r$-dimensional subspace of the corresponding fiber of End
$E$. (This notion of regularity agrees with the one in Example
\ref{example-coadjoint-orbits}.)

We conclude that the fiber of the Hitchin map $H : T^*U^s_\Sigma(r,d)
\rightarrow B_\omega$ over a characteristic polynomial $b \in B_\omega$ is
precisely the open subset of the Jacobian $J_C^{d+r(1-g_\Sigma)+g_C-1}$
consisting of
the line bundles $L$ whose direct image is a stable vector bundle.
Moreover, the construction of the characteristic polynomial map and a
similar description of its fibers applies to moduli spaces of pairs with
$K$-valued endomorphism where $K$ need not be the canonical line bundle
(Theorem \ref{thm-markman-botachin}).

The missing line bundles in the fibers of the Hitchin map indicate that we
need to relax the stability condition for the pair $(E,\varphi)$.

\smallskip
\noindent
{\bf Definition:} \ \ A pair $(E,\varphi)$ is stable (semistable) if the
slope of every $\varphi$-invariant subbundle of $E$ is less than (or equal)
to the slope of $E$.

As in the case of vector bundles we can define an equivalence relation for
semistable pairs, where two bundles are equivalent, roughly, if they admit
$\varphi$-invariant filtrations by semistable pairs (of constant slope)
with isomorphic graded pieces.  Two stable pairs are equivalent if and only
if they are isomorphic.

\begin{theorem} \cite{hitchin,simpson-moduli,Nit}
There exists an algebraic coarse moduli
scheme $\HiggsModuli_K := \HiggsModuli_\Sigma(r,d,K)$ parametrizing equivalence
classes of semistable $K$-valued pairs.
\end{theorem}

The characteristic polynomial map $H : \HiggsModuli_K \rightarrow B_K$
is a proper
algebraic morphism.

A deeper reason for working with the above definition of stability is
provided by the following theorem from nonabelian Hodge theory:

\begin{theorem} \label{thm-higgs-pairs-and-representations-of-pi1-for-curves}
\cite{hitchin,simpson-higgs-bundles-and-local-systems}
\ There is a canonical real analytic
diffeomorphism between
\begin{description}
\item[-] The moduli space of conjugacy classes of semisimple
representations of the fundamental group $\pi_1(\Sigma)$ in $GL(r,\C)$ and
\item[-] The moduli space of semistable $\omega$-valued (Higgs) pairs
$(E,\varphi)$ of rank $r$ and degree $0$.
\end{description}
\end{theorem}

In the case of Hitchin's system $(K = \omega_\Sigma)$, the symplectic
structure of the cotangent bundle extends to the stable locus of the moduli
space of Higgs pairs giving rise to an integrable system $H :
\HiggsModuli_\Sigma(r,d,\omega_\Sigma) \rightarrow B_\omega$ whose generic
fiber is a complete Jacobian of a spectral curve.

We will show in Chapter \ref{ch6} that the Hitchin system is, in fact, the
lowest
rank symplectic leaf of a natural infinite dimensional Poisson variety
$\HiggsModuli_\Sigma(r,d)$ obtained as an inductive limit of the moduli spaces
$\HiggsComponent_\Sigma(r,d,\omega(D))$  of $\omega(D)$-valued pairs as $D$
varies through all effective divisors on $\Sigma$.
The basic fact, generalizing
the results of \cite{hitchin-integrable-system,B-N-R,B} is:

\begin{theorem} \label{thm-markman-botachin} \cite{botachin,markman-higgs}  \
Let
$D$ be an effective divisor  (not necessarily reduced)
on a smooth algebraic curve $\Sigma$ of genus
$g$.  Assume that $[\omega(D)]^{\otimes r}$ is very ample and if $g=0$
assume further that $\deg(D) > \max(2,\rho)$ where $0 \le \rho < r$ is the
residue of $d$ mod $r$.  Then
\begin{itemize}
\begin{enumerate}
\item The moduli space $\HiggsModuli^s_\Sigma(r,d,\omega(D))$
of stable rank $r$ and
degree $d$ \ $\omega(D)$-valued Higgs pairs has a smooth component
$\HiggsComponent_\Sigma(r,d,\omega(D))$ of top dimension
 $r^2(2g - 2 + \deg(D)) + 1 + \epsilon_{D=0}$, where $\epsilon_{D=0}$ is $1$
if $D=0$ and zero if $D > 0$.
$\HiggsComponent_\Sigma(r,d,\omega(D))$ is the unique component which
contains Higgs pairs supported on irreducible and reduced spectral curves.
\item $\HiggsComponent_\Sigma(r,d,\omega(D))$ has a canonical Poisson
structure.
\item The characteristic polynomial map
$H : \HiggsComponent_\Sigma(r,d,\omega(D))
\rightarrow B_{\omega(D)}$ is an algebraically completely integrable
Hamiltonian system.  The generic (Lagrangian) fiber is a complete Jacobian
of a smooth spectral curve of genus $r^2(g-1) + 1 + (\deg D)( {r(r-1) \over
2})$.
\item The foliation of $\HiggsComponent_\Sigma(r,d,\omega(D))$
by closures of top dimensional symplectic leaves  is induced by the cosets of
$$
H^0 \left( \Sigma , \left[ \bigoplus^r_{i=1} \omega_\Sigma(D)^{\otimes i}
\right] (-D) \right) \; {\rm in} \; B_{\omega(D)}.
$$
\end{enumerate}
\end{itemize}
\end{theorem}

\begin{definition} \label{def-good-component-of-higgs-pairs}
As in the theorem, we will denote by $\HiggsComponent_\Sigma(r,d,\omega(D))$
the unique component which
contains Higgs pairs supported on irreducible and reduced spectral curves.
\end{definition}

In Chapter \ref{ch6} we will discuss the relationship of these integrable
systems
with flows of KdV type.  In the next section we will discuss the example of
geodesic flow on the ellipsoid as a Hamiltonian flow of a symplectic leaf
of one of these spaces. See \cite{B,markman-higgs} for more examples.

The Hitchin system has been useful in the study of the geometry of the
moduli space of vector bundles. The main technique is to reduce questions
about vector bundles to questions about spectral Jacobians.  Hitchin used
these ideas to compute the cohomology groups $H^i({\cal U}, S^kT)$, $i=0,1$,
of the symmetric products of the tangent bundle of the moduli space
$\cal U$ of rank $2$ and odd degree stable vector bundles.
In \cite{B-N-R} these techniques
provided the first mathematical proof that the dimensions of the space of
sections of the generalized theta line bundle are
\begin{eqnarray*}
h^0({\cal U}_\Sigma(n, \, n(g-1)), \; \Theta) &=& 1, \\
h^0({\cal SU}_\Sigma(n), \; \Theta) &=& n^g, \\
\end{eqnarray*}
where ${\cal SU}_\Sigma(n)$ denotes
the moduli space of vector bundles with trivial determinant line bundle.
(This of course is now subsumed in the Verlinde Formula for sections
of powers of theta bundles.)
These ideas were proven useful in
the proof of the existence of a projectively flat connection on the bundles
of level $k$ theta sections over the moduli space ${\cal M}_g$ of curve of
genus
$g$ \cite{hitchin-flat-connection}, an important fact in conformal field
theory.  Kouvidakis and Pantev applied these ideas to the study of
automorphisms of the moduli
space of vector bundles \cite{kouvidakis-pantev}.

\subsection{  Polynomial Matrices}
\label{sec-polynomial-matrices}

Theorem \ref{thm-markman-botachin} has a concrete
description when the base curve $\Sigma$ is
$\bP^1$.  Let $K$ be the line bundle ${\cal O}_{\bP^1}(d)$.  Consider the
moduli space $\HiggsModuli_K := \HiggsComponent_K(-r,r)$ of pairs
$(E,\varphi)$ consisting of a vector bundle $E$ of rank $r$ and degree $-r$
with a $K$-valued endomorphism $\varphi : E \rightarrow E \otimes K$
(we also follow the notation of definition
\ref{def-good-component-of-higgs-pairs} singling out a particular component).
Choose a coordinate $x$ on $\bP^1 - \{ \infty \}$.  The space $B_K$ of
characteristic polynomials becomes
$$
\{ P(x,y) = y^r + b_1(x)y^{r-1} + \cdots + b_r(x) \; | \; b_i(x) \; {\rm is
\; a \; polynomial \, in} \; x \; {\rm of \, degree} \; \le i \cdot d \}.
$$
The total space of the line bundle ${\cal O}_{\bP^1}(d \cdot \infty)$
restricted
to the affine line $\bP^1 - \{ \infty \}$ is isomorphic to the affine
plane, and under this isomorphism $P(x,y)$ becomes the equation of the
spectral curve as an affine plane curve.

Denote by $B^0 \subset B_K$ the subset of smooth spectral curves.  Let
$Q := Q_r(d)$ be the subset of $\HiggsModuli_K$ parametrizing pairs
$(E,\varphi)$ with a smooth spectral curve and a vector bundle $E$
isomorphic to $E_0 := \oplus^r {\cal O}_{\bP^1}(-1)$. $Q$ is a
Zariski open (dense) subset of $\HiggsModuli_K$ because:

\begin{description}
\item [i)]  by definition \ref{def-good-component-of-higgs-pairs}
$\HiggsComponent_K(r,-r)$ is irreducible,
\item [ii)]  $E_0$ is the unique semistable rank $r$ vector bundle of
degree $-r$ on $\bP^1$ and semistability is an open condition.
\end{description}
The bundle ${\rm End} \; E_0$ is the trivial Lie algebra bundle $\goth{gl}
_r(\C) \otimes {\cal O}_{\bP^1}$.  Hence, every point in $Q$ is
represented by an element $\varphi \in M_r(d) := H^0 (\bP^1,\goth{gl}
_r(\C) \otimes {\cal O}_{\bP^1}(d \cdot \infty))$, i.e., by an $r \times
r$ matrix $\varphi$ with polynomial entries of degree $\le d$.  Denote the
inverse image of $B^0$ in $M_r(d)$ by $M^0_r(d)$.  The subset $Q \subset
\HiggsModuli_K$ is simply the quotient of $M^0_r(d)$ by the conjugation
action of $PGL_r(\C)$.

\[
{\divide\dgARROWLENGTH by 2
\begin{diagram}
\node{M^{0}_r(d)}
\arrow{s}
\\
\node{Q}
\arrow{s}
\arrow{e}
\node{\HiggsModuli_{K}}
\arrow{s}
\\
\node{B^0}
\arrow{e}
\node{B_K}
\end{diagram}
}
\]

In this setting, Theorem \ref{thm-markman-botachin} specializes to the
following theorem of
Beauville and Adams-Harnad-Hurtubise-Previato generalizing results of
Mumford and Moser \cite{Mum} in rank $2$:

\begin{theorem}\label{thm-beauville} \cite{B,AHH}
\begin{itemize}
\begin{enumerate}
\item The quotient $Q$ of the action of $PGL_r(\C)$ by conjugation on
$M^0_r(d)$ is a smooth variety.
\item The fiber of the characteristic polynomial maps $H : Q
\rightarrow B^0$ over the polynomial of a spectral curve $C$ is the
complement $J^{g-1}_C - \Theta$ of the theta divisor in the Jacobian of
line bundles on $C$ of degree $g-1$ \ $(g = {\rm genus \; o}f \, C)$.
\item The choice of $d+2$ points $a_1, \cdots a_{d+2}$ on $\bP^1$
determines a Poisson structure on $Q$. The characteristic polynomial
map $H : Q \rightarrow
B^0$ is an algebraically completely  integrable Hamiltonian system with
respect to each of these Poisson structures.
\item The symplectic leaves of $Q$ are obtained by fixing the values
(of the coefficients) of the characteristic polynomials at the points $\{
a_i \}^{d+2}_{i=1}$.
\end{enumerate}
\end{itemize}
\end{theorem}
We note that in \cite{B} the Poisson structure on $Q$ was obtained as the
reduction of a Poisson structure on $M^0_r(d)$.  The latter was the
pullback of the Kostant-Kirillov Poisson structure via the embedding
$$
M_r(d) \hookrightarrow \goth{gl}_r(\C)^{d+2}
$$
by Lagrange interpolation at $a_1, \cdots , a_{d+2}$.
This embedding will be used in section
\ref{sec-geodesic-flow-via-polynomial-matrices} where geodesic flow on
ellipsoids is revisited.

A choice of a divisor $D = a_1 + \cdots + a_{d+2}$ of degree $d+2$ on
$\bP^1$ determines an isomorphism of ${\cal O}_{\bP^1}(d \cdot \infty)$
with $\omega_{\bP^1}(D)$.  For example, if $a_1, \cdots , a_i$ are finite,
$a_{i+1} = a_{a+2} = \cdots = a_{d+2} = \infty$ then we send a polynomial
$f(x)$ of degree $\le d$ to the meromorphic $1$-form
$$
{f(x) \over {\prod^i_{j=1}(x-a_j)}} dx.
$$
When the $d+2$ points are distinct, Lagrange interpolation translates to
the embedding
$$
{\rm Res}: M_r(d) = H^0(\bP^1, \goth{gl}_r(\C) \otimes \omega_{\bP^1}(D))
\hookrightarrow \goth{gl}_r(\C)^{d+2}
$$
via the residues of meromorphic $1$-form valued matrices at the points
$a_i$ (if $a_i$ has multiplicity $2$ or higher, we replace the $i$-th copy
of $\goth{gl}_r(\C)$ by its tangent bundle or higher order infinitesimal
germs at  $a_i$ of sections of the trivial bundle $\goth{gl}_r(\C) \otimes
{\cal O}_{\bP^1})$.

\subsubsection{Explicit Equations for Jacobians of Spectral Curves with a
Cyclic Ramification Point}
\label{sec-explicit-equations-for-jacobians}

A further simplification occurs for matrices with a nilpotent leading
coefficient (nilpotent at $\infty$).  The projection $M^0_r(d) \rightarrow
Q$ has a natural section over the image $N \subset Q$ of this locus.  So
$N$ can be described concretely as a space of polynomials (rather than as a
quotient of such a space).

As a consequence we obtain explicit equations in $M_r(d)$ for the
complement $J_C - \Theta$ of the theta divisor in the Jacobian of every
irreducible and reduced $r$-sheeted spectral curve over $\bP^1$ which is
totally ramified and smooth at $\infty$ (generalizing the equations for
hyperelliptic curve (case $r=2$) obtained in \cite{Mum}).

\begin{lem}\label{lemma-normal-form}  Let $A = A_dx^d + \cdots + A_1x +
A_0$ be an $r \times r$ traceless matrix with polynomial entries of degree
$\le d$
\begin{description}
\item[i)] whose spectral curve in ${\cal O}_{\bP^1}(d \cdot \infty)$ is
irreducible and reduced and smooth over $\infty$, and
\item[ii)] whose leading coefficient $A_d$ is a nilpotent (necessarily
regular) matrix.
\end{description}
Then there exists a unique element $g_0 \in PGL_r(\C)$ conjugating $A$ to a
matrix $A' = x^d \cdot J + \sum ^{d-1}_{i=0} A'_i x^i$ of the form:

\begin{equation}\label{eqn-normal-form}
A' =
x^d
\left(
\begin{array}{cccccc}
0 &          &          &          & 0 & 0 \\
1 &          &          &          & 0 & 0 \\
0 &    \cdot &          &          & 0 & 0 \\
0 &          &    \cdot &          & 0 & 0 \\
0 &          &          &    \cdot & 0 & 0 \\
0 &          &          &          & 1 & 0
\end{array}
\right)
+ x^{d-1}
\left(
\begin{array}{cccccc}
\star &          &  \dots   &          & \star & \beta_{r} \\
\star &          &  \dots   &          & \star & 0 \\
\\
\vdots &         &  \vdots   &          & \vdots & \vdots \\
\\
\star &          &  \dots   &          & \star & 0
\end{array} \right)
+ \sum_{i=0}^{d-2}x^i A'_{i}
\end{equation}

\noindent
where $(-1)^{r+1}\beta_r$ is the (leading) coefficient of $x^{dr-1}$ in
the determinant $b_r(x)$ of $A(x)$.
\end{lem}

\noindent
{\bf Remark:} \ Notice that the coefficients $b_i(x)$
in the characteristic polynomial $P(x,y) = y^r + b_1(x)y^{r-1} + \cdots +
b_r(x)$ of $A(x)$ satisfy degree $b_i(x) \le d \cdot i-1$ since $A$ is
nilpotent at $\infty$, and degree $b_r(x) = dr-1$ since the spectral curve
is smooth over $\infty$. Thus $\beta_r \ne 0$.

\smallskip
\noindent
{\bf Proof} (of lemma \ref{lemma-normal-form}): \
Let $J$ be the nilpotent regular constant matrix appearing as
the leading
coefficient of $A'(x)$ in the normalized form (\ref{eqn-normal-form}).  Let
$\C[J]$ be the algebra of polynomials in $J$ with constant coefficients.
The proof relies on the elementary fact that $\C^r$, as a left
$\C[J]$-module, is free.  Any vector with non zero first entry is a
generator.  $A_d$ is conjugate to $J$.  Thus we may assume that $A_d = J$
and it remains to show that there exists a unique element in the stabilizer
of $J$ in $PGL_r(\C)$ conjugating $A(x)$ to the normal form
(\ref{eqn-normal-form}).

Since $A_d = J$, the first entry in the right column $R$ of $A_{d-1}$ is
$\beta_r$.  Thus $R$ is a generator of $\C^r$ as a $\C[J]$-module.  Any
element $g \in PGL_r(\C)$ in the commutator subgroup of $J$ is an invertible
element
in $\C[J]$ and can be written (up to scalar multiple)
in the form $g = I + N$, $N$ nilpotent.
The right column of $gA_{d-1}g^{-1}$ is $R + NR$
and there exists a unique nilpotent $N \in \C[J]$ such
that $NR =
\left(
\begin{array}{c}
\beta_r \\
0 \\
\vdots\\
0
\end{array}
\right) - R.$
Thus $g$ is unique up to a scalar factor.
\EndProof

Denote the affine subvariety of $M^0_d(r)$ of matrices satisfying the
$r^2+r - 1$ equations (\ref{eqn-normal-form}) by ${\tilde N}$.
The subvariety ${\tilde N}$ is
a section of the principal $PGL_r(\C)$ bundle $M^0_r(d) \rightarrow Q$ over
the locus $N$ of conjugacy classes of polynomial matrices with a nilpotent
leading coefficient.  $N$ is a Poisson subvariety of $Q$ with respect to
any Poisson structure on $Q$ determined by a divisor $D$ as in
theorem \ref{thm-beauville}, provided that $D$ contains the point at
infinity $\infty \in \bP^1$.

Choose a characteristic polynomial $P(x,y) = y^r + b_1(x)y^{r-1} + \cdots +
b_r(x)$  in $B_{{\cal O}_{\bP^1}(d)}$ of a smooth spectral curve $C$ with
degree $b_i(x) \le id - 1$,  $b_r(x)$ of degree $rd - 1$ with leading
coefficient $(-1)^{r+1}\beta_r$.  Theorem \ref{thm-beauville} implies
that the equations

\begin{description}
\item[a)] $A_d = J,$
\item[b)] The $r$-th column of $A_{d-1}$ is
$
\left(
\begin{array}{c}
{\cal \beta}_r \\
0 \\
\vdots \\
0 \\
\end{array}
\right),
$
\item[c)] char $(A(x)) = P(x,y)$
\end{description}
define a subvariety of $M_r(d)$ isomorphic to the complement $J^{g-1}_C -
\Theta_C$ of the theta divisor in the Jacobian of $C$.

\subsubsection
{Geodesic Flow on Ellipsoids via $2 \times 2$ Polynomial Matrices}
\label{sec-geodesic-flow-via-polynomial-matrices}

We use polynomial matrices to retrieve the Jacobi-Moser-Mumford system
which arose in chapter \ref{ch3} out of the geodesic flow on an ellipsoid.  Our
presentation follows \cite{B}.

Consider a spectral polynomial $P(x,y)$  in $B^0_{{\cal O}_{\bP^1}(d \cdot
\infty)}$ of the form
\begin{description}
\item[(i)]  $P(x,y) = y^2 - f(x)$ where f(x) is monic of degree $2d - 1$.
\end{description}
\noindent
The corresponding spectral curve $C$ is smooth, hyperelliptic of genes $g
= d-1$ and ramified over $\infty$.  Theorem \ref{thm-beauville} implies that
the fiber
$$
H^{-1}(P(x,y)) =
\left\{
\left(
\begin{array}{cc}
V & U \\
W & -V \\
\end{array}
\right)
\; \mid \; V^2 + UW = f(x) \right\} / PGL_2(\C)
$$
of the characteristic polynomial map is isomorphic to the complement
$J^{g-1}_C - \Theta$ of the theta divisor.

Lemma \ \ref{lemma-normal-form} specializes in our case to the following
statement (note that $\beta_r=1$ since $f$ is taken to be monic):\\
{\em
The $PGL_2(\C)$ orbit of a matrix
$\left(
\begin{array}{cc}
V & U \\
W & -V \\
\end{array}
\right)
$
over $H^{-1}(P(x,y)) \cong J_C^{g-1} - \Theta$ contains a unique matrix
satisfying
\begin{description}
\item [(ii)]  $W$ is monic of degree $d$, \\
$U$ is monic of degree $d-1$ and \\
deg $V \le d-2$.
\end{description}
} 
\noindent
In other words, condition (ii) and
\begin{description}
\item [(iii)] $V^2 + UW = f(x)$
\end{description}
are the equations of $J_C^{g-1} - \Theta$ as an affine subvariety of the
subspace of traceless matrices in $M_2(d)$.  In fact, condition (ii)
defines a section $\varphi : N \rightarrow M_2(d)$ over the locus $N$ in $Q$
of conjugacy classes
with characteristic polynomial satisfying condition (i).

Recall that the Jacobi-Moser-Mumford integrable system linearizing the
geodesic flow of the ellipsoid $\sum^d_{i=1} \; a^{-1}_i x^2_i = 1$ is
supported on the tangent bundle $TS$ of the sphere $S \subset \RealNumbers^d$.
Our
discussion ended by describing the quotient of $TS$ by the group $G \simeq
(\Z/2\Z)^d$ of involutions.  We will describe in the next three steps an
isomorphism of this quotient with a symplectic leaf $X$ of $Q$.
\begin{description}
\item[\underline{Step I:}]  (Identification of the symplectic leaf  $X$).
Assume that the points $a_1, \cdots, a_d \in \bP^1 - \{ \infty \}$ are
distinct and let $a_{d+1} = a_{d+2} = \infty$.  Let  $X \subset Q$ be the
symplectic leaf over the subspace of characteristic polynomials $P(x,y) =
y^2 - f(x)$ satisfying
$$f(a_i) = 0, \; 1 \le i \le d, \ \ \deg f = 2d-1 \ \ {\rm
and}\; f \; {\rm is \; monic}.$$
The spectral curves of matrices in the leaf $X$
have genus $d-1$, and are branched over the fixed
$g+2$ points $a_1, \cdots , a_d,\infty$ and $g$ varying points.
\item[\underline{Step II:}]  (Embedding of $X$ in the product
${\cal N}^d$ of the regular nilpotent orbit).
The isomorphism
${\cal O}_{\bP^1}(d \cdot \infty) \stackrel{\sim}{\rightarrow}
\omega_{\bP^1} (\sum^d_{i=1} a_i + 2 \cdot \infty)$ sending $F(x)$ to
${{F(x)dx} \over {\prod^d_{i=1}}(x-a_i)}$ translates the matrix
$
\left(
\begin{array}{cc}
V(x) & U(x) \\
W(x) & -V(x) \\
\end{array}
\right)
$
to a matrix $\varphi$ of meromorphic $1$-forms.  The
residues of $\varphi$ satisfy:\\

\smallskip
${
R_\infty := Res_\infty(\varphi) =
\left(
\begin{array}{cc}
0 & -1 \\
s & 0
\end{array}
\right)
{\rm for \hspace{1ex} some}
 \ s \in \C \ \ \ ({\rm condition} \; (ii)),
}$

${
R_i := Res_{a_i}(\varphi) =
\left(
\begin{array}{cc}
V(a_i) & U(a_i) \\
W(a_i) & -V(a_i) \\
\end{array}
\right)
{{1} \over {\prod^d_{\stackrel{j=1}{j \ne i}}} (a_i-a_j)}.
}$

\smallskip
\noindent
The residues at the finite points $a_i$ can be calculated using
Lagrange interpolation of
polynomials of degree $d$ at the $d+1$ points $a_1, \cdots , a_d, \infty$
given by the formula
\begin{equation}
\label{eqn-lagrange-interpolation}
 F(x) = \sum^d_{i=1}F(a_i)
\frac
{\prod^d_{\stackrel{j=1}{j \ne i}}(x   - a_j)}
{\prod^d_{\stackrel{j=1}{j \ne i}}(a_i - a_j)}
+ F(\infty) \prod^d_{j=1} (x - a_j)
\end{equation}
where $F(\infty)$ is the leading coefficient of $F(x)$.

\smallskip
\noindent
The residues $R_\infty$, $R_i$ are nilpotent regular $2 \times 2$
matrices and the residue
theorem implies that $R_\infty = - \sum^d_{i=1} R_i$.  The residue map $Res
: X \rightarrow {\cal N}^d$ defines a symplectic embedding $\varphi \mapsto
(R_1,
\cdots, R_d)$ of the symplectic leaf $X$ of $Q$ in the Cartesian product of
$d$ copies of the regular nilpotent orbit ${\cal N}$ in $\goth{gl}_2(\C)$.

\item[\underline{Step III}]  (The $2^d$ covering $TS \rightarrow X$).
Endow $\C^2$ with the symplectic structure $2dx \wedge dy$.  The map
$\C^2 - \{(0,0)\} \rightarrow {\cal N}$ sending $(x,y)$ to
$
\left(
\begin{array}{cc}
xy  & -x^2 \\
y^2 & -xy \\
\end{array}
\right)
$ is a symplectic $SL_2(\C)$-equivariant double cover of the regular
nilpotent orbit $\cal N$ (where $SL_2(\C) \cong Sp_2(\C,2dx \wedge dy)$ acts
on $\C^2$ via the standard representation).
We obtain a $2^d$-covering $\tau : (\C^2 - \{(0,0)
\})^d \rightarrow {\cal N}^d$.  The residue theorem translates to the fact
that the image $Res(X) \subset {\cal N}^d$ is covered by
$$
\left\{ ({\bar x}, {\bar y}) = ((x_1,y_1), \cdots, (x_d,y_d)) \ | \ \sum \;
x_iy_i = 0 \; {\rm and} \sum^d_{i=1}x^2_i = 1 \right\}.
$$
This is exactly the tangent bundle $TS \subset (\C^2)^d$ of the sphere
$S \subset \C^d$.
\end{description}
\newpage

\section{Poisson structure via levels} \label{ch5}

We construct a Poisson structure on the moduli space of meromorphic
Higgs pairs in two steps (following \cite{markman-higgs}):

\smallskip
\noindent
- First we realize a dense open subset of moduli as the orbit space of
a Poisson action of a group on the cotangent bundle of the moduli
space of vector bundles with level structures (sections
\ref{sec-level-structures},
\ref{sec-the-cotangent-bundle} and
\ref{the-poisson-structure}).

\noindent
- Next we
exhibit a $2$-vector on the smooth locus of moduli, using
a cohomological construction (section \ref{sec-linearization}).
On the above dense open set this agrees with the Poisson structure,
so it is a Poisson structure everywhere.

\medskip
We summarize the construction in section
\ref{sec-hamiltonians-and-flows}
in a diagram whose rotational symmetry
relates dual pairs of realizations.

\subsection{Level structures} \label{sec-level-structures}

     Fix a curve $\Sigma$, an effective divisor $D = \sum p_i$
in $\Sigma$, and a rank $r$ vector bundle $E$ on $\Sigma$.  A
level $D$ structure on $E$ is an $\O_D$-isomorphism $\eta: E
\otimes \O_D \stackrel{\sim}{\rightarrow} \O^{\oplus r}_D$.
Seshadri \cite{seshadri-construction-moduli-vb}
constructs a smooth, quasi-projective moduli space
$\ModuliVB_\Sigma(r,d,D)$ parametrizing stable pairs $(E,\eta)$.  Here
stability means that for any subbundle $F \subset E$, $$\frac{deg
F - deg D}{rank F} < \frac{deg E - deg D}{rank E}.$$

     The level-$D$ group is the projectivized group of
$\O_D$-algebra automorphisms,
$$
G_D:= \P Aut_{\O_D}(\O^{\oplus
r}_{D}).
$$
(i.e. the automorphism group modulo complex scalars
$\C^*$.) It acts on $\ModuliVB_\Sigma(r,d,D)$: an element $g \in G_D$
sends $$[(E,\eta)] \mapsto [(E,\overline{g} \circ \eta)],$$ where
$\overline{g} \in Aut_{\O_D}(\O^{\oplus r}_{D})$ lifts $g$, and
$[\cdot]$ denotes the isomorphism class of a pair.  The open set
$\ModuliVB^\circ_\Sigma (r,d,D)$, parametrizing pairs $(E,\eta)$ where
$E$ itself is stable, is a principal $G_D$-bundle over
$\ModuliVB^s_\Sigma (r,d)$.  The Lie algebra $\LieAlg{g}_D$ of
$G_D$ is given by
$\LieAlg{gl}_r(\O_D)$/scalars.


\subsection{The cotangent bundle} \label{sec-the-cotangent-bundle}

     We compute deformations of a pair $$(E,\eta) \in \ModuliVB_D :=
\ModuliVB_\Sigma(r,d,D)$$ as we did for the single vector bundle
$$
E \in \ModuliVB := \ModuliVB_\Sigma (r,d).
$$
Namely, $E$ is given (in terms of an open
cover of $\Sigma$) by a $1$-cocycle with values in the sheaf of groups
$GL_r(\O_\Sigma)$.  Differentiating this cocycle with respect to
parameters gives a $1$-cocycle with values in the corresponding
sheaf of Lie algebras, so we obtain the identification
$$
T_E\ModuliVB
\approx H^1(End E).
$$
Similarly, the pair $(E,\eta)$ is given by
a $1$-cocycle with values in the subsheaf
$$
GL_{r,D}(\O_\Sigma)
:= \{f \in GL_r(\O_\Sigma)|f-1 \in \LieAlg{gl}_r({\cal I}_D)\}.
$$
Differentiating, we find the natural isomorphism
$$
T_{(E,\eta)}
\ModuliVB_D \approx H^1({\cal I}_D \otimes End E),
$$
so by Serre duality,
$$
T^*_{(E,\eta)} \ModuliVB_D \approx H^{0}(End E \otimes \omega(D)).
$$
(we
identify $End E$ with its dual via the trace.)  We will denote a
point of this cotangent bundle by a triple $(E,\varphi,\eta)$,
where $(E,\eta) \in \ModuliVB_D$ and $\varphi$ is a $D$-Higgs field,
$$\varphi: E \to E \otimes \omega(D).$$


\subsection{The Poisson structure} \label{the-poisson-structure}

     The action of the level group $G_D$ on $\ModuliVB_D$ lifts naturally
to an action of $G_D$ on $T^*\ModuliVB_D$.  Explicitly, an element $g \in
G_D$ with lift $\overline{g} \in GL_r(\O_D)$ sends
$$
(E,\varphi,\eta) \mapsto (E, \varphi, \overline{g} \circ
\eta).
$$
The lifted action has the following properties:
\begin{list}{{\rm(\arabic{bean})}}{\usecounter{bean}}
\item It is Poisson with respect to the standard symplectic
structure on $T^*\ModuliVB_D$ (holds for any lifted action, see
example \ref{example-coadjoint-orbits}).
\item The moment map
$$
\mu: T^*\ModuliVB_D \to \LieAlg{g}^*_D
$$
is given by
\begin{equation} \label{eq-moment-map-of-finite-dim-level-action}
\mu(E,\varphi,\eta): A \mapsto \mbox{Res\ Trace}\ (A\cdot
\varphi^\eta),
\end{equation}
where
$$
\begin{array}{lcl}
A & \in & \LieAlg{g}_D =
(\LieAlg{gl}_r(\O_D))/\mbox{(scalars)} \approx
(\LieAlg{gl}_r(\O_D))_{\mbox{traceless}},\\ \varphi^\eta & := & \eta \circ
\varphi \circ \eta^{-1} \in H^0(\LieAlg{gl}_r(\omega(D) \otimes \O_D)),\\
\end{array}
$$
and the residue map
$$
Res : H^0(\omega(D) \otimes \O_D) \to H^1(\omega)
\approx \C
$$
is the coboundary for the restriction sequence
$$
0 \to \omega \to \omega(D) \to \omega(D) \otimes \O_D \to 0
$$
(cf. \cite{markman-higgs} Proposition 6.12).
\item $G_D$ acts freely on the open subset $(T^*\ModuliVB_D)^{\circ}$
parametrizing triples $(E,\varphi,\eta)$ where $(E,\eta)$ is
stable and $(E,\varphi)$ is a stable Higgs bundle, since such
bundles are simple.  This makes $(T^*\ModuliVB_D)^{\circ}$ into a principal
$G_D$-bundle over an open subset $\HiggsModuli^{\circ}_D$ of
$\HiggsModuli^s_D$.
\end{list}

     We conclude that the symplectic structure on $(T^*\ModuliVB_D)^{\circ}$
induces a Poisson structure on $\HiggsModuli^{\circ}_D$.  The
symplectic leaves will then be the inverse images under $\mu$ of
coadjoint orbits in $\LieAlg{g}^*_D$.

\subsection{Linearization} \label{sec-linearization}

The main remaining task is to find a two-vector on the non-singular locus
$\HiggsModuli^{ns}_D$ whose restriction to $\HiggsModuli^{\circ}_D$
is the above Poisson structure. This two vector is then automatically
Poisson. The algebraic complete integrability
of the component $\HiggsComponent_D$
(see definition \ref{def-good-component-of-higgs-pairs})
would then follow:  The spectral curve $C_b$, for generic $b \in B_D$,
is non-singular, so its Jacobian $\J(C)$ is contained in
$\HiggsComponent_D$.  Thus any Hamiltonian vector field must be
constant on the generic fiber $\J(C)$, hence on all fibers.

     A natural two-vector defined over all of $\HiggsModuli^{ns}_D$
can be given in several ways.  One \cite{markman-higgs} is to identify the
tangent spaces to $\HiggsModuli_D$ (and related spaces) at their
smooth points as hypercohomologies, $\h^1$, of appropriate
complexes:
$$
\begin{array}{ccc}
\mbox{{\underline{space}}} &
\mbox{{\underline{at}}} &
\mbox{{\underline{complex}}}
\\
\ModuliVB & E & End E
\\
&&
\\
\ModuliVB_D & (E,\eta) & End
E(-D)
\\
&&
\\ \HiggsModuli_D & (E,\varphi) &
\overline{{\cal K}}= [End E
\stackrel{ad
\varphi}{\longrightarrow} End E \otimes
\omega(D)]
\\
&&
\\
T^*\ModuliVB_D &
(E,\varphi,\eta) & {\cal K} := [End E(-D) \stackrel{ad \varphi
\otimes i}{\longrightarrow} End E \otimes \omega(D)].\\
\end{array}
$$

\vspace{0.1in}
\noindent
where $i: \O(-D) \hookrightarrow \O$ is the natural inclusion.
These identifications are natural, and differentials of maps
between these spaces are realized by maps of complexes.  For
example, the fibration $T^*\ModuliVB_D \to \ModuliVB_D$, with fiber
$T^*_{(E,\eta)}\ModuliVB_D$, gives the sequence
\vspace{0.1in}
$$
\begin{array}{ccccccccc}
0 & \to & T^*_{(E,\eta)} \ModuliVB_D & \to &
T_{(E,\varphi,\eta)}(T^*\ModuliVB_D) &
\to & T_{(E,\eta)}\ModuliVB_D & \to & 0
\\
& & \parallel & & \parallel & & \parallel &
\\
0 & \to & H^0(End E \otimes
\omega(D)) & \to & \h^1({\cal K}) & \to & H^1(End E(-D)) & \to &
0\\
\end{array}
$$

\vspace{0.1in}
\noindent
derived from the short exact sequence of
complexes,
$$
0 \rightarrow End E \otimes \omega(D)[-1] \rightarrow K
\rightarrow End E (-D)
\rightarrow 0 ,
$$
while the (rational) map $T^*\ModuliVB_D \to
\HiggsModuli_D$
gives
\vspace{0.1in}
$$\begin{array}{ccccccccc}0 & \rightarrow & \LieAlg{g}_D &
\rightarrow &
T_{(E,\varphi,\eta)}(T^*\ModuliVB_D) & \rightarrow &
T_{(E,\varphi)}\HiggsModuli_D & \rightarrow & 0
\\
& & \parallel & & \parallel & & \parallel &
\\
0 & \rightarrow & \frac{H^0(End E
\otimes
\O_D)}{H^0(End E)} & \rightarrow & \h^1({\cal K}) &
\rightarrow &
\h^1(\overline{{\cal K}}) & \rightarrow & 0\\
\end{array}
$$

\vspace{0.1in}
\noindent
which
derives from: $$0 \to {\cal K} \to \overline{{\cal K}} \to End E
\otimes \O_D \to 0.$$

     The dual of a complex $\L: A \to B$ of vector bundles is the
complex
$$
\L^\Dual: B^* \otimes \omega \to A^* \otimes \omega.
$$
Grothendieck
duality in this case gives a natural isomorphism
$$
\h^1(\L)
\approx \h^1(\L^\Dual)^*.
$$
We note that ${\cal K}$ is
self-dual, in the sense that there is a natural isomorphism of
complexes,
$J: {\cal K}^\Dual \stackrel{\sim}{\rightarrow} {\cal K}$.
For $\overline{{\cal K}}$ we obtain a natural isomorphism
of complexes,
$\overline{{\cal K}}^\Dual \stackrel{\sim}{\rightarrow}
\overline{{\cal K}} \otimes \O(-D)$,
hence (composing with $i$)
a morphism
$$I: \overline{{\cal K}}^\Dual \to \overline{{\cal K}}.$$
Combining with duality, we get maps
$$
\h^1({\cal K})^* \approx
\h^1({\cal K}^\Dual) \stackrel{J}{\approx} \h^1({\cal K})
$$
and
$$
\h^1(\overline{{\cal K}})^* \approx
\h^1(\overline{{\cal K}}^\Dual) \RightArrowOf{I}
\h^1(\overline{{\cal
K}}).
$$
These give elements of $\otimes^2\h^1({\cal K})$
and $\otimes^2 \h^1(\overline{{\cal K}})$.  Both are skew
symmetric (since $ad_{\varphi}$, and hence $I,J$, are), so we
obtain global two-vectors on $T^*\ModuliVB^{ns}_D$ and
$\HiggsModuli^{ns}_D$.  At stable points these agree with (the
dual of) the symplectic form and its reduction modulo $G_D$,
which is what we need.

     Another way to find the two-vector on $\HiggsModuli_D$ is
based on the interpretation of $\HiggsModuli_D$ as a moduli space
of sheaves on the total space $S$ of $\omega(D)$.  At such a simple
sheaf ${\cal E}$, with support on some spectral curve $C$, Mukai
\cite{mukai}
identifies the tangent space to moduli with
$$Ext^1_{\O_S}(\E,\E),$$
and notes that any two-form $\sigma \in
H^0(\omega_S)$ determines an alternating bilinear map:
$$
Ext^1_{\O_S}(\E,\E) \times Ext^1_{\O_S}(\E,\E) \to
Ext^2_{\O_S}(\E,\E) \stackrel{tr}{\rightarrow} H^2(\O_S)
\stackrel{\sigma}{\rightarrow} H^2(\omega_S) \approx \C,
$$
hence a two-form
on moduli.  Mukai uses this argument to produce symplectic
structures on the moduli spaces of sheaves on $K3$ and abelian
surfaces.  The same argument works, of course, for sheaves on
$T^*\Sigma$; this reconstructs the symplectic form on Hitchin's
system.

     Our surface $S$ (the total space of $\omega(D)$) is related to
$T^*\Sigma$ by
a birational morphism $\alpha: T^*\Sigma \to S$.  The symplectic
form $\sigma$ does not descend to $S$, but its inverse $\sigma^{-1}$ does
give a two vector on $S$ which is non-degenerate away from $D$
and is closed there (since it is locally equivalent to the
Poisson structure on $T^*\Sigma$).

     Tyurin notes \cite{tyurin-symplectic}
that a variant of this argument produces a
two-vector on moduli from a two-vector on $S$.  Now the birational
morphism $T^*\Sigma \to S$ takes the Poisson structure
on $T^*\Sigma$ to one on $S$, so the Mukai-Tyurin argument gives
the desired two-vector on $\HiggsModuli_D$. In chapter \ref{ch8} this
approach is generalized to higher dimensional varieties.

     A third argument for the linearization is given by Bottacin
\cite{botachin}.
He produces an explicit two-vector at stable points using
a deformation argument as above, and then makes direct, local
computations to check closedness of the Poisson structures and
linearity of the flows.


\subsection{Hamiltonians and flows in $T^*\ModuliVB_D$}
\label{sec-hamiltonians-and-flows}

     We saw that the level group $G_D$ acts on $T^*\ModuliVB_D$, inducing
the Poisson structure on the quotient $\HiggsModuli_D$,  and that
the moment map is
$$
\begin{array}{lcc}
\mu: T^*\ModuliVB_D & \rightarrow & \LieAlg{g}^*_D
\\
\mu(E,\varphi, \eta)(A) & := & \mbox{Res\ Trace}\ (A
\cdot \varphi^{\eta}).\\
\end{array}
$$ The characteristic
polynomial map $\widetilde{h}: T^*\ModuliVB_D \to B_D$ is a composition
of
the Poisson map $T^*\ModuliVB_D \to \HiggsModuli_D$ with the Hamiltonian
map $h: \HiggsModuli_D \to B_D$.  Hence $\widetilde{h}$ is also
Hamiltonian.  Clearly, $\widetilde{h}$ is $G_D$-invariant.

     The composition $T^*\ModuliVB_D \stackrel{\mu}{\rightarrow}
\LieAlg{g}^*_D \rightarrow  \LieAlg{g}^*_D/G_D$ is a
$G_D$-invariant Hamiltonian
morphism and hence factors through $\HiggsModuli_D$.  It follows
that it factors also through $B_D$ since
$h: \HiggsModuli_D \to B_D$ is a Lagrangian
fibration whose generic fiber is connected
(see remark \ref{rem-acihs-implies-maximal-commutative-subalgebra}).
The conditions of
example \ref{example-diagram-hexagon-plus-realization}
in section \ref{subsec-moment-maps}
are satisfied and we get a diagram with
a $180^\circ$ rotational symmetry in which opposite spaces are
dual pairs of realizations.  The realization dual to $T^*\ModuliVB_D \to
\LieAlg{g}^*_D/G_D$ is the rational morphism $T^*\ModuliVB_D \to
G\HiggsModuli_D :=
\HiggsModuli_D \times_{(\LieAlg{g}^*_D/G_D)}\LieAlg{g}^*_D$
to the fiber product.
The one dual to $\widetilde{h}: T^*\ModuliVB_D \to B_D$ is the morphism
$T^*\ModuliVB_D \to G B_D:= B_D \times _{(\LieAlg{g}^*_D/G_D)}\LieAlg{g}^*_D$
to the fiber product.
We write down the spaces and typical elements in them:
\begin{equation} \label{diagram-hexagon-for-higgs-pairs}
\begin{array}{ccc}
{\divide\dgARROWLENGTH by 4
\divide\dgHORIZPAD by 2
\divide\dgCOLUMNWIDTH by 2
\begin{diagram}[TTT]
\node[3]{T^*\ModuliVB_D}
\arrow{s}
\\
\node[3]{G\HiggsModuli_D}
\arrow[2]{sw}
\arrow{se}
\\
\node[4]{GB_D}
\arrow[2]{sw}
\arrow{se}
\\
\node{\HiggsModuli_D}
\arrow{se}
\node[4]{\LieAlg{g}^*_D}
\arrow[2]{sw}
\\
\node[2]{B_D}
\arrow{se}
\\
\node[3]{\LieAlg{g}^*_D/G_D}
\arrow{s}
\\
\node[3]{(0)}
\end{diagram}
}
&
\hspace{3ex}
&
{\divide\dgARROWLENGTH by 4
\divide\dgHORIZPAD by 2
\divide\dgCOLUMNWIDTH by 2
\begin{diagram}[TTT]
\node[3]{(E,\varphi,\eta)}
\arrow{s}
\\
\node[3]{(E,\varphi,\varphi^{\eta})}
\arrow[2]{sw}
\arrow{se}
\\
\node[4]{({\rm char}\varphi,\varphi^{\eta})}
\arrow[2]{sw}
\arrow{se}
\\
\node{(E,\varphi)}
\arrow{se}
\node[4]{\varphi^{\eta}}
\arrow[2]{sw}
\\
\node[2]{{\rm char}\varphi}
\arrow{se}
\\
\node[3]{{\rm char}\varphi^{\eta}}
\arrow{s}
\\
\node[3]{(0)}
\end{diagram}
}

\end{array}
\end{equation}
\newpage

\section {Spectral flows and $KP$} \label{ch6}
\label{sec-spectral-flows-and-kp}

     Our aim in this section is to relate the general spectral
system which we have been considering to the $KP$ and multi-component
$KP$ hierarchies.  We start by reviewing these
hierarchies and their traditional relationship to curves and
bundles via the Krichever map.  We then reinterpret these flows
as coming from Hamiltonians on the limit $T^*U_\infty$ of our
previous symplectic spaces.  We show that $\mbox{Higgs}_\infty$
can be partitioned into a finite number of loci, each of which
maps naturally to one of the $mcKP$-spaces in a way which
intertwines isospectral flows with $KP$ flows.
As an example we consider the Elliptic solitons studied by Treibich and
Verdier.

\subsection{The hierarchies} \label{sec-the-heirarchies}

\noindent
\underline{KP}

     Following the modern custom (initiated by Sato, explained by
Segal-Wilson \cite{segal-wilson-loop-groups-and-kp},
 and presented elegantly in
\cite{AdC,mulase-cohomological-structure,li-mulase-category}
and elsewhere), we think of the $KP$ hierarchy
as given by the action of an infinite-dimensional group on an
infinite-dimensional Grassmannian: set
 $$\begin{array}{lrl}K & :=
& \C((z)) = \mbox{field\ of\ formal\ Laurent\ series\ in\ a\
variable}\ z\\
Gr & := & \{\mbox{subspaces}\ W \subset K|
\mbox{projection}\ W \to K/\C[[z]] z\ \mbox{is\ Fredholm}\}\\
 & = & \{\mbox{subspaces\ ``comparable\ to}\ \C[z^{-1}]"\}.\\
\end{array}$$ This can be given an algebraic structure which
allows us to talk about vector fields on $Gr$, finite-dimensional
algebraic subvarieties, etc.  Every $a \in K$ determines a
vector field $KP_a$ on $Gr$, whose value at $W \in Gr$ is the map
$$W \hookrightarrow K \stackrel{a}{\rightarrow} K \rightarrow
K/W,$$ considered as an element of $$Hom(W,K/W) \approx T_W Gr.$$

The (double) $KP$ hierarchy on $Gr$ is just this collection of
commuting vector fields.  For $a \in \C[[z]]$, this vector field
comes from the action  on $Gr$ of the one-parameter subgroup $exp(ta)$ in
$\C[[z]]^*$, which we consider trivial.  The $KP$
hierarchy itself thus consists of the vector fields $KP_a$, for $a \in
\C[z^{-1}]z^{-1}$, on the quotient $Gr/(\C[[z]]^*)$.  This quotient is
well-behaved:
the action of $\C ^*$ is trivial, and the unipotent part
$1 + z \C [[z]] $ acts freely and with transversal slices. One restricts
attention
to the open subset of this quotient ("the big cell") parametrizing $W$ of fixed
index (the index of $W$ is the index of the Fredholm projection) and satisfying
a general position condition with respect to the standard subspace
$W_0 := \C[z^{-1}] $.
Sato's construction identifies this subset with the space $\Psi$ of
pseudo differential operators of the form
$$\L = D + \sum^{\infty}_{i=1} u_iD^{-i}$$
where
$$u_i = u_i(t_1,t_2,\cdots)$$
and
$D = \partial/\partial t_1$.
The resulting flows on $\Psi$ have the familiar Lax form:
$$\frac{\partial \L}{\partial t_i} = [(\L^i) _{+}, \L],$$
where $(\L^i) _{+}$ is the differential operator part of $\L^i$.

\bigskip
\noindent
\underline{multi component KP}

    The $k^{th}$ multi-component $KP$ hierarchy ($mcKP$) is
obtained by considering instead the Grassmannian $Gr_k$ of
subspaces  of $K^{\oplus k}$ comparable to $(\C[z^{-
1}])^{\oplus k}$.  The entire ``loop algebra'' $gl(k,K)$ acts here, but
to obtain commuting flows we need to restrict to a commutative subalgebra.
For the k-th multi-component KP we take the simplest choice, of diagonal
matrices, i.e.
we consider the action of $(\C[z^{-1}]z^{-1})^{\oplus k}$ on the
quotient $Gr_k/(\C[[z]]^*)^k$.  There is a big cell
$\Psi_{k} \subset  Gr_k/(\C[[z]]^*)^k$, consisting as before of subspaces in
general position with respect to a reference subspace $W_0$,  on which the
flow is given by a Lax equation (for vector-valued operators). \\

\noindent
\underline{Heisenberg flows}

More generally, for a partition
$$\underline{n} = (n_1,\cdots,n_k)$$
of the positive integer $n$,
we can consider, following \cite{adams-bergvelt},
the maximal torus $Heis_{\underline{n}}$ of type
$\underline{n}$ in
$GL(n,K)$, as well as  $heis_{\underline{n}}$, the corresponding
Lie subalgebra  in
$gl(n,K)$.  These consist of matrices in block-diagonal form,
where the $i^{th}$ block is a formal power series in the $n_i
\times n_i$ matrix

\begin{equation} \label{eq-the-generator-of-the-ith-heisenberg-block}
P_{n_i}: =
\left(
\begin{array}{cccccc}
0 &          &          &          & 0 & z \\
1 &          &          &          & 0 & 0 \\
0 & \cdot &          &          & 0 & 0 \\
0 &          & \cdot &          & 0 & 0 \\
0 &          &          & \cdot & 0 & 0 \\
0 &          &          &          & 1 & 0
\end{array}
\right)
\end{equation}

\noindent
We recall
that this matrix arises naturally when we consider a
vector bundle which is the direct image of a line bundle, near a
point where $n_i$ sheets come together: in terms of a natural
local basis of the vector bundle, it expresses multiplication by
a coordinate upstairs (see (\ref{eq-ramification-matrix})).
The $\underline{n}^{th}$ $mcKP$ (or
``Heisenberg flows'' of type $\underline{n}$) lives on the
quotient of $Gr_n$ by the non-negative powers of the $P_{n_i}$,
and a basis for the surviving flows is indexed by $k$-tuples
$(d_1,\cdots,d_k)$, $d_i > 0$.  Again, this can all be realized
by Lax equations on an appropriate space $\Psi_{\underline{n}}$
of pseudo differential operators.  When $\underline{n} =
(1,\cdots,1)$ we recover the $n^{th}$ $mcKP$.  When
$\underline{n} = (n)$, the flows are pulled back from the
standard $KP$ flows on $Gr$, via the mixing map
$$m_n: Gr_n \to Gr$$
sending
$$\widetilde{W} \subset \C((\widetilde{z}))^{\oplus n}$$
 to
$$W := \{\sum^{n-1}_{i=0}
a_i(z^n)z^i|(a_0(\widetilde{z}),\cdots,a_{n-1}(\widetilde{z}))
\in
\widetilde{W}\} \subset \C((z)).$$
 An arbitrary $k$-part
partition $\underline{n}$ of $n$ determines a map
$$m_{\underline{n}}: Gr_n \to Gr_k,$$
and the
$\underline{n}^{th}$ Heisenberg flows are pullbacks of the
$k^{th}$ $mcKP$. The natural big cell in this situation is determined by the
cartesian diagram:

$$
\begin{array}{lcccc}
Gr_n  & \longrightarrow & Gr_n / Heis^+_{\underline{n}} & \hookleftarrow &
\Psi _{\underline{n}}\\
 \mbox{   } \downarrow m_{\underline{n}}  & & \downarrow & & \downarrow  \\
Gr_k & \longrightarrow & Gr_k / (C[[z]]^*)^k &  \hookleftarrow & \Psi _k
\end{array}
$$

\subsection{ Krichever maps} \label{sec-krichever-maps}

\noindent
\underline{The data}

    A basic Krichever datum (for the $KP$ hierarchy) consists of
a quintuple $$(C,p,z,L,\eta)$$ where:
\begin{tabbing}..................\= \kill
\> $C$ is a (compact, non-singular) algebraic curve \\
\> $p \in C$\\
\> $z$ is a  local (analytic or formal) coordinate at $p$\\
\> $L \in Pic C$\\
\> $\eta: L \otimes \hat{\O}_p \stackrel{\approx}{\rightarrow}
\hat{\O}_p \approx \C[[z]]$ is a (formal) trivialization of $L$
near $p$.\\
\end{tabbing}
If we fix $C,p$ and $z$, we think of the Krichever datum as
giving a point of
$${\cal U}_C(1,\infty p) :=
\displaystyle{\lim_{\stackrel{\leftarrow}{\ell}}} {\cal U}_C(1,\ell p).$$
The Krichever map
$$\{\mbox{Krichever\ data} \}\to Gr$$
sends
the above datum to the subspace $$W := \eta (H^0(C,L(\infty p)))
= \bigcup_k\eta (H^0(C,L(kp))) \subset \C((z)).$$
This subspace is comparable to $\C[z^{-1}]$, since it follows from
 Riemann-Roch that the dimension of
$H^0(L(kp))$ differs from $k$ by a bounded (and eventually
constant) quantity.\\

\noindent
\underline{The flows}

    Let's work with a coordinate $z$ which is analytic, i.e. it actually
converges on some disc. A line bundle $L$ on $C$ can be trivialized
(analytically)
on the Stein manifold $C\setminus p$.  We can think of
$(L,\eta)$ as being obtained from $\O_{C\setminus p}$ by glueing
it to $\hat{\O}_p$ via a $1$-cocycle, or transition function,
which should consist of an invertible function $g$ on a
punctured neighborhood of $p$ in $C$.  Conversely, we claim there is a  map:
$$   \exp{}  :  K \longrightarrow {\cal U}_C(1,\infty p)  ,  $$
$$                    f  \longmapsto            (L,\eta)                   . $$
For
$f \in \C(z)$,
this is defined by the above analytic gluing, using
$g:=\exp{f}$,
which is indeed analytic on a punctured neighborhood.  For
$f \in \C[[z]] \approx  \hat{\O}_p$,
on the other hand, we take
$ (L, \eta) := ({\O} , \exp{f} ).  $
These two versions agree on the intersection,
$ f \in \C[z]_{(0)} $,
so the map is uniquely defined as claimed.
(The bundles we get this way all have degree 0, but we can also obtain maps
$$   \exp{_{g_{0}}}  :  K \longrightarrow {\cal U}_C(1,d,\infty p)$$
to the moduli space of level-$\infty p$ line bundles of degree d, simply by
fixing a meromorphic function $g_0$ on a neighborhood of $p$ which has order
$d$ at $p$, and replacing the previous $g$ by $g_0 \exp{f}$.
We will continue to suppress the degree $d$ in our notation.)

Any $a \in K$ gives an additive flow on $K$, which at $f \in K$ is
$$ t \longmapsto f+ta .  $$
Under the composed map
$$  K \stackrel{\exp{_{g_{0}}}} {\longrightarrow}  {\cal U}_C(1,\infty p)
\stackrel{{Krichever}} {\longrightarrow}  Gr     ,                     $$
this is mapped to the double KP flow $KP_a$ on $Gr$.
For $a \in \C[[z]]$ this flow does not affect the isomorphism class of $L$,
and simply multiplies $\eta$ by $exp(ta)$.  On the other hand,
the $i^{th}$ $KP$ flow, given by $a = z^{-i}$, changes both $L$
and $\eta$ if $i > 0$.  The projection to $Pic\ C$ is a linear
flow, whose direction is the $i^{th}$ derivative at $p$, with
respect to the coordinate $z$, of the Abel-Jacobi map
$C \to Pic\ C$.  Dividing out the trivial flows corresponds to suppressing
$\eta$, so we obtain, for each $C,p,z$ and degree $d \in \Z$, a
finite-dimensional orbit of the $KP$ flows in $Gr / \C[[z]]^*$, isomorphic
to $Pic^dC$.\\

\noindent
\underline{Multi-Krichever data}

     Several natural generalizations of the Krichever map  to
the multi-component KP can be found in
\cite{adams-bergvelt,li-mulase-category} and elsewhere.
Here are some of the possibilities.  We can consider
 ``multi-Krichever'' data
$$(C,D,z_i,L,\eta)$$
involving a curve $C$ with a divisor $D$ consisting of $k$
distinct points $p_i(1 \leq i \leq k)$, a coordinate $z_i$
at each $p_i$, a line bundle $L$, and a formal trivialization
$\eta_i$ at each $p_i$.  Fixing $C,p_i$ and $z_i$, we have a
multi-Krichever map
$$\{\mbox{multi-Krichever\ data}\} \approx
{\cal U}_C(1,\infty D) \longrightarrow Gr_k$$
sending
$$(L,\eta_i) \mapsto W :=
(\eta_1,\cdots,\eta_k)(H^0(C,L(\infty D))) \subset
\C((z))^{\oplus k}.$$
The $k$-component KP flow  on the right hand side given by
$a=(a_1, \ldots, a_k) \in K^k$
restricts to the flow on the multi-Krichever data which multiplies the
transition function at $p_i$ (for an analytic trivialization of $L$ on
$C \setminus D$) by $\exp{a_i}$.

We can also consider "vector-Krichever" data
$(C,p,z,E,\eta)$
where the line bundle $L$ is replaced by a rank $n$ vector bundle $E$,  and
$$\eta: E \otimes \hat{\O}_p \stackrel{\approx}{\rightarrow}
(\hat{\O}_p)^n \approx (\C[[z]])^n$$
is now a (formal) trivialization of $E$ near $p$. Not too surprisingly,
the vector-Krichever map
$$\{\mbox{vector-Krichever\ data} \}\to Gr_n$$
sends the above datum to the subspace
$$W := \eta (H^0(C,E(\infty p))) \subset (\C((z)))^n.$$

In the next subsection we will see that the interesting interaction of these
two types of higher Krichever maps occurs not by extending further
(to objects such as $(C,D,z_i,E,{\eta}_i)$), but by restricting to those
vector data  on one curve which match some multi-data on another.\\

\noindent
\underline{KdV-type subhierarchies}

      Among the Krichever data one can restrict attention to those
quintuples where $z^{-n}$ (for some fixed $n$) happens to extend to a regular
function on $C \setminus p$, i.e. gives a morphism
$$f = z^{-n}: C \to \P^1$$
of degree $n$, such that the fiber $f^{-1}(\infty)$
is $n \cdot p_0$.  The Krichever map sends such data to the
$n^{th}$ $KdV$ hierarchy, the distinguished subvariety of $Gr$
(invariant under the (double) $KP$ flows) given by
$$ \mbox{KdV}_n := \{ W \in Gr \  |  \ z^{-n} W \subset W \} .$$
The corresponding subspace of $\Psi$ is
$$\{\L|\L^n = \L^n_+ \ \mbox{is\ a\ \underline{differential}\ operator}\}.$$

     Fixing a partition $\underline{n} = (n_1,\cdots,n_k)$ of
$n$, we can similarly consider the covering data of
type $\underline{n}$, consisting of the multi-Krichever data
$(C,p_i,z_i,L,\eta_i)$ plus a map $f: C \to \Sigma$ of degree $n$
to a curve $\Sigma$ with local coordinate $z$ at a point $\infty
\in \Sigma$, such that
$$f^{-1}(p) = \Sigma n_kp_i,\ \ f^{-1}(z) = z_i^{n_i}\ \mbox{at}\ p_i.$$

     Such a covering datum clearly gives a multi-Krichever datum on $C$, but
it also determines a vector-Krichever datum  $(E,\eta)$ on $\Sigma$:
The standard $m$-sheeted branched cover
$$\begin{array}{cccc}
f_m: & \C & \rightarrow & \C\\
&  \widetilde{z} & \mapsto & z = \widetilde{z}^m\\
\end{array}$$
of the $z$-line determines an isomorphism
 $$s_m: (f_m)_* \O \stackrel{\approx}{\rightarrow} \O^{\oplus m}$$
given by
$$\sum^{m-1}_{i=0} a_i(\widetilde{z}^m)\widetilde{z}^i \mapsto
(a_0(z),\cdots,a_{m-1}(z)).$$
To the covering datum above we
can then associate the rank-$n$ vector bundle $E := f_*L$ on $\Sigma$,
together with the trivialization at $p$ obtained by composing
$$\oplus_i f_{*,p_i}(\eta_i): (f_*L)_p
\stackrel{\sim}{\rightarrow} \oplus_i f_{*,p_i}(\O_{p_i})$$
with the isomorphisms
$$f_{*,p_i}(\O_{p_i}) \stackrel{\sim}{\rightarrow} \O^{\oplus n_i}_{p}$$
which are conjugates of the standard isomorphisms $s_{n_i}$ by the chosen
local coordinates $z,z_i$. Finally, we note that there are obvious geometric
flows on these covering data: $L$ and $\eta _i$ flow as before, while
everything else stays put.
The compatibility of the two types of higher Krichever data is expressed by
the commutativity of the diagram:

\begin{equation} \label{eq-diagram-of-krichever-maps}
\begin{array}{ccc}
\{\underline{n}-\mbox{covering\ data}\} & \approx
& \{ f:C \rightarrow \Sigma; \  p_i,z_i,L,\eta_i;  \  z \ | \ldots     \}\\
 \downarrow & & \downarrow\\
\{vector \ Krichever \ data \ on \ \Sigma \} & &
\{multi \ Krichever \ data \ on \ C \} \\
\parallel & & \parallel \\
\cup_{\Sigma,p,z}{\cal U}_\Sigma(n,\infty p) &
&  \cup_{C,D,z_i} {\cal U}_C(1,\infty D)\\
\downarrow & & \downarrow\\
Gr_n & \stackrel {m_{\underline{n}}} {\longrightarrow} & Gr_k\\
\downarrow & & \downarrow\\
Gr_n / Heis_{\underline{n}}^+ & \longrightarrow & Gr_k / (\C[[z]]^\times)^k.\\
\end{array}
\end{equation}

The mcKP flows on the bottom right pull back to the Heisenberg flows on the
bottom left, and to the geometric flows on the $\underline{n}$-covering data.\\

\subsection{Compatibility of hierarchies}
\label{sec-compatibility-of-heirarchies}

\bigskip
Fix a smooth algebraic curve $\Sigma$ of arbitrary genus and a point $P$
in it. The moduli space
$\HiggsModuli_D := \HiggsComponent_{\Sigma}(n,d,\omega(D))$
(see definition \ref{def-good-component-of-higgs-pairs})
can be partitioned into type loci. We
consider the Zariski dense subset consisting of
the union of  finitely many type loci
$\HiggsModuli_D^{\underline{n}}$ indexed
by partitions $\underline{n}$ of $n$. A Higgs pair in
$\HiggsModuli_D^{\underline{n}}$ has a spectral curve $C \rightarrow \Sigma$
whose ramification type over $P \in \Sigma$ is
$\underline{n}=(n_1,\dots,n_k)$.

Fix a formal local parameter $z$ on the base curve $\Sigma$ at $P$.
A Higgs pair in $\HiggsModuli_D^{\underline{n}}$
(or rather its spectral pair $(C,L)$, see proposition
\ref{prop-ordinary-spectral-construction-higgs-pairs})
can be completed
to an $\HeisN^+$-orbit of an $\underline{n}$-covering data
$(C,P_i,z_i,L,\eta_i) \rightarrow (\Sigma,P,z,E,\eta)$ in finitely many ways.
These extra choices form
a natural finite Galois cover $\CoverHiggsModuli_D^{\underline{n}}$
of each type locus $\HiggsModuli_D^{\underline{n}}$.
We obtain Krichever maps (see diagram
\ref{eq-diagram-of-krichever-maps})
from the Galois cover $\CoverHiggsModuli_D^{\underline{n}}$
to the quotients $Gr_n/\HeisN^+$ and $Gr_k/(\ComplexNumbers[[z]]^\times)^k$.
Both the mcKP and Heisenberg flows
pull back to the same geometric flow on the Galois cover.
It is natural to ask:

\begin{question} \label{question-compatibility}
{\bf (The compatibility question)}
Is the Heisenberg flow Poisson with respect to the
natural Poisson structure on $\HiggsModuli_D$?
\end{question}

The Compatibility Theorem \ref{thm-compatibility} and its
extension \ref{thm-compatibility-singular-case}
provide an affirmative answer. We factor the moment map
of the Heisenberg action
through natural finite Galois covers of the ramification
type loci in the space of characteristic polynomials
(equations (\ref{eq-the-jth-hamiltonian-of-the-ith-component}) and
(\ref{eq-the-moment-map-for-the-heisenberg-action})).

The compatibility naturally follows from the construction of the Poisson
structure via level structures. Recall the birational realization of
the moduli space $\HiggsComponent_{\Sigma}(n,d,\omega(lP+D))$
as a quotient of the cotangent bundle
$T^*\ModuliVB_{lP+D}$ of the moduli space
$\ModuliVB_{lP+D} := \ModuliVB_{\Sigma}(n,d,lP+D)$
of vector bundles with level structures
(Chapter \ref{ch5}). This realization is a finite dimensional approximation
of the limiting realization of the moduli space
\[
\HiggsModuli_{\infty{P}+D} :=
\lim_{l\rightarrow \infty} \HiggsModuli_{l{P}+D}
\]
as a quotient of (a subset of) the cotangent bundle
$T^*\ModuliVB_{\infty{P}+D}$.
The ramification type loci $\HiggsModuli_D^{\underline{n}}$,
their Galois covers $\CoverHiggsModuli_D^{\underline{n}}$ and the
infinitesimal $\HeisN$-action on $\CoverHiggsModuli_D^{\underline{n}}$
become special cases of those appearing in the construction of
section \ref{sec-type-loci}. The compatibility theorem follows from
corollary \ref{cor-hamiltonians-on-the-base}
accompanied by the concrete identification of the moment maps in our
particular example.

\bigskip
The rest of this section is organized as follows.
In section \ref{sec-galois-covers-and-relative-krichever-maps} we emphasize
the ubiquity of the setup of relative Krichever maps.
They can be constructed for any family ${\cal J} \rightarrow B$
of Jacobians of branched covers of a fixed base curve $\Sigma$.
The analogue of the compatibility question \ref{question-compatibility}
makes sense whenever the family ${\cal J} \rightarrow B$ is an
integrable system (see for example question
\ref{question-compatibility-for-mukai-system-on-elliptic-k3}).

Starting with section \ref{sec-compatibility-of-stratifications}
we concentrate on the moduli spaces of Higgs pairs.
Sections \ref{sec-compatibility-of-stratifications} and
\ref{sec-the-compatibility-thm-smooth-case} consider the case of
smooth spectral curves. Especially well behaved is the case where the point
$P\in\Sigma$ of the
$\underline{n}$-covering data is in the support of the polar divisor
$D$. In this case the symplectic leaves foliation of the moduli space of Higgs
pairs is a refinement of the type loci partition
(lemma \ref{lemma-compatibility-of-stratifications}).

Section \ref{sec-the-compatibility-thm-singular-case}
is a generalization to singular cases. As an example,
we consider in section
\ref{sec-elliptic-solitons} the Elliptic Solitons studied by
Treibich and Verdier. We conclude with an outline of the proof of the
compatibility theorem in section
\ref{sec-proof-of-compatibility-theorem}.

\noindent
{\em Note:} Type-$(1,1,\dots,1)$ relative Krichever maps
were independently considered by Y. Li and M. Mulase
in a recent preprint \cite{li-mulase-compatibility}.

\subsubsection{Galois covers and relative Krichever maps}
\label{sec-galois-covers-and-relative-krichever-maps}

Let $B_D := \oplus_{i=1}^{n}
H^{0}(\Sigma,(\omega_{\Sigma}(D))^{\otimes i})$ be the
space of characteristic polynomials. For simplicity, we restrict ourselves to
the Zariski open subset $Bsm_D$ of reduced and irreducible
$n$-sheeted spectral curves
$\pi:C\rightarrow \Sigma$ in $T^*_{\Sigma}(D)$ whose fiber over
$P\in \Sigma$ consists of
{\em smooth} points of $C$. Denote by $\Higgsm_{D}$ the corresponding open
subset of $\HiggsModuli_{D}$.
The ramification type stratification of $\Higgsm_{D}$ is induced by that
of $Bsm_D$
\[
Bsm_D = \DisjointUnion_{\underline{n}} Bsm_D^{\underline{n}}.
\]
Given a Higgs pair $(E,\varphi)$ in $\HiggsmN_{D}$
corresponding to a torsion free sheaf $L$ on a spectral cover
$C\rightarrow \Sigma$
we can complete it to an $\underline{n}$-covering data
\[
(C,P_i,z_i,L,\eta_i) \rightarrow (\Sigma,P,z,E,\eta)
\]
by choosing
i) a formal local parameter $z$ at $P$, ii) an $n_i$-th root $z_{i}$
of $\pi^*{z}$ at each point $P_i$ of $C$ over $P$ and
iii) formal trivializations $\eta_{i}$ of the sheaf $L$ at $P_i$.
The $\HeisN^+$-orbit of an $\underline{n}$-covering data
consists precisely of all possible choices of the $\eta_i$'s. Thus, for
fixed $P$ and $z$ only a finite choice is
needed in order to obtain the points of the quotients of the Grassmannians
$Gr_n/\HeisN^+$ and $Gr_k/(\ComplexNumbers[[z]]^\times)^k$
(see diagram \ref{eq-diagram-of-krichever-maps}).
These choices are independent of the sheaf $L$.
The choices are parametrized by the Galois cover
$\CoverBsmN_D \rightarrow \BsmN_D$
consisting of pairs $(C,\lambda)$ of a spectral curve $C$ in
$\BsmN$ and the discrete data $\lambda$ which amounts to:

\smallskip
\noindent
i) (Parabolic data)
An ordering $(P_1,P_2,\dots,P_k)$ of the points
(eigenvalues) in the fiber over
$P$ compatible with the fixed order of the ramification indices
$(n_1,\dots,n_k)$ (say, $n_1\leq n_2 \leq \dots \leq n_k$).

\noindent
ii) A choice of an $n_i$-th root $z_{i}$ of $\pi^*{z}$ at $P_i$.

\smallskip
Denote by $\CoverHiggsmN_{D}$ the corresponding Galois cover of $\HiggsmN_{D}$.
We get a canonical relative Krichever map
\begin{equation} \label{eq-relative-krichever-map}
\kappa_{\underline{n}}: \CoverHiggsmN_{D} \rightarrow Gr_n/\HeisN^+
\end{equation}
from the Galois cover to the quotient Grassmannian.

The Galois group of $\CoverBsmN_D \rightarrow \BsmN_D$
is the Weyl group $W_{\underline{n}} := N(\HeisN^+)/\HeisN^+$
of the maximal torus of the level infinity group $\LevelInfinityGroup$.
For example, $W_{(1,\dots,1)}$ is the symmetric group $S_n$, while
$W_{(n)}$ is the cyclic group of order $n$.
The discrete data $\lambda = [(P_1,P_2,\dots,P_k),(z_1,\dots,z_k)]$
of a point $(C,\lambda)$ in $\CoverBsmN_D$ is equivalent to a
commutative $\ComplexNumbers[[z]]$-algebras isomorphism
\begin{equation} \label{eq-isomorphism-heis-to-structure-sheaf}
\lambda:\heisN^+
\IsomRightArrow
\oplus_{i=1}^{k}\CompletedSheafOfAt{C}{(P_i)}
\end{equation}
of the torus algebra with the formal completion of
the structure sheaf $\StructureSheaf{C}$ at the fiber over $P$.
The inverse $\lambda^{-1}$ sends $z_i$ to the generator  of the $i$-th block
of the torus $\heisN^+$ given by
(\ref{eq-the-generator-of-the-ith-heisenberg-block}).
The finite Weyl
group $W_{\underline{n}}$ acts on $\heisN^+$, hence on $\lambda$, introducing
the $W_{\underline{n}}$-torsor structure on $\CoverBsmN_D$.
(See also lemma
\ref{lemma-two-type-loci-coincide} part
\ref{lemma-item-two-galois-covers-coincide}
for a group theoretic interpretation.)

\bigskip
We note that a Galois cover $\widetilde{B} \rightarrow B$ as above
can be defined for any family ${\cal J} \rightarrow B$
of Jacobians of a family ${\cal C} \rightarrow B$ of branched covers with a
fixed ramification type $\underline{n}$
of a fixed triple $(\Sigma,P,z)$. We obtain a relative Krichever map
\[
\kappa_{\underline{n}}: \widetilde{{\cal J}} \rightarrow Gr_n/\HeisN^+
\]
as above. The Heisenberg flow pulls back to an
infinitesimal action, i.e., a
Lie algebra homomorphism
\[
d\rho : \heisN \rightarrow V(\widetilde{{\cal J}})
\]
into a commutative algebra of vertical tangent vector fields. When
${\cal J}$ (and hence $\widetilde{{\cal J}}$) is an integrable system
we are led to ask the compatibility question
\ref{question-compatibility}: {\em is the action Poisson?}
{\em i.e., can $d\rho$ be lifted to a Lie algebra homomorphism}
\[
\mu^*_{\heisN}: \heisN \rightarrow
\Gamma(\StructureSheaf{\widetilde{B}}) \hookrightarrow
[\Gamma(\StructureSheaf{\widetilde{{\cal J}}}),\{,\}]?
\]
A priori, the vector fields $d\rho(a)$, $a\in\heisN$ may not even be
{\em locally} Hamiltonian.

Inherently nonlinear examples arise from
the Mukai-Tyurin integrable system of a family of Jacobians of
a linear system $B := {\bP}H^{0}(S,{\cal L})$ of curves  on a symplectic or
Poisson surface $S$ (see chapter \ref{ch8}). Consider for example the:

\begin{question} \label{question-compatibility-for-mukai-system-on-elliptic-k3}
Let $\pi : S \rightarrow \PiOne$  be an elliptic K3 surface and
${\cal L}$ a very ample line bundle on $S$.
Fix $P\in \PiOne$ and a local parameter $z$ and consider the Galois cover
$\widetilde{B}^{(1,1,\dots,1)}$ of the generic ramification type locus.
Is the Heisenberg action Poisson on
$\widetilde{{\cal J}} \rightarrow \widetilde{B}^{(1,1,\dots,1)}$
(globally over $\widetilde{B}^{(1,1,\dots,1)}$)?
\end{question}

The compatibility question has an intrinsic algebro-geometric formulation:
The $j$-th KP flow of the $P_i$-component
is the vector field whose direction along the fiber over
$(b,\lambda) \in \widetilde{B}^{(1,1,\dots,1)}$ is the $j$-th derivative
of the Abel-Jacobi map at $P_i$.
Using the methods of chapter \ref{ch8} it is easy to see that
the Heisenberg action is symplectic.
As we move the point $P(0):=P$ in ${\Bbb P}^1$ and its (Lagrangian) fiber
in $S$ we obtain an analytic (or formal) family of {\em Lagrangian} sections
${\cal AJ}(P_i(z))$ of
$\widetilde{{\cal J}} \rightarrow \widetilde{B}^{(1,1,\dots,1)}$
(see corollary \ref{cor-canonical-symplectic-str-on-albanese}).
Translations by the sections ${\cal AJ}(P_i(z))-{\cal AJ}(P_i(0))$ is a
family of symplectomorphisms of $\widetilde{{\cal J}}$.
Thus, the vector field corresponding to its $j$-th derivative with
respect to the local parameter $z$ is locally Hamiltonian.
%
%
%
%
It seems unlikely however that the Heisenberg flow integrates to a global
Poisson action for a general system as in question
\ref{question-compatibility-for-mukai-system-on-elliptic-k3}. It is
the {\em exactness} of the symplectic structure in a neighborhood of the
fiber over $P$ in $T^*\Sigma$ which lifts the infinitesimal symplectic
Heisenberg action to a Poisson action in the Hitchin's system case
(see equation (\ref{eq-the-moment-map-for-the-heisenberg-action})).

\subsubsection{Compatibility of stratifications}
\label{sec-compatibility-of-stratifications}

Prior to stating the compatibility theorem
\ref{thm-compatibility} we need to examine the Poisson nature of the
Galois covers $\CoverHiggsmN_{D}$.
$\Higgsm^{(1,1,\dots,1)}_D$ is a Zariski {\em open Poisson} subvariety
of $\Higgsm_D$. Hence, the unramified Galois cover
$\CoverHiggsm^{(1,1,\dots,1)}_D$ is endowed with the canonical pullback
Poisson structure.

Though non generic, the other type strata are as important.
The cyclic ramification type $(n)$, for example, corresponds to the single
component KP-hierarchy (see \cite{segal-wilson-loop-groups-and-kp}).
When the point $P$ of the $\underline{n}$-covering data
is in the support of the divisor $D$, we obtain a {\em strict compatibility}
between the $P$-ramification type stratification of
$\Higgsm_D$ and its symplectic leaves foliation. All Galois covers
$\CoverHiggsmN_D$, $P \in D$ are thus endowed with the canonical pullback
Poisson structure:

\begin{lem} \label{lemma-compatibility-of-stratifications}
(conditional compatibility of stratifications)
When the point $P \in \Sigma$ is in the support of $D$,
the symplectic leaves foliation is a {\em refinement} of the ramification type
stratification $\Higgsm_D = \DisjointUnion_{\underline{n}} \HiggsmN_D$.
\end{lem}

\noindent
{\bf Proof:} We need to show that the ramification type $\underline{n}$ of the
spectral cover $\pi:C\rightarrow \Sigma$ over $P\in D$ is fixed throughout the
symplectic leaf $\Higgsm_S \subset \Higgsm_D$ of a Higgs pair
$(E_{0},\varphi_{0})$.
The symplectic leaves of $\HiggsModuli_D$ are determined by
coadjoint orbits of the level $D$ algebra $\LieAlg{g}_D$.
The coadjoint orbit $S$ is determined by the residue
of the Higgs field, namely, the infinitesimal data $(\restricted{E}{D},
\restricted{\varphi}{D})$ encoded in the value of $\varphi$ at $D$ (see
\cite{markman-higgs} Remark 8.9 and Proposition 7.17)). Thus,
the Jordan type of the Higgs field $\varphi$ at $P$
is fixed throughout $\Higgsm_S$. In general,
(allowing singularities over $P$), the Jordan type
depends both on the ramification type of $C\rightarrow \Sigma$ and
on the sheaf $L$ on $C$ corresponding to the Higgs pair.
For Higgs pairs $(E,\varphi)$ in
$\Higgsm_D$, the spectral curve $C$ is smooth over $P$,
hence, its ramification type coincides with the Jordan type of $\varphi$ at
$P$.
\EndProof

\subsubsection{The compatibility theorem, the smooth case}
\label{sec-the-compatibility-thm-smooth-case}

We proceed to introduce the moment map
of the infinitesimal Poisson action
\[d\rho : \heisN \rightarrow V(\CoverHiggsmN_D).\]
Throughout the end of
this subsection
we will assume the
\begin{condition}\label{condition-ramification-type}
Ramification types $\underline{n}$ other than $(1,1,\dots,1)$
are considered only if $P$ is in the support of $D$.
\end{condition}
This condition will be relaxed later by conditions
\ref{condition-only-generic-type} or
\ref{condition-only-singularities-which-are-resolved-by-spectral-sheaf}.

Let $b \in \BsmN_D$ be the polynomial of the spectral curve
$\pi_b:C_b \rightarrow \Sigma$.
Recall that spectral curves are endowed with a tautological
meromorphic $1$-form $y_b \in H^0(C_b,\pi_b^*\omega_\Sigma(D))$ with poles
over $D \subset \Sigma$
(see section \ref{sec-spectral-curves-and-the-hitchin-system}).
Let $\phi^j_{P_i}$ be the function on $\CoverBsmN_D$ given at
a pair $(b,\lambda) \in \CoverBsmN_D$ by
%
%
%
%
\begin{equation} \label{eq-the-jth-hamiltonian-of-the-ith-component}
\phi^j_{P_i}(b,\lambda) := Res_{P_i}((z_i)^{-j}\cdot y_b).
\end{equation}
The  Lie algebra homomorphism $\mu_{\heisN}^*$ sends the inverse
(in $\heisN$) of the generator of the $i$-th block
of the torus $\heisN^+$ given in
(\ref{eq-the-generator-of-the-ith-heisenberg-block})
to the function $\phi^1_{P_i}\circ \widetilde{{\rm char}}$ on
$\CoverHiggsmN_D$. In other words,
$\mu_{\heisN}^* :\heisN \rightarrow
[\Gamma(\StructureSheaf{\CoverHiggsmN_D}),\{,\}]$
factors
as a composition $\phi \circ \widetilde{{\rm char}}^*$ through
$\Gamma(\StructureSheaf{\CoverBsmN_D})$.
If we regard $\lambda$ also as an isomorphism from $\heisN$ to
$\oplus_{i=1}^{k}\CompletedSheafOfAt{C}{(P_i)}$
(eq. (\ref{eq-isomorphism-heis-to-structure-sheaf})), then
$\phi:\heisN \rightarrow \Gamma(\StructureSheaf{\CoverBsmN_D})$
is given by
%
%
%
%
\begin{equation} \label{eq-the-moment-map-for-the-heisenberg-action}
(\phi(a))(b,\lambda) =
\sum_{\{P_i\}} Res_{P_i}(\lambda(a)\cdot y_{b}),
\ \ \ a\in\heisN.
\end{equation}

\begin{theorem} \label{thm-compatibility}
({\bf The Compatibility Theorem, smooth case})
(Assuming condition \ref{condition-ramification-type})
The relative Krichever map
\[
\kappa_{\underline{n}} : \CoverHiggsmN_D \rightarrow Gr_n/\HeisN^+
\]
intertwines the Heisenberg flow on $Gr_n/\HeisN^+$ (and the mcKP flow
on $Gr_k/(\ComplexNumbers[[z]]^\times)^k$)
with an infinitesimal Poisson action of the maximal torus $\heisN$
on $\CoverHiggsmN_D$. The latter is induced by the Lie algebra homomorphism
\[
\mu_{\heisN}^* = \phi \circ \widetilde{{\rm char}}^* :\heisN
\RightArrowOf{\phi}
\Gamma(\StructureSheaf{\CoverBsmN_D})
\HookRightArrowOf{\widetilde{{\rm char}}^*}
[\Gamma(\StructureSheaf{\CoverHiggsmN_D}),\{,\}]
\]
which factors through the homomorphism $\phi$
given by (\ref{eq-the-moment-map-for-the-heisenberg-action}).
\end{theorem}

\begin{rems}
{\rm
\begin{enumerate}
\item \label{rems-item-nonvanishing-casimirs}
The subalgebra $\heisN^+$ acts trivially on $Gr_n/\HeisN^+$ hence
also on $\CoverHiggsmN_D$. This corresponds to the fact that the functions
$\phi^j_{P_i}$, $j\leq 0$ are Casimir. Indeed, if $P$ is not contained in $D$,
the $1$-form $y_b$ is holomorphic at the fiber over $P$ and
$\phi^j_{P_i}$ is identically zero for $j\leq 0$.
If $P$ is in $D$ then the finite set of
non-zero $\phi^j_{P_i}$, indexed by finitely many non-positive integers $j$,
are among the Casimirs that induce the highest rank symplectic leaves foliation
(see \cite{markman-higgs} Proposition 8.8).

\item
The multi-Krichever map $\kappa_{\underline{n}}$ depends
on auxiliary parameters $P$ and $z$.
In other words, it lives naturally on an infinite dimensional space
$\cup_{P,z} \ \CoverHiggsmN_{D,P,z}$ in which $P$ and $z$ are allowed to
vary. This is not as bad as it might seem, since
the $j$-th flow on $\CoverHiggsmN_D$ really depends
only on our finite dimensional choices of $P$ and the $j$-th order germ
of $z$ there. Similarly, our
Hamiltonians $\phi^j_{P_i}$ depend at most on the $(j+n\deg D)$-th germ of $z$,
the shift arising,
as in part \ref{rems-item-nonvanishing-casimirs} of this remark,
from the possible poles above $P$ of the tautological
$1$-form on the spectral curve.
So we may think of $\cup_{P,z} \ \CoverHiggsmN_{D,P,z}$
as the inverse limit of a family
of finite dimensional moduli spaces, indexed by the level. Each KP flow
or Hamiltonian is defined for sufficiently high level.
\end{enumerate}
}
\end{rems}

\subsubsection{The compatibility theorem, singular cases}
\label{sec-the-compatibility-thm-singular-case}

The condition that the fiber over $P$ of the embedded spectral curve
be smooth is too restrictive. The embedded spectral data
$(\bar{C} \subset T^*_{\Sigma}(D),\bar{L})$
of a Higgs pair $(E,\varphi)$ may have singularities
over $P$ which are canonically resolvable. The point is
that the rank $1$ torsion free sheaf $\bar{L}$ on $\bar{C}$ determines
a partial normalization $\nu: C \rightarrow \bar{C}$
and a unique rank $1$ torsion free sheaf $L$ on $C$ such that
i) $\bar{L}$ is isomorphic to the direct image $\nu_{*}L$, and
ii) $L$ is locally free at the fiber over $P$.
We are interested
in those Higgs pairs for which the fiber of $C$ over $P$ is smooth.
Such data may also be completed in finitely many
ways to an $\underline{n}$-covering data as in section
\ref{sec-galois-covers-and-relative-krichever-maps}.

\begin{definition}\label{def-non-essential-singularities}
{\rm
The singularities over $P$ of
a spectral pair $(\bar{C} \subset T^*_{\Sigma}(D),\bar{L})$
are said to be {\em resolved by the spectral sheaf} $\bar{L}$
if i) $\bar{C}$ is irreducible and reduced. ii) The sheaf $\bar{L}$
is the direct image $\nu_{*}L$ of a rank $1$ torsion free sheaf
$L$ on the normalization $\nu:C\rightarrow \bar{C}$ of the
fiber of $\bar{C}$ over $P$.
}
\end{definition}

Fixing a symplectic leaf $\HiggsModuli_S$ we may consider the type sub-loci
in the locus of Higgs pairs whose spectral curve has at worst singularities
over $P$ which are resolved by the spectral sheaf.
The topology of these type loci is quite complicated. As a result,
the Galois covers of these type loci may not have a symplectic structure.
Nevertheless, the construction of section \ref{sec-type-loci}, as used in
section \ref{sec-proof-of-compatibility-theorem},
provides  {\em canonical embeddings} of the Galois covers of these
type loci in (finite dimensional) symplectic varieties. These embeddings
realize the Heisenberg flow as a Hamiltonian flow.

Control over the topology is regained below by restraining the singularities.
In the smooth case (section \ref{sec-the-compatibility-thm-smooth-case})
it is the smoothness which assures that the generic ramification type locus
in a symplectic leaf $\HiggsModuli_S$ is open (rather than only
Zariski dense). The point is that
degenarations from a ramification type $\underline{n}$ through other types
back to type $\underline{n}$ must end with a singular fiber over $P$
(and are thus excluded).
If $P \in D$,
there are symplectic leaves $\HiggsModuli_S$, $S\subset \LieAlg{g}_D^*$ of
$\HiggsModuli_D$ for which
the singularities over $P$
are encoded in the infinitesimal
data associated to the coadjoint orbit $S$ and shared by the generic Higgs
pair in $\HiggsModuli_S$.
Often, this is an indication that the Poisson
surface $T^*_{\Sigma}(D)$ is not the best to work with. Moreover,
a birational transformation
$T^*_{\Sigma}(D) \DotRightArrow X_S$, centered at points of the fiber over $P$
and encoded in $S$, can simultaneously resolve the singularities (over $P$) of
the spectral curve of the generic point in $\HiggsModuli_S$
(see example \ref{example-the-coadjoint-orbit-of-elliptic-solitons}). In such a
case,
smoothness of the proper transform of the spectral curve in $X_S$ at points of
the fiber over $P$ is an {\em open} condition and the corresponding locus
in $\HiggsModuli_S$ with generic ramification type is {\em symplectic}.
When the multiplicity of $P$ in $D$ is greater than $1$ the correspondence
between the coadloint orbits $S\subset \LieAlg{g}_D^*$ and their surfaces $X_S$
can be quite complicated. Instead of working the correspondence out,
we will use the following notion of $S$-smoothness to assure
(see condition \ref{condition-only-generic-type}) that the
generic ramification type locus
in a symplectic leaf $\HiggsModuli_S$ is open (rather than only
Zariski dense).

\begin{definition}\label{def-S-smoothness-over-D}
{\rm
Let $S$ be a coadjoint orbit of $\LieAlg{g}_D$.
An irreducible and reduced spectral curve $\pi:\bar{C}\rightarrow\Sigma$ is
{\em $S$-smooth over $D$} if
i) a line bundle $L$ on the resolution of the singularities
$\nu: C \rightarrow \bar{C}$ of the fibers over $D$ results in a Higgs pair
$(E,\varphi) := \pi\circ\nu_*(L,\otimes \nu^*(y))$ in $\HiggsModuli_S$.
and
ii) the arithmetic genus
of the normalization $C$ above is equal to half
the dimension of the symplectic leaf $\HiggsModuli_S$.
}
\end{definition}

If $\bar{C}$ is an irreducible and reduced
spectral curve then, by the construction of \cite{simpson-moduli},
the fiber of the characteristic polynomial map
in $\HiggsComponent_{\Sigma}(n,d,\omega(D))$ is
its compactified Jacobian, the latter being the moduli space of all rank $1$
torsion free sheaves on $\bar{C}$ with a fixed Euler characteristic.
The compactified Jacobian is known to be irreducible for irreducible and
reduced curve on a surface (i.e., with planar singularities,
\cite{a-i-k}). Moreover, the symplectic leaf $\HiggsModuli_S$ intersects
the compactified Jacobian of $\bar{C}$
in a union of strata determined by partial normalizations. If
$(E,\varphi)$ in $\HiggsModuli_S$ corresponds to $(\bar{C},\bar{L})$ and
$\nu:(C,L)\rightarrow (\bar{C},\bar{L})$ is a partial normalization
where $L$ is a locally free sheaf on $C$,
then any twist $F$ of $L$ by a  locally free sheaf in $Pic^0(C)$
(the component of $\StructureSheaf{C}$) will result in a
Higgs pair
$(E',\varphi') := (\pi\circ\nu_*(F),\pi_*(\otimes y))$ in $\HiggsModuli_S$.
The point is that the residue of $(E',\varphi')$
(with respect to any level-$D$ structure)
will be in the same coadjoint orbit as that of $(E,\varphi)$.
$S$-smoothness of $\bar{C}$ over $D$ is thus
equivalent to the geometric condition:


\bigskip
\noindent
{\em
The fiber of the characteristic polynomial map over $\bar{C}$
intersects $\HiggsModuli_S$ in a
Lagrangian subvariety
isomorphic to
the compactified Jacobian of the resolution $C$ of the singularities of
$\bar{C}$ over $D$.
}

The following example will be used in section \ref{sec-elliptic-solitons}
to describe a symplectic leaf which parametrizes Elliptic solitons.

\begin{example} \label{example-the-coadjoint-orbit-of-elliptic-solitons}
{\rm
Let $D = P$, $S \subset \LieAlg{g}^*_D \cong \LieAlg{gl(n)}^*
\cong \LieAlg{gl(n)}$ be the coadjoint orbit containing the diagonal
matrix
\[ A = \left(
\begin{array}{cccccc}
-1     &     0     &       & \dots    &    & 0 \\
0      &    -1     &       &          &    &   \\
       &           & \cdot &          &    & \vdots \\
\vdots &           &       & \cdot    &    &   \\
       &           &       &          & -1 & 0 \\
0      &           & \dots &          & 0  & n-1
\end{array} \right)
\]
$S$ is the coadjoint orbit of lowest dimension with characteristic
polynomial $(x+1)^{n-1}(x-(n-1))$. Its dimension $2n-2$ is
$(n-2)(n-1)$ less than the
generic rank of the Poisson structure of $\LieAlg{gl(n)}^*$.
If non-empty, each component of
$\HiggsModuli_S$ is a smooth symplectic variety of dimension
$\dim \HiggsModuli_{P}-(n-1)-(n-2)(n-1) = [n^2(2g-1)+1]-(n-1)^2$
(see theorem \ref{thm-markman-botachin}
and \cite{markman-higgs} proposition 7.17).
The  spectral curves $\bar{C} \subset T^*_{\Sigma}(P)$ which are
$S$-smooth over $P$
will have two points in the fiber over $P$,
one smooth at the eigenvalue with residue $n-1$
and one (singular if $n\geq 3$) at the eigenvalue with residue $-1$.
Assume $n \geq 3$.
The resolution $C$ of the singularity of a typical (though not
all) such $\bar{C}$ will be unramified over $P$ with $n-1$ points
collapsed to one in $\bar{C}$. The sheaf $\bar{L}$
(of a Higgs pair in $\HiggsModuli_S$ with spectral curve $\bar{C}$)
will be a pushforward of a
torsion free sheaf $L$ on $C$. In contrast, line bundles on that $\bar{C}$
will result in Higgs pairs in another symplectic leaf
$\HiggsModuli_{S_{reg}}$ corresponding to the {\em regular} coadjoint orbit
$S_{reg}$ in $\LieAlg{gl(n)}^*\cong \LieAlg{gl(n)}$ with characteristic
polynomial $(x+1)^{n-1}(x-(n-1))$.
These Higgs pairs will {\em not} be $S_{reg}$-smooth.
$S_{reg}$-smoothness coincides with usual smoothness of
the embedded spectral curve which is necessarily ramified
with ramification index $n-1$ at the point with residue $-1$ over $P$.

The $S$-smooth spectral curves will be smooth on the blowup $\widehat{X_S}$
of $T^*_{\Sigma}(P)$ at residue $-1$ over $P$. If we blow down
in $\widehat{X_S}$ the proper transform of the fiber of  $T^*_{\Sigma}(P)$
we get a surface $X_S$ with a marked point $x_{n-1}$ over $P$.
An $S$-smooth spectral curve $\bar{C}$
will correspond to a curve on $X_S$ through $x_{n-1}$. It will be smooth at
$x_{n-1}$ if in addition it is of ramification type $(1,1,\dots,1)$ over $P$.
Consider the compactification
${\Bbb P}(T^*_{\Sigma}(P)\oplus \StructureSheaf{\Sigma})$ of
$T^*_{\Sigma}(P)$.
Blowing it up and down as above we get a ruled
surface $\bar{X}_S$ over $\Sigma$ which is isomorphic to the projectivization
${\Bbb P}W$ of the unique nontrivial extension
\[
0 \rightarrow \omega_{\Sigma} \rightarrow W \rightarrow
\StructureSheaf{\Sigma} \rightarrow 0.
\]
In particular, the surface $\bar{X}_S$ is
{\em independent of the point $P$}.
The point is that blowing up and down the ruled surface ${\Bbb P}V :=
{\Bbb P}(T^*_{\Sigma}(P)\oplus \StructureSheaf{\Sigma})$
at residue $-1$ over $P$ results with
the ruled surface of a rank $2$ vector bundle
$W$ whose sheaf of sections is a subsheaf of $V :=
T^*_{\Sigma}(P)\oplus \StructureSheaf{\Sigma}$.
This subsheaf
consists of all sections which restrict at $P$ to the
subspace of the fiber $\restricted{V}{P}$ spanned by $(-1,1)$ (i.e., $W$ is
a Hecke transform of $V$, see \cite{tyurin-vb-survey}).
Clearly, $\omega_{\Sigma}$
is a subsheaf of $W$ and the  quotient $W / \omega_{\Sigma}$ is
isomorphic to $\StructureSheaf{\Sigma}$.
The resulting extension is non-trivial because $H^0(\Sigma,W)$ and
$H^0(\Sigma,\omega_{\Sigma})$ are equal as subspaces of $H^0(\Sigma,V)$.
}
\end{example}

\begin{rem} \label{rem-resolved-by-the-spectral-sheaf-vs-S-smoothness}
{\rm
$S$-smoothness over $D$ is stronger than having singularities over $D$
which are resolved by the spectral sheaf. They differ when the singularity
appears at an infinitesimal germ of too high an order to be detected by $S$.
E.g.,
take $n=3$ in example \ref{example-the-coadjoint-orbit-of-elliptic-solitons}
and consider a pair $(\bar{C},\bar{L})$ with a tacnode at
residue $-1$ over $P$ (two branches meet with a common tangent).
The arithmetic genus of the normalization
$\nu:C\rightarrow \bar{C}$ of the fiber over $P$ will drop by $2$ while
the pushforward $\bar{L}:=\nu_*(L)$ of a line bundle $L$ on $C$
will belong to a symplectic leaf $\HiggsModuli_{S}$ whose rank is $2$ less
than the maximal rank (rather than $4$). Hence $(\bar{C},\bar{L})$
is $S$-singular.
}
\end{rem}

We denote by
$\Higgsm_{S/D}$ the subset of $\HiggsModuli_S \subset \HiggsModuli_D$
parametrizing Higgs pairs whose spectral curve is $S$-smooth over $D$.

Unfortunately, the compatibility of stratifications (lemma
\ref{lemma-compatibility-of-stratifications})
does not extend to the $S$-smooth case. To overcome this
inconvenience we may either
assume condition \ref{condition-only-generic-type} or condition
\ref{condition-only-singularities-which-are-resolved-by-spectral-sheaf}.

\begin{condition}\label{condition-only-generic-type}
Consider the ramification locus $\Higgsm_{S/D}^{\underline{n}}$
in a component of $\Higgsm_{S/D}$ only if it is the generic ramification type
in this component.
\end{condition}

Note that $\Higgsm_{S/D}$ parametrizes
only Higgs pairs whose spectral curve is $S$-smooth over $D$.
$S$-smoothness over $D$ assures that if the type $\underline{n}$ is a
generic ramification type in a component of
$\Higgsm_{S/D}$, then the corresponding component of $\HiggsmN_{S/D}$
is an {\em open} subset of $\Higgsm_{S/D}$ (i.e., it excludes degenerations
of type $\underline{n}$ Higgs pairs through other types back to
type $\underline{n}$). In particular, these components of
$\HiggsmN_{S/D}$ are {\em symplectic}. Alternatively, we may
relax condition \ref{condition-only-generic-type}
even further at the expense of losing the symplectic nature of the loci and
having to resort to convention \ref{convention-abused-hamiltonian-language}:

\begin{condition}
\label{condition-only-singularities-which-are-resolved-by-spectral-sheaf}
Consider only Higgs pairs with a spectral curve whose singularities
over $P$ are resolved by its spectral sheaf. Adopt convention
\ref{convention-abused-hamiltonian-language}.
\end{condition}

\begin{theorem} \label{thm-compatibility-singular-case}
The compatibility theorem \ref{thm-compatibility} holds
for:
i) The $W_{\underline{n}}$-Galois covers of the type loci in $\Higgsm_{S/D}$
satisfying the genericity condition
\ref{condition-only-generic-type} (instead of condition
\ref{condition-ramification-type}).
ii) The $W_{\underline{n}}$-Galois covers of the locus in
$\HiggsModuli^{\underline{n}}_D$ consisting of Higgs pairs
satisfying condition
\ref{condition-only-singularities-which-are-resolved-by-spectral-sheaf}.
In ii) however we adopt convention
\ref{convention-abused-hamiltonian-language}.
\end{theorem}

\subsubsection{Elliptic Solitons}
\label{sec-elliptic-solitons}
In this subsection we illustrate the possibilities in the singular case
with a specific example.
Let $\Sigma$ be a smooth  elliptic curve. A {\em $\Sigma$-periodic
Elliptic KP soliton} is
a finite dimensional solution to the KP hierarchy, in which the orbit of
the first KP equation is isomorphic to $\Sigma$.
Its Krichever data
$(C,\tilde{P},\frac{\partial}{\partial z},L)$
consists of a reduced and irreducible curve $C$,
a smooth point $\tilde{P}$, a nonvanishing tangent vector at $\tilde{P}$
and a rank $1$ torsion free sheaf $L$ on $C$ of Euler characteristic $0$
(we suppress the non essential formal trivialization $\eta$ and consider
only the first order germ of $z$ which is equivalent to choosing a
nonzero tangent vector at $P$). We will denote
the global vector field extending $\frac{\partial}{\partial z}$ also
by $\frac{\partial}{\partial z}$.
The periodicity implies that
the image of the
Abel Jacobi map $AJ: C  \hookrightarrow J_C$, $Q \mapsto Q-\tilde{P}$
is tangent at $0$
to a subtorus isomorphic to $\Sigma$. Composing the Abel Jacobi map
with projection to $\Sigma$ we get a
{\em tangential morphism}
$\pi:(C,\tilde{P}) \rightarrow (\Sigma,P)$.
Its degree $n$ is called the {\em order} of the Elliptic soliton.
In general, a tangential morphism  $\pi:(C,\tilde{P}) \rightarrow (\Sigma,P)$
is a morphism with the property that $AJ(C)$ is tangent at $AJ(\tilde{P})$ to
$\pi^*J^0_\Sigma$.  Notice that composing a tangential morphism
$\pi:(C,\tilde{P}) \rightarrow (\Sigma,P)$ with a normalization
$\nu:\tilde{C} \rightarrow C$ results in a tangential morphism.
A tangential morphism is called {\em minimal} if it does not factor through
another tangential morphism.

The KP elliptic solitons enjoyed a careful and detailed study by
A. Treibich and J.-L. Verdier in a series of beautiful papers
(e.g., \cite
{treibich-verdier-elliptic-solitons,treibich-verdier-krichever-variety}).
Their results fit nicely with our picture:

\begin{theorem} \label{thm-variety-of-kp-solitons}
The variety of Krichever data of Elliptic KP solitons of order $n$
with a fixed pointed elliptic curve and a tangent vector
$(\Sigma,P,\frac{\partial}{\partial z})$ is canonically birational
to the divisor of traceless Higgs pairs
in the symplectic leaf $\HiggsModuli_S$  of
$\HiggsComponent_{\Sigma}(n,0,\omega_{\Sigma}(P))$
corresponding to the coadjoint orbit $S$ of example
\ref{example-the-coadjoint-orbit-of-elliptic-solitons}.
The KP flows are well defined on $\Higgsm_{S/P}\subset \HiggsModuli_S$
as the Hamiltonian vector fields of the functions
$\phi^j_{\tilde{P}}$ given in
(\ref{eq-the-jth-hamiltonian-of-the-ith-component}).
\end{theorem}

\medskip
Note that
if non-empty (which is the case)
$\HiggsModuli_S$ is $2n$-dimensional (see example
\ref{example-the-coadjoint-orbit-of-elliptic-solitons}).
The correspondence between tangential covers and
spectral covers is a corollary of
the following characterization of tangential covers due to
I. M. Krichever and A. Treibich.
For simplicity we consider only the case
in which the tangency point $\tilde{P}$ is not a ramification point of $\pi$.

\begin{theorem} \label{thm-characterization-of-tangential-covers}
\cite{treibich-verdier-krichever-variety}
Assume that $\pi:C\rightarrow \Sigma$ is unramified at $\tilde{P}$.
Then $\pi$ is tangential if and only if there exists a section
$y \in H^0(C,\pi^*[\omega_{\Sigma}(P)])$ satisfying:

\smallskip
\noindent
a) Near a point of $\pi^{-1}(P) - \tilde{P}$
(away from the tangency point $\tilde{P}$), $y- \pi^*(dz/z)$ is a holomorphic
multiple of
$\pi^*(dz)$, where $z$ is a local parameter at $P$.
(If $\pi:C\rightarrow \Sigma$ is unramified over $P$, this is equivalent to
saying that the residues $Res_{P_i}(y)$ are the same at all  $P_i$
other than $\tilde{P}$ in the fiber
over $P \in \Sigma$).

\noindent
b) The residue $Res_{\tilde{P}}(y)$ at $\tilde{P}$ does not vanish
if $n \geq 2$.
\end{theorem}

It follows by the residue theorem that there is a unique such section
which has residue $n-1$ at $\tilde{P}$ and which is moreover
traceless $tr(y) = 0 \in H^0(\Sigma,\omega_{\Sigma}(P))$.
Let $dz$ be a global holomorphic non zero $1$-form on $\Sigma$.
The function $k := y/\pi^*(dz)$ is called a {\em tangential function}
in  \cite{treibich-verdier-krichever-variety}. It was also proven
that a tangential morphism of order $n$ has arithmetic genus  $\leq n$
and is minimal if and only if its arithmetic genus is $n$.
(\cite{treibich-verdier-elliptic-solitons} Corollaire 3.10).

\bigskip
\noindent
{\bf Sketch of proof of theorem
\ref{thm-characterization-of-tangential-covers}:}
(for $C$ smooth, $\pi:C\rightarrow \Sigma$ unramified over $P$.)

\noindent
\underline{Step 1}: (Cohomological identification of the differential of
the Abel-Jacobi map)

The differential
$dAJ : T_{Q}C \cong H^{0}(Q,\StructureSheaf{Q}(Q)) \rightarrow
H^{1}(C,\StructureSheaf{C})$
of the Abel-Jacobi map at $Q \in C$ is identified
as the connecting homomorphism of the short exact sequence
\[
0 \rightarrow
\StructureSheaf{C} \rightarrow
\StructureSheaf{C}(Q) \rightarrow
\StructureSheaf{Q}(Q) \rightarrow 0.
\]
Similarly, the differential $d(AJ \circ \pi^{-1})$
of the composition
$
\Sigma \HookRightArrowOf{\pi^{-1}} \Sym^{n}C \RightArrowOf{AJ} J^{n}_C
$
is given by
\[
T_{P}\Sigma \cong H^{0}(P,\StructureSheaf{P}(P))
\HookRightArrowOf{\pi^*}
H^{0}(\pi^{-1}(P),\StructureSheaf{\pi^{-1}(P)}(\pi^{-1}(P))) \IsomRightArrow
T_{[\pi^{-1}(P)]}\Sym^{n}C \RightArrowOf{dAJ}
H^{1}(C,\StructureSheaf{C}),
\]
where the composition of the last two arrows is the
connecting homomorphism of the short exact sequence
\[
0 \rightarrow
\StructureSheaf{C} \rightarrow
\StructureSheaf{C}(\pi^{-1}(P)) \rightarrow
\StructureSheaf{\pi^{-1}(P)}(\pi^{-1}(P)) \rightarrow 0.
\]

\smallskip
\noindent
\underline{Step 2}: (residues as coefficients in a linear dependency
of tangent lines)

Clearly, the tangent line to $(AJ \circ \pi^{-1})(\Sigma)$ at the
image of $P$ is in the span of the tangent lines to $AJ(C)$ at
the points $P_i$. If, in addition, $\pi$ is tangential with tangency point
$\tilde{P}\in C$, then
the tangent lines to $AJ(C)$ at the points in the fiber over $P$
are linearly dependent.
If $\pi:C\rightarrow \Sigma$ is unramified over $P$ we can write
these observations in the form of two linear equations:
\begin{equation} \label{eq-linear-condition-tangentiality}
\sum dAJ_{P_i}(\frac{\partial}{\partial z_i}) =
d(AJ\circ \pi^{-1})_{P}(\frac{\partial}{\partial z}),
\end{equation}
and
\begin{equation} \label{eq-linear-dependency}
\sum a_i dAJ_{P_i}(\frac{\partial}{\partial z_i}) = 0
\ \ \ \mbox{linear} \ \  \mbox{dependency}.
\end{equation}

Above, $\frac{\partial}{\partial z_i}$ is the lift of
$\frac{\partial}{\partial z}$, i.e.,
$d\pi(\frac{\partial}{\partial z_i})=\frac{\partial}{\partial z}$.
We claim that the coefficients $a_i$ in (\ref{eq-linear-dependency})
are residues of a meromorphic $1$-form $y$ at the points of the fiber.
More precisely we have:
\begin{lem} \label{lemma-residues-vs-differential-of-abel-jacobi}
Assume that $\pi:C\rightarrow \Sigma$ is unramified over $P$.
There exists a section $y \in H^0(C,\pi^*\omega_{\Sigma}(P))$
with residues $(a_1,a_2,\dots,a_n)$ at the fiber over $P$ if
and only if $(a_1,a_2,\dots,a_n)$ satisfy equation
(\ref{eq-linear-dependency}).
\end{lem}

\noindent
{\bf Proof:}
The global tangent vector field
$\frac{\partial}{\partial z} \in H^{0}(\Sigma,T\Sigma)$
gives rise to the commutative diagram:
\begin{equation}\label{diag-residues-and-differential-of-abel-jacobi}
{\divide\dgARROWLENGTH by 4
\begin{diagram}
\node{H^0(C,\pi^*\omega_{\Sigma}(P))}
\arrow{s,lr}{\cong}{\rfloor\frac{\partial}{\partial z}}
\arrow{e}
\node{H^0(\pi^{-1}(P),\restricted{\pi^{*}\omega_{\Sigma}(P)}{\pi^{-1}(P)})}
\arrow{s,lr}{\cong}{\rfloor\frac{\partial}{\partial z}}
\arrow{e}
\node{H^1(C,\pi^*\omega_{\Sigma})}
\arrow{s,lr}{\cong}{\rfloor\frac{\partial}{\partial z}}
\\
\node{H^0(C,\StructureSheaf{C}(\pi^{-1}(P)))}
\arrow{e}
\node{H^0(\pi^{-1}(P),\StructureSheaf{\pi^{-1}(P)}(\pi^{-1}(P))}
\arrow{e,t}{dAJ}
\node{H^1(C,\StructureSheaf{C}).}
\end{diagram}
}
\end{equation}
The middle contraction $\rfloor\frac{\partial}{\partial z}$ maps
residues $(a_1,a_2,\dots,a_n) \in
H^0(\pi^{-1}(P),\restricted{\pi^{*}\omega_{\Sigma}(P)}{\pi^{-1}(P)})$ to
$(a_1\frac{\partial}{\partial z_1},
a_2\frac{\partial}{\partial z_2},\dots,
a_n\frac{\partial}{\partial z_n})$.
The lemma follows
by the exactness of the horizontal sequences in the diagram.
\EndProof


\smallskip
\noindent
\underline{Step 3}:
We conclude that $\pi:C\rightarrow \Sigma$ is tangential if
and only if there is a $1$-form $y$ as in the theorem.
If $\pi:C\rightarrow \Sigma$ is tangential then
using lemma \ref{lemma-residues-vs-differential-of-abel-jacobi}
we see that
equation (\ref{eq-linear-condition-tangentiality})
gives rise
to a $1$-form $y$ with residues $(-1,-1,\dots,n-1)$ as required.
Conversely, given a $1$-form $y$ with residues $(-1,-1,\dots,n-1)$  lemma
\ref{lemma-residues-vs-differential-of-abel-jacobi}
and equations (\ref{eq-linear-dependency}) and
(\ref{eq-linear-condition-tangentiality}) imply the tangentiality.
This completes the proof of theorem
\ref{thm-characterization-of-tangential-covers} in the generic case considered.
\EndProof

Theorem \ref{thm-variety-of-kp-solitons} would follow once the existence
of either $n$-sheeted spectral covers
$S$-smooth over $P$, or degree $n$ tangential covers
smooth and unramified over $P$ is established for every choice of
$(\Sigma,P)$.
This was done in (\cite{treibich-verdier-elliptic-solitons} Theorem 3.11)
by studying the linear system of transferred spectral
curves on the surface $\bar{X}_{S}$ of example
\ref{example-the-coadjoint-orbit-of-elliptic-solitons} and applying
Bertini's theorem to show that the generic transferred spectral curve is smooth
in $X_{S}$.
It follows that the blow up of the point with residue $-1$ over $P$ resolves
the generic $n$-sheeted spectral curve of Higgs pairs in $\HiggsModuli_{S}$
to a smooth curve of genus $n$ unramified over $P$.

Any $S$-smooth spectral curve $\bar{C}$
in $T^*\Sigma(P)$ of a Higgs pair $(E,\varphi)$ in
$\HiggsModuli_S$ admits
a unique partial normalization $C$ of arithmetic genus $n$
by the spectral sheaf $\bar{L}$ corresponding to $(E,\varphi)$.
The tautological $1$-form $\bar{y}$ pulls back to a $1$-form on $C$
of the type which characterize tangential covers by theorem
\ref{thm-characterization-of-tangential-covers}.

Conversely, a degree $n$ tangential cover
$\pi:(C,\tilde{P})\rightarrow (\Sigma,P)$
of arithmetic genus $n$
which is smooth and unramified over $P$ is sent to the
spectral curve  $\bar{C}$ in $T^*\Sigma(P)$  of the Higgs pair
\[
(E,\varphi):=
(\pi_*(L),[\otimes y :
\pi_*(L) \rightarrow
\pi_*(L)\otimes\omega_{\Sigma}(P)])
\]
for some, hence every,
choice of a line bundle $L$ on $C$.
Note that $\bar{C}$ is reduced since it is irreducible and
the branch through residue $n-1$ over $P$ is reduced. The
canonical morphism $\nu : C \rightarrow \bar{C}$ is the resolution
by the spectral sheaf $\nu_*L$. Hence $(\bar{C},\nu_*L)$ is $S$-smooth.
(The arithmetic genus of $\bar{C}$ is $\frac{1}{2}(n^2 -n+2)$,
the common arithmetic genus to all
$n$-sheeted spectral curves in $T^*\Sigma(P)$.)

\medskip
Finally we note that, as the tangency point $\tilde{P}$ over $P$ is
marked by having residue $n-1$, it does not have monodromy and
all the KP flows corresponding to it are well define on $\Higgsm_{S/P}$
as the Hamiltonian vector fields of the functions
$\phi^j_{\tilde{P}}$ given in
(\ref{eq-the-jth-hamiltonian-of-the-ith-component}).

\subsubsection{Outline of the proof of the compatibility theorem}
\label{sec-proof-of-compatibility-theorem}

For simplicity we assume that $D = lP$, $l \geq 0$. The general case is
similar replacing  $\ModuliVBLevels_{\Sigma}(n,d,\infty P)$ by
$\ModuliVBLevels_{\Sigma}(n,d,\infty P + D)$.
Let $G_{\infty} :=
GL(n,K)$ be the loop group and $G_{\infty}^{+}$ the level infinity subgroup.
(More canonically,
$K \cong \ComplexNumbers((z))$ should be thought of as the completion of
the function field of $\Sigma$ at $P$, and we may postpone the choice
of a coordinate $z$ until we need to choose generators for a
maximal torus $\HeisN$ of $G_{\infty}$.)

Denote by $M_{l,k}$, $k\geq l$ the pullback of
$T^*\ModuliVBLevels_{\Sigma}(n,d,lP)$ to $\ModuliVBLevels_{\Sigma}(n,d,kP)$
via the rational forgetful morphism.
$M:=T^*\ModuliVBLevels_{\Sigma}(n,d,\infty P)$
is defined as the limit of finite dimensional approximations
(see \ref{sec-finite-dim-approaximations} for the terminology)
\[
T^*\ModuliVBLevels_{\Sigma}(n,d,\infty P) :=
\lim_{l \rightarrow l}\lim_{\infty \leftarrow k} M_{l,k}.
\]
Denote by $M^s$ (resp. $M_{l,k}^s$)
the subset of $T^*\ModuliVBLevels_{\Sigma}(n,d,\infty P)$
(resp. $M_{l,k}$)
consisting of triples $(E,\varphi,\eta)$ with a {\em stable} Higgs pair
$(E,\varphi)$.
We arrive at the setup of section \ref{sec-finite-dim-approaximations}
with the Poisson quotient $Q_{\infty} := M^s/G^{+}$ being the direct limit
\[
Higgs_{\infty} := \lim_{l \rightarrow \infty}
\HiggsComponent_{\Sigma}(n,d,lp).
\]
We emphasize that the stability condition is used here for the morphism
\[M_{l,k}^s \rightarrow \HiggsComponent_{\Sigma}(n,d,lp)\]
between the two
{\em existing} coarse moduli spaces to be well defined, and {\em not} to
define the quotient.

The {\em infinitesimal} loop group action  on
$\ModuliVBLevels_{\infty P} := \ModuliVBLevels_{\Sigma}(n,d,\infty P)$
(the derivative of the action defined in section \ref{sec-krichever-maps} on
the level of \v{C}ech $1$-cocycles)
may be lifted to an infinitesimal action on its cotangent bundle. The
point is that the infinitesimal action of
$a \in \LieAlg{g}_{\infty}$, with poles of order
$\leq l_0$, is well defined on the finite dimensional approximation
$\ModuliVBLevels_{\Sigma}(n,d,l P)$ for $l \geq l_0$. Thus, it lifts
to all cotangent bundles $M_{l,l}$, $l\geq l_0$.
$M_{l,k}$ embeds naturally in $M_{k,k}$
as an invariant subvariety. This defines the action on the limit
$T^*\ModuliVBLevels_{\Sigma}(n,d,\infty P)$.

As a lifted action it is automatically Poisson. Its moment map
\[\mu_{\infty}^* :
\LieAlg{g}_{\infty} \rightarrow
\Gamma[\StructureSheaf{T^*\ModuliVBLevels_{\infty P}},\{,\}]
\]
(the limit of the moment maps for the finite dimensional approximations)
is given by the same formula that we have already encountered for the level
groups
(see \ref{eq-moment-map-of-finite-dim-level-action}):
\begin{equation} \label{eq-moment-map-lifted-loop-group-action}
(\mu_{\infty}^*(a))(E,\varphi,\eta) = Res_{P}trace(a\cdot(\varphi)^{\eta}), \
\ \ a \in \LieAlg{g}_{\infty}.
\end{equation}

Choosing a maximal torus $\heisN \subset \LieAlg{g}_{\infty}$ of
type $\underline{n}$ we arrive at the setup of section
\ref{sec-type-loci}. In particular, we obtain the type locus
$\HiggsModuli^{\underline{n}}_{\Sigma}(n,d,lp)$ in
$\HiggsComponent_{\Sigma}(n,d,lp)$.

\begin{lem} \label{lemma-two-type-loci-coincide}
\begin{enumerate}
\item \label{lemma-item-two-type-loci-coincide}
 The algebro-geometric definition of the ramification type loci
coincides with the group theoretic definition
\ref{def-group-theoretic-definition-type-loci} when
$char(E,\varphi)$ is an integral (irreducible and reduced) spectral curve.
\item \label{lemma-item-two-galois-covers-coincide}
A choice of generators for a maximal torus $\heisN$ as in
(\ref{eq-the-generator-of-the-ith-heisenberg-block})
determines a canonical isomorphism between the group theoretic and
the algebro-geometric $W_{\underline{n}}$-Galois covers.
\end{enumerate}
\end{lem}

\noindent
{\bf Proof:} \ref{lemma-item-two-type-loci-coincide})
The stabilizer $\LieAlg{t} \subset \LieAlg{g}_{\infty}$
of $(\varphi)^{\eta} \in
\LieAlg{g}_{\infty}\otimes_{\StructureSheaf{(P)}} \omega_{\Sigma,(P)}
\cong \LieAlg{g}_{\infty}^*$ with spectral curve
$\pi : C = char(E,\varphi) \rightarrow \Sigma$
is precisely $U^{\eta}$ where $U$ is the
stalk of
$Ker[ad\varphi : \End E \rightarrow \End E\otimes \omega_{\Sigma}(lP)]$
at the formal punctured neighborhood of $P$.
In addition, $U$ is canonically isomorphic to the
stalk of the structure sheaf
at the formal punctured neighborhood
of the fiber of
$C$ over $P$ via the completion of the canonical embedding:
\[
\pi_*\StructureSheaf{C} \hookrightarrow \pi_*\End L
\hookrightarrow \End E.
\]
Hence, the level infinity structure $\eta$ provides a {\em canonical}
isomorphism
\begin{equation} \label{eq-isomorphism-of-tori-induced-by-level-structure}
\lambda : \LieAlg{t} \IsomRightArrow U
\end{equation}
from the stabilizer algebra $\LieAlg{t}$ to the
structure sheaf at the  formal punctured neighborhood of the fiber of
$C$ over $P$.

\smallskip
\noindent
\ref{lemma-item-two-galois-covers-coincide})
As the types coincide, we may choose the level infinity structure
$\eta$ so that the stabilizer $\LieAlg{t}$ coincides with the
fixed torus $\heisN$.
We may further require that the isomorphism
$\lambda$ given by (\ref{eq-isomorphism-of-tori-induced-by-level-structure})
coincides with the one in (\ref{eq-isomorphism-heis-to-structure-sheaf}).
This determines the $\HeisN^+$-orbit of $\eta$ uniquely, i.e.,
a point in the group theoretic Galois cover.
\EndProof

Theorems
\ref{thm-compatibility}
and \ref{thm-compatibility-singular-case}
would now follow from corollary \ref{cor-hamiltonians-on-the-base}
provided that we prove that the homomorphism $\phi$ of the theorems
(given by \ref{eq-the-moment-map-for-the-heisenberg-action}) is
indeed the factorization of the $\heisN$-moment map through the
characteristic polynomial map. (Note that the existence of this factorization
follows from diagram (\ref{diag-factoring-the-moment-map-through-base})).
In other words, we need to prove the identity
\begin{equation} \label{eq-equality-of-two-moment-maps}
\sum_{\{P_i\}} Res_{P_i}(\lambda(a)\cdot y_b) =
Res_{P}trace(a\cdot(\varphi)^{\eta})
\ \ \ a\in\heisN, \ b={\rm char}(E,\varphi)
\end{equation}
as functions on the set of all Higgs pairs $(E,\varphi,\lambda)$ in
$\CoverHiggsModuli_{lp}^{\underline{n}}$
for which the spectral sheaf resolves the singularities
of their spectral curve over $P$
(see definition \ref{def-non-essential-singularities}).
Above, $\eta$ is any level infinity structure in the $\HeisN^+$-orbit
as in the proof of lemma \ref{lemma-two-type-loci-coincide}
or, equivalently, $\lambda(a) = \eta^{-1} \circ a \circ \eta$
where we identify the structure sheaf of the formal punctured neighborhood
with $U$ of that lemma.
If the embedded spectral curve $\bar{C_b}$ is singular,
the $P_i$ are the points over $P$ of its resolution
$\nu: C_b \rightarrow \bar{C_b}$,
and the tautological meromorphic $1$-form
$y_b$ should be replaced by the pullback $\nu^*(y)$
of the tautological $1$-form $y$ on the surface $T^*_{\Sigma}(lP)$.

Conjugating the right hand side of (\ref{eq-equality-of-two-moment-maps})
by $\eta$, we get
\[
\sum_{\{P_i\}} Res_{P_i}(A \cdot \nu^*(y)) =
Res_{P}trace((\pi\circ \nu)_*[A \cdot \nu^*(y)])
\]
for $A$ a (formal) meromorphic function at the fiber over $P$.
Working formally, we can consider only the ``parts'' with first order
pole $r_i dlog{z_i}$ of $A \cdot \nu^*(y)$ at $P_i$.
The equality follows from the identity
$dlog{z_i} = (\pi\circ\nu)^*[\frac{1}{n_i}dlog{z}]$ which imply (projection
formula) that $(\pi\circ\nu)_*(\otimes dlog{z_i})$ acts as
$\frac{1}{n_i}e_{P_i}\otimes dlog{z}$ were $e_{P_i}$ is the
projection onto the eigenspace
of the point $P_i$.
\EndProof
\newpage
\section{The Cubic Condition and Calabi-Yau threefolds}
\label{ch7}

We pose in section \ref{sec-families-of-tori} the general question:
when does a family of polarized abelian varieties or complex tori
support a completely integrable system?
In section \ref{subsec-cubic-condition}
we describe a general necessary infinitesimal symmetry condition
on the periods of the family
(the cubic condition of lemmas \ref{lemma-weak-cubic-cond-poisson} and
\ref{lemma-weak-cubic-cond-symplectic}) and a sufficient local condition
(lemmas \ref{lemma-strong-cubic-cond-local-coordinates}
and \ref{lemma-strong-cubic-cond-coordinate-free}).

In section \ref{sec-cy-threefolds} we use the Yukawa cubic
to construct a symplectic structure (and an ACIHS)
on the relative intermediate Jacobian over the moduli space
of gauged Calabi-Yau threefolds (theorem \ref{thm-cy-acihs}).
The symplectic structure extends to the bundle of Deligne cohomologies
and we show that the image of the relative cycle map
as well as bundles of sub-Hodge-structures are isotropic
(corollary \ref{cor-contact-structure-extends-to-deligne-coho}).

\subsection{Families of Tori} \label{sec-families-of-tori}

Consider a Poisson manifold $(X,\psi)$ together with a Lagrangian fibration
$$
\pi : {\cal X} \longrightarrow B
$$
over a base $B$, whose fibers
$$X_b := \pi^{-1}(b), \quad \quad b \in B $$
are tori.  (We say $\pi$ is Lagrangian if each fiber $X_b$ is a Lagrangian
submanifold of some symplectic leaf in $\cal X$.)  All these objects may be
$C^\infty$, or may be equipped with a complex analytic or algebraic
structure.

On $B$ we have the tangent bundle $\T_B$ as well as the vertical
bundle $\cal V$, whose sections are vector fields along the fibers of $\pi$
which are constant on each torus.  The pullback $\pi^*{\cal V}$ is the
relative tangent bundle $\T_{{\cal X}/B}$; in the analytic or
algebraic situations, we can define $\cal V$ as $\pi_*\T_{{\cal
X}/B}$.  The data $\pi$ and $\psi$ determine an injection
$$
i : {\cal V}^* \hookrightarrow \T_B
$$
or, equivalently, a surjection
$$
i' : \T^*_B \twoheadrightarrow {\cal V}
$$
sending a 1-form $\alpha$ on $B$ to the vertical vector field
$$
i'(\alpha) := \pi^* \alpha \, \rfloor \, \psi.
$$
The image $i({\cal V}^*) \subset \T_B$ is an integrable distribution
on $B$.  Its integral manifolds are the images in $B$ of symplectic leaves
in $\cal X$.

In this section we start with a family of tori $\pi : {\cal X} \rightarrow
B$ and ask whether there is a {\it Lagrangian structure} for $\pi$, i.e. a
Poisson structure on $\cal X$ making the map $\pi$ Lagrangian.  More
precisely, we fix $\pi : {\cal X} \rightarrow B$ and an injection $i :
{\cal V}^* \hookrightarrow \T_B$ with integrable image, and ask for
existence of a Lagrangian structure $\psi$ on $\cal X$ inducing the given $i$.

In the $C^\infty$ category there are no local obstructions to existence of
a Lagrangian structure:  the fibration $\pi$ is locally trivial, so one can
always find action-angle coordinates near each fiber, and use them to
define $\psi$.  In the analytic or algebraic categories, on the other hand,
the fibers $X_b$ (complex tori, or abelian varieties) have invariants, given
essentially by their {\it period matrix} $p(X_b)$, so the fibration
may not be analytically locally trivial.  We will see that there is an
obstruction to existence of a Lagrangian structure for $\pi : {\cal X}
\rightarrow B$, which we formulate as a symmetry condition on the
derivatives of the period map $p$.  These derivatives can be considered as
a linear system of quadrics, and the condition is, roughly, that they be
the polars of some cubic (= section of $\Sym^3{\cal V})$.

Let $X$ be a $g$-dimensional complex torus, and $\gamma_1, \cdots ,
\gamma_{2g}$ a basis of the integral homology $H_1(X, \Z)$.  There is a
unique basis $\alpha_1 , \cdots, \alpha_g$ for the holomorphic
differentials $H^0(X,\Omega^1_X)$ satisfying
$$
\int_{\gamma_{g+i}} \alpha_j = \delta_{ij}, \quad \quad 1 \le i,j \le g,
$$
so we define the period matrix  $P= p(X,\gamma)$ by
$$
p_{ij} := \int_{\gamma_i} \alpha_j, \quad \quad  1 \le i,j \le g.
$$
Riemann's first and second bilinear relations say that $X$ is a principally
polarized abelian variety (PPAV) if and only if $P$ is in {\it Siegel's
half space}:
$$
\h_g := \{ {\rm symmetric} \; g \times g \; {\rm complex \, matrices \,
whose \, imaginary \, part \, is \, positive \, definite} \}.
$$
In terms of a dual basis $\gamma_1^*, \cdots , \gamma_{2g}^*$
of $H^1(X,\Integers)$, the integral class
$\omega := \sum_{i=1}^{g}\gamma_i^*\wedge\gamma_{g+i}^* \in
H^2(X,\Integers)$ is a K\"{a}hler class if and only if $P$ is in Siegel's
half space. In this case we call $\omega$ a {\em principal polarization}.

Given a family $\pi : {\cal X} \rightarrow B$ of PPAVs together with a
continuously varying family of symplectic bases $\gamma_1 , \cdots ,
\gamma_{2g}$ for the fiber homologies, we then get a period map
$$
p : B \longrightarrow \h_g.
$$
If we change the basis $\gamma$ by a symplectic transformation
$$
\left(
\begin{array}{cc}
A & B \\
C & D \\
\end{array}
\right) \in Sp(2g, \Z),
$$
the period matrix $P$ goes to $(AP + B)$ $(CP + D)^{-1}$.  So given a
family $\pi$ without the choice of $\gamma$, we get a multi-valued map of $B$
to $\h_g$, or a map
$$   p : B \longrightarrow {\cal A}_g $$
to the moduli space of PPAV.  The latter is a quasi projective variety,
which can be described analytically as the quotient
$$ {\cal A}_g = \h_g/\Gamma$$
of $\h_g$ by the action of the modular group
$$ \Gamma := Sp(2g, \Z) / (\pm 1).$$

A PPAV $X$ determines a point $[X]$ (or carelessly, $X$) of ${\cal A}_g$.
This point is non-singular if $X$ has no automorphisms other than $\pm 1$,
and then the tangent space $T_{[X]}{\cal A}_g$ can be identified with
$\Sym^2V_X$, where $V_X$ is the tangent space (at $0 \in X$) to $X$.  This
can be seen by identifying $T_{[X]}{\cal A}_g$ with $T_{[X]} \h_g$ and
recalling that $\h_g$ is an open subset of $\Sym^2V_X$.  More algebraically,
this follows from elementary deformation theory:  all first-order
deformations of $X$ are given by
$$
H^1(X,\T_X) \approx H^1(X,V_X \, \otimes_{\C} {\cal O}_X) \approx V_X
\otimes H^1(X, {\cal O}_X) \approx \otimes^2 V_X,
$$
and in there the deformations as abelian variety, i.e., the deformations
preserving the polarization bilinear form on $H_1(X,\Z)$, are given by the
symmetric tensors $\Sym^2V_X$.


\subsection{The Cubic Condition} \label{subsec-cubic-condition}

Our condition for an analytic or algebraic family $\pi : {\cal X}
\rightarrow B$ of PPAVs, given by a period map $p : B \rightarrow {\cal
A}_g$, to have a Lagrangian structure $\psi$ inducing a given $i : {\cal
V}^* \hookrightarrow \T_B$, can now be stated as follows.  The
differential of $p$ is a map of bundles:
$$ dp : \T_B \longrightarrow \Sym^2{\cal V},  $$
so the composite
$$ dp \circ  i \; : \; {\cal V}^* \longrightarrow \Sym^2 {\cal V}$$
can be considered as a section of ${\cal V} \otimes \Sym^2{\cal V}$, and the
condition is that it should come from the subbundle $\Sym^3{\cal V}$.  In
other words, there should exist a cubic $c \in H^0(B,\; Sym^3{\cal V})$
whose polar quadrics give the directional derivatives of the period map:  if
the tangent vector $\partial / \partial b \in T_bB$ equals $i(\beta)$ for
some
$\beta \in {\cal V}^*$, then:
$$ {\partial p \over \partial b} = \beta \, \rfloor \, c.$$
We give two versions of this cubic condition.  In the first, we check the
existence of a two vector $\psi$, not necessarily satisfying the Jacobi
identity, for which $\pi$ is Lagrangian, and which induces a given
injection $i : {\cal V}^* \hookrightarrow \T_B$.  (Note that neither
the definition of the map $i$ induced by the two-vector $\psi$, nor the
notion of $\pi$ being Lagrangian, require $\psi$ to be Poisson.)
\begin{lem}\label{lemma-weak-cubic-cond-poisson} \ (Weak cubic condition,
Poisson form).  \ A family $\pi : {\cal X} \rightarrow B$ of polarized
abelian varieties has a two vector $\psi$ satisfying

\smallskip
\noindent
a) $\pi : {\cal X} \rightarrow B$ is Lagrangian\\
b) $\psi$ induces a given $i : {\cal V}^* \hookrightarrow \T_B$ \\

\smallskip
\noindent
if and only if
$$  dp \circ i \in {\rm Hom}({\cal V}^*, Sym^2 {\cal V}) $$
comes from a cubic
$$c \in H^0(B, \; Sym^3{\cal V}).$$
Moreover, in this case there is a unique such 2-vector $\psi$ which
satisfies also

\smallskip
\noindent
c) The zero section $z : B \rightarrow {\cal X}$ is
Lagrangian, i.e., $(T^*_{{\cal X}/B})|_{z(B)}$ is $\psi$-isotropic
(here we identify the conormal bundle of the zero section with
$(T^*_{{\cal X}/B})|_{z(B)}$.)
\end{lem}
{\bf Proof:}
(Note: we refer below to the vertical bundle
$\T_{{\cal X}/B}$ by its, somewhat indirect, realization as the pullback
$\pi^*{\cal V}$.)
The short exact sequence of sheaves on $\cal X$:
$$
0 \rightarrow \pi^*{\cal V} \rightarrow \T_{\cal X} \rightarrow
 \pi^*\T_B \rightarrow 0
$$
determines a subsheaf $\cal F$ of $\Wedge{2}\T_{\cal X}$ which fits
in the exact sequences:
$$
\begin{array}{ccccccccc}
0 &\rightarrow& {\cal F} &\rightarrow& \Wedge{2}\T_{\cal X}
&\rightarrow& \pi^*\Wedge{2}\T_B &\rightarrow& 0 \\

0 &\rightarrow& \pi^*\Wedge{2}{\cal V} &\rightarrow& {\cal F} &\rightarrow&
\pi^*({\cal V} \otimes \T_B) &\rightarrow& 0. \\
\end{array}
$$
The map $\pi$ is Hamiltonian with respect to the two-vector $\psi \in
H^0(B_{\cal X}, \Wedge{2}\T_{\cal X})$ if and only if $\psi$ goes to
$0$ in $\Wedge{2}\T_B$, i.e., if and only if it comes from $H^0({\cal
F})$.  The question is therefore whether $i \in H^0(B, {\cal V} \otimes
\T_B) \subset H^0({\cal X}, \pi^*({\cal V} \otimes \T_B))$
is in the image of $H^0({\cal X},{\cal F})$.  Locally in $B$, this happens
if and only if $i$ goes to $0$ under the coboundary map
\[
\begin{array}{ccc}
\pi_*\pi^*({\cal V} \otimes \T_B) &\longrightarrow& R^1 \pi_*\pi^*
\Wedge{2}{\cal V} \\
\parallel & & \parallel \\
{\cal V} \otimes \T_B &\longrightarrow& \Wedge{2} {\cal V} \otimes
{\cal V}. \\
\end{array}
\]
This latter map factors through the period map
$$
1 \otimes dp \; : \; {\cal V} \otimes \T_B \longrightarrow {\cal V}
\otimes \Sym^2{\cal V}
$$
and a Koszul map
$$
{\cal V} \otimes \Sym^2{\cal V} \longrightarrow \Wedge{2}{\cal V} \otimes
{\cal V}.
$$
Now exactness of the Koszul sequence
$$
0 \rightarrow {\rm \Sym}^3 {\cal V} \rightarrow {\cal V} \otimes {\rm \Sym}^2
{\cal V} \rightarrow \Wedge{2} {\cal V} \otimes {\cal V}
$$
shows that the desired $\psi$ exists if and only if
$$
dp \circ i \; = \; (1 \otimes dp)(i) \in {\cal V} \otimes \Sym^2 {\cal V}
$$
is in the subspace $\Sym^3{\cal V}$.  (The Hamiltonian map $\pi$ will
automatically be Lagrangian, since $i$ is injective.)

We conclude that, locally on $B$, $i$ lifts to a $2$-vector $\psi$
satisfying conditions a), b), if and only if  $dp \circ i$ is a cubic.  If
$\psi_1,\psi_2$ are two such lifts then
$\psi_1 - \psi_2 \in H^0({\cal X},\Wedge{2} \pi^*{\cal V})$.
Moreover, $\psi_1 - \psi_2$ is
determined by its restriction to the zero section because
$\Wedge{2}\pi^*{\cal V}$ restricts to a trivial bundle on each fiber.
The zero section induces a splitting
$\T_{{\cal X}_{|z(B)}} \simeq \pi^*\T_B
\oplus (\pi^*{\cal V})_{|z(B)}$
and hence a well defined pullback
$z^*(\psi) \in H^0(B, \stackrel{2}{\wedge} {\cal V})$ (locally on $B$).
The normalizations $\psi - \pi^*(z^*(\psi))$ patch to a unique
global section satisfying a), b) c).
\EndProof

\bigskip

The symplectic version of this lemma is:
\begin{lem}: \label{lemma-weak-cubic-cond-symplectic} (Weak cubic
condition, quasi-symplectic form).  A family $\pi : {\cal X} \rightarrow B$ of
principally polarized abelian varieties has a $2$-form $\sigma$ satisfying

\smallskip
\noindent
a) $\pi : {\cal X} \rightarrow B$ has isotropic fibers,\\
b) $\sigma$ induces a given (injective) homomorphism $j : \T_B
\hookrightarrow {\cal V}^*,$

\smallskip
\noindent
if and only if
$$ (1\otimes j^*)\circ dp \in \Hom(\T_B,\T^*_B\otimes {\cal V}) \cong
\T^*_B\otimes \T^*_B\otimes {\cal V}
\ \ \mbox{is} \ \mbox{in} \ \ \Sym^2\T^*_B\otimes {\cal V}.
$$
Moreover, in this case, there exists a unique 2-form $\sigma$ satisfying
a), b), and the additional condition

\smallskip
\noindent
c) the zero section is isotropic $(z^*\sigma = 0)$.
\end{lem}

\begin{rem}
{\rm Riemann's first bilinear condition implies further that
$(1\otimes j^*)\circ dp$ maps to  $\Sym^3\T^*_B$,
 i.e., $(\Sym^{2}j^*)\circ dp
\in {\rm Hom}(\T_B,\Sym^2\T^*_B)$
comes from a cubic $c \in H^0(B, \Sym^3\T^*_B)$.}
\end{rem}


The cubic condition for an embedding $j : \T_B \hookrightarrow {\cal V}^*$
does not
guarantee that the induced 2-form $\sigma$ on $\cal X$ is closed.  In that
sense, the cubic condition is a necessary condition for $j$ to induce a
symplectic structure while the following condition is necessary and
sufficient (but, in general, harder to verify).

\bigskip
\noindent
\underline{Closedness Criterion for a Symplectic Form:}
{\em
Given a family $\pi : {\cal X} \rightarrow B$ of polarized abelian
varieties and a surjective $j' : {\cal V} \rightarrow \T^*_B$,
there exists a {\it
closed} $2$-form $\sigma$ on $\cal X$ satisfying conditions a), b), c) of
Lemma \ref{lemma-weak-cubic-cond-symplectic} if and only if $j'({\cal
H}_1({\cal X}/B, \Z)) \subset \T^*_B$ is a Lagrangian lattice in $T^*B$,
i.e., if locally on $B$ it consists of closed $1$-forms.
Moreover, the $2$-form $\sigma$ is uniquely determined by $j'$.
}

\smallskip
\noindent
{\bf Proof:} \ $j'({\cal H}_1({\cal X}/B, \Z))$ is Lagrangian
 $\Longleftrightarrow$ the canonical symplectic structure
$\tilde{\sigma}$
on $T^*B$ is translation invariant under $j'({\cal H}_1({\cal X}/B, \Z))$
$\Longleftrightarrow$  $(j')^*(\tilde{\sigma})$ descends to the unique
$2$-form
$\sigma$ on ${\cal X} = {\cal V}/{\cal H}_1({\cal X}/B, \Z)$  satisfying
conditions a), b),
c) of Lemma \ref{lemma-weak-cubic-cond-symplectic}.
\EndProof

\smallskip
It is instructive to relate the cubic condition to the above criterion.
This is done in lemma \ref{lemma-strong-cubic-cond-local-coordinates}
in a down to earth manner and is reformulated in lemma
\ref{lemma-strong-cubic-cond-coordinate-free}
as a coordinate free criterion.

\noindent
\begin{lem} \label{lemma-strong-cubic-cond-local-coordinates}
\ (``Strong Cubic Condition'')
\noindent
Let $V$ be a $g$-dimensional vector space, $\{e_1$ ,$\cdots$, $e_g \}$ a basis,
$B \subset V^*$ an open subset,
$p:B \rightarrow \h_g \hookrightarrow \Sym^2V$
a holomorphic map ($\h_g$ is embedded in $\Sym^2V$ via the basis $\{e_j\}$),
$\pi : {\cal X} \rightarrow B$ the corresponding family of principally
polarized abelian varieties.  Then the following are equivalent:

\smallskip
\noindent
(i)  There exists a symplectic structure $\sigma$ on ${\cal X}$
such that $\pi : ({\cal X},\sigma) \rightarrow B$ is a Lagrangian fibration and
$\sigma$ induces the identity isomorphism
$$
{\rm id} \in \Hom(\T_{{\cal X}/B}, \pi^*\T^*_B) \simeq
\Hom(\pi^*{\cal V},\pi^*{\cal V}).
$$
(ii) $p : B \rightarrow \Sym^2V$ is, locally in $B$, the Hessian of a
function on $B$,\\
(iii) $dp \in {\rm Hom}(\T_B, \Sym^2{\cal V}) \simeq
({\cal V} \otimes \Sym^2{\cal V})$
is a section of $\Sym^3{\cal V}.$
\end{lem}

\noindent
{\bf Proof:} \ Let $\{e^*_j\}$ be the dual basis of $V$.

\noindent
\underline{(i) $\Leftrightarrow$ (iii):}
By the closedness criterion above, there exists $\sigma$ as in (i)
if and only if the subsheaf of lattices
${\cal H}^1({\cal X}/B, \Z) \subset T^*B$ is Lagrangian, i.e.,
if and only if its basis
\[
\{e_1,\cdots,e_g, (p \; \rfloor \; e^*_1), \cdots,
(p \; \rfloor \; e^*_g) \}
\]
consists of closed $1$-forms. The $e_i$'s are automatically closed.
If we regard the differential $dp$
as a section of $\T^*_B \otimes \Sym^2{\cal V}$, then
the two-form $d(p \; \rfloor \; e^*_j)$ is equal to
the anti-symmetric part of the contraction
$dp \rfloor \; e^*_j \in \T^*_B \otimes {\cal V} \cong
{\cal V}\otimes {\cal V}$. Hence,
closedness of $(p \; \rfloor \; e^*_j)$, $1 \le j \le g$,
is equivalent to the symmetry of $dp \in {\cal V}\otimes \Sym^2{\cal V}$
also with respect to the first two factors, i.e., to $dp$
being a section of $\Sym^3{\cal V}$.

\noindent
\underline{(ii) $\Rightarrow$ (iii)}. \ \ Clear.

\noindent
\underline{(iii) $\Rightarrow$ (ii)}. \ \ Follows from the Poincare lemma.
\EndProof

\medskip

The additional information contained in the ``Strong Cubic Condition'' and
lacking in the ``Weak Cubic Condition'' is that a Lagrangian sublattice
(with respect to the polarization) ${\cal L} \subset {\cal H}_1({\cal X}/B,
\Z)$ is mapped via $j' : {\cal V} \tilde{\rightarrow} \T^*_B$ to a
sublattice $j'({\cal L}) \subset T^*B$ Lagrangian with respect to the
holomorphic symplectic structure on $T^*B$.  (In the above lemma, ${\cal L}
= {\rm Sp} \{ e_1,\cdots,e_g \}$). The coordinate free reformulation of
lemma \ref{lemma-strong-cubic-cond-local-coordinates} is:

\begin{lem} \label{lemma-strong-cubic-cond-coordinate-free}
(``Strong cubic condition'')
\ Let $j' : {\cal V} \tilde{\rightarrow} \T^*_B$ be an
isomorphism of the vertical bundle
${\cal V} = R^0_{\pi_*}(T_{{\cal X}/B})$ of
the family $\pi : {\cal X} \rightarrow B$ of polarized abelian varieties
with the cotangent bundle of the base.
Assume only that $j'$ maps a sublattice
${\cal L} \subset {\cal H}_1({\cal X}/B, Z)$
Lagrangian with respect to the polarization  to a sublattice
$j'({\cal L}) \subset T^*B$ Lagrangian with respect to the holomorphic
symplectic structure on $T^*B$.  Then there exists a symplectic structure
$\sigma$ on $\cal X$ s.t. $\pi : {\cal X} \rightarrow B$ is a Lagrangian
fibration and inducing $j'$ if and only if $j'$ satisfies the weak cubic
condition, i.e.
$$
dp \circ i \in H^0(B,\Sym^3{\cal V}) \quad \quad \mbox{where}
\ \ i = (j')^{*^{-1}}.
$$
\end{lem}

\begin{rem}
{\rm In most cases however, $j'({\cal L})$ being Lagrangian implies
$j'({\cal H}_1({\cal X}/B, \Z))$ being Lagrangian via the global monodromy
action and without reference to the weak cubic condition.}
\end{rem}

Finally we remark that the above discussion applies verbatim to the case of
polarized complex tori (not necessarily algebraic) since only the first
Riemann bilinear condition was used.


\subsection{An Integrable System for Calabi-Yau Threefolds}
\label{sec-cy-threefolds}

The {\it Hodge group} $H^{p,q}$ of an $n$-dimensional compact K\"{a}hler
manifold $X$ is defined as the space of harmonic forms on $X$ of type
$(p,q)$, i.e. involving $p$ holomorphic and $q$ antiholomorphic
differentials.  Equivalently, $H^{p,q}$ is isomorphic to the $q-th$
cohomology $H^q(X, \Omega^p)$ of the sheaf of holomorphic $p$-forms on $X$.
The {\it Hodge theorem} gives a natural decomposition of the complex
cohomology,
$$
H^k(X,\C) \approx \oplus_{p+q=k} \; H^{p,q}
\approx \oplus_{p+q=k} \; H^q(\Omega^p).
$$
The {\it Hodge number} $h^{p,q}$ is the complex dimension of $H^{p,q}$.
The {\it Hodge filtration} of $H^k(X,\C)$  is defined by
$$
F^iH^k(X,\C) := \oplus_{\stackrel{p+q=k}{p \geq i}} \; H^{p,q}.
$$
The {\it k-th intermediate Jacobian} of $X$
\cite{c-g}
is:
$$J^k(X) := H^{2k-1}(X,\C)/(F^kH^{2k-1}(X,\C) + H^{2k-1}(X,\Z))$$
$$\approx (F^{n-k+1} H^{2n-2k+1}(X,\C))^*/H_{2n-2k+1}(X,\Z).$$
Elementary properties of the Hodge filtration imply that this is a complex
torus, but generally not an abelian variety unless $k=1$ or $k=n$:  it
satisfies Riemann's first bilinear condition (which expresses the skew
symmetry of the cup product on $H^{2k-1}$), but not the second, since the
sign of the product (on primitive pieces) will vary with the parity of $p$.
The extreme cases correspond to the connected component of the Picard
$(k=1)$ and Albanese $(k=n)$ varieties.

The Hodge decomposition does not depend holomorphically on parameters, since
both holomorphic and antiholomorphic differentials are involved.  The
advantage of the Hodge filtration is that it does vary holomorphically and
even algebraically when $X$ is algebraic. The $F^p$ can be defined
algebraically, as the hypercohomology of the complex
\begin{equation} \label{eq-the-quotient-of-the-algebraic-derham-complex}
0 \rightarrow
\Omega^p \rightarrow \Omega^{p+1} \rightarrow \cdots \rightarrow \Omega^n
\rightarrow 0.
\end{equation}
In particular, the intermediate Jacobian $J^k(X)$
varies holomorphically with  $X$.  This means that a smooth analytic family
${\cal X} \rightarrow B$ of compact K\"{a}hler manifolds gives rise to
analytic vector bundles $F^i{\cal H}^k({\cal X}/B)$ and to smooth analytic
families ${\cal J}^k({\cal X}/B) \longrightarrow B$ of intermediate
Jacobians of the fibers.

The bundle ${\cal H}^k({\cal X}/B)$ is the complexification of a bundle
${\cal H}^k({\cal X}/B, \Z)$ of discrete groups.  In particular, it has a
natural local trivialization.  In other words, it admits a natural flat
connection, called the {\it Gauss-Manin} connection.  The holomorphic
subbundles $F^i{\cal H}^k({\cal X}/B)$ are in general not invariant with
respect to this connection, since the Hodge decomposition and filtration do
change from point to point.  {\it Griffiths' transversality} says that when
a holomorphic section of $F^i{\cal H}^k$ is differentiated, it can move at
most one step:
$$ \nabla(F^i{\cal H}^k) \subset F^{i-1}{\cal H}^k \otimes \Omega^1_B.
$$

An $n$-dimensional compact K\"{a}hler manifold $X$ is called {\it
Calabi-Yau} if it has trivial canonical bundle,
$$
\omega_X = \Omega_X^n \; \approx \; {\cal O}_X,
$$
and satisfies
$$
h^{p,0} = 0 \  {\rm for} \ 0 < p < n.
$$
A {\it gauged Calabi-Yau} is a pair $(X,s)$ consisting of a Calabi-Yau
manifold $X$ together with a non-zero volume form
$$
s : {\cal O}_X \stackrel{\approx}{\longrightarrow} \omega_X.
$$
A theorem of Bogomolov, Tian and Todorov
\cite{bogomolov,Ti,To} says that $X$ has
a smooth (local analytic)  universal deformation space $M_X$.
We say that a family $\chi: {\cal X} \rightarrow {\cal M}$ of Calabi-Yaus
$X_t, t\in {\cal M}$, is {\em complete} if the local classifying map
${\cal M} \supset U_t \rightarrow M_{X_t}$ is an isomorphism
for some neighborhood of every point $t \in {\cal M}$.
It follows that ${\cal M}$ is smooth and that the tangent space at $t$ to
${\cal M}$ is naturally isomorphic to
$H^1(X, \T_X)$. Typically, such families might consist of all Calabi-Yaus
in some open subset of moduli, together with some ``level'' structure.

The choice of gauge $s$ gives an isomorphism
$$
\rfloor \;s : \T_X \longrightarrow \Omega^{n-1}_X,
$$
hence an isomorphism
$$
T_X{\cal M} \; \approx \; H^{n-1,1}(X).
$$
Starting with a complete family $\chi : {\cal X} \rightarrow {\cal M}$, we
can construct
\begin{itemize}
\item The bundle ${\cal J}^k \rightarrow {\cal M}$ of intermediate
Jacobians of the Calabi-Yau fibers.
\item The space $\tilde{\cal M}$ of gauged Calabi-Yaus, a
$\C^*$-bundle over ${\cal M}$ obtained by removing the $0$-section from the
line bundle $\chi_*(\omega_{{\cal X}/{\cal M}}).$
\item The fiber product
$$\tilde{\cal J}^k := {\cal J}^k \; \times_{\cal M} \; \tilde{\cal M} ,
$$
which is an analytic family of complex tori $\pi : \tilde{\cal J}^k
\rightarrow \tilde{\cal M}$.
\end{itemize}

\begin{theorem} \label{thm-cy-acihs}
Let ${\cal X} \rightarrow {\cal M}$ be a complete family of
Calabi-Yau manifolds of odd
dimension $n=2k-1 \geq 3$.  Then there exists a canonical closed
holomorphic $2$-form $\sigma$ on
the relative $k$-th intermediate Jacobian $\pi : \tilde{\cal J} \rightarrow
\tilde{\cal M}$ with respect to which $\pi$ has maximal isotropic
fibers. When $n=3$, the $2$-form $\sigma$ is a symplectic structure and
$\pi : \tilde{\cal J} \rightarrow \tilde{\cal M}$
is an analytically completely integrable Hamiltonian system.
\end{theorem}

\noindent
{\bf Proof.}

\underline{Step I.} \ There is a canonical isomorphism
$$T_{(X,s)} \tilde{\cal M} \; \approx \; F^{n-1} H^n(X,\C).$$
Indeed, the natural map $p : \tilde{\cal M} \rightarrow {\cal M}$ gives a
short exact sequence
$$
0 \rightarrow T_{(X,s)}(\tilde{\cal M}/{\cal M}) \rightarrow T_{(X,s)}
\tilde{\cal M} \rightarrow T_X{\cal M} \rightarrow 0,
$$
in which the subspace can be naturally identified with $H^0(\omega_X) =
H^0(\Omega^n_X)$, and the quotient with $H^1(\T_X)$, which goes
isomorphically to $H^1(\Omega^{n-1}_X)$ by $\rfloor \; s$.  What we are
claiming is that this sequence can be naturally identified with the one
defining $F^{n-1}H^n$:
$$
0 \rightarrow H^0 (\omega_X) \rightarrow F^{n-1} \rightarrow H^1
(\Omega^{n-1}_X) \rightarrow 0,$$
i.e., that the extension data match, globally over $\tilde{\cal M}$.  To see
this we need a natural map $T_{(X,x)} \tilde{\cal M} \rightarrow F^{n-1}
H^n$ inducing the identity on the sub and quotient spaces.

Over $\tilde{\cal M}$ there is a tautological section $s$ of $F^n{\cal
H}^n(\tilde{\cal X}/\tilde{\cal M},\C)$.  The Gauss-Manin connection
defines an embedding
$$
\nabla_{(\cdot)}s : T_{(X,s)} \tilde{\cal M} \longrightarrow H^n(X,\C).
$$
Griffiths' transversality implies that the image is in $F^{n-1}H^n(X,\C)$.
Clearly $\nabla_{(\cdot)}s$ has the required properties.

We will need also a description of the isomorphism in terms of Dolbeault
cohomology.  We think of a 1-parameter family $(X_t,s_t) \in
\tilde{\cal M}$, depending on the parameter $t$, as living on a fixed
topological model $X$ on which there are families $\bar{\partial}_t$ of
complex structures (given by their $\bar{\partial}$-operator) and $s_t$ of
$C^\infty$ $n$-forms, such that $s_t$ is of type $(n,0)$ with respect to
$\bar{\partial}_t$, all $t$.  Since the $s_t$ are now on a fixed underlying
$X$, we can differentiate with respect to $t$:
$$
s_t = s_0 + ta \quad \quad ({\rm mod}\;t^2).
$$
Griffiths transversality now says that $a$ is in $F^{n-1}H^n(X_0)$.  It
clearly depends only on the tangent vector to $\tilde{\cal M}$ along
$(X_t,s_t)$ at $t = 0$, so we get a map $T_{(X,s)} \tilde{\cal M}
\longrightarrow F^{n-1}H^n$ with the desired properties.

\medskip

\underline{Step II}.  Let ${\cal V}$ be the vertical bundle on
$\tilde{\cal M}$ coming from $\pi:\tilde{\cal J} \rightarrow \tilde{\cal M}$.
It is
isomorphic to
$$
F^k{\cal H}^n(\tilde{\cal X}/\tilde{\cal M})^*
$$
(recall $n = 2k-1$).  Combining with Step I, we get a natural injection
$$j : \T_{\tilde{\cal M}} \hookrightarrow {\cal V}^*, $$
which above a given $(X,s)$ is the inclusion of $F^{n-1}H^n(X)$ into
$F^kH^n(X)$.  Its transpose
$$
j' : {\cal V} \twoheadrightarrow  T^*\tilde{\cal M}
$$
determines a closed $2$-form $\sigma := (j')^* \tilde{\sigma}$ on $\cal V$,
where $\tilde{\sigma}$ is the standard symplectic form on $T^* \tilde{\cal
M}$ (see example \ref{examples-symplectic-varieties}).
By construction, the fibers of $\cal V$ over $\tilde{\cal M}$
are maximal isotropic with respect to this form.
\medskip

\underline{Step III}.  We need to verify that $\tilde{\sigma}$ descends to
$${\cal J}^k(\tilde{\cal X} /\tilde{\cal M} ) = {\cal V}/{\cal H}_n
(\tilde{\cal X}/\tilde{\cal M}, \Z).$$
Equivalently, a locally constant integral cycle
$$
\gamma \in \Gamma(B,{\cal H}_n(\tilde{\cal X} / \tilde{\cal M},\Z)),
$$
defined over some open subset $B$ of $\tilde{\cal M}$, gives a section of
$\cal V$ on $B$; hence through $j'$, a $1$-form $\xi$ on $B$, and we need
this $1$-form to be closed.  Explicitly, if $a$ is a section of
$\T_{\tilde{\cal M}}$ over $B$, we have
$$
a \, \rfloor s \in \Gamma (B, F^{n-1}{\cal H}^n(\tilde{\cal X} /\tilde{\cal
M}))
\subset \Gamma(B, {\cal H}^n(\tilde{\cal X} /\tilde{\cal M}))
$$
and $\xi$ is defined by:
$$
\langle \xi,a \rangle := \int_\gamma \; (a \rfloor s).
$$
Consider the function
$$
\begin{array}{c}
g : B \longrightarrow \C \\
 \\
g(X,s) := \int_\gamma \; s. \\
\end{array}
$$
If we set
$$
a = \left. {\partial \over {\partial t}} \right|_{t=0}  (X_t,s_t)
$$
as in Step I, we get:
$$
\langle dg,a \rangle = \left.
{\partial \over {\partial t}} \right|_{t=0} \; g(X_t,s_t) =
\left. {\partial \over {\partial t}} \right|_{t=0} \; \int_\gamma \; s_t =
\int_\gamma \; (a\, \rfloor s) = \langle \xi, a \rangle,
$$
so $\xi = dg$ is closed.
\EndProof

\bigskip
\begin{rem} \label{rem-to-thm-cy-acihs}
{\rm
\begin{enumerate}
\item[(1)]  The most interesting case is clearly $n=3$, when
$\tilde{{\cal J}}$ has
an honest symplectic structure.  The cubic field on $\tilde{\cal M}$
corresponding to this structure by
lemma \ref{lemma-weak-cubic-cond-poisson}
made its first appearance in \cite{BG}
and is essentially the {\it Yukawa coupling}, popular among physicists and
mirror-symmetry enthusiasts.  At $(X,x) \in \tilde{\cal M}$ there is a
natural cubic form on $H^1(\T_x)$:
$$
c : \otimes^3H^{1}(\T_X) \rightarrow H^3(\Wedge{3}\T_X) = H^3(\omega^{-1}_X)
\stackrel{\cdot s^2}{\rightarrow} H^3(\omega_X) \RightArrowOf{\int} \C,
$$
which pulls back to the required cubic on $\T_{(X,s)}\tilde{{\cal M}}$.  Hodge
theoretically, this cubic can be interpreted as the third iterate of the
infinitesimal variation of the periods, or the Hodge structure, of $X$ c.f.
\cite{IVHS-I} and \cite{BG}.  By Griffiths transversality, each tangent
direction on ${\cal M}$, $\theta \in H^1(\T_X)$, determines a linear map
$$
\theta_i : H^{i,3-i} \longrightarrow H^{i-1,4-i} \quad \quad i = 3,2,1,
$$
and clearly the composition
$$
\theta_1 \circ \theta_2 \circ \theta_3 : H^{3,0} \longrightarrow H^{0,3}
$$
becomes  $c(\theta)$ when we use $s$ to identify $H^{3,0}$ and its dual
$H^{0,3}$ with $\C$.

\item[(2)] For $n=2k-1 \geq 5$, we get a closed $2$-form on $\tilde{\cal J}$
which is in general not of maximal rank.  The corresponding cubic is
identically $0$.  Hodge theoretically, the ``cubic'' multiplies the gauge
$s \in H^0(\omega_X)$ by two elements of $H^1(\T_X)$ (landing in
$H^{n-2,2}$) and then with an element of $F^kH^n$.  When $k > 2$ there are
too many $dz's$, so the product vanishes.

\item[(3)]  The symplectic form $\sigma$ which we constructed on
$\tilde{\cal J}$ is actually exact.  Recall that the natural symplectic
form $\tilde{\sigma}$ on $T^*\tilde{\cal M}$ is exact:  $\tilde{\sigma} =
d\tilde{\alpha}$, where $\tilde{\alpha}$ is the action $1$-form.  We
obtained $\sigma$ by pulling $\tilde{\sigma}$ back to
$(j')^*\tilde{\sigma}$ on ${\cal V}$, and observing that the latter is
invariant under translation by locally constant integral cycles $\gamma$,
hence descends to $\tilde{{\cal J}}$.
Now a first guess for the
anti-differential of $\sigma$ would be the $1$-form $(j')^*
\tilde{\alpha}$; but this is {\it not} invariant under translation:  if the
cycle $\gamma$ corresponds, as in Step III of the proof, to a $1$-form
$\xi$ on $\tilde{\cal M}$, then the translation by $\gamma$ changes $(j')^*
\tilde{\alpha}$ by $\pi^*\xi$, where $\pi : {\cal V} \rightarrow
\tilde{\cal M}$ is the projection.  To fix this discrepancy, we consider
the tautological function $f \in \Gamma({\cal O}_{\cal V})$ whose value at
a point $(X,s,v) \in {\cal V}$ (where $(X,s) \in \tilde{\cal M}$ and $v \in
F^kH^n(X)^*)$ is given by
\begin{equation} \label{eq-moment-map}
f(X,s,v) = v(s).
\end{equation}
This $f$ is linear on the fibers of $\pi$, so $df$ is constant on these
fibers, and therefore translation by $\gamma$ changes $df$ by $\pi^*$ of a
$1$-form on the base $\tilde{\cal M}$.  This $1$-form is clearly $\xi$, so
we conclude that
\begin{equation} \label{eq-contact-structure}
(j')^* \tilde{\alpha} - df
\end{equation}
is a global 1-form on ${\cal V}$ which is invariant under translation by
each $\gamma$, hence descends to a $1$-form $\alpha$ on $\tilde{{\cal J}}$.
It satisfies $d\alpha = \sigma$, as claimed.

\item[(4)]
Another way to see the exactness of $\sigma$ on $\tilde{{\cal J}}$
is to note that it comes from a {\it quasi-contact structure} $\kappa$ on
${\cal J}$.
By a quasi-contact structure we mean a  line subbundle
$\kappa$ of $T^*{\cal J}$. It determines a tautological $1$-form
on the $\C^*$-bundle $\tilde{\cal J}$ obtain from $\kappa$ by omitting its
zero section. Hence, it determines also an exact $2$-form $\sigma$
on $\tilde{\cal J}$.
We refer to the pair $(\tilde{\cal J},\sigma)$ as the
{\em quasi-symplectification} of $({\cal J},\kappa)$.
Conversely, according to \cite{AG}, page
78, a $2$-form $\sigma$ on a manifold $\tilde{\cal J}$ with a $\C^*$-action
$\rho$ is the quasi-symplectification of a line subbundle of the cotangent
bundle of the quotient ${\cal J}$ if and only if $\sigma$ is homogeneous of
degree $1$ with respect to $\rho$ (and the contraction of $\sigma$ with the
vector field generating $\rho$ is nowhere vanishing).

In our case, there
are two independent $\C^*$-actions on the total space of $T^*\tilde{{\cal M}}
\simeq [F^{n-1}{\cal H}^n(\tilde{\cal X} / \tilde{\cal M} , \C)]^*$: the
$\C^*$-action on $\tilde{\cal M}$ lifts to an action $\bar{\rho}'$ on
$T^*\tilde{\cal M}$, and there is also the action $\bar{\rho}''$ which
commutes with the projection to $\tilde{\cal M}$ and is linear on the
fibers.  The symplectic form $\tilde{\sigma}$ is homogeneous of weight $0$
with respect to $\bar{\rho}'$ and of weight $1$ with respect to
$\bar{\rho}''$, hence of weight $1$ with respect to $\bar{\rho} :=
\bar{\rho}' \cdot \bar{\rho}''$.  Hence, $\tilde{\sigma}$ is the
symplectification of a contact structure on $T^*\tilde{\cal M}/\bar{\rho}
\simeq [F^{n-1}{\cal H}^n({\cal X} / {\cal M}, \C)]^*$
(suppressing the gauge).

Denote a point in $[F^{n-1}{\cal H}^n(\tilde{\cal X} / \tilde{\cal M},
\C)]^*$ by $(X,s,\xi)$.  The actions, for  $t \in \C^*$, are given by:
\[
\begin{array}{lcl}
\bar{\rho}' &:& (X,s,\xi) \longmapsto (X, ts, t^{-1} \xi) \\
\bar{\rho}'' &:& (X,s,\xi) \longmapsto (X, s, t\xi) \\
\bar{\rho}  &:& (X,s,\xi) \longmapsto (X, ts, \xi). \\
\end{array}
\]
The function $f$ on ${\cal V}$, given by (\ref{eq-moment-map}),
is the pullback
$(j')^*(\bar{f})$ of the function $\bar{f}$ on $T^*\tilde{{\cal M}}$
given by
$$
\bar{f}(X,s,\xi) = \xi(s).
$$
The symplectic structure $\tilde{\sigma}$ on $T^*\tilde{\cal M}$ takes the
vector fields generating the actions $\bar{\rho}', \bar{\rho}''$, and
$\bar{\rho}$ to the $1$-forms $-d\bar{f}, \tilde{\alpha}$ and
$\tilde{\alpha} - d \bar{f}$, respectively.
The $1$-form $\tilde{\alpha} - d \bar{f}$,
which is homogeneous of degree $1$ with
respect to $\bar{\rho}$, is the $1$-form canonically associated to the
contact structure on $T^*\tilde{\cal M}/\bar{\rho}
\simeq [F^{n-1}{\cal H}^n({\cal X} / {\cal M}, \C)]^*$
(namely, the contraction of $\tilde{\sigma}$ with the vector field of
$\bar{\rho}$.)

Similarly, we have three action $\rho', \rho''$ and $\rho =
\rho' \cdot \rho'' $ on the total space of ${\cal V} \simeq
[F^k {\cal H}^n( \tilde{\cal X}/ \tilde{\cal M}, \C)]^*$.  The surjective
homomorphism $j' : {\cal V} \rightarrow T^* \tilde{\cal M}$ is $(\rho',
\bar{\rho}' ), (\rho'', \bar{\rho}'')$, and $(\rho ,\bar{\rho})$-equivariant.
The $1$-form $\tilde{\alpha} - d \bar{f}$ pulls back to the $1$-form
$(j')^*(\tilde{\alpha}) - df$ given by
(\ref{eq-contact-structure}).  Clearly, the
action $\rho$ commutes with translations by ${\cal H}_n(\tilde{\cal X} /
\tilde{\cal M}, \Z)$. Since the $2$-form $(j')^*\tilde{\sigma}$ is also
${\cal H}_n(\tilde{\cal X} / \tilde{\cal M}, \Z)$-equivariant,
$(j')^*(\alpha) - df$ descends to a $1$-form $\alpha$ on $\tilde{\cal J}$.
Clearly, $d\alpha = \sigma$ and $\alpha$ is homogeneous of degree $1$ with
respect to the $\C^*$-action on $\tilde{\cal J}$.  Hence $\alpha$ comes
from a quasi-contact structure $\kappa$ on ${\cal J}$.
\end{enumerate}
}
\end{rem}

The Abel-Jacobi map of a curve to its Jacobian has an analogue for
intermediate Jacobians.  Let $Z$ be a codimensional-$k$ cycle in $X$, i.e.
a formal linear combination $Z = \sum m_i Z_i$, with integer
coefficients, of codimension  $k$  subvarieties $Z_i \subset X$.  If $Z$ is
homologous to $0$, we can associate to it a point $\mu(Z) \in {\cal
J}^k(X)$, as follows.  Choose a real $(2n-2k+1)$-chain $\Gamma$ in $X$
whose boundary is $Z$, and let $\mu(Z)$ be the image in
$$
{\cal J}^k(X) \approx (F^{n-k+1} H^{2n-2k+1}(X,\C))^* / H_{2n-2k+1}(X,\Z)
$$
of the linear functional
$$
\int_\Gamma \in (F^{n-k+1} H^{2n-2k+1}(X,\C))^*
$$
sending a cohomology class represented by a harmonic form $\alpha$ to
$\int_\Gamma  \alpha$.  Changing the choice of $\Gamma$ changes
$\int_\Gamma$ by an integral class, so $\mu(Z)$ depends only on $Z$.  This
construction depends continuously on its parameters:  given a family
$\pi:{\cal X} \rightarrow B$ and a family ${\cal Z} \rightarrow B$ of
codimension-$k$ cycles in the fibers which are homologous to $0$ in the
fibers, we get the {\it normal function}, or Abel-Jacobi map
$$
\mu : B \longrightarrow {\cal J}^k({\cal X}/B)
$$
to the family of intermediate Jacobians of the fibers.

Abstractly, a normal function $\nu:B \rightarrow {\cal J}^k({\cal X}/B)$ is a
section satisfying the infinitesimal condition:

\smallskip
\noindent
{\em Any lift
$$
\tilde{\nu} : B \longrightarrow {\cal H}^n({\cal X}/B, \C)
$$
of
$$
\nu:B \rightarrow J^k({\cal X}/B) \simeq {\cal H}^n({\cal X}/B, \C)
/[F^k{\cal H}^n({\cal X}/B, \C) + {\cal H}^n({\cal X}/B, \Z)]
$$
satisfies
\begin{equation} \label{eq-infinitesimal-cond-normal-fn}
\nabla \tilde{\nu} \in F^{k-1}{\cal H}^n({\cal X}/B, \C) \otimes \Omega^1_B
\end{equation}
or equivalently
$$
(\nabla \tilde{\nu}, s) = 0 \; {\rm for \; any \; section} \; s \; {\rm of} \;
F^{k+1} {\cal H}^n({\cal X}/B, \C)
$$
where $\nabla \tilde{\nu}$ is the Gauss-Manin derivative of $\tilde{\nu}$.
}

\smallskip
\noindent
This condition is independent of the choice of the lift $\tilde{\nu}$ by
Griffiths' transversality.  It is satisfied by the Abel-Jacobi image of a
relative codimension $k$-cycle (see \cite{griffith-normal-functions}).
More generally, we can consider maps
$$
{\divide\dgARROWLENGTH by 2
\begin{diagram}
\node{B}
\arrow[2]{e,t}{\mu}
\arrow{se,b}{q}
\node[2]{{\cal J}^k}
\arrow{sw}
\\
\node[2]{{\cal M}}
\end{diagram}
}
$$
The pullback ${\cal J}^k({\cal X}/B) \rightarrow B$ of the relative
intermediate Jacobian to $B$ has a canonical section $\nu:B \rightarrow
{\cal J}^k({\cal X}/B)$.  We will refer to the subvariety $\mu(B)$ as a {\it
multivalued normal function} if
$\nu : B \rightarrow {\cal J}^k({\cal X}/B)$
is a normal function.

\begin{theorem} \label{thm-integrality-of-normal-fn}
Let
$\tilde{{\cal X}} \rightarrow \tilde{{\cal M}}$ be a complete family of
gauged Calabi-Yau manifolds of dimension $n = 2k - 1 \geq 3$,
$\tilde{{\cal J}} \rightarrow \tilde{{\cal M}}$ the relative intermediate
Jacobian, $B \RightArrowOf{q} \tilde{{\cal M}}$ a base of a family
${\cal Z} \rightarrow B$
of codimension-$k$ cycles homologous to $0$ in the fibers of
$q^*\tilde{{\cal X}} \rightarrow B$. Then
%
i) the
Abel-Jacobi image in $\tilde{\cal J}$ of $B$ is isotropic with respect to
the quasi-symplectic form $\sigma$ of theorem \ref{thm-cy-acihs}.
ii) Moreover, the Abel-Jacobi image is also integral with respect to the
$1$-form $\alpha$ given by (\ref{eq-contact-structure}).
\end{theorem}
\medskip
\noindent
{\bf Proof.} i) We follow Step III of the proof of theorem
\ref{thm-cy-acihs}.  We thus
think, locally in $B$, of $X$ as being a fixed $C^\infty$ manifold with
variable complex structure $\bar{\partial}_b$, $n$-form $s_b$, and cycle
$Z_b$, subject to the obvious compatibility. We choose a family $\Gamma_b$,
$b \in B$ of $n$-chains whose boundary is $Z_b$, and consider the $1$-form
$\xi$ on $B$ given at $b \in B$ by $\int_{\Gamma_b}$; we need to show that
$\xi$ is closed.  (The new feature here is that instead of the cycles
$\gamma_b \in H_n(X_b,\Z)$ we have chains, or relative cycles $\Gamma_b \in
H_n(X_b, |Z_b|, \Z)$, where $|Z_b|$ is the support of $Z_b$, which varies
with $b$.)

As before, we consider the function
\begin{eqnarray*}
g:B &\longrightarrow& \C \\
g(X,s,Z,\Gamma) &:=& \int_\Gamma \, s, \\
\end{eqnarray*}
and we claim $\xi = dg$.  This time, in the integral $\int_{\Gamma_b}
\; s_b$, both the integrand and the chain depend on $b$.  So if we take a
normal vector $v$ to the supports $|Z_b|$ along $\Gamma_b$, we obtain two
terms:
$$
{\partial \over {\partial b}} \; \int_{\Gamma_b} \; s_b = \int_{\Gamma_b}
{{\partial s} \over {\partial b}} + \int_{\partial\Gamma_b} \; (v \, \rfloor
s_b).
$$
In the second term, however, $s_b$ is of type $(n,0)$ with respect to the
complex structure $\bar{\partial}_b$, so the contraction $v \, \rfloor s_b$
is of type $(n-1,0)$ regardless of the type of $v$.  Since $\partial
\Gamma_b = Z_b$ is of the type $(k-1, k-1)$, the second term vanishes
identically, so we have $dg = \xi$ as desired.

\smallskip
\noindent
ii)
Integration $\int_\Gamma(\cdot)$ defines a section  of ${\cal V} \simeq
[F^{k} {\cal H}^n]^*$.  The function $g$ on $B$ is the pullback via
$\int_\Gamma(\cdot)$ of the function $f$ on ${\cal V}$ given by the formula
(\ref{eq-moment-map}).  Similarly, integration $\int_\Gamma(\cdot)$ defines
the section $\xi$ of $T^*\tilde{{\cal M}} \simeq [F^{n-1}{\cal H}^n]^*$.  The
pullback of the tautological $1$-form $\tilde{{\alpha}}$
by $\xi$ is $\xi$ itself.  The
equation $\xi - dg = 0$ translates to the statement that the $1$-form
$(j')^*\tilde{{\alpha}}-df$
vanishes on the section $\int_\Gamma(\cdot)$ of ${\cal V}$
(see formula (\ref{eq-contact-structure})).
In particular, its descent $\alpha$ vanishes on the Abel-Jacobi image of
${\cal Z} \rightarrow B$.
\EndProof
\medskip

Again, the most interesting case is $n=3$.  When $B$ dominates the moduli
space $\tilde{\cal M}$, i.e. for a multivalued choice of cycles on the
general gauged Calabi-Yau of a given type, the normal function produces a
Lagrangian subvariety of the symplectic $\tilde{{\cal J}}$, generically
transversal to the fibers of the completely integrable system.

\begin{rem} \label{rem-integrality-of-normal-fn}

\begin{enumerate}
{\rm
\item[(1)] The result of Theorem \ref{thm-integrality-of-normal-fn} holds
for every multi-valued normal function $\mu : B \rightarrow
\tilde{\cal J}^k(\tilde{\cal X}/\tilde{\cal M})$ , not only for those coming
from cycles.  Given a vector field ${\partial \over {\partial b}}$ on $B$, a
lift $\tilde{\nu} : B \rightarrow {\cal H}^n$, and any section $s$ of $F^{k+1}
{\cal H}^n(\tilde{\cal X}/\tilde{\cal M},\C)$ , the infinitesimal condition
for normal functions (\ref{eq-infinitesimal-cond-normal-fn})
becomes
\begin{equation} \label{eq-derivative-of-normal-fn}
0 = \left( \nabla_{\frac{\partial}{\partial b}}\tilde{\nu}, s \right) =
{\partial \over {\partial b}}(\tilde{\nu},s) -
\left( \tilde{\nu}, \nabla_{\frac{\partial}{\partial b}}s \right).
\end{equation}
When $s$ is the tautological gauge, ${\partial \over {\partial b}}
(\tilde{\nu}, s)$ is the pullback of $df$ by the projection of $\tilde{\nu}$
to ${\cal V} \cong {\cal H}^n / F^k{\cal H}^n$ (where $f$ is defined by the
equation (\ref{eq-moment-map})).  Similarly,
$\left(\tilde{\nu},
\nabla_{\frac{\partial}{\partial b}} s \right)$  is the contraction $\xi \;
\rfloor \; {\partial \over {\partial b}}$ of the pullback $\xi$ of the
tautological $1$-form $\tilde{\alpha}$ on $T^*\tilde{\cal M}$ by the
composition
$$
\tilde{\mu} : B \rightarrow {\cal H}^n(\tilde{\cal X}/\tilde{\cal M})
\rightarrow {\cal H}^n / F^2 \simeq [F^{n-1}{\cal H}^n ]^* \simeq
T^*\tilde{\cal M}.
$$
Thus, the infinitesimal condition for a normal function
(\ref{eq-derivative-of-normal-fn}) implies that the image $\mu(B) \subset
\tilde{{\cal J}}^k(\tilde{\cal X}/\tilde{\cal M})$ is integral with
respect to the
1-form $\alpha$ (defined in (\ref{eq-contact-structure})).

In the case of $CY$ $3$-folds $(n=3, k = 2)$ we see that the
Legendre subvarieties of ${\cal J}^2 \rightarrow {\cal M}$
(i.e., the $\kappa$-integral subvarieties of maximal dimension $h^{2,1}$
where $\kappa$ is the contact structure of Remark
\ref{rem-to-thm-cy-acihs}(4))
are precisely the
multivalued normal functions.

\item[(2)] Both the infinitesimal condition for a normal function
(\ref{eq-infinitesimal-cond-normal-fn}) and the (quasi) contact structure
$\kappa$ on the relative Jacobian ${\cal J}^k \rightarrow {\cal M}$
(see  Remark \ref{rem-to-thm-cy-acihs}(4)) are special cases of a more general
filtration of Pfaffian exterior differential systems on the relative
intermediate Jacobian ${\cal J}^k \rightarrow {\cal M}$ of any family
${\cal X} \rightarrow {\cal M}$
of $n = 2k-1$ dimensional projective algebraic
varieties.

The tangent bundle $T{\cal J}^k$ has a canonical decreasing filtration
(defined by (\ref{eq-filtration}) below)
$$
T{\cal J}^k = F^0T{\cal J}^k \supset F^1T{\cal J}^k \supset \cdots
\supset F^{k-1}T{\cal J}^k \supset
0.
$$
The quotient
$T{\cal J}^k/F^iT{\cal J}^k$ is canonically isomorphic to the
pullback of the Hodge bundle
${\cal H}^n/F^i{\cal H}^n$.  The $F^{k-1}T{\cal J}^k$ integral subvarieties
are precisely the multi-valued normal functions.

When ${\cal J}^k$ is the relative intermediate Jacobian of a family of $CY$
$n$-folds, the subbundle $F^1T{\cal J}^k$ is a hyperplane distribution on
${\cal J}^k$
which defines the (quasi) contact structure $\kappa$ of Remark
\ref{rem-to-thm-cy-acihs}(4).

When $n=3$, $k=2$, the filtration is a two step filtration
$$
T{\cal J}^2 \supset F^1T{\cal J}^2 \supset 0
$$
and the $F^1T{\cal J}^2$-integral subvarieties are precisely the normal
functions.

The filtration $F^iT{\cal J}^k$, $0 \leq i \leq k-1$ is defined at a point
$(b,y) \in {\cal J}^k$ over $b \in {\cal M}$ as follows:
Choose a section $\tilde{\nu} :
{\cal M} \rightarrow {\cal J}^k$ through $(b,y)$ with the property that
any lift
$\tilde{\nu} : {\cal M} \rightarrow {\cal X}^n({\cal H}/{\cal M}, \C)$ of
$\nu$ satisfies the horizontality condition
\begin{equation} \label{eq-generalized-normal-fn}
\nabla \tilde{\nu} \in F^{i} {\cal H}^n({\cal X}/{\cal M}, \C) \otimes
\Omega^1_{\cal M}
\end{equation}
The section $\nu$ defines a splitting
$$
T_{(b,y)}{\cal J}^k = T_b{\cal M} \oplus
\left[ H^n(X_b,\C) / F^k H^n(X_b,\C) \right]
$$
and the $i$-th piece of the filtration is defined by
\begin{equation} \label{eq-filtration}
F^i T_{(b,y)}{\cal J}^k := T_b{\cal M} \oplus \left[ F^i H^n(X_b,\C) /
F^k H^n(X_b,\C) \right].
\end{equation}
The horizontality condition (\ref{eq-generalized-normal-fn}) implies that the
subspace $F^{i}T{\cal J}^k_{(b,y)}$  is independent of the choice of the
section
$\nu$ through $(b,y)$.  Moreover, the subbundle $F^iT{\cal J}^k$ is invariant
under translations by its integral sections, namely, by sections $\nu :
{\cal M} \rightarrow {\cal J}^k$ satisfying the i-th horizontality condition.
}
\end{enumerate}
\end{rem}

We noted above that when $k=1$ the intermediate Jacobian ${\cal J}^k(X)$
becomes the connected component $\Pic^{0}(X)$ of the Picard variety.  The
generalization of the Picard variety itself is the {\it Deligne cohomology}
group $D^k(X)$, cf. \cite{EZ}.  This fits in an exact sequence
\begin{equation} \label{eq-exact-seq-deligne-coho}
0 \rightarrow {\cal J}^k(X) \rightarrow D^k(X)
\stackrel{p}{\rightarrow} H^{k,k}(X,\Z) \rightarrow 0,
\end{equation}
where the quotient is the group of Hodge $(k,k)$-classes,
$$
H^{k,k}(X,\Z) := H^{k,k}(X,\C) \cap H^{2k}(X,.\Z).
$$
Any codimension-$k$ cycle $Z$ in $X$ has an Abel-Jacobi image, or cycle
class $\mu(Z)$ in $D^k(X)$.  Its image $p(\mu(Z))$ is the
cycle class of $Z$ in ordinary cohomology.

Formally, $D^k(X)$ is defined as the hypercohomology $\Hyper^{2k}$ of the
following complex of sheaves on $X$ starting in degree $0$.
$$
0 \rightarrow \Z \rightarrow {\cal O}_X \rightarrow \Omega^1_X \rightarrow
\cdots \rightarrow \Omega^{k-1}_X \rightarrow 0.
$$
The forgetful map to $\Z$ is a map of complexes, with kernel the complex
$$
0 \rightarrow {\cal O}_X \rightarrow \Omega^1_X \rightarrow \cdots
\rightarrow \Omega^{k-1}_X \rightarrow 0.
$$
The resulting long exact sequence of hyper cohomologies gives
(\ref{eq-exact-seq-deligne-coho}).

Let $H^{k,k}_{alg}$ be the subgroup of $H^{k,k}(X,\Z)$ of classes of
algebraic cycles.  (The Hodge conjecture asserts that $H^{k,k}_{alg}$ is of
finite index in $H^{k,k}(X,\Z)$.)  The inverse image
$$
D^k_{alg}(X) := p^{-1} (H^{k,k}_{alg})
$$
has an elementary description:  it is the quotient of
$$
{\cal J}^k(X) \times \; \{ {\rm codimension-}k \; {\rm algebraic \;
cycles} \}
$$
by the subgroup of codimension-$k$ cycles homologous to $0$, embedded
naturally in the second component and mapped to the first by Abel-Jacobi.

As $X$ varies in a family, the rank of $H^{k,k}(X,\Z)$ can jump up (at
those $X$ for which the variable vector subspace $H^{k,k}(X,\C)$  happens
to be in special position with respect to the ``fixed'' lattice
$H^{2k}(X,\Z)$).  To obtain a well-behaved family of Deligne cohomology
groups, we require that
$$
H^{k,k}(X,\C)  = H^{2k}(X,\C)  .
$$
For example, this holds for $k=1$ or $k=n-1$ if $h^{2,0} = 0$.  In this
case we also have $H^{k,k}_{alg} = H^{k,k}(X,\Z)$  and hence $D^k_{alg}(X)
= D^k(X)$, by the Lefschetz theorem on $(1,1)$-classes
\cite{griffiths-harris}.
\begin{corollary} \label{cor-contact-structure-extends-to-deligne-coho}
Let ${\cal X} \rightarrow {\cal M}$ be a complete family of $3$-dimensional
Calabi-Yau manifolds, $\tilde{\cal X} \rightarrow \tilde{\cal M}$ the
corresponding gauged family.
Let ${\cal D} \rightarrow {\cal M}$, $\tilde{\cal D} \rightarrow
\tilde{\cal M}$, be their families of (second) Deligne cohomology groups,
${\cal J}$, $\tilde{\cal J}$ their relative intermediate Jacobians.  Then
there is a natural contact structure $\kappa$ on ${\cal D}$ with
symplectification $\sigma = d \alpha$ on $\tilde{\cal D}$ with the
following properties:

\begin{enumerate}
\item[(a)] $\sigma$, $\alpha$, and $\kappa$ restrict to the previously
constructed structures on $\tilde{\cal J}$ and ${\cal J}$.

\item[(b)] The fibration $\tilde{\cal D} \rightarrow \tilde{\cal M}$
is Lagrangian.

\item[(c)] The multivalued normal functions of ${\cal D}$ (resp.
$\tilde{\cal D}$) are precisely the $\kappa$-integral (resp.
$\alpha$-integral) subvarieties.  In particular, all multi-valued normal
functions in $\tilde{\cal D}$ are isotropic.
\end{enumerate}
\end{corollary}
\medskip
\noindent
{\bf Proof:} The contact structure $\kappa$ on ${\cal J}$ defines one on
${\cal J} \times \{ cycles \}$, which descends to ${\cal D}$ since the
equivalence relation is $\kappa$-integral by remark
\ref{rem-integrality-of-normal-fn}.
\EndProof

\bigskip
The mirror conjecture of conformal field theory predicts that to a family
${\cal X} \rightarrow {\cal M}$ of Calabi-Yau three folds, with some extra
data, corresponds a ``mirror'' family ${\cal X}' \rightarrow {\cal M}'$, cf.
\cite{morrison-guide} for the details.
A first property of the conjectural symmetry
is that for $X \in {\cal M}$, $X' \in {\cal M}'$,
$$
h^{2,1}(X) = h^{1,1}(X'), h^{1,1} = h^{2,1}(X').
$$
The conjecture goes much deeper, predicting a relation between the Yukawa
cubic of ${\cal M}$ and the numbers of rational curves of various homology
classes in a typical $X' \in {\cal M}'$.  This has been used spectacularly
in \cite{candelas} and subsequent works, to predict those numbers on a
non-singular quintic hypersurface in $\bP^4$ and in a number of other
families.

We wonder whether the conjecture could be reformulated and understood as a
type of Fourier transform between the integrable systems on the universal
Deligne cohomologies  $\tilde{\cal D}$ and $\tilde{\cal D}'$ of the mirror
families $\tilde{\cal M}$ and $\tilde{\cal M}'$.  Note that the dimensions
$h^{2,1}$ and $h^{1,1}$ which are supposed to be interchanged by the mirror,
can be read
off the continuous and discrete parts of the fibers of $\pi : \tilde{\cal
D} \rightarrow \tilde{\cal M}$, respectively.  One may try to imagine the
mirror as a transform, taking these Lagrangian fibers over $\tilde{\cal M}$
(which encode the Yukawa cubic, as in Section
\ref{subsec-cubic-condition}) to Lagrangian sections
over $\tilde{\cal M}'$, which should somehow encode the numbers of
curves in $X'$ via
their Abel-Jacobi images.
\newpage
%


\section{The Lagrangian Hilbert scheme and its relative Picard} \label{ch8}
\label{sec-lagrangian-sheaves}
\subsection{Introduction}

The Lagrangian Hilbert scheme of a symplectic variety $X$ parametrizes
Lagrangian subvarieties of $X$. Its relative Picard parametrizes
pairs $(Z,L)$ consisting of a line bundle $L$ on a Lagrangian
subvariety $Z$.
We use the cubic condition of chapter \ref{ch7} to construct an integrable
system
structure on components of the relative Picard bundle over the
Lagrangian Hilbert scheme.

We interpret the generalized Hitchin integrable system, supported by
the moduli space of Higgs pairs over an algebraic curve (see Ch V), as a
special case of this construction.  Other examples discussed include:
\begin{description}
\item [a)] Higgs pairs over higher dimensional base varieties (example
\ref{moduli-higgs-pairs-as-lagrangian-sheaves}), and
\item [b)] Fano varieties of lines on hyperplane sections of a cubic
fourfold (example \ref{subsec-fanos-of-cubics}).
\end{description}


Understanding the {\em global} geometry of such an
integrable system requires a compactification and a study of
its boundary. Our compactifications  of the relative Picard
are moduli spaces of sheaves and
we study the symplectic structure at (smooth, stable) points of
the boundary.

\bigskip
Let $X$ be a smooth projective symplectic algebraic variety, $\sigma$ an
everywhere non degenerate algebraic 2-form on $X$.  A smooth projective
Lagrangian subvariety $Z_0$ of $X$ determines a component $\bar{B}$ of the
Hilbert scheme parametrizing deformations of $Z_0$ in $X$.  The component
$\bar{B}$ consists entirely of Lagrangian subschemes.  Its dense open
subset $B$, parametrizing smooth deformations of $Z_0$, is a smooth
quasi-projective variety \cite{ziv-ran-lifting,voisin}.

Choose a very ample line bundle ${\cal O}_X(1)$ on $X$ and a Hilbert
polynomial $p$.  The relative Picard $h:{\cal M}^{p} \rightarrow B$,
parametrizing line bundles with Hilbert polynomial $p$ which are
supported on Lagrangian subvarieties of $X$, is a quasi-projective variety
(see \cite{simpson-moduli}).
If the Chern class $c_1(L_0) \in H^2(Z_0,\Integers)$ of a
line bundle on $Z_0$ deforms as a $(1,1)$-class over the whole of $B$, then
$L_0$ belongs to a component ${\cal M}$ of ${\cal M}^p$  which {\it
dominates} the Hilbert scheme $B$.
(By Griffiths' and Deligne's Theorem of the Fixed Part,
\cite{schmid-vhs-the-singularities} Corollary 7.23,
this is the case for example, if
$c_1(L_0)$ belongs to the image of $H^2(X,\Q)$).  Such components ${\cal M}$
are integrable systems, in other words:

\begin{theorem} \label{thm-symplectic-structure-on-relative-picard}
There exists a canonical symplectic structure $\sigma_{\cal M}$ on the
relative Picard bundle ${\cal M} \stackrel{h}{\rightarrow} B$ over the open
subset $B$ of the Hilbert scheme of smooth projective Lagrangian
subvarieties of $X$.  The support map $h:{\cal M} \rightarrow B$ is a
Lagrangian fibration.
\end{theorem}

The relative Picard over the Hilbert scheme of curves on a $K3$ or abelian
surface is an example \cite{mukai}.  In example
\ref{subsec-fanos-of-cubics}, $X$ is a symplectic fourfold.
\begin{rem} \label{rem-conditions-for-thm-symp-case}
{\rm
Theorem \ref{thm-symplectic-structure-on-relative-picard} holds in a more
general setting where $X$ is a smooth projective
algebraic variety, $\sigma$ is a meromorphic, generically non degenerate
closed 2-form on $X$.  We let $D_0$ denote its degeneracy divisor,
$D_\infty$ its polar divisor, and set $D = D_0 \cup D_\infty$.  Let $Z_0$
be a smooth projective Lagrangian subvariety of $X$ which does not
intersect $D$.  Denote by $B$ the open subset of a component of the Hilbert
scheme parametrizing smooth deformations of $Z_0$ which stay in $X-D$.
Then $B$ is smooth and Theorem 1 holds.
A special case is when $X$ has a
generically non-degenerate Poisson structure $\psi$. In this case
$D_{\infty}$, the polar divisor of the inverse symplectic
structure, is just the degeneracy locus of $\psi$, while $D_{0}$ is empty.
The case where the subvariety $Z_0$  does intersect the degeneracy
locus $D_{\infty}$ of the Poisson structure is also of interest.
It is discussed below
under the category of Poisson integrable systems.
 }
\end{rem}

The moduli space of 1-form valued Higgs pairs is related to the case where
$X = \bP(\Omega^1_Y \oplus {\cal O}_Y)$ is the compactification of the
cotangent bundle of a smooth projective algebraic variety $Y$, and $D = \bP
\Omega^1_Y$ is the divisor at infinity
(see example \ref{moduli-higgs-pairs-as-lagrangian-sheaves}).

\smallskip
The relative Picard bundle ${\cal M}$ is in fact also a Zariski open subset
of a component of the moduli space of stable coherent sheaves on $X$ (see
\cite{simpson-moduli} for the construction of the moduli space).
Viewed in this
way, Theorem \ref{thm-symplectic-structure-on-relative-picard}
 extends a result of Mukai \cite{mukai} for sheaves on a $K3$ or
abelian surface.

\noindent
{\bf Theorem} \cite{mukai}:
{\it Any component of the moduli space of simple
sheaves on $X$ is smooth and has a canonical symplectic structure.}

Kobayashi \cite{kobayashi} generalized the above theorem to the case of
simple vector
bundles on a (higher dimensional) compact complex symplectic manifold $(X,
\sigma)$:\\
{\it The smooth part of the moduli space has a canonical symplectic
structure.}

In view of Theorem \ref{thm-symplectic-structure-on-relative-picard}
and Kobayashi's result one might be tempted to
speculate that every component of the moduli space of (simple) sheaves on a
symplectic algebraic variety has a symplectic structure.  This is {\em false}.
In fact, some components are odd dimensional
(see example \ref{example-odd-dimensional-moduli-spaces}).

\bigskip
Returning to our symplectic relative Picard $\M$, it is natural to ask
whether its {\em compactification} is symplectic. More precisely:

\begin{description}
{\it \item [(i)]  Does the symplectic structure extend to the smooth locus of
the closure of the relative Picard ${\cal M}$ in the moduli
space of stable (Lagrangian) sheaves?
\item[(ii)] Which of these components $\bar{\M}$ admits a smooth
projective birational model which is  symplectic? }
\end{description}
A partial answer to (i) is provided in Theorem
\ref{thm-extension-of-symplectic-str}.
We provide a cohomological identification of the
symplectic structure which extends as a 2-form $\sigma_{\bar{\M}}$ over the
smooth locus of $\bar{\M}$.  We do not know at the moment if the 2-form
$\sigma_{\bar{\M}}$ is {\em non-degenerate} at every smooth point of
$\bar{\M}$.
The cohomological identification of $\sigma_{\bar{\M}}$ involves a surprisingly
rich {\em polarized Hodge-like structure} on the algebra
$\Ext^*_X(L,L)$ of extensions of a Lagrangian line bundle $L$ by itself
as an $\StructureSheaf{X}$-module.

Much of the above generalizes to Poisson integrable systems.
Tyurin showed in \cite{tyurin-symplectic} that Mukai's theorem generalizes to
Poisson surfaces:

\noindent
{\it The smooth part of any component of the moduli space of simple sheaves
on a Poisson surface has a canonical Poisson structure.  }

\smallskip
\noindent
When the sheaves are supported as line bundles on curves in the surface, we
get an integrable system. More precisely:

\begin{theorem} \label{thm-poisson-structure-on-relative-picard}
Let $(X,\psi)$ be a Poisson surface, $D_\infty$ the degeneracy divisor of
$\psi$.  Let $B$ be the Zariski open subset of a component of the Hilbert
scheme of $X$ parametrizing smooth irreducible curves on $X$ which are not
contained in $D_\infty$.  Then
\begin{description}
\item [i)] $B$ is smooth,
\item [ii)] the relative Picard bundle $h:{\cal M} \rightarrow B$ has a
canonical Poisson structure $\psi_M$,
\item [iii)] The bundle map $h : {\cal M} \rightarrow B$ is a Lagrangian
fibration and
\item[iv)] The symplectic leaf foliation of ${\cal M}$ is induced by the
canonical morphism $B \rightarrow \; {\rm Hilb}_{D_\infty}$ sending a curve
$Z$ to the subscheme $Z \cap D_\infty$ of $D_\infty$.
\end{description}
\end{theorem}
The generalization to higher dimensional Poisson varieties
is treated here under rather restrictive conditions on the component of the
Lagrangian Hilbert scheme (see condition \ref{setup-for-poisson-case}).
These restrictions will be relaxed in \cite{markman-lagrangian-sheaves}.

\smallskip
The rest of this chapter is organized as follows: In section
\ref{subsec-hilbert-schemes} we review the deformation theory of
Lagrangian subvarieties.
The construction of the symplectic structure is
carried out in section \ref{subsec-construction}
where we prove Theorems
\ref{thm-symplectic-structure-on-relative-picard} and
\ref{thm-poisson-structure-on-relative-picard}.
In section \ref{subsec-extension-to-singular-lagrangian-sheaves}
we outline the extension of the symplectic structure to the
smooth locus of the moduli space of Lagrangian sheaves (Theorem
\ref{thm-extension-of-symplectic-str}).
We discuss the examples of Higgs pairs and of Fano varieties
of lines on cubics in section
\ref{subsec-examples-of-lagrangian-sheaves}.

\subsection{Lagrangian Hilbert Schemes}\label{subsec-hilbert-schemes}

Let $X$ be a smooth $n$-dimensional projective algebraic variety, $Z
\subset X$ a codimension $q$ subvariety and ${\cal O}_X(1)$ a very ample
line bundle.  The Hilbert polynomial $p$ of $Z$ is defined to be
$$
p(n) := \chi\Big( {\cal O}_Z(n) \Big) := \sum (-1)^i \dim H^i\Big( Z,
{\cal O}_Z(n) \Big).
$$

Grothendieck proved in \cite{grothendieck-existence}
that there is a projective scheme ${\rm
Hilb}^p_X$ parametrizing all algebraic subschemes of $X$ with Hilbert
polynomial $p$ and having all the expected functoriality and naturality
properties.

The Zariski tangent space $T_{[Z]}{\rm Hilb}^p_X$ at the point $[Z]$
parametrizing a subvariety $Z$ is canonically identified with the space of
sections $H^0(Z,N_{X/Z})$ of the normal bundle (normal sheaf if $Z$ is
singular).

The scheme ${\rm Hilb}^p_X$ may, in general, involve pathologies.  In
particular, it may be non-reduced.  A general criterion for the smoothness
of the Hilbert scheme at a point $[Z]$ parameterizing a locally complete
intersection subscheme $Z$ is provided by:
\begin{definition}
The semi-regularity map $\pi : H^1(Z,N_{X/Z}) \longrightarrow
H^{q+1}(X,\Omega^{q-1}_X)$ is the dual of the natural homomorphism
$$
\pi^* : H^{n-q-1} \Big( X,\Omega^{n-q+1}_X \Big) \longrightarrow H^{n-q-1}
\Big(Z, \omega_Z \otimes N^*_{Z/X} \Big).
$$
Here $\omega_Z \simeq \; \stackrel{q}{\wedge} N_{Z/X} \otimes \omega_X$ is the
dualizing sheaf of $Z$ and the homomorphism $\pi^*$ is induced by the sheaf
homomorphism
\begin{equation} \label{eq-sheaf-homomorphism-inducing-the-semiregularity-map}
\Omega^{n-q+1}_X \simeq \omega_X \otimes \stackrel{q-1}{\wedge}T_X
\longrightarrow  \omega_X \otimes \Wedge{q-1}N_{Z/X} \cong
\omega_{Z}\otimes\Normal{Z}{X}^{*}.
\end{equation}
\end{definition}

\smallskip
\noindent
{\bf Theorem}
{\em
(Severi-Kodaira-Spencer-Bloch \cite{kawamata})
If the semi-regularity map $\pi$ is injective, then the Hilbert scheme is
smooth at $[Z]$.
}

\smallskip
Together with a result of Ran it implies:

\begin{corollary}
Let $(X,\psi)$ be a Poisson surface with a degeneracy divisor $D_\infty$
(possibly empty).  Let $Z \subset X$ be a smooth irreducible curve which is
{\it not} contained in $D_\infty$.  Then the Hilbert scheme ${\rm
Hilb}^p_X$ is smooth at $[Z]$.
\end{corollary}

{\bf Proof:}  The Poisson structure induces an injective homomorphism
$   \phi : N^*_{Z/X} \  \hookrightarrow \  T_Z. $
If $Z$ intersects $D_\infty$ non-trivially then
$N_{Z/X} \simeq \omega_Z(Z \cap D_\infty)$
and hence $H^1(Z, N_{X/Z}) = (0)$ and the semi-regularity map
is trivially injective.

Note that in our case $n=2$, $q=1$ and
the dual of the semi-regularity map
$$
\pi^* : H^{0} \Big(X,\omega_{X} \Big) \longrightarrow H^{0}
\Big(Z, \omega_Z \otimes N^*_{Z/X} \Big)
$$
is induced by the sheaf homomorphism
\[
\omega_{X} \rightarrow   \omega_{\restricted{X}{Z}}
           \rightarrow \omega_Z \otimes N^*_{Z/X}
\]
given by (\ref{eq-sheaf-homomorphism-inducing-the-semiregularity-map}).
If $D_\infty = \emptyset$ $(X$ is symplectic) then
$\omega_{X}$, $\omega_{\restricted{X}{Z}}$ and $\omega_Z \otimes N^*_{Z/X}$
are all trivial line bundles and hence both $\pi^*$ and
the semi-regularity map are isomorphisms.  If $D_\infty \cap Z =
\emptyset$ but $D_\infty \not= \emptyset$ then $\pi$ fails to be injective
but the result nevertheless holds by a theorem of Ran which we recall below
(Theorem \ref{thm-voisin-ziv-ran}).
\EndProof

\bigskip
The condition that the curve $Z$ is not contained in $D_\infty$ is
necessary as can be seen by the following counterexample due to Severi and
Zappa:

\begin{example}
{\rm (\cite{mumford-curves-on-surface} Section 22)
Let $C$ be an elliptic curve, $E$ a nontrivial extension $0 \rightarrow E_1
\rightarrow E \rightarrow E_2 \rightarrow 0$, \ $E_i \simeq {\cal O}_C$ and
$\pi : X = \bP(E) \rightarrow C$ the corresponding ruled surface over $C$.
Denote by $Z$ the section $s:C \rightarrow X$ given by the line subbundle
$E_1 \subset E$.  Let ${\cal O}_X(-1)$ be the tautological subbundle of
$\pi^*E$. Then ${\cal O}_X(1)$ is isomorphic to the
line bundle ${\cal O}_X(Z)$ and
the canonical bundle $\omega_X$ is isomorphic to
$\pi^*\Big(\omega_C\Big) \otimes {\cal O}_X(-2) \simeq {\cal O}_X(-2)$.
$H^0(X, \stackrel{2}{\wedge} T_X)$ is thus isomorphic to
$H^0(C,\Sym^2 E^*)$
which is one dimensional.
It follows that $X$ has a
unique Poisson structure $\psi$ up to a scalar factor.  The divisor
$D_\infty = 2Z$ is the degeneracy divisor of $\psi$.

Clearly, $N_{Z/X} \simeq \pi^* T_C \simeq T_Z$ and hence $H^0(Z, N_{Z/X})$
is one dimensional.  On the other hand, $Z$ has no deformations in $X$ (its
self intersection is $0$ and a deformation $Z'$ of $Z$ will contradict the
nontriviality of the extension $0 \rightarrow E_1
\rightarrow E \rightarrow E_2 \rightarrow 0$).
\EndProof
}
\end{example}
A curve $Z$ on a symplectic surface $X$ is automatically Lagrangian.
In the higher dimensional case we replace the curve $Z$ by a Lagrangian
subvariety. Lagrangian subvarieties of symplectic varieties have two
pleasant properties:
\begin{description}
{\it \item [i)] The condition of being Lagrangian is both open and closed,
\item [ii)] Their deformations are unobstructed.}
\end{description}
More precisely, we have:

\begin{theorem} (Voisin \cite{voisin}, Ran \cite{ziv-ran-lifting})
\label{thm-voisin-ziv-ran}
Let $X$ be a smooth projective algebraic variety, $\sigma$ a generically non
degenerate meromorphic closed $2$-form, $D_\infty$ its polar divisor, $D_0$
its degeneracy divisor.  Assume that $Z_0 \subset X - D_\infty - D_0$ is a
smooth {\it projective} Lagrangian subvariety.  Then
\begin{description}
\item [(i)] The subset of the Hilbert scheme ${\rm Hilb}^p_X$ parametrizing
deformations of $Z_0$ in $X - D_\infty$ consists entirely of Lagrangian
subvarieties.
\item [(ii)] The Hilbert scheme is smooth at $[Z_0]$.
\end{description}
\end{theorem}
\noindent
{\bf Sketch of Proof:} \ \  (i) \ \ The Lagrangian condition is closed.
Thus, it suffices to prove that the open subset of smooth deformations of
$Z_0$ is Lagrangian.  If $Z \subset X-D_\infty$ then $\sigma_{|_Z}$ is a
closed holomorphic $2$-form and the cohomology class $[\sigma_{|_Z}]$ in
$H^{2,0}(Z)$ vanishes if and only if $\sigma_{|_Z}$ is identically zero.  Since
$\sigma$ induces a {\it flat} section of the Hodge bundle of relative
cohomology with $\C$-coefficients, then $[\sigma_{|_Z}] = 0$ is an open and
closed condition.

(ii).  The symplectic structure $\sigma$ induces a canonical isomorphism
$N_{Z/X} \simeq \Omega^1_Z$ for any Lagrangian projective smooth subvariety
$Z \subset X - D_\infty - D_0$. Ran proved a criterion for
unobstructedness of deformations:  the $T^1$-lifting property (see
\cite{ziv-ran-lifting,kawamata}).  Let $S_n = {\rm Spec}(\C[t]/t^{n+1})$.  Any
flat $(n+1)$-st order infinitesimal embedded deformation $Z_{n+1} \rightarrow
S_{n+1}$ of $Z_0 = Z$ restricts canonically to an n-th order deformation
$Z_n \rightarrow S_n$.  In our context, the $T^1$-lifting property amounts
to the following criterion:

\noindent
{
\it Given any $(n+1)$-st order flat embedded deformation $Z_{n+1}
\rightarrow S_{n+1}$, every extension \\
(a) of $Z_n \rightarrow S_n$ to a flat embedded deformation
$\tilde{Z}_n \rightarrow S_n \times_{\ComplexNumbers} S_1$

\noindent
\smallskip
lifts to an extension

\noindent
\smallskip
(b) of $Z_{n+1} \rightarrow S_{n+1}$ to $\tilde{Z}_{n+1}
\rightarrow S_{n+1} \times_{\ComplexNumbers} S_1$.
}

\smallskip
\noindent
Extensions in (a) and (b) are classified by $T^1(Z_i / S_i) \cong H^0 (Z_i,
{\cal N}_{\varphi_i/S_i})$ where $\varphi_i : Z_i \rightarrow S_i \times X$
is the canonical morphism and
${\cal N}_{\varphi_i/S_i} $ is the relative normal sheaf.
Recall that the De Rham cohomology and its Hodge filtration
can be computed using the algebraic De Rham complex
(\ref{eq-the-quotient-of-the-algebraic-derham-complex}).
Consequently, the discussion of part (i) applies in the infinitesimal setting
to show that $T_{{\cal Z}_i/S_i}$ is {\em Lagrangian} as a subbundle of
the pullback $(\varphi^i)^*T_X$ with respect to the
non-degenerate $2$-form $(\varphi^i)^*(\sigma)$ on $(\varphi^i)^*T_X$.
The relative normal sheaf is the quotient
\[
0 \rightarrow T_{{\cal Z}_i/S_i} \rightarrow
(\varphi^i)^*T_X \rightarrow
{\cal N}_{\varphi_i/S_i} \rightarrow 0.
\]
Hence the symplectic structure induces an
isomorphism ${\cal N}_{\varphi_i/S_i} \simeq \Omega^1_{Z_i/S_i}$.  By a
theorem of Deligne, $H^0\Big( \Omega^1_{Z_i/S_i}\Big)$, and hence also
$H^0\Big({\cal N}_{\varphi_i/S_i}\Big)$, is a free ${\cal O}_{S_i}$-module
\cite{deligne-leray-degenerates}.  Thus, $H^0\Big({\cal
N}_{\varphi_{n+1}/S_{n+1}}\Big) \longrightarrow H^0\Big( {\cal N}
_{\varphi_n/S_n}\Big)$ is surjective and the $T^1$-lifting property holds.
\EndProof

\bigskip

Note that the naive analogue of the above theorem fails for Poisson
varieties.  In general, deformations of Lagrangian subvarieties need not
stay Lagrangian.  Consider for example $(\bP^{2n},\psi)$ where the Poisson
structure $\psi$ is the extension of the standard (non degenerate)
symplectic structure on $\A^{2n}\subset \bP^{2n}$.  The Lagrangian
Grassmannian has positive codimension in $Gr(n+1, \; 2n+1)$.
\bigskip


\subsection{The construction of the symplectic structure}
\label{subsec-construction}

The construction of the symplectic structure on the relative Picard bundle
is carried out in three steps:

In Step I we reduce it to the construction of the symplectic structure on
the relative $Pic^0$-bundle.

In Step II we verify the cubic condition and thus construct the 2-form (or
the 2-tensor in the Poisson case).

In Step III we prove the closedness of the 2-form.

\bigskip
\noindent
{\bf Step I:}\ \  \underline{Reduction to the $Pic^0$-Bundle Case}:

The construction of a 2-form on the relative Picard bundle
${\M} \stackrel{h}{\rightarrow} B$ reduces to constructing it on its zero
component ${\M}^0 \stackrel{h}{\rightarrow} B$,
namely the $Pic^0$-bundle, by the following:
\begin{proposition}
Any closed 2-form $\sigma_{{\M}^0}$ on ${\M}^0$, with respect to
which the zero section of ${\M}^0$ is Lagrangian, extends to a closed
2-form $\sigma_{\M}$ on the whole Picard bundle $h :{\M} \rightarrow B$.
The extension $\sigma_{\M}$ depends canonically on $\sigma_{{\M}^0}$
and the polarization ${\cal O}_X(1)$ of $X$.
\end{proposition}

{\bf Proof:}  The point is that Picard bundles are rationally split.  For
any polarized projective variety $(Z, {\cal O}_Z(1))$, we have the Lefschetz
map
$$
Lef: Pic \   Z \longrightarrow  Alb \ Z
$$
$$
[D] \mapsto  \Big[ D \cap [{\cal O}_Z(1)]^{n-1} \Big]
$$
inducing an isogeny
$$   Lef^0 : Pic^0 \; Z \ \longrightarrow  Alb^0 \; Z.   $$
We can set
$$
L_Z := \left\{ s \in Pic \; Z\; | \; \exists \ \ell,m, \ \ \ell \not= 0,
\ {\rm such\; that}\; \ell \cdot Lef(s) = m \cdot Lef({\cal O}(1)) \right\}.
$$
This is an extension of $H^{1,1}_\Z(Z)$ by the torsion subgroup $L_Z^{tor}$ of
$Pic(Z)$.  In a family ${\cal Z} \rightarrow B$, these groups form a
subsheaf ${\cal L}$ of ${\M} := Pic({\cal Z}/B)$, intersecting ${\M}^0
:= Pic^0({\cal Z}/B)$ in its torsion subsheaf ${\cal L}^0$.  In our
situation, the 2-form $\sigma_{{\M}^0}$ is ${\cal L}^0$-invariant, so
it extends uniquely to an ${\cal L}$-invariant closed 2-form $\sigma_{\M}$
on ${\M}$.
\EndProof

\bigskip
\noindent
{\bf Step II:} \ \ \underline{Verification of the Cubic Condition:}

In this step we construct the $2$-form (or $2$-vector)
on the relative Picard bundle.
In the next step we will prove that it is closed
(respectively, a Poisson structure).

Let $(X,\psi)$ be a smooth projective variety, $\psi$ a generically
non-degenerate holomorphic Poisson structure.  Denote by $D_\infty$ the
degeneracy divisor of $\psi$. We will assume throughout this step that
$Z \subset X$ is a smooth
subvariety,  $Z \cap (X - D_\infty)$ is non empty and Lagrangian, and

\begin{condition} \label{setup-for-poisson-case}
\begin{description}
{\it \item [i)] $[Z]$ is a smooth point of the Hilbert scheme, and
\item [ii)]  all deformations of $Z$ in $X$ are Lagrangian.}
\end{description}
\end{condition}
\noindent

%

As we saw in the previous section, conditions i) and ii) hold in case
$\psi$ is
everywhere non-degenerate $((X, \psi^{-1})$ is a symplectic projective
algebraic variety), and also in case $X$ is a surface.  Such $[Z]$ vary in a
smooth Zariski open subset $B$ of the Hilbert scheme and we denote by
$h:{\M} \rightarrow B$ the relative $Pic^0$-bundle.
Condition \ref{setup-for-poisson-case} can be relaxed considerably
(see \cite{markman-lagrangian-sheaves}).

Let
\begin{equation} \label{eq-sheaf-homomorphism-induced-by-poisson-str}
\phi :N^*_{Z/X} \hookrightarrow T_Z
\end{equation}
be the injective
homomorphism induced by the Poisson structure $\psi$.  Its dual $\phi^* :
T^*_Z \rightarrow N_{Z/X}$ induces an injective homomorphism.
\begin{equation}
\label{eq-global-sections-homomorphism-induced-by-poisson-str}
i:H^0 \Big(Z,T^*_Z \Big) \hookrightarrow H^0 \Big( Z,{\cal N}_{Z/X}
\Big).
\end{equation}
The vertical tangent bundle $V := h_* {\cal T}_{{\M}/B}$ is isomorphic to the
Hodge bundle ${\cal H}^{0,1}({\cal Z}/B)$.  The polarization induces an
isomorphism $V^* \simeq {\cal H}^{1,0}$.  We get a global injective
homomorphism $i : V^* \hookrightarrow T_B$.

\begin{proposition} \label{prop-verification-of-cubic-condition}
The homomorphism $i$ is induced by a canonical
$2$-vector $\psi_{\M} \in \\
H^0({\M}, \; \stackrel{2}{\wedge} \; T_{\M})$ with respect to which
$h:{\M}\rightarrow B$ is a Lagrangian fibration.  (We do not assert yet
that $\psi_{\M}$ is a Poisson structure).
\end{proposition}

{\bf Proof:}  It suffices to show that $i$ satisfies the (weak) cubic
condition, namely, that $dp \circ i$ comes from a cubic.  The derivative of
the period map
\begin{equation} \label{eq-differential-of-period-map}
dp:H^0 \Big(Z,N_{Z/X}\Big) \longrightarrow  \ {\rm Hom} \Big( H^{1,0}(Z),
H^{0,1}(Z) \Big) \simeq \Big[ H^{1,0}(Z)^* \Big]^{\otimes 2}
\end{equation}
is identified by the composition
$$
H^0 \Big(Z,N_{Z/X}\Big) \stackrel{{\rm K-S}}{\longrightarrow}
H^1(Z,T_Z) \stackrel{VHS}{\longrightarrow} \Sym^2 H^{1,0}(Z)^*
$$
where K-S is the Kodaira-Spencer map given by cup product with the
extension class of $T_{\restricted{X}{Z}}$:
\begin{equation} \label{eq-extension-class-related-to-kodaira-spencer-map}
\tau \in \; {\rm Ext}^1 \Big(N_{Z/X}, T_Z \Big) \simeq H^1 \Big(Z, N^*_{Z/X}
\otimes T_Z \Big),
\end{equation}
and the variation of Hodge structure map
VHS is given by cup product and contraction
$$
H^1(Z, T_Z) \otimes H^0(Z,T^*_Z) \longrightarrow H^1(Z,{\cal O}_Z).
$$
The composition
$({\rm K\!\!-\!\!S})\; \circ \; i:H^0(Z,T^*_Z) \rightarrow H^1(Z,T_Z)$
is then given by cup product with the class $(\phi \; \otimes \;{\rm id})(\tau)
\in H^1(Z, \; T_Z \otimes T_Z)$. We will
show that $(\phi \; \otimes \;{\rm id})(\tau)$ is symmetric, that is,
an element of  $H^1(Z, \Sym^2T_Z)$. This would imply that $dp$, regarded
as a section of
$H^0 \Big(Z,N_{Z/X}\Big)^* \otimes \Sym^2 H^{1,0}(Z)^* \stackrel{i}{\cong}
H^{1,0}(Z)^* \otimes \Sym^2 H^{1,0}(Z)^*$, is symmetric also with respect to
the first two factors.
The cubic condition will follow.

Lemma
\ref{lemma-alternating-two-forms-and-symmetric-extensions} below
implies that $(\phi \otimes \;{\rm id})(\tau)$
is in  $H^1(Z, \Sym^2T_Z)$ if and only if $\phi$ is induced by a section $\psi$
in $H^0\Big(Z, \; \stackrel{2}{\wedge}  T_{X_{|Z}}\Big)$ with respect
to which $Z$ is Lagrangian (i.e., $N^*_{Z/X}$ is isotropic).  This is indeed
the way $\phi$ was defined.
\EndProof

\begin{lem} \label{lemma-alternating-two-forms-and-symmetric-extensions}
Let $T$ be an extension
\begin{equation} \label{eq-exact-seq-with-symmetric-extension-class}
0 \rightarrow Z \rightarrow T \rightarrow N \rightarrow 0
\end{equation}
of a vector bundle $N$ by a vector bundle $Z$.
Then the following are equivalent for any homomorphism
$\phi : N^* \rightarrow Z$.

\smallskip
\noindent
i) The homomorphism
$\phi$ is induced by a section $\psi \in H^0(\Wedge{2}{T})$
with respect to which $N^*$ is isotropic.

\smallskip
\noindent
ii) The homomorphism
$\phi_* := H^1(\phi\otimes 1) : H^1(N^*\otimes Z) \rightarrow
H^1(Z \otimes Z)$
maps the extension class
$\tau \in H^1(N^*\otimes Z)$ of $T$ to a symmetric class
$\phi_*(\tau)\in H^1(\Sym^2 Z) \subset H^1(Z \otimes Z)$.
\end{lem}

\noindent
{\bf Proof:}
We argue as in the proof of the cubic
condition (lemma \ref{lemma-weak-cubic-cond-poisson}).
The extension (\ref{eq-exact-seq-with-symmetric-extension-class})
induces an extension
\[
0 \rightarrow \Wedge{2}{Z} \rightarrow F \rightarrow Z\otimes N \rightarrow 0,
\]
where $F$ is the subsheaf of $\Wedge{2}{T}$ of sections with respect to which
$N^*$ is isotropic. The homomorphism $\phi$, regarded as a section of
$Z\otimes N$, lifts to a section
$\psi$ of $F$ if and only if it is in the kernel of
the connecting homomorphism
\[
\delta: H^0(Z\otimes N)  \rightarrow H^1(\Wedge{2}{Z}).
\]
The latter is given by a) pairing with the extension class $\tau$
\[
(\cdot)_{*}\tau : H^0(Z\otimes N) \rightarrow H^1(Z\otimes Z),
\]
followed by b) wedge product
\[
H^1(Z\otimes  Z) \RightArrowOf{\wedge} H^1(\Wedge{2}Z).
\]
Thus, $\delta(\phi)$ vanishes if and only if $\phi_{*}(\tau)$ is in the
kernel of $\wedge$, i.e., in $H^1(\Sym^2 Z)$.
\EndProof

\bigskip
The identification of the cubic is particularly simple in the case of a
curve $Z$ on a surface $X$.  In that case Serre's duality identifies VHS
with the dual of the multiplication map
$$
\Sym^2 H^0(Z,\omega_Z) \stackrel{VHS^*}{\longrightarrow} H^0(Z,
\omega_Z^{\otimes 2}).
$$
The cubic $c \in \Sym^3H^0(Z,\omega_Z)^*$ is given by composing the
multiplication
$$
\Sym^3 H^0(Z, \omega_Z) \longrightarrow H^0(Z, \omega^{\otimes 3}_Z)
$$
with the linear functional
$$
(\phi \otimes \;{\rm id})(\tau) \in  H^1(Z,T^{\otimes 2}_Z) \simeq H^0(Z,
\omega_Z^{\otimes 3})^*
$$
corresponding to the extension class $\tau$.

In higher dimension (say $n$), the cubic depends on the choice of a
polarization $\alpha \in H^{1,1}(X)$:
\[
\Sym^3 H^0(Z, \Omega^1_Z) \rightarrow H^0(Z, \Sym^3\Omega^1_Z)
\LongRightArrowOf{\phi_*(\tau)} H^1(Z,\Omega^1_Z)
\LongRightArrowOf{\restricted{\alpha}{Z}^{n-1}} H^{n,n}(Z)
\cong \ComplexNumbers.
\]
The choice of $\alpha$
is implicitly made in the proof of proposition
\ref{prop-verification-of-cubic-condition}
when we identify $H^{0,1}(Z)$ with $H^{1,0}(Z)^*$ via the Lefschetz
isomorphism (see (\ref{eq-differential-of-period-map}))
{}.

\bigskip
\noindent
{\bf Step III:} \ \ \underline{Closedness}:

In this step we prove that the canonical $2$-vector $\psi_{\M}$
constructed in the previous step is a Poisson structure.  We first prove
it in the symplectic case and later indicate the modifications needed for
the Poisson case (assuming condition
\ref{setup-for-poisson-case} of the previous step).
This completes the proof of Theorems
\ref{thm-symplectic-structure-on-relative-picard} and
\ref{thm-poisson-structure-on-relative-picard} stated in the introduction to
this chapter.

\noindent
\underline{Symplectic Case:}

We assume, for simplicity of exposition, that $(X,\sigma)$ is a smooth
projective symplectic algebraic variety.  The arguments apply verbatim to
the more general setup involving a smooth projective algebraic variety $X$,
a closed generically non-degenerate meromorphic $2$-form $\sigma$ on $X$
with degeneracy divisor $D_0$ and polar divisor $D_\infty$, and Lagrangian
smooth projective subvarieties which do not intersect $D_0 \cup D_\infty$.

We then have a non-degenerate $2$-tensor $\psi_{\M}$ on $h: \M
\rightarrow B$ and hence a $2$-form $\sigma_{\M}$.  The closedness of
$\sigma_{\M}$ follows from that of $\sigma_X$ as we now show.
A polarization of X induces a
relative polarization on the universal Lagrangian subvariety
\[
{\divide\dgARROWLENGTH by 2
\begin{diagram}
\node{{\cal Z}}
\arrow{s,l}{\pi}
\arrow{e}
\node{B\times X}
\arrow{sw,r}{p_{B}}
\arrow{se,r}{p_{X}}
\\
\node{B}
\node[2]{X.}
\end{diagram}
}
\]
The relative polarization induces an isogeny
 \[
{\divide\dgARROWLENGTH by 2
\begin{diagram}[B]
\node{{\cal M}}
\arrow{se} \arrow[2]{e}
\node[2]{{\cal A}}
\arrow{sw}
\\
\node[2]{B}
\end{diagram}
}
\]
between the relative $Pic^0$-bundle and the relative Albanese
$h : {\cal A} \rightarrow B$.
Hence, a $2$-form $\sigma_{\cal A}$ on ${\cal A}$.
Clearly, closedness of $\sigma_{\M}$ is equivalent to that of $\sigma_{\cal
A}$.  Since the question is local, we may assume that we have a section
$\xi : B \rightarrow {{\cal Z}}$.
We then get for each positive integer $t$ a relative Albanese map
 \[
{\divide\dgARROWLENGTH by 2
\begin{diagram}
\node{{\cal Z}^{t}}
\arrow[2]{e,t}{a_{t}}
\arrow{se,l}{\pi}
\node[2]{{\cal A}}
\arrow{sw,l}{h}
\\
\node[2]{B}
\end{diagram}
}
\]
from the fiber product over $B$ of $t$ copies of the universal Lagrangian
subvariety ${\cal Z} \rightarrow B$.  For a fixed subvariety $Z_b$ and points
$(z_1, \dots, z_t) \in Z^t_b, \ \ a_t$ is given by integration
$$
\sum^t_{i=1} \ \int^{z_i}_{\xi(b)}(\cdot) \ \ ({\rm modulo} \ H_1(Z_b,
\Z)) \in H^{1,0}(Z)^* \Big/ H_1(Z_b,\Z).
$$
We may assume, by choosing $t$ large enough, that $a_t$ is surjective.
Thus, closedness of $\sigma_{\cal A}$ is equivalent to closedness of
$a^*_t(\sigma_{\cal A})$.
The closedness of $a^*_t(\sigma_{\cal A})$ now follows from
that of $\sigma_{X}$ by lemma
\ref{lemma-pullback-of-symplectic-structure-from-albanese}.

\begin{lem} \label{lemma-pullback-of-symplectic-structure-from-albanese}
Let $\ell : {\cal Z}^t \rightarrow X^t$ be the natural morphism;
$\sigma_{X^t}$
the product symplectic structure on $X^t$.  Then,
\begin{equation}\label{eq-pulled-back-symplectic-str}
 a^*_t(\sigma_{\cal A}) = \ell^*(\sigma_{X^t}) -
\pi^*(\xi^t)^* \; \ell^* (\sigma_{X^t}) .
\end{equation}
\end{lem}

{\bf Proof:}  The fibers of $\pi : {\cal Z}^t \rightarrow B$
are isotropic with respect to the
$2$-forms on both sides of the equation (\ref{eq-pulled-back-symplectic-str}).
Hence, these $2$-forms induce (by contraction) homomorphisms
$$
f_{\cal A}, f_X : T_{Z^t_b} \longrightarrow N^*_{Z^t_b/{\cal Z}^t}.
$$
The section $\xi^t (B) \subset {\cal Z}^t$ is also isotropic with respect to
the $2$-forms on both sides of equation
(\ref{eq-pulled-back-symplectic-str}).
Thus, equality in (\ref{eq-pulled-back-symplectic-str})
will follow from equality of the induced homomorphisms $f_{\cal A}, f_X$.
Proving the equality $f_{\cal A}= f_X$ is a straightforward,
though lengthy, unwinding of cohomological identifications.

The relative normal bundle is identified as the pullback of the
tangent bundle of the Hilbert scheme
\[
N_{Z^t_b/{\cal Z}^t} \simeq
{\cal O}_{Z^t_b} \otimes (T_b B) \simeq {\cal O}_{Z^t_b} \otimes H^0
\Big(Z_b, N_{Z_b/X} \Big).
\]
We will show that the duals of
both $f_{\cal A}$ and $f_X$
\[
f_{\cal A}^*,f_X^*: {\cal O}_{Z^t_b} \otimes H^0
\Big(Z_b, N_{Z_b/X} \Big) \rightarrow
T_{Z^t_b}^*
\]
are identified as the composition of

\noindent
i) the diagonal homomorphism
$$
{\cal O}_{Z^t_b} \otimes H^0(Z_b, N_{Z_b/X}) \stackrel{\Delta}{\hookrightarrow}
{\cal O}_{Z^t_b} \otimes \Big[ H^0(Z_b, N_{Z_b/X})\Big]^t
\ \ \ \ \mbox{followed} \ \mbox{by}
$$

\noindent
ii) the evaluation map
$$
e_t :
{\cal O}_{Z^t_b} \otimes \Big[ H^0(Z_b,N_{Z_b/X})\Big]^t
\simeq
{\cal O}_{Z^t_b} \otimes H^0 \left( Z^t_b, N_{Z^t_b/X^t} \right)
\longrightarrow N_{Z^t_b/X^t}
\ \ \ \ \mbox{followed} \ \mbox{by}
$$

\noindent
iii) contraction with the $2$-form $\sigma_{X^t}$
\[
(\phi^{-1^*})^t : N_{Z^t_b/X^t} \stackrel{\sim}{\rightarrow} T_{Z^t_b}^*
\]
($\phi$ is given by contraction with the Poisson structure
(\ref{eq-sheaf-homomorphism-induced-by-poisson-str})).

\smallskip
\noindent
\underline{Identification of $f_{\cal A}$}:  (for simplicity assume t=1).

The 2-form $\sigma_{\cal A}$ is characterized as the unique 2-form with
respect to which the three conditions of lemma
\ref{lemma-weak-cubic-cond-poisson} hold, i.e.,
i) ${\cal A} \rightarrow B$ is a Lagrangian fibration, ii) the zero
section is Lagrangian, and iii) $\sigma_{\cal A}$ induces the
homomorphism
$$
H^0(\phi^{-1^*}) = i^{-1} : H^0 (Z, N_{Z/X})
\stackrel{\sim}{\longrightarrow} H^0(Z,T^*_Z).
$$
Thus, $a^*(\sigma_{\cal A})$ induces
$$
f_{\cal A}^* =
\left(
{\cal O}_{Z_b} \otimes H^0 \Big(Z_b,N_{Z_b/X} \Big)
\stackrel{(i^{-1})}{\longrightarrow}
{\cal O}_{Z_b} \otimes H^0 \Big( Z_b,T^*_{Z_b} \Big)
\stackrel{da^*}{\longrightarrow}
T_{Z_b}^*
\right)
$$
and the codifferential $da^*$ of the Albanese map is the evaluation map.

\noindent
\underline{Identification of $f_X$}: $(t=1)$

Both 2-forms $\ell^*(\sigma_X)$ and $\ell^*(\sigma_X)-\pi^* \xi^*
\ell^*(\sigma_X)$ induce the same homomorphism
$f_X^*: {\cal O}_{Z_b} \otimes T_bB \rightarrow T^*_{Z_b}$.
This homomorphism is the composition
$\phi^{-1^*} \circ  \overline{(d\ell)}$,
where $\overline{d\ell}$ is the
homomorphism $N_{Z_b/{\cal Z}} \rightarrow N_{Z_b/X}$ induced by the
differential of $\ell : {\cal Z} \rightarrow X$:
\[
{\divide\dgARROWLENGTH by 2
\begin{diagram}
\node{0} \arrow{e}
\node{T_{Z_b}}
\arrow{e}
\arrow{s,l}{=}
\node{(T_{\cal Z})_{|Z_b}}
\arrow{e}
\arrow{s,l}{d\ell}
\node{H^0( Z_b, N_{Z_b/X}) \otimes {\cal O}_{Z_b}}
\arrow{s,l}{\overline{d\ell}}
\arrow{e} \node{0}
\\
\node{0} \arrow{e}
\node{T_{Z_b}}
\arrow{e}
\node{(\ell^*TX)_{|Z_b}}
\arrow{e}
\node{N_{Z_b/X}}
\arrow{e} \node{0.}
\end{diagram}
}
\]
Clearly $\overline{d\ell}$ is given by evaluation. This completes the
proof of lemma \ref{lemma-pullback-of-symplectic-structure-from-albanese}.
\EndProof

\bigskip

As a simple corollary of lemma
\ref{lemma-pullback-of-symplectic-structure-from-albanese} we have:

\begin{corollary} \label{cor-canonical-symplectic-str-on-albanese}
There exists a canonical symplectic structure $\sigma_{{\cal A}^t}$ on the
relative Albanese of degree $t\in \Integers$,
depending canonically on the symplectic
structure $\sigma_X$ (independent of the polarization ${\cal O}_X(1)$!) and
satisfying, for $t \geq 1$,
$$
a_t^* (\sigma_{{\cal A}^t}) = \ell^*(\sigma_{X^t}).
$$
(the pullback to the fiber product
${\cal Z}^t := \times^t_B {\cal Z}$ via the Albanese map
coincides with the pullback of the symplectic structure $\sigma_{X^t}$ on
$X^t$).
\end{corollary}
{\bf Proof:}  The $t=0$ case is proven. We sketch the proof of the
$t \geq 1$ case. The $t \leq -1$ case is similar.
Let $\xi$ be a local section of ${\cal Z}^t \rightarrow B$.
Translation by the section $-a_t(\xi)$ of ${\cal A}^{-t}$ defines a local
isomorphism
$$
\tau_\xi : {\cal A}^t \longrightarrow {\cal A}^0.
$$
Let
$$
\sigma_{{\cal A}^t} \ := \ \tau^*_\xi (\sigma_{{\cal A}^0}) +
h^*\xi^*\ell^*(\sigma_{X^t}).
$$
We claim that $\sigma_{{\cal A}^t}$ is independent of $\xi$.  This amounts
to the identity
\[
\tau^*_{(\xi_1 - \xi_2)}(\sigma_{{\cal A}^0}) = \sigma_{{\cal A}^0}
- h^*[a_0(\xi_1-\xi_2)]^* \sigma_{{\cal A}^0}
\]
for any two sections $\xi_1,\xi_2$ of ${\cal Z}^t \rightarrow B$.
\EndProof

\bigskip
\noindent
\underline{Poisson Case}: (assuming condition \ref{setup-for-poisson-case})

Showing that the $2$-vector $\psi_{\M}$ constructed in  step II
is a Poisson structure, amounts to showing that
\begin{lem} \label{lemma-involutive-distribution}
 $\psi_{\M}(T^*_{\M}) \subset T_{\M}$ is an involutive distribution,
\end{lem}
and
\begin{lem} \label{lemma-closedness-poisson-case}
the induced $2$-form on each symplectic leaf is closed.
\end{lem}

{\bf Sketch of Proof of Lemma \ref{lemma-involutive-distribution}}:  Since
$h : {\M} \rightarrow B$ is a Lagrangian fibration with respect to
$\psi_{\M}$
(by proposition \ref{prop-verification-of-cubic-condition}),
the distribution is the pullback
of the distribution on the base $B$. The latter is
induced by the image of the injective homomorphism
$i: V^* \hookrightarrow T_B$ identified by
(\ref{eq-global-sections-homomorphism-induced-by-poisson-str})
$$
i = H^0(\phi^*) :
H^0(Z,T^*_Z) \hookrightarrow H^0(Z,N_{Z/X}).
$$
Recall (\ref{eq-sheaf-homomorphism-induced-by-poisson-str})
that $\phi$, in turn, is induced by the Poisson structure $\psi_X$
on $X$.
The involutivity now follows from that of
$\psi_X(T^*_X) \subset T_X$ by a deformation theoretic argument.  The
details are omitted.
\EndProof

\bigskip

In case $X$ is a surface, the degeneracy divisor $D_\infty$ of $\psi_X$ is
a curve and $iH^0(Z, T^*_Z) \subset H^0(Z, N_{Z/X})$ is the subspace of all
infinitesimal deformations of $Z$ which {\it fix} the divisor $Z \cap
D_\infty$.  Thus, the distribution $i(V^*)$ on the Hilbert scheme $B$
(as in the proof of lemma \ref{lemma-involutive-distribution})
corresponds to the foliation by
level sets of the algebraic morphism
\begin{eqnarray*}
R : B &\longrightarrow& \; {\rm Hilb}(D_\infty) \\
Z &\longmapsto& Z \cap D_\infty.
\end{eqnarray*}
The higher dimensional case is analogous.  The degeneracy divisor
$D_\infty$ has an algebraic rank stratification
$$
D_\infty = \bigcup^{n-1}_{r=0} \ D_\infty[2r] \qquad \qquad (\dim X = 2n).
$$
Each rank stratum is foliated, local analytically, by symplectic leaves.
The subspace
$$
i:H^0(Z,T^*_Z) \subset H^0(Z,N_{Z/X})
$$
is characterized as the subspace of all infinitesimal deformations of $Z$
which deform the subscheme $Z \cap D_\infty[2r]$ fixing the image $f(Z \cap
D_\infty[2r])$ with respect to any Casimir function $f$ on $D_\infty[2r]$.
As an illustration, consider the case where $X$ is the logarithmic
cotangent bundle $T^*_M (log(D))$ and $Z$
is a $1$-form with logarithmic poles along a divisor $D$ with normal
crossing.  In this case the residues induce the symplectic leaves foliation.

\noindent
{\bf Sketch of Proof of Lemma \ref{lemma-closedness-poisson-case}:}
The proof is essentially the same as in the symplectic case.  We consider
an open (analytic) subset $B_1$ of a leaf in $B$, the universal Lagrangian
subvariety
$
{\divide\dgARROWLENGTH by 2
\begin{diagram}
\node{{\cal Z}_1}
\arrow{e,t}{\ell}
\arrow{s}
\node{X}
\\
\node{B_1}
\end{diagram}
}
$,
and the relative Albanese
$\begin{array}{c}
{\cal A}\\
\downarrow \\
B_1
\end{array}.
$
One has to choose the section $\xi:B_1 \rightarrow {\cal Z}_1$
outside $\ell^{-1}(D_\infty)$ and
notice that the identity (\ref{eq-pulled-back-symplectic-str})
implies that the pullback
$\ell^*(\sigma^t_X)$ of the {\it meromorphic} closed $2$-form $\sigma^t_X$
(inverse of the generically non-degenerate Poisson structure
on the product of $t$ copies of $X$) is a {\it
holomorphic} $2$-form on ${\cal Z}^t_1$ (because $a^*_t(\sigma_{\cal A})$ is)
and that $a^*_t(\sigma_{\cal A})$ is {\it closed} (because
$\ell^*(\sigma_{X^t})$ is).
\EndProof

\subsection
{Partial compactifications: a symplectic structure on the moduli space of
Lagrangian sheaves}
\label{subsec-extension-to-singular-lagrangian-sheaves}

We describe briefly in this section the extension of the symplectic
structure on the relative Picard ${\M}$ to an algebraic $2$-form on the
smooth locus of a partial compactification.  For details see
\cite{markman-lagrangian-sheaves}.  For simplicity, we assume that
$(X,\sigma)$ is a smooth $2n$-dimensional projective symplectic variety.
We note that with obvious modifications, the extension of the
$2$-form will hold in the setup $(X,\sigma,D_0,D_\infty)$ as in remark
\ref{rem-conditions-for-thm-symp-case} allowing $\sigma$ to degenerate and
have poles away from the support of the sheaves.

When $X$ is a symplectic surface, some of these extensions give rise to
smooth projective symplectic compactifications \cite{mukai}.
These projective symplectic compactifications
appear also in the higher dimensional case:
\begin{example}
{\rm
A somewhat trivial reincarnation of a relative Picard of a linear
system on a K3 surface $S$ as a birational model of a relative Picard of a
Lagrangian Hilbert scheme over a higher dimensional symplectic
variety $X$ is realized as follows. Let $X$ be the Beauville
variety $S^{[n]}$ which is the resolution of the $n$-th symmetric
product of $S$ provided by the Hilbert scheme of zero cycles of length $n$
\cite{beauville-zero-first-chern-class}. The symmetric powers $C^{[n]}$ of
smooth curves on $S$ are smooth Lagrangian subvarieties of $S^{[n]}$.
Components of the relative Picard over the smooth locus in the
linear system $|C|$ are isomorphic to Zariski open subsets of
components of the relative Picard over the Lagrangian Hilbert scheme
of $S^{[n]}$.
}
\end{example}

This leads us to speculate that genuinely new examples
of smooth symplectic {\em projective} varieties
will arise as birational models of moduli spaces of Lagrangian line bundles.
(see
section \ref{subsec-fanos-of-cubics} for new {\em quasiprojective} examples).

We worked so far with a component ${\M} \rightarrow B$ of the relative
Picard of the universal smooth Lagrangian subvariety ${\cal Z} \rightarrow
B$ which dominates the corresponding component $\bar{B}$ of the Lagrangian
Hilbert scheme (i.e., if $L$ is supported on $Z$, $c_1(L) \in H^{1,1}_\Z(Z)$
remains of type $(1,1)$ over $B$).
Let $p(n) := \chi\Big(L \otimes_{{\cal O}_X} {\cal O}_X(n)\Big)$
be the Hilbert polynomial of a
Lagrangian line bundle $L$ parametrized by ${\M}$.  A construction of C.
Simpson enables us to compactify ${\M}$ as an open subset of a component
${\M}^{ss}$ of the moduli space of equivalence classes of coherent
semistable sheaves on $X$ with Hilbert polynomial $p$
\cite{simpson-moduli}.
Denote by ${\M}^s$ the open subset of ${\M}^{ss}$ parametrizing
isomorphism classes of stable sheaves, ${\M}^{s,sm}$ the smooth locus of
${\M}^s$. Then  $\M \subseteq {\M}^{s,sm} \subseteq {\M}^s \subseteq
{\M}^{ss}$.  In addition, the moduli space ${\M}^s$ embeds as a Zariski
open subset of the moduli space of simple sheaves
\cite{altman-kleiman-compactifying}.
The Zariski tangent space $T_{[L]}{\M}^s$ at a stable sheaf $L$ is thus
canonically isomorphic to  the Zariski tangent space of the moduli space
of simple sheaves. The latter is identified as
the group ${\rm Ext}^1_{{\cal O}_X}(L,L)$ of
extensions $0 \rightarrow L \rightarrow E \rightarrow L \rightarrow 0$ of
$L$ by $L$ as an ${\cal O}_X$-module.  When $X$ is a $K3$ or abelian
surface, Mukai's symplectic structure is given by the pairing
$$
{\rm Ext}^1_{{\cal O}_X}(L,L) \otimes {\rm Ext}^1_{{\cal O}_X}(L,L)
\stackrel{{\rm Yoneda}}{\longrightarrow} {\rm Ext}^2_{{\cal O}_X}(L,L)
\stackrel{{\rm S.D.}}{\longrightarrow} \ {\rm Hom}_X(L,L\otimes \omega_X)^*
\stackrel{id\otimes \sigma}{\longrightarrow} \C
$$
(Composition of the Yoneda pairing, Serre Duality, and evaluation at
$$id \otimes \sigma \in {\rm Hom}_X(L,L\otimes \omega_X)).$$

The generalization of Mukai's pairing requires the construction of a
homomorphism, depending linearly on the Poisson structure $\psi$,
\begin{equation}
\label{eq-homomorphism-lifting-the-polarization-to-a-two-extension-class}
y : H^{1,1}(X)  \rightarrow {\rm Ext}^2_{{\cal O}_X}(L,L).
\end{equation}
It sends the Kahler
class $\alpha := c_1 ({\cal O}_X(1)) \in H^{1,1}(X)$ to a
$2$-extension class $y(\alpha) \in {\rm Ext}^2_{{\cal O}_X}(L,L)$.
Once this is achieved, the $2$-form $\sigma_{\M}$ will become:
\begin{equation} \label{eq-generalized-mukai-pairing}
\begin{array}{l}
{\rm Ext}^1_{{\cal O}_X}(L,L) \otimes {\rm Ext}^1_{{\cal O}_X}(L,L)
\LongRightArrowOf{{\rm Yoneda}} {\rm Ext}^2_{{\cal O}_X}(L,L)
\LongRightArrowOf{y(\alpha)^{n-1}} {\rm Ext}^{2n}_{{\cal O}_X}(L,L)
\stackrel{S.D.}{\rightarrow}
\\
{\rm Hom}_X(L,L\otimes \omega_X)^*
\LongRightArrowOf{id\otimes \sigma^n} \C.
\end{array}
\end{equation}

\begin{rem} \label{rem-polarized-hodge-like-structure}
{\rm
When $L$ is a line bundle on a smooth Lagrangian subvariety $Z$
the construction involves a surprisingly rich
polarized Hodge-like structure on the algebra
$$
{\rm Ext}^*_{{\cal O}_X}(L,L) := \bigoplus^{2n}_{k=0} {\rm Ext}^k
_{{\cal O}_X}(L,L).
$$
Since
${\rm Ext}^k_{{\cal O}_X}(L,\cdot)$ is the right derived functor of the
composition $\Gamma \circ \; {\cal H}{\rm om}_{{\cal O}_X}(L,\cdot)$ of the
Sheaf Hom and the global sections functors, there
is a spectral sequence converging to ${\rm Ext}^k_{{\cal O}_X}(L,L)$ with
$$
E^{p,q}_2 = H^p\Big( Z, {\cal E}xt^q_{{\cal O}_X}(L,L) \Big)
$$
(see \cite{hilton-stammbach}).
The sheaf of $q$-extensions ${\cal E}xt^q_{{\cal O}_X}(L,L)$ is canonically
isomorphic to $\stackrel{q}{\wedge}  N_{Z/X}$ and thus, via the
symplectic structure, to $\Omega^q_Z$.  We obtain a canonical isomorphism
$E^{p,q}_2 \simeq H^{q,p}(Z)$ with the Dolbeault groups of $Z$.  Notice
however, that the Dolbeault groups appear in {\it reversed order} compared
to their order in the graded pieces of the Hodge filtration on the
cohomology ring $H^*(Z,\C)$.
}
\end{rem}

The construction of the $2$-extension class $y(\alpha)$ and
hence of the generalized
Mukai pairing (\ref{eq-generalized-mukai-pairing}) can be carried out
for all coherent sheaves parametrized by ${\M}^{s,sm}$.
We obtain:

\begin{theorem}
\label{thm-extension-of-symplectic-str}
\cite{markman-lagrangian-sheaves}
The symplectic structure $\sigma_{\M}$ on the relative Picard ${\M}$
extends to an algebraic $2$-form over the smooth locus ${\M}^{s,sm}$ of the
closure of ${\M}$ in the moduli space of stable sheaves on $X$.  It is
identified by the pairing (\ref{eq-generalized-mukai-pairing}).
\end{theorem}

The non-degeneracy of $\sigma_{\M}$ at a point $[L] \in {\M}$
parametrizing a line bundle on a smooth Lagrangian subvariety $Z$ follows
from the Hard Lefschetz theorem.  We expect $\sigma_{\M}$ to be non degenerate
everywhere on ${\M}^{s,sm}$.

Finally we remark that the pairing (\ref{eq-generalized-mukai-pairing}) can be
used to define a $2$-form on other components of the moduli space of stable
sheaves on $X$.  This $2$-form will, in general, be degenerate.
In fact, some components are odd dimensional:
\begin{example} \label{example-odd-dimensional-moduli-spaces}
{\rm
Consider an odd
dimensional complete linear system $|Z|$ whose
generic element is a smooth ample divisor on an abelian variety $X$ of even
dimension $\ge 4$, with a symplectic structure $\sigma$.  The dimension of
the component of the Hilbert scheme parameterizing deformations of $Z$ is
$\dim ({\rm Pic} \; X) + \dim|Z| = \dim X + \dim|Z|$.  Since  $h^{1,0}(Z) =
h^{1,0}(X)$, the component of the moduli space of sheaves parameterizing
deformations of the structure sheaf ${\cal O}_Z$, as an ${\cal O}_X$-module,
is of dimension $2 \cdot \dim X + \dim |Z|$ which is odd.
\EndProof
}
\end{example}
It is the
Hodge theoretic interpretation of the graded pieces of the spectral
sequence of ${\rm Ext}^k_{{\cal O}_X}(L,L)$ for {\em Lagrangian} line bundles
which assures the non degeneracy of $\sigma_{\M}$.

\subsection{Examples} \label{subsec-examples-of-lagrangian-sheaves}
\subsubsection{Higgs Pairs}
\label{moduli-higgs-pairs-as-lagrangian-sheaves}
In chapter \ref{ch9} we define the notion of a 1-form
valued Higgs pair $(E,\varphi)$ over a smooth n-dimensional projective
algebraic variety $X$.  It consists of a torsion free sheaf $E$ over $X$
and a homomorphism $\varphi : E \rightarrow E \otimes \Omega^1_X$
satisfying the symmetry condition $\varphi \wedge \varphi = 0$.

The moduli space Higgs$_X$ of semistable Higgs pairs of rank $r$ with
vanishing first and second Chern classes may be viewed
as the Dolbeault non-abelian first $GL_r(\C)$-cohomology group
of $X$ (cf. \cite{simpson-higgs-bundles-and-local-systems}
and theorem
\ref{thm-higgs-pairs-and-representations-of-pi1-for-curves}
when $X$ is a curve):

Non-abelian Hodge theory introduces a hyperkahler structure on the smooth
locus of the space ${\M}_{{\rm Betti}}$ of isomorphism classes of semisimple
$GL_r(\C)$-representations of the fundamental group $\pi_1(X)$ of $X$
\cite{deligne-twistors,hitchin,simpson-internetional-congress}.  The
hyperkahler structure consists of a Riemannian metric and an action of the
quaternion algebra $\h$ on the real tangent bundle with respect to which
\begin{description}
\item [(i)] the (purely imaginary)
unit vectors $\{a|a \bar{a} = 1\}$ in $\h$ correspond to a
(holomorphic) $\bP^1$-family of integrable complex structures,
\item [(ii)] the metric is Kahler with respect to these complex structures.
\end{description}
All but two of the complex structures are isomorphic to that of ${\M}_{{\rm
Betti}}$,
the two special ones are that  of $\HiggsModuli_X$ and its conjugate
(${\M}_{{\rm Betti}}$ and Higgs$_X$ are diffeomorphic).

The hyperkahler structure introduces a holomorphic symplectic structure
$\sigma$ on the smooth locus of Higgs$_X$.  In case $X$ is a Riemann
surface, that symplectic structure is the one giving rise to the Hitchin
integrable system of spectral Jacobians.

Our aim is to interpret the symplectic structure on Higgs$_X$ as an example
of a Lagrangian structure over the relative Picard of a Lagrangian
component of the Hilbert scheme of the cotangent bundle $T^*_X$ of $X$.
This interpretation will apply to the Hitchin system (where $\dim X = 1$).
For higher dimensional base varieties $X$ it will apply only to certain
particularly nice cases. See also \cite{biswas-a-remark} for a deformation
theoretic study of the holomorphic symplectic structure.

The spectral construction (proposition
\ref{prop-ordinary-spectral-construction-higgs-pairs})
can be carried out also for Higgs pairs over a higher dimensional
smooth projective variety $X$
(cf. \cite{simpson-moduli}).
We have a one to one correspondence between
\begin{description}
\item [(i)] (Stable) Higgs pairs $(E,\varphi)$ on $X$ (allowing $E$ to
be a rank $r$ torsion free sheaf) and
\item [(ii)] (Stable) sheaves $F$ on the cotangent bundle $T^*_X$ which are
supported on (pure) $n$-dimensional projective subschemes of $T^*_X$ which
are finite, degree $r$, branched coverings (in a scheme theoretic sense) of
$X$.
\end{description}

Projective subvarieties of $T^*_X$ which are finite over $X$ are called
{\it spectral coverings}.  Spectral coverings $\tilde{X}$ are necessarily
Lagrangian since the symplectic form $\sigma$ on $T^*_X$, which restricts
to a global exact 2-form on $\tilde{X}$, must vanish on $\tilde{X}$.

Let $B$ be the open subset of a component of the Hilbert scheme of
$\bP(T^*_X \oplus {\cal O}_X)$ parametrizing degree $r$ smooth
spectral coverings (closed subvarieties of $\bP(T^*_X \oplus {\cal O}_X)$
which are contained in $T^*_X$).  The above
correspondence embeds components of the relative Picard ${\cal M}
 \rightarrow B$
as open subsets of components of the moduli spaces of stable rank $r$ Higgs
pairs over $X$.  Theorem \ref{thm-symplectic-structure-on-relative-picard}
of this chapter implies
\begin{corollary}
\begin{description}
\item [(i)] The open subset ${\M}$ of the moduli Higgs$_X$  of Higgs pairs
over $X$ which, under the spectral construction, parametrizes line bundles
on smooth spectral covers, has a canonical symplectic structure $\sigma_\M$
(we do not require the Chern classes of the Higgs pairs to vanish).
\item [(ii)] The support morphism $h:{\cal M} \rightarrow B$ is a Lagrangian
fibration.
\end{description}
\end{corollary}

\begin{rem} \label{rem-bad-components}
{\rm
In general, when $\dim X > 1$, there could be components of the
moduli spaces of Higgs pairs for which the open set ${\M}$ above is empty,
i.e.,
\begin{enumerate}
\item the spectral coverings of all Higgs pairs in this
component are singular, or
\item the corresponding sheaves on the spectral coverings are
torsion free but not locally free.
\end{enumerate}
}
\end{rem}

\subsubsection{Fano Varieties of Lines on Cubic Fourfolds}
\label{subsec-fanos-of-cubics}

We will use theorem \ref{thm-symplectic-structure-on-relative-picard}
 to prove:
\begin{example} \label{fano-varieties-of-cubics}
Let $Y$ be a smooth cubic hypersurface in $\bP^5$.  The relative
intermediate Jacobian ${\cal J} \rightarrow B$ over the family
$B \subset | {\cal O}_{\bP^5} (1)|$ of smooth cubic hyperplane sections of
$Y$ is an algebraically completely integrable Hamiltonian system.
\end{example}

The statement follows from a description of the family
${\cal J} \rightarrow B$ as
an open subset of the moduli space of Lagrangian sheaves on the Fano
variety $X$ of lines on $Y$.

A. Beauville and R. Donagi proved \cite{beau-donagi} that $X$ is symplectic
(fourfold).
Clemens and Griffiths proved in \cite{c-g}
that the intermediate jacobian $J_b$ of
a smooth hyperplane section $Y \cap H_b$ is isomorphic to the Picard
$Pic^0 Z_b$ of the $2$-dimensional Fano variety $Z_b$ of lines on the
cubic $3$-fold $Y \cap H_b$.  C. Voisin observed that $Z_b$
is a Lagrangian subvariety of $X$ \cite{voisin}.
Since $h^{1,0}(Y \cap H_b) = 5$, $B$ is
isomorphic to a dense open subset of a component of the Hilbert scheme.
In fact, using results of Altman and Kleiman, one can show that the
corresponding component is isomorphic to $| {\cal O}_{\bP^5} (1)|$
(see \cite{altman-kleiman-fano} Theorem 3.3 (iv)).
Theorem \ref{thm-symplectic-structure-on-relative-picard} implies that
the relative Picard ${\M} \rightarrow B$ has a completely integrable
Hamiltonian system structure.

The symplectic structure $\sigma_\M$ is defined also at the fiber
of  the
relative Picard corresponding to a Fano variety $Z_b$ of lines on a
hyperplane section $Y\cap H_b$
with an ordinary double point $x_{b}\in Y\cap H_b$
(Theorem \ref{thm-extension-of-symplectic-str}). In that case, we have
a genus $4$ curve $C_b$ in $Z_b$ parametrizing lines through $x_{b}$. $Z_b$
is isomorphic to the quotient
$S^{2}C_b/(C_{1}\sim C_{2})$ of the second symmetric product of $C_b$ modulo
the identification of two disjoint copies of $C_b$ \cite{c-g}.
It is not difficult to check that $\sigma_{\M}$ is
{\em non-degenerate} also on the fiber $Pic^0(Z_b)$ of  the
relative Picard which is a ${\Bbb C}^{\times}$-extension of
the Jacobian of genus $4$.
The non-degeneracy of the symplectic structure implies that we get an
{\em induced boundary integrable system} on the relative Picard of
the family of genus $4$ curves
\begin{equation} \label{equation-integrable-system-of-genus-four-curves}
\RelPic({\cal C}) \rightarrow (Y^{*}-\Delta)
\end{equation}
over the complement of the singular locus $\Delta$ of the dual variety of
the cubic fourfold.

It is interesting to note that the boundary integrable system
(\ref{equation-integrable-system-of-genus-four-curves})
can not be realized as the relative Picard of a family of curves on a
symplectic surface. If this were the case, the generic rank of the pullback
$a^{*}(\sigma_{\RelPic^1({\cal C})})$
of the symplectic
structure from $\RelPic^1({\cal C}) \rightarrow (Y^{*}-\Delta)$ to
${\cal C} \rightarrow (Y^{*}-\Delta)$ via the Abel-Jacobi map would be $2$.
On the other hand, $a^{*}(\sigma_{\RelPic^1({\cal C})})$
is equal to the pullback of the
symplectic structure $\sigma_X$ on $X$ via the natural dominant map
${\cal C} \rightarrow X$
(corollary \ref{cor-canonical-symplectic-str-on-albanese}).
Thus, its generic rank is $4$.
The importance of this rank as an invariant of integrable systems supported by
families of Jacobians is illustrated in an interesting recent study
of J. Hurtubise \cite{hurtubise-local-geometry}.

More examples of nonrigid Lagrangian subvarieties can be found in
\cite{voisin,ye}.
\newpage





\section{Spectral covers} \label{ch9}
\subsection{Algebraic extensions} \label{intro}\
\indent  We have seen that Hitchin's system, the geodesic flow on an ellipsoid,
the polynomial matrices system of Chapter \ref{ch4} , the elliptic solitons,
and so on, all fit as special cases of the  spectral system on a curve. In this
final chapter, we consider some algebraic properties of the general spectral
system. We are still considering families of Higgs pairs
$ (E  \;,  \;  \varphi:E \longrightarrow E \otimes K) $,
but we  generalize in three separate directions:
\begin{enumerate}
\item The base curve $C$ is replaced by an arbitrary complex algebraic variety
$S$. The spectral curve $\widetilde{C}$ then becomes a spectral cover
$\widetilde{S}\longrightarrow S$.
\item The line bundle $K$ in which the endomorphism $\varphi$ takes its values
is replaced by a vector bundle, which we still denote by $K$. (this requires an
integrability condition on  $\varphi$.) Equivalently, $\widetilde{S}$ is now
contained in the total space $\Bbb{K}$ of a vector bundle over $S$.

\item Instead of the vector bundle $E$ we consider a principal $G$-bundle $\cal
G$, for an arbitrary complex reductive group $G$. The $G$-vector bundle $E$ is
then recovered as
$E := {\cal G} \times^{G} V$,
given the choice of a representation
$\rho :  G \longrightarrow Aut(V)$. The twisted endomorphism $\varphi$ is
replaced by a section of
$K \otimes \bdl{ad}(\cal G)$.
Even in the original case of $G = GL(n)$ one encounters interesting phenomena
in studying the dependence of  $\widetilde{S} := \widetilde{S}_V$, for a given
$ ({\cal G} , \varphi)$, on the representation $V$ of $G$.

\end{enumerate}

We will see that essentially all {\em algebraic} properties (but not the {\em
symplectic} structure) of  the Hitchin system, or of the (line-bundle valued,
$G=GL(n)$) spectral system on a curve, survive in this new context. In fact,
the added generality forces the discovery of some symmetries which were not
apparent in the original:
\begin{itemize}
\item Spectral curves are replaced by spectral covers. These come in several
flavors:
$\widetilde{S}_V,\widetilde{S}_{\lambda},\widetilde{S}_P$,
indexed by representations of $G$, weights, and parabolic subgroups.  The most
basic object is clearly the {\em cameral} cover $\widetilde{S}$; all the others
can be considered as associated objects.  In case $G=GL(n)$, the cameral cover
specializes not to our previous spectral cover, which has degree $n$ over $S$,
but roughly to its Galois closure, of degree $n!$ over $S$.
\item The spectral Picards, $Pic(\widetilde{S}_V)$ etc., can all be written
directly in terms of the decomposition of  $Pic(\widetilde{S})$ into Prym-type
components under the action of the Weyl group $W$. In particular, there is a
distinguished Prym component common to all the nontrivial
$Pic(\widetilde{S}_V)$. The identification of this component combines and
unifies many interesting constructions in Prym theory.
\item The Higgs bundle too can be relieved of its excess baggage. Stripping
away the representation $V$ as well as the values bundle $K$, one arrives (in
subsection \ref{abstract_objects}) at the notion of abstract, principal Higgs
bundle. The abelianization procedure assigns to this a spectral datum,
consisting of  a cameral cover with an equivariant bundle on it.
\item There is a Hitchin map (\ref{BigHitchin}) which is algebraically
completely integrable in the sense that its fibers can be naturally identified,
up to a "shift" and a "twist", with the distinguished Pryms (Theorem \ref{main}
).
\item The "shift" is a property of the group $G$, and is often nonzero even
when $\widetilde{S}$ is etale over $S$, cf. Proposition
\ref{reg.ss.equivalence}. The "twist", on the other hand, arises from the
ramification of $\widetilde{S}$ over $S$, cf. formula (\ref{G twist}).
\item The resulting abelianization procedure is local in the base $S$, and does
not require particularly nice behavior near the ramification, cf. example
\ref{nilpo}. It does require that $\varphi$ be  regular (this means that its
centralizer has the smallest possible dimension; for $GL(n)$, this means that
each eigenvalue may have arbitrary multiplicity, but the eigen-{\em space} must
be 1-dimensional), at least over the generic point of $S$. At present we can
only guess at the situation for irregular Higgs bundles.

\end{itemize}
The consideration of general spectral systems is motivated in part by work of
Hitchin \cite{hitchin-integrable-system} and Simpson \cite{simpson-moduli}. In
the remainder of this section we briefly recall those works.  Our exposition in
the following sections closely follows that of \cite{MSRI},  which in turn is
based  on \cite{D2}, for  the group-theoretic approach to spectral
decomposition used in section (\ref{deco}), and on \cite{D3} for the
Abelianization procedure, or equivalence of Higgs and spectral data, in
section(\ref{abelianization}). Some of these results, especially in the case of
a base curve, can also be found in \cite{AvM, BK, F, K, Me, MS, Sc}. \\

\noindent \underline{\bf Reductive groups} \nopagebreak

\noindent  Principal $G$-bundles {\cal G} for arbitrary reductive $G$ were
considered already in Hitchin's original paper
\cite{hitchin-integrable-system}.  Fix a curve $C$ and a line bundle $K$. There
is a moduli space ${\cal M}_{G,K}$ parametrizing equivalence classes of
semistable $K$-valued $G$-Higgs bundles, i.e. pairs
$({\cal G}, \varphi)$  with $\varphi \in K \otimes \bdl{ad}(\cal G)$. The
Hitchin map goes to $$B:=\oplus_{i} H^0(K^{\otimes d_i}),$$ where the $d_i$ are
the degrees of the $f_i$, a basis for the $G$-invariant polynomials on the Lie
algebra $\frak g$. It is:
\[ h: ({\cal G}, \varphi) \longrightarrow (f_i (\varphi))_{i}.
\]
When $K=\omega_C$, Hitchin showed \cite{hitchin-integrable-system} that one
still gets a completely integrable system, and that this system is
algebraically completely integrable for the classical groups $GL(n), SL(n),
SP(n), SO(n).$ The generic fibers are in each case  (not quite canonically; one
must choose various square roots! cf. sections \ref{reg.ss} and \ref{reg})
isomorphic to abelian varieties given in terms of the spectral curves
$\widetilde{C}$:

\begin{center}
\begin{equation}
\begin{array}{cl}                                      \label{Pryms for groups}
  GL(n)&    \widetilde{C}
                       \mbox{ has degree n over C, the AV is Jac(}
                       \widetilde{C}).    \\
  SL(n)&    \widetilde{C}
                      \mbox{ has degree n over C, the AV is Prym(}
                      \widetilde{C} / C).  \\
  SP(n)&    \widetilde{C}
                      \mbox{ has degree  2n over C and an involution }
                       x  \mapsto -x.  \\
             &    \mbox{ The map factors through the quotient }
                       \overline{C}.   \nonumber \\
             &    \mbox{ The AV is }
                       Prym( \widetilde{C} / \overline{C}).  \nonumber \\
  SO(n)&   \widetilde{C}  \mbox{ has degree  n and an involution , with: }  \\
             &  \bullet \mbox{ a fixed  component,  when n is odd.} \\
             &  \bullet \mbox{ some fixed double points, when n is even.} \\
             &   \mbox{ One must desingularize }
                      \widetilde{C}
                      \mbox{ and the quotient }
                      \overline{C}, \\
             &    \mbox{and ends up with the Prym  of the} \\
             &    \mbox{desingularized double cover.}  \
\end{array}
\end{equation}
\end{center}

For the exceptional group $G_2$, the algebraic complete integrability was
verified in \cite{KP1}.
A sketch of the argument for any reductive $G$ is in \cite{BK},  and a complete
proof was given in \cite{F}. We will outline a proof in section
\ref{abelianization} below.   \\

\noindent \underline{\bf Higher dimensions} \nopagebreak

\noindent  A sweeping extension of the notion of Higgs bundle is suggested by
the work of Simpson  \cite{simpson-moduli}, which was already discussed in
Chapter \ref{ch8}.  To him, a Higgs bundle on a projective variety S is a
vector bundle (or principal $G$-bundle \ldots) $E$  with a  {\em symmetric},
$\Omega^1_S$-valued  endomorphism
\[ \varphi : E \longrightarrow E \otimes \Omega^1_S.
\]
Here {\em symmetric} means the vanishing of:
\[ \varphi\wedge\varphi : E \longrightarrow E \otimes \Omega^2_S,
\]
a condition which is obviously vacuous on curves. Simpson  constructs a moduli
space for such Higgs bundles (satisfying appropriate stability conditions), and
establishes diffeomorphisms to corresponding moduli spaces  of connections and
of representations of $\pi_1(S)$ .

In our approach, the $\Omega^1$-valued Higgs bundle will be considered as a
particular realization of an abstract Higgs bundle, given by a subalgebra of
$ad(\cal{G})$. The symmetry condition will be expressed in the definition
\ref{princHiggs}  of an abstract Higgs bundle by requiring the abelian
subalgebras of  $ad(\cal{G})$ to be abelian.

\subsection {Decomposition of spectral Picards}              \label{deco}
\subsubsection{The question}\
\indent Throughout this section we fix a vector bundle $K$ on a complex variety
$S$, and a pair
 $({\cal G},\varphi)$ where  ${\cal G}$ is a principal $G$-bundle on $S$ and
$\varphi$  is a regular section of  $K \otimes \bdl{ad}(\cal G)$. (This data is
equivalent to the regular case of what we call in section
\ref{abstract_objects} a $K$-valued principal Higgs bundle.) Each
representation
\[   \rho :  G \longrightarrow Aut(V)
\]
determines an associated $K$-valued Higgs (vector) bundle
\[   ( {\cal V} := {\cal G} \times^{G} V, \qquad{\rho}(\varphi)\ ),
\]
which in turn determines a spectral cover $\widetilde{S}_V \longrightarrow S$.

The question, raised first in \cite{AvM} when $S={\bf P}^1$, is to relate the
Picard varieties of the
$\widetilde{S}_V$ as $V$ varies, and in particular to find pieces common to all
of them. For Adler and van Moerbeke, the motivation was that  many evolution
DEs (of Lax type) can be {\em linearized} on the Jacobians of spectral curves.
This means that the "Liouville tori", which live naturally in the complexified
domain of the DE (and hence are independent of the representation $V$) are
mapped isogenously to their image in $\mbox{Pic}(\widetilde{S}_V)$ for each
nontrivial $V$ ; so one should be able to locate these tori among the pieces
which occur in an isogeny decomposition of each of the
$\mbox{Pic}(\widetilde{S}_V)$. There are many specific examples where a pair of
 abelian varieties constructed from related covers of curves are known to be
isomorphic or isogenous, and some of these lead to important identities among
theta functions.

\begin{eg}
\begin{em}
Take $G=SL(4)$ . The standard representation $V$ gives a branched cover
$\widetilde{S}_V \longrightarrow S$ of degree 4. On the other hand, the
6-dimensional representation $\wedge ^2 V$ (=the standard representation of the
isogenous group $SO(6)$)  gives a cover
$ \stackrel{\approx}{S} \longrightarrow S$  of degree 6, which factors through
an involution:
\[  \stackrel{\approx}{S} \longrightarrow \overline{S} \longrightarrow S.
\]
One has the isogeny decompositions:
\[ Pic \, (\widetilde{S}) \sim Prym(\widetilde{S} / S) \oplus Pic \,(S)
\]
\[ Pic \,(\stackrel{\approx}{S}) \sim
Prym(\stackrel{\approx}{S} / \overline{S})        \oplus
Prym(\overline{S} / S)              \oplus                     Pic \,(S).
\]
It turns out that
\[  Prym(\widetilde{S} / S)            \sim        Prym(\stackrel{\approx}{S} /
\overline{S}) .
\]
For  $S={\bf P}^1$, this is Recillas' {\em trigonal construction} \cite{R}. It
says that every Jacobian of a trigonal curve is the Prym of a double cover of a
tetragonal curve, and vice versa.
\end{em}
\end{eg}

\begin{eg}
\begin{em}
Take $G=SO(8)$ with its standard 8-dimensional representation $V$.  The
spectral cover has degree 8 and factors through an involution,
 $ \stackrel{\approx}{S} \longrightarrow \overline{S} \longrightarrow S.$
The two half-spin representations $V_1, V_2$ yield similar covers
\[ \stackrel{\approx}{S} _i \longrightarrow \overline{S} _i \longrightarrow S,
\qquad i=1,2.
\]
The {\em tetragonal construction} \cite{D1} says that the three Pryms of the
double covers are isomorphic. (These examples, as well as Pantazis' {\em
bigonal construction} and constructions based on some exceptional groups, are
discussed in the context of spectral covers in \cite{K} and \cite{D2}.)
\end{em}
\end{eg}

It turns out in general that there is indeed a distinguished, Prym-like isogeny
component common to all the spectral Picards, on which the solutions to
Lax-type DEs evolve linearly.  This was noticed in some cases already in
\cite{AvM}, and was greatly extended by Kanev's construction of Prym-Tyurin
varieties. (He still needs $S$ to be ${\bf P}^1$ and the spectral cover to have
generic ramification; some of his results apply only to {\em minuscule
representations}.)
Various parts of the general story have been worked out recently by a number of
authors, based on either of two approaches: one, pursued in  \cite{D2,Me,MS},
is to decompose everything according to the action of the Weyl group $W$ and to
look for common pieces; the other, used in \cite{BK,D3,F,Sc}, relies on the
correspondence of spectral data and Higgs bundles . The group-theoretic
approach is described in the rest of this section. We take up the second
method, known as {\em abelianization}, in section~\ref{abelianization}.

\subsubsection{Decomposition of spectral covers}
\label{decomp covers}
\indent  The decomposition of spectral Picards arises from three sources.
First, the spectral cover for a sum of representations is the union of the
individual covers $\widetilde{S}_V$.  Next, the cover  $\widetilde{S}_V$ for an
irreducible representation is still the union of subcovers
$\widetilde{S}_{\lambda}$  indexed by weight orbits.  And finally, the Picard
of  $\widetilde{S}_{\lambda}$ decomposes into Pryms.
We start with a few observations about the dependence of the covers themselves
on the representation.  The decomposition of the Picards is taken up in the
next subsection. \\

\noindent \underline{\bf Spectral covers} \nopagebreak

\noindent Whenever a representation space $V$ of $G$ decomposes,
 $$V = \oplus V_i,$$
there is a corresponding decomposition
$$ \widetilde{S}_V =\cup \widetilde{S}_{V_i}, $$
so we may restrict attention to irreducible representations $V$.
There is an \mbox{\em infinite} collection (of  irreducible representations $V
:= V_{\mu}$, hence) of  spectral covers $\widetilde{S}_V$, which can be
parametrized by their highest weights $\mu$  in the dominant  Weyl chamber
$\overline{C}$ ,  or equivalently by the $W$-orbit of extremal weights, in
$\Lambda  / W$.  Here $T$ is a maximal torus in $G$, $\Lambda := Hom(T, {\bf
C}^*)$ is the {\em weight lattice } (also called {\em character lattice })  for
$G$, and $W$ is the Weyl group.  Now $V_{\mu}$ decomposes as the sum of its
weight subspaces $V_{\mu}^{\lambda}$, indexed by certain weights $\lambda$ in
the convex hull in $\Lambda$ of the $W$-orbit of $\mu$. We conclude that each
$\widetilde{S}_{V_{\mu}}$ itself decomposes  as the union of its subcovers
$\widetilde{S}_{\lambda}$, each of which involves eigenvalues in a given
$W$-orbit   $W{\lambda}$ . ($\lambda$ runs over the weight-orbits in
$V_{\mu}$.) \\

\noindent \underline{\bf Parabolic covers} \nopagebreak

\noindent There is a {\em finite} collection of covers $\widetilde{S}_P$,
parametrized by the conjugacy classes in $G$ of parabolic subgroups (or
equivalently by arbitrary dimensional faces $F_P$ of the chamber
$\overline{C}$) such that (for general $S$) each  eigenvalue cover
$\widetilde{S}_{\lambda}$ is birational to some  parabolic cover
$\widetilde{S}_{P}$, the one whose open face $F_P$ contains ${\lambda}$.  \\

\noindent \underline{\bf The cameral cover} \nopagebreak

\noindent There is a $W$-Galois cover   $\widetilde{S} \longrightarrow S$ such
that each
 $\widetilde{S}_{P}$ is isomorphic to  $\widetilde{S} / W_P$, where $W_P$ is
the Weyl subgroup of $W$ which stabilizes $F_P$. We call $\widetilde{S}$ the
{\em cameral cover} ,  since, at least generically,  it parametrizes the
chambers determined by $\varphi$ (in the duals of the Cartans).
Informally, we think of  $\widetilde{S} \longrightarrow S$ as the cover which
associates to a point $s \in S$ the set of Borel subalgebras of  $ad({\cal
G})_s$ containing $\phi(s)$. More carefully, this is constructed as follows:
There is a morphism
${\frak g}\longrightarrow {\frak t}/W$ sending $g \in {\frak g}$ to the
conjugacy class of its semisimple part $g_{ss}$. (More precisely, this is
$Spec$ of the composed ring homomorphism
${\bf C} [ {\frak t} ] ^{W}
{ \stackrel{\simeq}{\leftarrow}}
 {\bf C}[{\frak g}]^{G}       \label{t/W}
\hookrightarrow
{\bf C}[{\frak g}]$.)
Taking fiber product with the quotient map ${\frak t}\longrightarrow {\frak
t}/W$, we get the cameral cover ${\tilde{\frak g}}$ of  ${\frak g}$. The
cameral cover $\widetilde{S} \longrightarrow S$ of a $K$-valued principal Higgs
bundle on  $S$  is glued from covers of open subsets in $S$ (on which $K$ and
$\cal G$ are trivialized) which in turn are pullbacks by $\varphi$ of
${\tilde{\frak g}}  \longrightarrow  {\frak g} $.

\subsubsection{Decomposition of spectral Picards}\
\indent The decomposition of  the Picard varieties of spectral covers can be
described as follows:\\

\noindent \underline{\bf The cameral Picard} \nopagebreak

\noindent  From each isomorphism class of irreducible $W$-representations,
choose an integral representative $\Lambda _i$. (This can always be done, for
Weyl groups.)  The group ring
${\bf Z} [W]$  which acts on $Pic(\widetilde{S}) $ has an isogeny
decomposition:

\begin{equation}\label{regular rep}
{\bf Z} [W] \sim \oplus _i \Lambda _i \otimes_{\bf Z} \Lambda _i^{*},
\end{equation}

\noindent
which is just the decomposition of the regular representation. There is a
corresponding isotypic decomposition:

\begin{equation}\label{cameral Pic decomposition}
Pic(\widetilde{S}) \sim \oplus _i \Lambda _i \otimes_{\bf Z} Prym_{\Lambda
_i}(\widetilde{S}),
\end{equation}

\noindent
where

\begin{equation}\label{def of Prym_lambda}
Prym_{\Lambda _i}(\widetilde{S} ):= Hom_W (\Lambda _i , Pic(\widetilde{S})).
\end{equation}\\

\noindent \underline{\bf Parabolic Picards} \nopagebreak

\noindent There are at least three reasonable ways of obtaining an isogeny
decomposition of $Pic(\widetilde{S}_P) $, for a parabolic subgroup $P \subset
G$:
\begin{itemize}

\item The `Hecke' ring $Corr_P$ of correspondences on $\widetilde{S}_P$ over
$S$ acts on  $Pic(\widetilde{S}_P) $, so every irreducible integral
representation $M$ of $Corr_P$ determines a generalized Prym
$$   Hom_{Corr_P} (M, Pic(\widetilde{S}_P)),   $$
and we obtain an isotypic decomposition of $Pic(\widetilde{S}_P)$ as before.

\item $Pic(\widetilde{S}_P)$ maps, with torsion kernel, to
$Pic(\widetilde{S})$, so we obtain a decomposition of the former by
intersecting its image with the isotypic components
$\Lambda _i \otimes_{\bf Z} Prym_{\Lambda _i}(\widetilde{S})$ of the latter.

\item Since $\widetilde{S}_P$ is the cover of $S$ {\em associated} to the
$W$-cover $\widetilde{S}$ via the permutation representation ${\bf Z} [W_P
\backslash W]$ of $W$, we get an isogeny decomposition of
$Pic(\widetilde{S}_P)$ indexed by the irreducible representations in
${\bf Z} [W_P \backslash W]$.

\end{itemize}

It turns out (\cite{D2},section 6) that all three decompositions agree
and can be given explicitly as

\begin{equation}
\label{multiplicity spaces}
\oplus _i M _i \otimes Prym_{\Lambda _i}(\widetilde{S}) \subset
\oplus _i \Lambda _i \otimes Prym_{\Lambda _i}(\widetilde{S}),\qquad
  M_i := (\Lambda_i)^{W_P}.
\end{equation}

\noindent \underline{\bf Spectral Picards} \nopagebreak

\noindent  To obtain the decomposition of the Picards of the original  covers
$\widetilde{S}_V$ or
 $\widetilde{S}_{\lambda}$, we need, in addition to the decomposition of
$Pic(\widetilde{S}_P)$, some information on the singularities. These can arise
from two separate sources:
\begin{description}
\item[Accidental singularities of the $\widetilde{S}_{\lambda}$. ]
For a sufficiently general Higgs bundle, and for a weight $\lambda$ in the
interior of the face $F_P$ of the Weyl chamber  $\overline{C}$, the natural
map:
$$    i_{\lambda}: \widetilde{S}_P\longrightarrow  \widetilde{S}_{\lambda}   $$
is birational. For the {\em standard} representations of the classical groups
of types
 $A_n, B_n$ or $C_n$, this {\em is} an isomorphism. But for general ${\lambda}$
it is {\em not}: In order for $i_{\lambda}$ to be an isomorphism, ${\lambda}$
must be a multiple of a fundamental weight, cf. \cite{D2}, lemma 4.2. In fact,
the list of fundamental weights for which this happens is quite short; for the
classical groups we have only:  $\omega_1$ for  $A_n, B_n$ and $C_n$,
$\omega_n$  (the dual representation) for $A_n$, and $\omega_2$ for $B_2$. Note
that for $D_n$ the list is {\em empty}. In particular, the covers produced by
the standard representation of $SO(2n)$ are singular; this fact, noticed by
Hitchin In \cite{hitchin-integrable-system},  explains the need for
desingularization in his result~(\ref{Pryms for groups}).

\item[Gluing the $\widetilde{S}_{V}$. ]
In addition to the singularities of each $i_{\lambda}$, there are the
singularities created by the gluing map $\amalg_{\lambda}
\widetilde{S}_{\lambda} \longrightarrow  \widetilde{S}_V$. This makes explicit
formulas somewhat simpler in the case, studied by Kanev \cite{K}, of {\em
minuscule} representations, i.e. representations whose weights form a single
$W$-orbit.  These singularities account, for instance,  for the
desingularization required in the $SO(2n+1)$ case in
(\ref{Pryms for groups}).

\end{description}

\subsubsection{The distinguished Prym}   \label{distinguished}\
\indent Combining much of the above, the Adler--van Moerbeke problem of finding
a component common to  the $Pic(\widetilde{S}_V)$ for all non-trivial $V$
translates into: \\

\begin{em}
Find the non trivial irreducible representations
$\Lambda_i $  of $W$   which occur in  ${\bf Z} [W_P \backslash W] $
with positive multiplicity for all proper Weyl subgroups
$W_P \subsetneqq W.$
\end{em} \\

It is easy to see that for arbitrary finite groups $W$, or even for Weyl groups
$W$ if we allow arbitrary rather than Weyl subgroups $W_P$, there may be no
common factors \cite{D2}.  For example, when $W$ is the symmetric group $S_3$
(=the Weyl group of $GL(3)$) and $W_P$ is  $S_2$ or $A_3$, the representations
${\bf Z} [W_P \backslash W] $  are 3 or 2 dimensional, respectively, and have
only the trivial representation as common component. In any case, our problem
is equivalent (by Frobenius reciprocity or (\ref{multiplicity spaces})) to \\

\begin{em}
Find the irreducible representations
$\Lambda_i $
  of  W such that  for every  proper Weyl subgroup
$W_P \subsetneqq W, $
 the space of invariants
$M_i := (\Lambda_i)^{W_P} $
 is non-zero.
\end{em} \\

One solution is now obvious: the {\em{reflection representation}} of $W$ acting
on the weight lattice $\Lambda$ has this property. In fact,
$\Lambda^{W_P}$ in this case is just the face $F_P$ of $\overline{C}$. The
corresponding component  $Prym_{\Lambda }(\widetilde{S})$ , is called {\em{the
distinguished Prym}.} We will see in section \ref{abelianization} that its
points correspond, modulo some corrections,  to Higgs bundles.

For the classical groups, this turns out to be the only common component. For
$G_2$ and $E_6$ it turns out (\cite{D2}, section 6) that a second common
component exists. The geometric significance of points in these extra
components is not known. As far as we know, the only component other than the
distinguished Prym which has arisen `in nature' is the one associated to the
1-dimensional sign representation of $W$, cf.  \cite{KP2}.

\subsection {Abelianization}\label{abelianization}
\subsubsection{Abstract vs. $K$-valued objects}\label{abstract_objects}\
\indent We want to describe the abelianization procedure in a somewhat abstract
setting, as an equivalence between {\em{principal Higgs bundles}} and certain
{\em spectral data}.
Once we fix a {\em{values}} vector bundle $K$, we obtain an equivalence between
{\em $K$-valued principal Higgs bundles} and {\em K-valued spectral data}.
Similarly,
the choice of a representation $V$ of $G$ will determine an equivalence of
{\em $K$-valued Higgs bundles} (of a given representation type) with $K$-valued
spectral data.

As our model of a $W$-cover we take the natural quotient map
$$G/T \longrightarrow G/N
$$
and its partial compactification

\begin{equation}
\overline{G/T} \longrightarrow \overline{G/N}.      \label{partial
compactification}
\end{equation}

Here $T \subset G$ is a maximal torus, and $N$ is its normalizer in $G$.
The quotient $G/N$ parametrizes maximal tori (=Cartan subalgebras) $\frak{t}$
in $\frak{g}$,
while $G/T$ parametrizes pairs ${\frak t \subset \frak b}$
with ${\frak b \subset \frak g}$ a Borel subalgebra.
An element $x \in {\frak g}$ is {\em regular} if the dimension of its
centralizer
${\frak c \subset \frak g}$ equals $\dim{T}$ (=the rank of $\frak{g}$). The
partial compactifications
$ \overline{G/N}$ and $ \overline{G/T}$  parametrize regular centralizers
${\frak c }$ and pairs ${\frak c \subset \frak b}$, respectively.

In constructing the cameral cover in section  \ref{t/W}, we used the $W$-cover
$\frak t \longrightarrow \frak t / W$ and its  pullback cover ${
\widetilde{\frak g} \longrightarrow \frak g}$.
Over the open subset $\frak g_{reg}$ of regular elements, the same cover is
obtained by pulling back (\ref{partial compactification}) via the map
$\alpha : \frak g_{reg}  \longrightarrow \overline{G/N}$ sending an element to
its centralizer:

\begin{equation}
                       \label{commutes}
\begin{array}{lccccc}
\frak t & \longleftarrow & \widetilde{\frak g}_{{reg}} & \longrightarrow &
\overline{G/T}  & \\
\downarrow & &\downarrow & & \downarrow  & \\
\frak t  /W & \longleftarrow & {\frak g}_{{reg}} &
\stackrel{\alpha}{\longrightarrow} & \overline{G/N} &.
\end{array}
\end{equation}

When working with $K$-valued objects, it is usually more convenient to work
with the left hand side of  (\ref{commutes}), i.e. with eigen{\em values}. When
working with the abstract objects, this is unavailable,  so we are forced to
work with the eigen{\em vectors},
or the right hand side of  (\ref{commutes}). Thus:

\begin{defn}
An abstract {\em cameral cover} of $S$ is a finite morphism $\widetilde{S}
\longrightarrow S$
with $W$-action, which locally (etale) in $S$ is a pullback of (\ref{partial
compactification}). \\
\end{defn}

\begin{defn}
A {\em $K$-valued cameral cover}  ($K$ is a vector bundle on $S$) consists of a
cameral cover  $\pi :  \widetilde{S} \longrightarrow S$  together with an
$S$-morphism
\begin{equation}
   \widetilde{S} \times \Lambda  \longrightarrow \Bbb{K}   \label{K-values}
\end{equation}
which is $W$-invariant ($W$ acts on $ \widetilde{S} , \Lambda,$ hence
diagonally on
$\widetilde{S} \times \Lambda $ ) and linear in $\Lambda$. \\
\end{defn}

We note that a cameral cover $\widetilde{S}$ determines quotients
$\widetilde{S}_P$ for parabolic subgroups $P \subset G$. A $K$-valued cameral
cover determines additionally the $\widetilde{S}_{\lambda}$ for $\lambda \in
\Lambda$, as images in $\Bbb{K}$ of
$\widetilde{S} \times \{ \lambda \}$.  The data of (\ref{K-values}) is
equivalent to a $W$-equivariant map $\widetilde{S} \longrightarrow
\frak{t}\otimes_{\bf C} K.$

\begin{defn}                    \label{princHiggs}
A $G$-principal Higgs bundle on $S$ is a pair ($\cal{G}, \bdl{c})$ with
$\cal{G}$ a principal $G$-bundle and  $\bdl{c} \subset ad(\cal{G})$ a subbundle
of regular centralizers.
\\
\end{defn}

\begin{defn}
A $K$-valued $G$-principal Higgs bundle consists of  $( \cal{G}, \bdl{c} )$
as above together with a section $\varphi$ of $\bdl{c} \otimes K$.
\end{defn}

A principal Higgs bundle $(\cal{G}, \bdl{c})$  determines a cameral cover
$\widetilde{S}\longrightarrow S$ and a homomorphism $\Lambda \longrightarrow
\mbox{Pic}(\widetilde{S}).$ Let $F$ be a parameter space for Higgs bundles with
a given $\widetilde{S}$. Each non-zero $\lambda \in \Lambda$ gives a
non-trivial map
$F\longrightarrow \mbox{Pic}(\widetilde{S})$. For $\lambda$ in a face $F_P$ of
$\overline{C}$, this factors through $\mbox{Pic}(\widetilde{S}_P)$. The
discussion in section \ref{distinguished} now suggests that $F$ itself should
be given roughly by the distinguished Prym,
$$ Hom_W (\Lambda , \mbox{Pic}(\widetilde{S})).
$$
It turns out that this guess needs two corrections. The first correction
involves restricting to a coset of a subgroup; the need for this is visible
even in the simplest case where
$\widetilde{S}$
 is etale over
$S$,
so
$(\cal{G}, \bdl{c})$
is everywhere regular and semisimple
(i.e.
$ \bdl{c}$
 is a bundle of Cartans.)
The second correction involves a twist along the ramification of
$\widetilde{S}$
over
$S$.
We explain these in the next two subsections.

\subsubsection{The regular semisimple case: the shift}       \label{reg.ss}

\begin{eg}   \label{unramified}
\begin{em}
Fix a smooth projective curve
$C$
and a line bundle
$K \in \mbox{Pic}(C)$
such that
$K^{\otimes 2} \approx  \cal{O}_C.$
This determines an etale double cover
$\pi : \widetilde{C} \longrightarrow C$
with involution
$i$,
and homomorphisms

\begin{center}
$\begin{array}{cccccc}
 \pi^{*}          &:&  \mbox{Pic}(C)                     &\longrightarrow
&\mbox{Pic}(\widetilde{C}) &, \\
 \mbox{Nm} &:&  \mbox{Pic}(\widetilde{C}) &\longrightarrow &\mbox{Pic}(C)
              &, \\
   i^{*}            &:& \mbox{Pic}(\widetilde{C}) &\longrightarrow
&\mbox{Pic}(\widetilde{C}) &,
\end{array}$
\end{center}

satisfying
$$   1+i^{*} = \pi^* \circ \mbox{Nm}.
$$

\begin{itemize}
\item For
$G = GL(2)$
we have
$\Lambda = \bf{Z} \oplus \bf{Z}$,
and
$W = {\cal{S}}_{2}$
permutes the summands, so
$$    Hom_W (\Lambda , \mbox{Pic}(\widetilde{C}))  \approx
\mbox{Pic}(\widetilde{C}).
$$
And indeed, the Higgs bundles corresponding to
$\widetilde{C}$
are parametrized by
$\mbox{Pic}(\widetilde{C})$:
send
$L \in \mbox{Pic}(\widetilde{C})$
to
$(\cal{G}, \bdl{c})$,
where
$\cal{G}$
has associated rank-2 vector bundle
${\cal V} := \pi_* L$,
and
$ \bdl{c} \subset \UnderlinedEnd{{\cal{V}}}$
is
$\pi_* {\cal O}_{\widetilde{C}}.$
\item On the other hand, for
$G=SL(2)$
we have
$\Lambda=\bf{Z}$
and
$W={\cal{S}}_2$
acts by
$\pm 1$,
so
$$   Hom_W (\Lambda , \mbox{Pic}(\widetilde{S}))  \approx
\{L \in \mbox{Pic}(\widetilde{C})\  | \  i^*L \approx L^{-1} \}
= \mbox{ker}(1+i^*).
$$
This group has 4 connected components. The subgroup
$\mbox{ker(Nm)}$
consists of 2 of these. The connected component of 0 is the classical Prym
variety, cf. \cite{MuPrym}. Now the Higgs bundles correspond, via the above
bijection
$L\mapsto \pi_*L$,
to
$$\{L \in \mbox{Pic}(\widetilde{C}) \ |\ \det (\pi_*L) \approx {\cal O}_C \} =
{\mbox{Nm}}^{-1}(K).
$$
Thus they form the {\em non-zero} coset of the subgroup
$\mbox{ker(Nm)}$.
(If we return to a higher dimensional  $S$, there is no change in the  $GL(2)$
story,
 but it is possible for $K$ not to be in the image of
$\mbox{Nm}$,
so there may be {\em no}
$SL(2)$-Higgs bundles corresponding to such a cover.)
\end{itemize}
\end{em}
\end{eg}

This example generalizes to all
$G$,
as follows. The equivalence classes of extensions
$$1 \longrightarrow  T \longrightarrow  N'  \longrightarrow  W \longrightarrow
1
$$
(in which the action of $W$ on $T$ is the standard one) are parametrized by the
group cohomology
$H^2(W,T)$.
Here the 0 element corresponds to the semidirect product . The class
$[N] \in H^2(W,T)$
of the normalizer $N$ of $T$ in $G$ may be 0, as it is for
$G=GL(n)  ,  {\bf P}GL(n)  ,  SL(2n+1)  $;
or not, as for
$G=SL(2n)$.

Assume first, for simplicity, that
$S,\widetilde{S}$
are connected and projective. There is then a natural group homomorphism
\begin{equation}
                                   \label{c}
c: Hom_W (\Lambda , \mbox{Pic}(\widetilde{S}))\longrightarrow H^2(W,T).
\end{equation}
Algebraically, this is an edge homomorphism for the Grothendieck spectral
sequence of equivariant cohomology, which gives the exact sequence

\begin{equation}
                 \label{c-edge}\qquad
0 \longrightarrow                            H^1(W,T)
\longrightarrow                                H^1(S,{\cal{C}})
\longrightarrow                                Hom_W (\Lambda ,
\mbox{Pic}(\widetilde{S}))
\stackrel{c}{\longrightarrow}         H^2(W,T).
\end{equation}
where
 ${\cal{C}} := \widetilde{S} \times _W T.$
Geometrically, this expresses a {\em Mumford group} construction: giving
${\cal{L}} \in \mbox{Hom}(\Lambda,\mbox{Pic}(\widetilde{S}))$
is equivalent to giving a principal $T$-bundle
$\cal T$
over
$\widetilde{S}$;
for
${\cal{L}} \in \mbox{Hom}_W(\Lambda,\mbox{Pic}(\widetilde{S}))$,
$c({\cal{L}})$
is the class in
$H^2(W,T)$
of the group
$N'$
of automorphisms of
$\cal T$
which commute with the action on
$\widetilde{S}$
of some
$w \in W$.

To remove the restriction on
$S, \widetilde{S}$,
we need to replace each occurrence of $T$ in (\ref{c},\ref{c-edge}) by
$\Gamma (\widetilde{S}, T)$,
the global sections of the trivial bundle on
$\widetilde{S}$
with fiber $T$. The natural map
$H^2(W,T) \longrightarrow H^2(W,\Gamma (\widetilde{S}, T))$
allows us to think of
$[N]$
as an element of
$H^2(W,\Gamma (\widetilde{S}, T))$.

\begin{prop} \cite{D3}                \label{reg.ss.equivalence}
Fix an etale $W$-cover
$\pi: \widetilde{S}\longrightarrow S$.
The following data are equivalent:

\begin{enumerate}
\item Principal $G$-Higgs bundles
$(\cal{G}, \bdl{c})$
with cameral cover
$\widetilde{S}$.
\item Principal $N$-bundles
$\cal N$
over $S$ whose quotient by $T$ is
$\widetilde{S}.$
\item $W$-equivariant homomorphisms
${\cal{L}} : \Lambda \longrightarrow \mbox{Pic}(\widetilde{S})$
with
$c({\cal L}) = [N] \in H^2(W,\Gamma (\widetilde{S}, T))$.
\end{enumerate}

\end{prop}

We observe that while the shifted objects correspond to Higgs bundles,
the unshifted objects
$$
{\cal{L}} \in \mbox{Hom}_W(\Lambda,\mbox{Pic}(\widetilde{S})),  \qquad  c({\cal
L})=0
$$

\noindent
come from the $\cal C$-torsers in $H^1(S, {\cal C} ).$

\subsubsection{The regular case: the twist along the ramification}
\label{reg}

\begin{eg}                               \label{ramified}
\begin{em}
Modify example \ref{unramified} by letting
$K \in \mbox{Pic}(C) $
be arbitrary, and choose a section $b$ of
$K ^{\otimes 2}$
which vanishes on a simple divisor
$B \subset C$.
We get a double cover
$\pi : \widetilde{C} \longrightarrow C$
branched along $B$, ramified along a divisor
$$
R \subset \widetilde{C}, \quad  \pi(R)=B.
$$
Via
$L\mapsto \pi_*L$,
the $SL(2)$-Higgs bundles still correspond to
$$\{L \in \mbox{Pic}(\widetilde{C}) \ |\ \det (\pi_*L) \approx {\cal O}_C \} =
{\mbox{Nm}}^{-1}(K).
$$
But this is no longer in
$   Hom_W (\Lambda , \mbox{Pic}(\widetilde{S}))$;
rather, the line bundles in question satisfy

\begin{equation}
                   \label{SL(2) twist}
i^*L \approx L^{-1}(R).
\end{equation}

\end{em}
\end{eg}

For arbitrary $G$, let
$\Phi$
denote the root system and
$\Phi^+$
the set of positive roots. There is a decomposition
$$   \overline{G/T} \  \smallsetminus \  G/T   = \bigcup _{\alpha \in
\Phi^+}R_{\alpha}
$$
of the boundary into components, with
$R_{\alpha}$
the fixed locus of the reflection
$\sigma_{\alpha}$
in
$\alpha$.
(Via (\ref{commutes}),  these correspond to the complexified walls in
$\frak t$.)
Thus each cameral cover
$\widetilde{S} \longrightarrow S$
comes with a natural set of (Cartier)  {\em ramification divisors}, which we
still denote
$R_{\alpha}, \quad  \alpha \in \Phi^+.$

For
$w \in W$,
set
$$  F_w := \left\{ \alpha \in \Phi^+ \ | \ w^{-1} \alpha \in \Phi^- \right\}
= \Phi^+ \cap w \Phi^-,
$$
and choose a $W$-invariant form
$\langle , \rangle$
  on
$\Lambda$.
We consider the variety
$$   Hom_{W,R} (\Lambda , \mbox{Pic}(\widetilde{S}))
$$
of $R$-twisted $W$-equivariant homomorphisms, i.e. homomorphisms
$\cal L$
satisfying

\begin{equation}      \qquad
                                \label{G twist}
w^*{\cal L}(\lambda) \approx
 {\cal L}(w\lambda)\left( \sum_{\alpha \in F_w}{
{\langle-2\alpha,w\lambda \rangle \over \langle \alpha ,\alpha \rangle}
R_{\alpha}
} \right) , \qquad  \lambda \in \Lambda, \quad w \in W.
\end{equation}

This turns out to be the correct analogue of  (\ref{SL(2) twist}). (E.g. for a
reflection
$w=\sigma_{\alpha}$,
\quad $F_w$
is
$\left\{ \alpha \right\}$,
so this gives
$ w^*{\cal L}(\lambda) \approx
 {\cal L}(w\lambda)\left(
{{\langle\alpha,2\lambda \rangle \over \langle \alpha,\alpha \rangle}
R_{\alpha}}
 \right),$
which specializes to (\ref{SL(2) twist}).)  As before, there is a class map

\begin{equation}
                                   \label{c,R}
c: Hom_{W,R} (\Lambda , \mbox{Pic}(\widetilde{S}))\longrightarrow
H^2(W,\  \Gamma (\widetilde{S}, T))
\end{equation}

\noindent
which can be described via a Mumford-group construction.

To understand this twist, consider the formal object

\begin{center}
$\begin{array}{cccc}
{1 \over 2} \mbox{Ram}: & \Lambda & \longrightarrow & {\bf Q}\otimes
\mbox{Pic}\widetilde{S},  \\
                               & \lambda & \longmapsto &
\sum_{ ( \alpha \in {\Phi^+} ) }{{\langle\alpha,\lambda \rangle \over \langle
\alpha,\alpha \rangle} R_{\alpha}}.
\end{array}$
\end{center}
In an obvious sense, a principal $T$-bundle
$\cal T$
on
$\widetilde{S}$
(or a homomorphism
${\cal L}: \Lambda \longrightarrow \mbox{Pic}(\widetilde{S})$)
is $R$-twisted $W$-equivariant if and only if
${\cal T} (-{1 \over 2} Ram)$
is $W$-equivariant, i.e. if
${\cal T}$
 and
${1 \over 2} Ram$
transform the same way under $W$.
The problem with this is that
${1 \over 2} Ram$
itself does not make sense as a $T$-bundle, because the coefficients
${\langle\alpha,\lambda\rangle \over \langle\alpha,\alpha\rangle} $
are not integers. (This argument shows that if
$Hom_{W,R} (\Lambda , \mbox{Pic}(\widetilde{S}))$
is non-empty, it is a torser over the untwisted
$Hom_{W} (\Lambda , \mbox{Pic}(\widetilde{S}))$.)

\begin{thm} \cite{D3}
  \label{main}
For a cameral cover
$\widetilde{S} \longrightarrow S$,
the following data are equivalent: \\
(1) $G$-principal Higgs bundles with cameral cover
$\widetilde{S}$. \\
(2) $R$-twisted $W$-equivariant homomorphisms
${\cal L} \in c^{-1}([N]).$
\end{thm}

The theorem has an essentially local nature, as there is no requirement that
$S$ be, say, projective. We also do not need the condition of generic behavior
near the ramification, which appears in \cite{F, Me, Sc}. Thus we may consider
an extreme case, where
$\widetilde{S}$
is `everywhere ramified':

\begin{eg}\begin{em}                                       \label{nilpo}
In example \ref{ramified}, take the section
$b=0$.
The resulting cover
$\widetilde{C}$
is a `ribbon', or length-2 non-reduced structure on $C$: it is the length-2
neighborhood of $C$ in
$\Bbb{K}$.
The SL(2)-Higgs bundles
$({\cal G},\bdl{c})$
for this
$\widetilde{C}$
have an everywhere nilpotent
$\bdl{c}$,
so the vector bundle
${\cal V} := {\cal G} \times^{SL(2)} V  \approx \pi_* L$
 (where $V$ is the standard 2-dimensional representation) fits in an exact
sequence
$$   0 \longrightarrow  {\cal S} \longrightarrow  {\cal V} \longrightarrow
{\cal Q} \longrightarrow 0
$$
with
${\cal S} \otimes K \approx {\cal Q}.$
Such data are specified by the line bundle
${\cal Q}$,
satisfying
${\cal Q}^{\otimes 2} \approx K$,
and an extension class in
$\mbox{Ext}^1({\cal Q}, {\cal S}) \approx H^1(K^{-1})$.
The kernel of the restriction map
$ \mbox{Pic}(\widetilde{C}) \longrightarrow \mbox{Pic}(C) $
is also given by
$H^1(K^{-1})$
(use the exact sequence
$0 \longrightarrow K^{-1} \longrightarrow  \pi_*{\cal
O}_{\widetilde{C}}^{\times}
   \longrightarrow {\cal O}_C^{\times} \longrightarrow  0$),
and the $R$-twist produces the required square roots of $K$. (For more details
on the nilpotent locus, cf. \cite{L} and \cite{DEL}.)
\end{em}\end{eg}

\subsubsection{Adding values and representations}\
\indent  Fix a vector bundle $K$, and consider the moduli space $ {\cal
M}_{S,G,K} $ of   $K$-valued $G$-principal Higgs bundles on $S$. (It can be
constructed as in Simpson's \cite{simpson-moduli}, even though the objects we
need to parametrize are slightly different than his. In this subsection we
sketch a direct construction.)
It comes with a Hitchin map:

\begin{equation}
  \label{BigHitchin}
  h:   {\cal M}_{S,G,K}  \longrightarrow  B_K
\end{equation}

\noindent where $B := B_K$ parametrizes all possible Hitchin data. Theorem
\ref{main} gives a precise description of the fibers of this map, independent
of the values bundle $K$. This leaves us with the relatively minor task of
describing, for each $K$, the corresponding  base, i.e. the closed subvariety
$B_s$ of $B$ parametrizing {\em split} Hitchin data, or $K$-valued cameral
covers.   The point is that  Higgs bundles satisfy a symmetry condition, which
in Simpson's setup is
$$    \varphi  \wedge \varphi = 0,
$$
and is built into our definition \ref{princHiggs} through the assumption that
\bdl{c} is a bundle of regular centralizers, hence is abelian. Since commuting
operators have common eigenvectors, this is translated into a splitness
condition on the Hitchin data, which we  describe below.  (When $K$ is a line
bundle, the condition is vacuous, $B_s = B$.) The upshot is:

\begin{lem}
           \label{parametrization}
The following data are equivalent: \\
(a) A $K$-valued cameral cover of $S$. \\
(b) A split, graded homomorphism
     $R{\bf \dot{\ }} \longrightarrow {Sym}{\bf \dot{\ }}K.$ \\
(c) A split Hitchin datum
$b \in B_s$.
\end{lem}

Here $R{\bf \dot{\ }}$ is the graded ring of $W$-invariant polynomials on
$\frak t$:

\begin{equation}
R{\bf \dot{\ }} := (\mbox{Sym}{\bf \dot{\ }} {\frak t}^*)^W
\approx {\bf C}[\sigma_1,\ldots,\sigma_l],  \qquad  \deg (\sigma_i) = d_i
\end{equation}

\noindent
where
$l := \mbox{Rank}({\frak g})$
and the
$\sigma_i$
form a basis for the $W$-invariant polynomials. The Hitchin base is the vector
space
$$  B :=  B_K  :=  \oplus _{i=1}^l H^0(S, {Sym}^{d_i}K)
\approx \mbox{Hom}(R{\bf \dot{\ }},\mbox{Sym}{\bf \dot{\ }}K).
$$

\noindent For each
$\lambda \in \Lambda$
(or
$\lambda \in {\frak t}^*$,
for that matter), the expression in an indeterminate $x$:

\begin{equation}
               \label{rep poly}
q_{\lambda}(x,t) := \prod_{w \in W}{(x-w\lambda(t))}, \qquad  t \in {\frak t},
\end{equation}

\noindent is $W$-invariant (as a function of $t$), so it defines an element
$q_{\lambda}(x) \in R{\bf \dot{\ }}[x].$
A Hitchin datum
$b \in B \approx \mbox{Hom}(R{\bf \dot{\ }},\mbox{Sym}{\bf \dot{\ }}K)$
sends this to
$$  q_{\lambda,b}(x) \in \mbox{Sym}\dot{\ }(K)[x].
$$
We say that $b$ is {\em split} if, at each point of $S$ and for each
$\lambda$,
the polynomial
$q_{\lambda,b}(x)$
factors completely, into terms linear in $x$.

We note that, for $\lambda$ in the interior of $C$ (the positive Weyl chamber),
$q_{\lambda,b}$
gives the equation in
$\Bbb K$
of the spectral cover
$\widetilde{S}_{\lambda}$
of section (\ref{decomp covers}):
$q_{\lambda,b}$ gives a morphism
$\Bbb K  \longrightarrow \mbox{Sym}^N \Bbb K$,
where $N:=\#W$, and $\widetilde{S}_{\lambda}$ is the inverse image of the
zero-section.
 (When
$\lambda$
is in a face
$F_P$
of
$\overline{C}$,
we define analogous polynomials
$q_{\lambda}^P(x,t)$
and
$q_{\lambda,b}^P(x)$
by taking the product in (\ref{rep poly}) to be over
$w \in W_P \backslash W.$
These give the reduced equations in this case, and
$q_{\lambda}$
is an appropriate power.)

Over $B_s$ there is a universal $K$-valued cameral cover
$$   \widetilde{\cal S}  \longrightarrow  B_s
$$
with ramification divisor $R \subset \widetilde{\cal S}$. From the relative
Picard,
$$   \mbox{Pic}( \widetilde{\cal S}  /  B_s)
$$
we concoct the relative $N$-shifted, $R$-twisted  Prym
$$   \mbox{Prym}_{\Lambda ,R}( \widetilde{\cal S}  /  B_s).
$$
By Theorem (\ref{main}), this can then be considered as a parameter space
$ {\cal M}_{S,G,K} $
for all $K$-valued $G$-principal Higgs bundles on $S$. (Recall that our objects
are assumed to be everywhere {\em regular}!)  It comes with a `Hitchin map',
namely the projection to $B_s$, and the fibers corresponding to smooth
projective $\widetilde{S}$ are abelian varieties. When $S$ is a smooth,
projective curve, we recover this way the algebraic complete integrability of
Hitchin's system and its generalizations. More generally, for any $S$, one
obtains an ACIHS (with symplectic, respectively Poisson structures) when the
values bundle has the same (symplectic, respectively Poisson) structure, by a
slight modification of the construction in Chapter \ref{ch8}. One considers
only Lagrangian supports which retain a $W$-action, and only equivariant
sheaves on them (with the numerical invariants of a line bundle). These two
restrictions are symplecticly dual, so the moduli space of Lagrangian sheaves
with these invariance properties is a symplectic (respectively, Poisson) subsp!
ace of
 the total  moduli space, and the
fibers of the Hitchin map are Lagrangian as expected.

\subsubsection{Irregulars?}    \nopagebreak

\noindent The Higgs bundles we consider in this survey are assumed to be
everywhere regular. This is a reasonable assumption for line-bundle valued
Higgs bundles on a curve or surface, but {\em not} in $\dim \geq 3$. This is
because the complement of ${\frak g}_{{reg}}$ has codimension 3 in ${\frak g}$.
The source of the difficulty is that the analogue of (\ref{commutes}) fails
over
${\frak g}$.  There are two candidates for the universal cameral cover:
$\widetilde{\frak g}$, defined by the left hand side of (\ref{commutes}), is
finite over ${\frak g}$ with $W$ action, but does not have a family of line
bundles parametrized by $\Lambda$.
These live instead on   $\stackrel{\approx}{\frak g}$, the object defined by
the right hand side, which parametrizes pairs
$(x,{\frak b}),    \qquad          x \in {\frak b} \subset {\frak g}$ .
This suggests that the right way to analyze irregular Higgs bundles may involve
spectral data consisting of a tower
$$  \stackrel{\approx}{S}  \stackrel{\sigma}{\longrightarrow}  \widetilde{S}
\longrightarrow  S
$$
together with a homomorphism
$ {\cal L} : \Lambda \longrightarrow  \mbox{Pic}(\stackrel{\approx}{S})$
such that  the collection of sheaves
$$  \sigma_*({\cal L}(\lambda)),   \qquad         \lambda \in \Lambda
$$
 on
$\widetilde{S}$
is
$R$-twisted $W$-equivariant in an appropriate sense.  As a first step, one may
wish to understand the direct images
$ R^i \sigma_*({\cal L}(\lambda))  $
and in particular the cohomologies
$H^i(F, {\cal L}(\lambda))$
where $F$, usually called a {\em Springer fiber}, is a fiber of $\sigma$. For
regular $x$, this fiber is a single point. For $x=0$, the fiber is all of
$G/B$, so the fiber cohomology is given by the Borel-Weil-Bott theorem. The
question may thus be considered as a desired extension of BWB to general
Springer fibers.
\newpage

\end{document}